%% file: main.tex
\documentclass{ldr-article}

\title{Persistent Homology of Configuration Spaces of Trees}

\author{
Wenwen Li\\
}

\date{}

\begin{document}

    \maketitle

    \begin{abstract}
       Multiparameter persistence modules come up naturally in topological data analysis and topological robotics. Given a metric graph $(X,\delta)$, the second configuration space of $(X,\delta)$ with proximity parameters (for example, the minimum distance allowed between each pair of robots) can be interpreted as a collection of all possible configurations of two robots moving in $(X,\delta)$. In this project, we study the $2$-parameter persistence modules associated with the second configuration spaces of the star graph ($\mathsf{Star}_k$, $k\geq 3$) and the generalized H-graph ($\mathcal{H}_{m,n}$, $m,n\geq 3$) with the edge length parameter $L_e$ and the restraint parameter $r$. Moreover, we provide the indecomposable direct summands for each persistence module. 
    \end{abstract}

    \tableofcontents
    
   \section{Introduction}
   \input{chapters/intro}

   \section{Preliminaries}\label{prelim}
   \input{chapters/section1.2}
   \input{chapters/section1.3}

   \input{chapters/section1.6}

\input{chapters/simplification}   
\input{chapters/chapter4-star-nofac}

\input{chapters/chapter4-H}

\printbibliography    
\begin{appendix}
\input{chapters/section1.2.5}

\end{appendix}
\end{document}

%% file: chapters/intro.tex
Configuration space is one of the key notions used in robot motion planning. Given a topological space $X$, the $n$-th configuration space of $X$ is $X^{\underline{n}}=\{(x_{1},\dots, x_{n})\in X^{n}\mid x_{i}\neq x_{j} \mbox{ if } i\neq j\}$. The $n$-th configuration space of a topological space $X$ can be interpreted as the collection of all possible arrangements of $n$ distinct robots (treated as points) in $X$ such that they can appear without colliding with others.\\

In the real-world scenario, however, there is a minimal distance allowed between each pair of robots. Let $(X,\delta)$ be a metric space and let $\mathbf{r}=(r_{ij})_{i<j}\in\mathbb{R}_{>0}^{\binom{n}{2}}$. The \textbf{$n$-th configuration space with restraint parameter $\mathbf{r}=(r_{ij})_{i<j}$} is: $$X^{n}_{\mathbf{r}}=\{(x_{1},\dots, x_{n})\in X^{n}\mid \delta(x_{i},x_{j})\geq r_{ij}\}$$ 

The topology of the configuration spaces with restraint parameters has been studied by various researchers. Deeley \cite{Dee} defines configuration spaces of thick particles on a graph $X$ and studies the second configuration space of thick particles. In \cite{carlsson2012computational}\cite{baryshnikov2014min}
\cite{alpert2020configuration}\cite{alpert2021configuration}, the configuration spaces of a given topological space $X$ with the critical values of the radius, assuming that all robots are disks with the same radius are considered. Alpert and Manin \cite{alpert2021configuration} investigate the topology of configuration space of $n$ distinct unit disks in an infinite strip.\\

In this paper, we study the second configuration space with restraint parameter, where its underlying space is a metric graph. Let $\Gamma=(V,E)$ be a finite connected graph\footnote{$\Gamma$ is allowed to have multiple edges. If $\Gamma$ has loops, we subdivide each loop of $\Gamma$ into two edges. We are slightly abusing the notation by denoting the new graph by $\Gamma$.} and $X=\lVert\Gamma\rVert$ be the (canonical) geometric realization of $\Gamma$. $X$ is a regular CW complex. Given a vector $\mathbf{L}=(L_e)_{e\in E}\in \mathbb{R}_{>0}^{\lvert E\rvert}$, $X$ is called a \textbf{metric graph} if a positive number $L_{e}$ is assigned to each edge $e$ for all $e\in E$, which represents its length. $X$ has the path metric $\delta$: for any $x,y\in X$, $\delta(x,y)$ is the length of a shortest path from $x$ to $y$. We use $X_{\mathbf{L}}$ to denote this metric space $(X,\delta)$. The \textbf{$n$-th configuration space} of a metric graph $(X,\delta)$ with restraint parameter $\mathbf{r}\in \mathbb{R}_{>0}^{\binom{n}{2}}$ and edge length vector $\mathbf{L}$ is 
$$X^{n}_{\mathbf{r},\mathbf{L}}=\{(x_{1},\dots, x_{n})\in (X_{\mathbf{L}})^{n}\mid \delta(x_{i},x_{j})\geq r_{ij}, \mbox{ for } i,j=1,\dots, n\}$$

Although the metric graph $X$ is a $1$-dimensional regular CW-complex which is less complicated than an arbitrary topological space, nevertheless, the homology groups of $X^{n}_{\mathbf{r},\mathbf{L}}$ are still hard to compute. Dover and \"Ozayd\i n\cite{dover2013homeomorphism} investigate the homeomorphism types of $X^n_{\mathbf{r}}$ over the space of parameters $r$ and provide a polynomial (of the number of the edges of $X$) upper bound for the number of isotopy types of $X^n_{\mathbf{r}}$.\\

Notice that there is a natural filtration $X^n_{-,-}:((\mathbb{R}_{>0},\leq)^{\binom{n}{2}})^{\op}\times (\mathbb{R}_{>0},\leq)^{\lvert E\rvert}\rightarrow \cat{Top}$, and it comes with $\binom{n}{2}+\lvert E\rvert$ parameters: $\mathbf{r}=(r_{ij})_{i<j}$ and $\mathbf{L}=(L_e)_{e\in E}$. When all robots have the same size and all edges have fixed lengths, in particular, the filtration reduces to a single-parameter filtration. Persistent homology is one of the standard tools to study the topological features of this filtration $X^n_{-,-}$. It is developed in the area of Topological Data Analysis (TDA) to analyze the shape of data. Let $\mathbb{F}$ be a field and $(P,\leq)$ be the set $P$ with a partial order. A persistence module over $(P,\leq)$ is a family of $\mathbb{F}$-vector spaces $\{Mt\mid t\in P\}$ and a doubly-indexed family of linear maps $\{\rho_{st}: Ms\rightarrow Mt\mid s\leq t\}$ where $\rho_{tu}\rho_{st}=\rho_{su}$ for any $s\leq t\leq u$ in $P$ and $\rho_{ss}=\mathrm{id}_{Ms}$ for all $s\in P$. \\

For single-parameter pointwise finite-dimensional persistence modules (or, more generally, persistence modules over totally ordered sets), Botman and Crawley-Boevey \cite{botnan2020decomposition} prove that such modules can be decomposed into a direct sum of interval modules. Therefore, one can obtain a barcode for a single-parameter pointwise finite-dimensional persistence module by assigning an interval (of the indexing totally ordered set) to each direct summand. A barcode can be viewed as the persistence analog of a Betti number. \cite{ghrist2008barcodes}\\

In the multiparameter persistence theory, however, there is no canonical way to define the barcodes (which represent indecomposable representations) analogous to the barcodes of the single-parameter persistence modules because there is no structure theorem available in $\cat{vect}_{\mathbb{F}}^{(P,\leq)}$ when the poset $P$ is not a totally ordered set. The indecomposable submodules of a multiparameter persistence module can be very complicated. A sequence of papers by Buchet and Escolar\cite{BE-2020}\cite{BE-2019} show that any $n$-parameter persistence module can be embedded as a slice of an indecomposable (n+1)-parameter persistence module. In this paper, we investigate the multiparameter persistence modules given by filtration of the second configuration spaces of metric tree graph $X$ with restraint parameter $r$ and edge length vector $L$.\\

\paragraph{Main results} We now discuss the main results of this paper. Since every tree is a union of star graphs and generalized H graphs, the $i$-th persistent homology group of the second configuration space of a tree with parameters $r$ and $L$ (denoted by $T^2_{r,L}$) can be calculated inductively using the Mayer-Vietoris spectral sequence, where $(\mathsf{Star}_k)^2_{r,L}$ and $(\mathcal{H}_{m,n})^2_{r,L}$ are a part of the cover of $T^2_{r,L}$. Note that, for a given metric graph $X$, at the object level, $PH_i(X^2_{r,L})$ is equal to the $i$-th homology group of $X^2_{r,L}$ for all parameters $r$ and $L$. In this paper, we focus on calculating $PH_i(X^2_{-,-})$ when $X$ is a metric star graph or generalized H graph.\\

\begin{minipage}{0.33\textwidth}
\centering
\resizebox{3cm}{!}{

\tikzset{every picture/.style={line width=0.75pt}} 

\begin{tikzpicture}[x=0.75pt,y=0.75pt,yscale=-1,xscale=1]

\draw    (100,113) -- (66,61.5) ;
\draw [shift={(66,61.5)}, rotate = 236.57] [color={rgb, 255:red, 0; green, 0; blue, 0 }  ][fill={rgb, 255:red, 0; green, 0; blue, 0 }  ][line width=0.75]      (0, 0) circle [x radius= 3.35, y radius= 3.35]   ;
\draw [shift={(100,113)}, rotate = 236.57] [color={rgb, 255:red, 0; green, 0; blue, 0 }  ][fill={rgb, 255:red, 0; green, 0; blue, 0 }  ][line width=0.75]      (0, 0) circle [x radius= 3.35, y radius= 3.35]   ;
\draw    (100,113) -- (127,55.5) ;
\draw [shift={(127,55.5)}, rotate = 295.15] [color={rgb, 255:red, 0; green, 0; blue, 0 }  ][fill={rgb, 255:red, 0; green, 0; blue, 0 }  ][line width=0.75]      (0, 0) circle [x radius= 3.35, y radius= 3.35]   ;
\draw    (100,113) -- (163,113) ;
\draw [shift={(163,113)}, rotate = 0] [color={rgb, 255:red, 0; green, 0; blue, 0 }  ][fill={rgb, 255:red, 0; green, 0; blue, 0 }  ][line width=0.75]      (0, 0) circle [x radius= 3.35, y radius= 3.35]   ;
\draw [shift={(100,113)}, rotate = 0] [color={rgb, 255:red, 0; green, 0; blue, 0 }  ][fill={rgb, 255:red, 0; green, 0; blue, 0 }  ][line width=0.75]      (0, 0) circle [x radius= 3.35, y radius= 3.35]   ;
\draw    (43,127.5) -- (100,113) ;
\draw [shift={(100,113)}, rotate = 345.73] [color={rgb, 255:red, 0; green, 0; blue, 0 }  ][fill={rgb, 255:red, 0; green, 0; blue, 0 }  ][line width=0.75]      (0, 0) circle [x radius= 3.35, y radius= 3.35]   ;
\draw [shift={(43,127.5)}, rotate = 345.73] [color={rgb, 255:red, 0; green, 0; blue, 0 }  ][fill={rgb, 255:red, 0; green, 0; blue, 0 }  ][line width=0.75]      (0, 0) circle [x radius= 3.35, y radius= 3.35]   ;
\draw    (72,175) -- (100,113) ;
\draw [shift={(100,113)}, rotate = 294.3] [color={rgb, 255:red, 0; green, 0; blue, 0 }  ][fill={rgb, 255:red, 0; green, 0; blue, 0 }  ][line width=0.75]      (0, 0) circle [x radius= 3.35, y radius= 3.35]   ;
\draw [shift={(72,175)}, rotate = 294.3] [color={rgb, 255:red, 0; green, 0; blue, 0 }  ][fill={rgb, 255:red, 0; green, 0; blue, 0 }  ][line width=0.75]      (0, 0) circle [x radius= 3.35, y radius= 3.35]   ;
\draw    (146,158.5) -- (100,113) ;
\draw [shift={(100,113)}, rotate = 224.69] [color={rgb, 255:red, 0; green, 0; blue, 0 }  ][fill={rgb, 255:red, 0; green, 0; blue, 0 }  ][line width=0.75]      (0, 0) circle [x radius= 3.35, y radius= 3.35]   ;
\draw [shift={(146,158.5)}, rotate = 224.69] [color={rgb, 255:red, 0; green, 0; blue, 0 }  ][fill={rgb, 255:red, 0; green, 0; blue, 0 }  ][line width=0.75]      (0, 0) circle [x radius= 3.35, y radius= 3.35]   ;

\draw (91.64,175.93) node [anchor=north west][inner sep=0.75pt]  [rotate=-346.72]  {$\cdots $};
\draw (128,97.4) node [anchor=north west][inner sep=0.75pt]  [color={rgb, 255:red, 15; green, 8; blue, 248 }  ,opacity=1 ]  {$e_{1}$};

\end{tikzpicture}
}

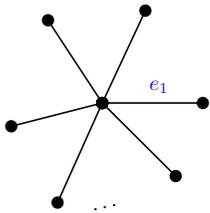
\captionof{figure}{Star graph}
\label{fig:star-intro}  
\end{minipage}%
\begin{minipage}{0.33\textwidth}
\centering
\resizebox{3.5cm}{!}{

\tikzset{every picture/.style={line width=0.75pt}} 

\begin{tikzpicture}[x=0.75pt,y=0.75pt,yscale=-1,xscale=1]

\draw    (40,32.7) -- (83.99,76.69) ;
\draw [shift={(83.99,76.69)}, rotate = 45] [color={rgb, 255:red, 0; green, 0; blue, 0 }  ][fill={rgb, 255:red, 0; green, 0; blue, 0 }  ][line width=0.75]      (0, 0) circle [x radius= 3.35, y radius= 3.35]   ;
\draw [shift={(40,32.7)}, rotate = 45] [color={rgb, 255:red, 0; green, 0; blue, 0 }  ][fill={rgb, 255:red, 0; green, 0; blue, 0 }  ][line width=0.75]      (0, 0) circle [x radius= 3.35, y radius= 3.35]   ;
\draw    (83.99,76.69) -- (41.63,119.05) ;
\draw [shift={(41.63,119.05)}, rotate = 135] [color={rgb, 255:red, 0; green, 0; blue, 0 }  ][fill={rgb, 255:red, 0; green, 0; blue, 0 }  ][line width=0.75]      (0, 0) circle [x radius= 3.35, y radius= 3.35]   ;
\draw [shift={(83.99,76.69)}, rotate = 135] [color={rgb, 255:red, 0; green, 0; blue, 0 }  ][fill={rgb, 255:red, 0; green, 0; blue, 0 }  ][line width=0.75]      (0, 0) circle [x radius= 3.35, y radius= 3.35]   ;
\draw    (145.9,76.69) -- (189.4,33.19) ;
\draw [shift={(189.4,33.19)}, rotate = 315] [color={rgb, 255:red, 0; green, 0; blue, 0 }  ][fill={rgb, 255:red, 0; green, 0; blue, 0 }  ][line width=0.75]      (0, 0) circle [x radius= 3.35, y radius= 3.35]   ;
\draw [shift={(145.9,76.69)}, rotate = 315] [color={rgb, 255:red, 0; green, 0; blue, 0 }  ][fill={rgb, 255:red, 0; green, 0; blue, 0 }  ][line width=0.75]      (0, 0) circle [x radius= 3.35, y radius= 3.35]   ;
\draw    (145.9,76.69) -- (189.89,120.68) ;
\draw [shift={(189.89,120.68)}, rotate = 45] [color={rgb, 255:red, 0; green, 0; blue, 0 }  ][fill={rgb, 255:red, 0; green, 0; blue, 0 }  ][line width=0.75]      (0, 0) circle [x radius= 3.35, y radius= 3.35]   ;
\draw [shift={(145.9,76.69)}, rotate = 45] [color={rgb, 255:red, 0; green, 0; blue, 0 }  ][fill={rgb, 255:red, 0; green, 0; blue, 0 }  ][line width=0.75]      (0, 0) circle [x radius= 3.35, y radius= 3.35]   ;
\draw    (83.99,76.69) -- (146.01,76.69) ;
\draw [shift={(146.01,76.69)}, rotate = 0] [color={rgb, 255:red, 0; green, 0; blue, 0 }  ][fill={rgb, 255:red, 0; green, 0; blue, 0 }  ][line width=0.75]      (0, 0) circle [x radius= 3.35, y radius= 3.35]   ;
\draw [shift={(83.99,76.69)}, rotate = 0] [color={rgb, 255:red, 0; green, 0; blue, 0 }  ][fill={rgb, 255:red, 0; green, 0; blue, 0 }  ][line width=0.75]      (0, 0) circle [x radius= 3.35, y radius= 3.35]   ;

\draw (109,81.4) node [anchor=north west][inner sep=0.75pt][color={rgb, 255:red, 15; green, 8; blue, 248 }  ,opacity=1 ]    {$e_{1}$};

\end{tikzpicture}
}
\captionof{figure}{H graph}
\label{fig:H-intro}  
\end{minipage}%
\begin{minipage}{0.33\textwidth}
\centering
\resizebox{4cm}{!}{

\tikzset{every picture/.style={line width=0.75pt}} 

\begin{tikzpicture}[x=0.75pt,y=0.75pt,yscale=-1,xscale=1]

\draw    (40,32.7) -- (83.99,76.69) ;
\draw [shift={(83.99,76.69)}, rotate = 45] [color={rgb, 255:red, 0; green, 0; blue, 0 }  ][fill={rgb, 255:red, 0; green, 0; blue, 0 }  ][line width=0.75]      (0, 0) circle [x radius= 3.35, y radius= 3.35]   ;
\draw [shift={(40,32.7)}, rotate = 45] [color={rgb, 255:red, 0; green, 0; blue, 0 }  ][fill={rgb, 255:red, 0; green, 0; blue, 0 }  ][line width=0.75]      (0, 0) circle [x radius= 3.35, y radius= 3.35]   ;
\draw    (83.99,76.69) -- (41.63,119.05) ;
\draw [shift={(41.63,119.05)}, rotate = 135] [color={rgb, 255:red, 0; green, 0; blue, 0 }  ][fill={rgb, 255:red, 0; green, 0; blue, 0 }  ][line width=0.75]      (0, 0) circle [x radius= 3.35, y radius= 3.35]   ;
\draw [shift={(83.99,76.69)}, rotate = 135] [color={rgb, 255:red, 0; green, 0; blue, 0 }  ][fill={rgb, 255:red, 0; green, 0; blue, 0 }  ][line width=0.75]      (0, 0) circle [x radius= 3.35, y radius= 3.35]   ;
\draw    (145.9,76.69) -- (189.4,33.19) ;
\draw [shift={(189.4,33.19)}, rotate = 315] [color={rgb, 255:red, 0; green, 0; blue, 0 }  ][fill={rgb, 255:red, 0; green, 0; blue, 0 }  ][line width=0.75]      (0, 0) circle [x radius= 3.35, y radius= 3.35]   ;
\draw [shift={(145.9,76.69)}, rotate = 315] [color={rgb, 255:red, 0; green, 0; blue, 0 }  ][fill={rgb, 255:red, 0; green, 0; blue, 0 }  ][line width=0.75]      (0, 0) circle [x radius= 3.35, y radius= 3.35]   ;
\draw    (145.9,76.69) -- (189.89,120.68) ;
\draw [shift={(189.89,120.68)}, rotate = 45] [color={rgb, 255:red, 0; green, 0; blue, 0 }  ][fill={rgb, 255:red, 0; green, 0; blue, 0 }  ][line width=0.75]      (0, 0) circle [x radius= 3.35, y radius= 3.35]   ;
\draw [shift={(145.9,76.69)}, rotate = 45] [color={rgb, 255:red, 0; green, 0; blue, 0 }  ][fill={rgb, 255:red, 0; green, 0; blue, 0 }  ][line width=0.75]      (0, 0) circle [x radius= 3.35, y radius= 3.35]   ;
\draw    (83.99,76.69) -- (146.01,76.69) ;
\draw [shift={(146.01,76.69)}, rotate = 0] [color={rgb, 255:red, 0; green, 0; blue, 0 }  ][fill={rgb, 255:red, 0; green, 0; blue, 0 }  ][line width=0.75]      (0, 0) circle [x radius= 3.35, y radius= 3.35]   ;
\draw [shift={(83.99,76.69)}, rotate = 0] [color={rgb, 255:red, 0; green, 0; blue, 0 }  ][fill={rgb, 255:red, 0; green, 0; blue, 0 }  ][line width=0.75]      (0, 0) circle [x radius= 3.35, y radius= 3.35]   ;
\draw    (27,57.5) -- (83.99,76.69) ;
\draw [shift={(83.99,76.69)}, rotate = 18.61] [color={rgb, 255:red, 0; green, 0; blue, 0 }  ][fill={rgb, 255:red, 0; green, 0; blue, 0 }  ][line width=0.75]      (0, 0) circle [x radius= 3.35, y radius= 3.35]   ;
\draw [shift={(27,57.5)}, rotate = 18.61] [color={rgb, 255:red, 0; green, 0; blue, 0 }  ][fill={rgb, 255:red, 0; green, 0; blue, 0 }  ][line width=0.75]      (0, 0) circle [x radius= 3.35, y radius= 3.35]   ;
\draw    (146.01,76.69) -- (203,55.5) ;
\draw [shift={(203,55.5)}, rotate = 339.6] [color={rgb, 255:red, 0; green, 0; blue, 0 }  ][fill={rgb, 255:red, 0; green, 0; blue, 0 }  ][line width=0.75]      (0, 0) circle [x radius= 3.35, y radius= 3.35]   ;
\draw [shift={(146.01,76.69)}, rotate = 339.6] [color={rgb, 255:red, 0; green, 0; blue, 0 }  ][fill={rgb, 255:red, 0; green, 0; blue, 0 }  ][line width=0.75]      (0, 0) circle [x radius= 3.35, y radius= 3.35]   ;

\draw (37.86,78.63) node [anchor=north west][inner sep=0.75pt]  [rotate=-63.61][color={rgb, 255:red, 15; green, 8; blue, 248 }  ,opacity=1 ]  {$\cdots $};
\draw (182.61,105.22) node [anchor=north west][inner sep=0.75pt] [rotate=-284.03][color={rgb, 255:red, 15; green, 8; blue, 248 }  ,opacity=1 ]  {$\cdots $};
\draw (107,81.4) node [anchor=north west][inner sep=0.75pt][color={rgb, 255:red, 15; green, 8; blue, 248 }  ,opacity=1 ]    {$e_{1}$};

\end{tikzpicture}
}

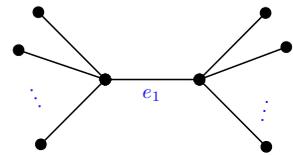
\captionof{figure}{Generalized H graph}
\label{fig:Hmn-intro}  
\end{minipage}

First, we calculated the rank of $H_i(Y^2_{r,L_{e_1}})$ for $i=0,1$ and $Y$ is the geometric realization of the Y-shaped graph, where every edge but $e_1$ has length $1$. We then construct a persistence module $PH_i(Y^2_{-,-})$ for $i=0,1$ and prove the following theorems regarding the interval decomposability of $PH_0(Y^2_{-,-};\mathbb{F})$ and $PH_1(Y^2_{-,-};\mathbb{F})$.

\begin{theorem}\label{Thm2'}
$PH_1(Y^2_{-,-};\mathbb{F})$ is an interval module, hence it is indecomposable.
\end{theorem}

\begin{theorem}\label{Thm3'}
$PH_0(Y^2_{-,-};\mathbb{F})$ is not interval decomposable. The indecomposables of $PH_0(Y^2_{-,-};\mathbb{F})$ is given in Figure \ref{H0summand-0} up to isomorphism. 
\end{theorem}
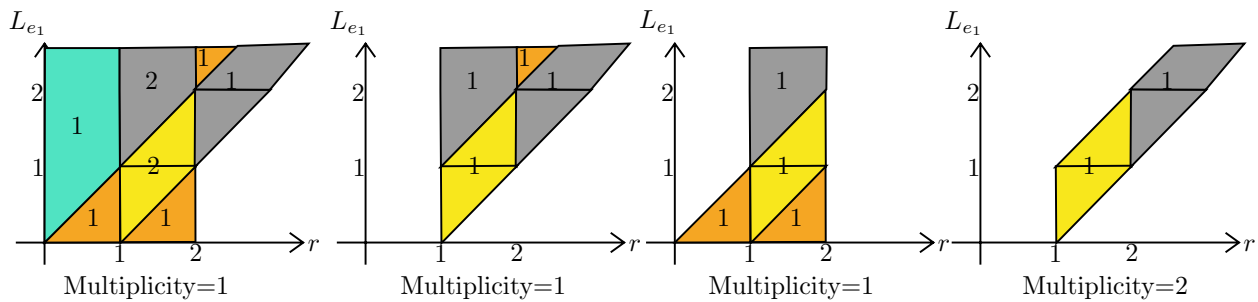
\begin{figure}[htbp!]
\centering
\resizebox{17cm}{!}{

\tikzset{every picture/.style={line width=0.75pt}} 

\begin{tikzpicture}[x=0.75pt,y=0.75pt,yscale=-1,xscale=1]

\draw  [draw opacity=0][fill={rgb, 255:red, 155; green, 155; blue, 155 }  ,fill opacity=1 ] (59.54,18.69) -- (99.37,18.22) -- (99.16,39.9) -- (59.81,79.8) -- (59.54,40.36) -- cycle ;
\draw  [draw opacity=0][fill={rgb, 255:red, 248; green, 231; blue, 28 }  ,fill opacity=1 ] (59.49,79.94) -- (98.7,40.36) -- (98.86,79.78) -- cycle ;
\draw  (6.06,118.47) -- (154.47,118.47)(20.9,16.37) -- (20.9,129.82) (147.47,113.47) -- (154.47,118.47) -- (147.47,123.47) (15.9,23.37) -- (20.9,16.37) -- (25.9,23.37)  ;
\draw  [draw opacity=0][fill={rgb, 255:red, 245; green, 166; blue, 35 }  ,fill opacity=1 ] (20.9,118.47) -- (59.81,79.8) -- (59.91,118.37) -- cycle ;
\draw  [draw opacity=0][fill={rgb, 255:red, 245; green, 166; blue, 35 }  ,fill opacity=1 ] (60.09,118.24) -- (98.73,79.4) -- (98.83,118.14) -- cycle ;
\draw  [draw opacity=0][fill={rgb, 255:red, 245; green, 166; blue, 35 }  ,fill opacity=1 ] (120.57,18.03) -- (99.16,39.9) -- (98.94,18.25) -- cycle ;
\draw  [draw opacity=0][fill={rgb, 255:red, 155; green, 155; blue, 155 }  ,fill opacity=1 ] (137.56,40.08) -- (98.4,79.78) -- (98.63,39.86) -- cycle ;
\draw  [draw opacity=0][fill={rgb, 255:red, 155; green, 155; blue, 155 }  ,fill opacity=1 ] (120.57,17.57) -- (157.25,16.37) -- (137.26,40.04) -- (137.56,40.08) -- (127.78,40.04) -- (98.23,39.9) -- cycle ;
\draw  [draw opacity=0][fill={rgb, 255:red, 80; green, 227; blue, 194 }  ,fill opacity=1 ] (59.54,18.69) -- (59.54,40.36) -- (59.81,79.8) -- (20.9,118.47) -- (21.11,18.69) -- cycle ;
\draw  [draw opacity=0][fill={rgb, 255:red, 248; green, 231; blue, 28 }  ,fill opacity=1 ] (99.19,79.08) -- (60.08,118.75) -- (59.81,79.34) -- cycle ;
\draw  [draw opacity=0][fill={rgb, 255:red, 155; green, 155; blue, 155 }  ,fill opacity=1 ] (225.28,18.67) -- (265.1,18.2) -- (264.9,39.88) -- (225.55,79.78) -- (225.28,40.34) -- cycle ;
\draw  [draw opacity=0][fill={rgb, 255:red, 248; green, 231; blue, 28 }  ,fill opacity=1 ] (225.22,79.92) -- (264.44,40.34) -- (264.6,79.76) -- cycle ;
\draw  (171.8,118.45) -- (320.21,118.45)(186.64,16.35) -- (186.64,129.8) (313.21,113.45) -- (320.21,118.45) -- (313.21,123.45) (181.64,23.35) -- (186.64,16.35) -- (191.64,23.35)  ;
\draw  [draw opacity=0][fill={rgb, 255:red, 245; green, 166; blue, 35 }  ,fill opacity=1 ] (286.31,18.01) -- (264.9,39.88) -- (264.68,18.23) -- cycle ;
\draw  [draw opacity=0][fill={rgb, 255:red, 155; green, 155; blue, 155 }  ,fill opacity=1 ] (303.3,40.06) -- (264.14,79.76) -- (264.37,39.84) -- cycle ;
\draw  [draw opacity=0][fill={rgb, 255:red, 155; green, 155; blue, 155 }  ,fill opacity=1 ] (286.31,17.55) -- (322.99,16.35) -- (303,40.02) -- (303.3,40.06) -- (263.97,39.88) -- cycle ;
\draw  [draw opacity=0][fill={rgb, 255:red, 248; green, 231; blue, 28 }  ,fill opacity=1 ] (264.93,79.05) -- (225.81,118.73) -- (225.55,79.32) -- cycle ;
\draw  [draw opacity=0][fill={rgb, 255:red, 155; green, 155; blue, 155 }  ,fill opacity=1 ] (384.94,18.67) -- (424.76,18.2) -- (424.55,39.88) -- (385.21,79.78) -- (384.94,40.34) -- cycle ;
\draw  [draw opacity=0][fill={rgb, 255:red, 248; green, 231; blue, 28 }  ,fill opacity=1 ] (384.88,79.92) -- (424.09,40.34) -- (424.26,79.76) -- cycle ;
\draw  (331.45,118.45) -- (479.86,118.45)(346.3,16.35) -- (346.3,129.8) (472.86,113.45) -- (479.86,118.45) -- (472.86,123.45) (341.3,23.35) -- (346.3,16.35) -- (351.3,23.35)  ;
\draw  [draw opacity=0][fill={rgb, 255:red, 245; green, 166; blue, 35 }  ,fill opacity=1 ] (346.3,118.45) -- (385.21,79.78) -- (385.31,118.35) -- cycle ;
\draw  [draw opacity=0][fill={rgb, 255:red, 245; green, 166; blue, 35 }  ,fill opacity=1 ] (385.48,118.22) -- (424.12,79.38) -- (424.22,118.11) -- cycle ;
\draw  [draw opacity=0][fill={rgb, 255:red, 248; green, 231; blue, 28 }  ,fill opacity=1 ] (424.58,79.05) -- (385.47,118.73) -- (385.21,79.32) -- cycle ;
\draw  [draw opacity=0][fill={rgb, 255:red, 248; green, 231; blue, 28 }  ,fill opacity=1 ] (542.66,79.91) -- (581.88,40.33) -- (582.04,79.75) -- cycle ;
\draw  (489.24,118.44) -- (637.65,118.44)(504.08,16.34) -- (504.08,129.79) (630.65,113.44) -- (637.65,118.44) -- (630.65,123.44) (499.08,23.34) -- (504.08,16.34) -- (509.08,23.34)  ;
\draw  [draw opacity=0][fill={rgb, 255:red, 155; green, 155; blue, 155 }  ,fill opacity=1 ] (620.74,40.05) -- (581.58,79.75) -- (581.81,39.83) -- cycle ;
\draw  [draw opacity=0][fill={rgb, 255:red, 155; green, 155; blue, 155 }  ,fill opacity=1 ] (603.75,17.53) -- (640.43,16.34) -- (620.44,40) -- (620.74,40.05) -- (581.41,39.87) -- cycle ;
\draw  [draw opacity=0][fill={rgb, 255:red, 248; green, 231; blue, 28 }  ,fill opacity=1 ] (582.37,79.04) -- (543.25,118.72) -- (542.99,79.3) -- cycle ;

\draw (59.78,124.26) node   [align=left] {1};
\draw (99.14,123.8) node   [align=left] {2};
\draw (17.64,79.81) node   [align=left] {1};
\draw (17.17,41.38) node   [align=left] {2};
\draw (32.94,52.76) node [anchor=north west][inner sep=0.75pt]    {$1$};
\draw (41.27,99.99) node [anchor=north west][inner sep=0.75pt]    {$1$};
\draw (78.78,99.99) node [anchor=north west][inner sep=0.75pt]    {$1$};
\draw (72.3,73.14) node [anchor=north west][inner sep=0.75pt]    {$2$};
\draw (71.37,29.61) node [anchor=north west][inner sep=0.75pt]    {$2$};
\draw (98.23,18.03) node [anchor=north west][inner sep=0.75pt]    {$1$};
\draw (112.58,29.61) node [anchor=north west][inner sep=0.75pt]    {$1$};
\draw (0.94,-1.61) node [anchor=north west][inner sep=0.75pt]    {$L_{e_{1}}{}$};
\draw (155.26,115.58) node [anchor=north west][inner sep=0.75pt]    {$r$};
\draw (225.51,124.24) node   [align=left] {1};
\draw (264.87,123.78) node   [align=left] {2};
\draw (183.38,79.79) node   [align=left] {1};
\draw (182.91,41.36) node   [align=left] {2};
\draw (238.04,73.12) node [anchor=north west][inner sep=0.75pt]    {$1$};
\draw (237.11,29.59) node [anchor=north west][inner sep=0.75pt]    {$1$};
\draw (263.97,18.01) node [anchor=north west][inner sep=0.75pt]    {$1$};
\draw (278.32,29.59) node [anchor=north west][inner sep=0.75pt]    {$1$};
\draw (166.68,-1.63) node [anchor=north west][inner sep=0.75pt]    {$L_{e_{1}}{}$};
\draw (321,115.56) node [anchor=north west][inner sep=0.75pt]    {$r$};
\draw (385.17,124.24) node   [align=left] {1};
\draw (424.53,123.78) node   [align=left] {2};
\draw (343.03,79.79) node   [align=left] {1};
\draw (342.57,41.36) node   [align=left] {2};
\draw (366.67,99.97) node [anchor=north west][inner sep=0.75pt]    {$1$};
\draw (404.17,99.97) node [anchor=north west][inner sep=0.75pt]    {$1$};
\draw (397.69,73.12) node [anchor=north west][inner sep=0.75pt]    {$1$};
\draw (396.76,29.59) node [anchor=north west][inner sep=0.75pt]    {$1$};
\draw (326.33,-2.33) node [anchor=north west][inner sep=0.75pt]    {$L_{e_{1}}{}$};
\draw (480.65,115.56) node [anchor=north west][inner sep=0.75pt]    {$r$};
\draw (542.95,124.23) node   [align=left] {1};
\draw (582.31,123.77) node   [align=left] {2};
\draw (500.82,79.78) node   [align=left] {1};
\draw (500.35,41.34) node   [align=left] {2};
\draw (555.48,73.1) node [anchor=north west][inner sep=0.75pt]    {$1$};
\draw (595.76,29.58) node [anchor=north west][inner sep=0.75pt]    {$1$};
\draw (484.12,-1.64) node [anchor=north west][inner sep=0.75pt]    {$L_{e_{1}}{}$};
\draw (638.44,115.55) node [anchor=north west][inner sep=0.75pt]    {$r$};
\draw (524.22,133.92) node [anchor=north west][inner sep=0.75pt]   [align=left] {Multiplicity=$2$};
\draw (29.22,133.92) node [anchor=north west][inner sep=0.75pt]   [align=left] {Multiplicity=$1$};
\draw (362.22,133.92) node [anchor=north west][inner sep=0.75pt]   [align=left] {Multiplicity=$1$};
\draw (203.22,133.92) node [anchor=north west][inner sep=0.75pt]   [align=left] {Multiplicity=$1$};

\end{tikzpicture}
}
\caption{Indecomposable direct summands of $PH_0(Y^2_{-,-};\mathbb{F})$. The number on the colored block indicates the dimension of the vector space in the representation.}
\label{H0summand-0}
\end{figure}

Next, we calculated the rank of $H_i((\hat{\mathsf{Star}}_{k})^2_r)$ for $k\geq 4, i=0,1$, where $\hat{\mathsf{Star}}_{k}$ is a metric star graph such that every edge but $e_1$ has length $1$.
\begin{theorem}\label{Thm1}
Let $k\geq 4$. Then 
\begin{equation}
\rank(H_0((\hat{\mathsf{Star}}_{k})^2_r))=\begin{cases} 1, & \mbox{if } r\leq 1 \\ 
k^2-3k+4, & \mbox{if } 1<r\leq 2 \mbox{ and } r\leq L_{e_1} \\
k^2-k, & \mbox{if } 1<r\leq 2 \mbox{ and } L_{e_1}<r\leq L_{e_1}+1 \\
k^2-3k+2, & \mbox{if } 1<r\leq 2 \mbox{ and } L_{e_1}+1<r \\
2, & \mbox{if } 2<r \mbox{ and } r\leq L_{e_1} \\
2k-2, & \mbox{if } 2<r \mbox{ and } L_{e_1}<r\leq L_{e_1}+1\\
0, & \mbox{else}
\end{cases}
\end{equation}
and
\begin{equation}
\rank(H_1((\hat{\mathsf{Star}}_{k})^2_r))=\begin{cases} k^2-3k+1, & \mbox{if } r\leq 1 \mbox{ and } r\leq L_{e_1}\\ 
k^2-5k+5, & \mbox{if } r\leq 1 \mbox{ and } L_{e_1}<r\leq L_{e_1}+1\\ 
0, & \mbox{else}
\end{cases}
\end{equation}
\end{theorem}

When $k\geq 4$, the rank of $H_0((\hat{\mathsf{Star}}_{k})^2_r)$ for all $r>0$ and $L_{e_1}>0$ is shown in Figure \ref{fig:h0star}. When $k\geq 4$, the rank of $H_1((\hat{\mathsf{Star}}_{k})^2_r)$ for all $r>0$ and $L_{e_1}>0$ is shown in Figure \ref{fig:h1star}.\\

\begin{minipage}{0.5\textwidth}
\centering
\resizebox{7cm}{!}{

\tikzset{every picture/.style={line width=0.75pt}} 

\begin{tikzpicture}[x=0.75pt,y=0.75pt,yscale=-1,xscale=1]

\draw  [draw opacity=0][fill={rgb, 255:red, 155; green, 155; blue, 155 }  ,fill opacity=0.6 ] (138.9,52.47) -- (237.3,51.32) -- (236.79,104.89) -- (139.57,203.47) -- (138.9,106.03) -- cycle ;
\draw  [draw opacity=0][fill={rgb, 255:red, 248; green, 231; blue, 28 }  ,fill opacity=0.8 ] (138.96,202.98) -- (235.84,105.18) -- (236.25,202.57) -- cycle ;
\draw  (6.75,299.04) -- (373.46,299.04)(43.42,46.75) -- (43.42,327.07) (366.46,294.04) -- (373.46,299.04) -- (366.46,304.04) (38.42,53.75) -- (43.42,46.75) -- (48.42,53.75)  ;
\draw  [draw opacity=0][fill={rgb, 255:red, 245; green, 166; blue, 35 }  ,fill opacity=0.6 ] (43.42,299.04) -- (139.57,203.47) -- (139.82,298.78) -- cycle ;
\draw  [draw opacity=0][fill={rgb, 255:red, 184; green, 233; blue, 134 }  ,fill opacity=0.6 ] (139.82,298.78) -- (235.28,202.82) -- (235.53,298.53) -- cycle ;
\draw  [draw opacity=0][fill={rgb, 255:red, 242; green, 145; blue, 197 }  ,fill opacity=0.6 ] (289.7,50.84) -- (236.79,104.89) -- (236.25,51.39) -- cycle ;
\draw  [draw opacity=0][fill={rgb, 255:red, 80; green, 227; blue, 194 }  ,fill opacity=0.59 ] (332.05,105.74) -- (235.28,203.82) -- (235.84,105.18) -- cycle ;
\draw  [draw opacity=0][fill={rgb, 255:red, 80; green, 227; blue, 194 }  ,fill opacity=0.59 ] (289.7,50.7) -- (380.32,47.75) -- (330.94,106.22) -- (331.68,106.33) -- (307.51,106.22) -- (234.5,105.89) -- cycle ;
\draw  [draw opacity=0][fill={rgb, 255:red, 245; green, 166; blue, 35 }  ,fill opacity=0.6 ] (138.9,52.47) -- (138.9,106.03) -- (139.57,203.47) -- (43.42,299.04) -- (43.93,52.47) -- cycle ;
\draw  [draw opacity=0][fill={rgb, 255:red, 248; green, 231; blue, 28 }  ,fill opacity=0.8 ] (235.87,201.79) -- (139.22,299.82) -- (138.58,202.43) -- cycle ;
\draw [color={rgb, 255:red, 0; green, 0; blue, 0 }  ,draw opacity=0.14 ] [dash pattern={on 4.5pt off 4.5pt}]  (37.32,106.05) -- (338.5,106.05) ;
\draw [color={rgb, 255:red, 0; green, 0; blue, 0 }  ,draw opacity=0.14 ] [dash pattern={on 4.5pt off 4.5pt}]  (37.5,202) -- (271.19,202) ;

\draw (139.47,313.34) node   [align=left] {1};
\draw (236.73,312.19) node   [align=left] {2};
\draw (30.35,202.5) node   [align=left] {1};
\draw (30.21,105.53) node   [align=left] {2};
\draw (86.99,185.84) node [anchor=north west][inner sep=0.75pt]    {$1$};
\draw (159.25,276.45) node [anchor=north west][inner sep=0.75pt]    {$k^{2} -3k+2$};
\draw (166.41,192.09) node [anchor=north west][inner sep=0.75pt]    {$k^{2} -k$};
\draw (150.12,86.54) node [anchor=north west][inner sep=0.75pt]    {$k^{2} -3k+4$};
\draw (243.31,62.03) node [anchor=north west][inner sep=0.75pt]    {$2$};
\draw (262.22,90.63) node [anchor=north west][inner sep=0.75pt]    {$2k-2$};
\draw (17.07,17.96) node [anchor=north west][inner sep=0.75pt]    {$L_{e_{1}}{}$};
\draw (376.76,302.52) node [anchor=north west][inner sep=0.75pt]    {$r_{12}$};

\end{tikzpicture}
}
\captionof{figure}{The rank of $H_0((\hat{\mathsf{Star}}_{k})^2_r)$}
\label{fig:h0star}  
\end{minipage}%
\begin{minipage}{0.5\textwidth}
\centering
\resizebox{7cm}{!}{

\tikzset{every picture/.style={line width=0.75pt}} 

\begin{tikzpicture}[x=0.75pt,y=0.75pt,yscale=-1,xscale=1]

\draw  (6.75,299.04) -- (373.46,299.04)(43.42,46.75) -- (43.42,327.07) (366.46,294.04) -- (373.46,299.04) -- (366.46,304.04) (38.42,53.75) -- (43.42,46.75) -- (48.42,53.75)  ;
\draw  [draw opacity=0][fill={rgb, 255:red, 126; green, 211; blue, 33 }  ,fill opacity=0.62 ] (43.42,299.04) -- (139.57,203.47) -- (139.82,298.78) -- cycle ;
\draw  [draw opacity=0][fill={rgb, 255:red, 248; green, 231; blue, 28 }  ,fill opacity=0.8 ] (138.9,52.47) -- (138.9,106.03) -- (139.57,203.47) -- (43.42,299.04) -- (43.93,52.47) -- cycle ;

\draw (139.47,313.34) node   [align=left] {1};
\draw (236.73,312.19) node   [align=left] {2};
\draw (35.35,203.5) node   [align=left] {1};
\draw (34.21,108.53) node   [align=left] {2};
\draw (17.07,17.96) node [anchor=north west][inner sep=0.75pt]    {$L_{e_{1}}{}$};
\draw (376.76,302.52) node [anchor=north west][inner sep=0.75pt]    {$r_{12}$};
\draw (51,146.4) node [anchor=north west][inner sep=0.75pt]    {$k( k-3) +1$};
\draw (57,279.4) node [anchor=north west][inner sep=0.75pt]    {$k( k-5) +5$};

\end{tikzpicture}
}
\captionof{figure}{The rank of $H_1((\hat{\mathsf{Star}}_{k})^2_r)$}
\label{fig:h1star}  
\end{minipage}

\vspace{0.8cm}
Next, we construct a persistence module $PH_i((\mathsf{Star_k})^2_{-,-})$ for $i=0,1$. We prove the following theorems regarding the interval decomposability of $PH_0((\mathsf{Star_k})^2_{-,-};\mathbb{F})$ and $PH_1((\mathsf{Star_k})^2_{-,-};\mathbb{F})$, where $(\mathsf{Star}_{k})^2_{r,L_{e_1}}:=(\hat{\mathsf{Star}}_{k})^2_r$.

\begin{theorem}\label{Thm2}
$k\geq 4$, $PH_1((\mathsf{Star_k})^2_{-,-};\mathbb{F})$ is interval decomposable.
\end{theorem}

\begin{theorem}\label{Thm3}
$k\geq 4$, $PH_0((\mathsf{Star_k})^2_{-,-};\mathbb{F})$ is not interval decomposable, and the indecomposable direct summands of $PH_0((\mathsf{Star_k})^2_{-,-};\mathbb{F})$ are given in Figure \ref{H0summand-1}.
\end{theorem}
\begin{figure}[htbp!]
\centering
\resizebox{18cm}{!}{

\tikzset{every picture/.style={line width=0.75pt}} 

\begin{tikzpicture}[x=0.75pt,y=0.75pt,yscale=-1,xscale=1]

\draw  [draw opacity=0][fill={rgb, 255:red, 155; green, 155; blue, 155 }  ,fill opacity=0.6 ] (58.11,18.42) -- (96.85,17.97) -- (96.65,39.05) -- (58.38,77.86) -- (58.11,39.5) -- cycle ;
\draw  [draw opacity=0][fill={rgb, 255:red, 248; green, 231; blue, 28 }  ,fill opacity=0.8 ] (58.06,78) -- (96.2,39.5) -- (96.36,77.84) -- cycle ;
\draw  (6.1,115.48) -- (150.44,115.48)(20.53,16.17) -- (20.53,126.51) (143.44,110.48) -- (150.44,115.48) -- (143.44,120.48) (15.53,23.17) -- (20.53,16.17) -- (25.53,23.17)  ;
\draw  [draw opacity=0][fill={rgb, 255:red, 245; green, 166; blue, 35 }  ,fill opacity=0.6 ] (20.53,115.48) -- (58.38,77.86) -- (58.48,115.38) -- cycle ;
\draw  [draw opacity=0][fill={rgb, 255:red, 184; green, 233; blue, 134 }  ,fill opacity=0.6 ] (58.35,115.68) -- (95.93,77.9) -- (96.03,115.58) -- cycle ;
\draw  [draw opacity=0][fill={rgb, 255:red, 242; green, 145; blue, 197 }  ,fill opacity=0.6 ] (117.47,17.78) -- (96.65,39.05) -- (96.43,17.99) -- cycle ;
\draw  [draw opacity=0][fill={rgb, 255:red, 80; green, 227; blue, 194 }  ,fill opacity=0.6 ] (134,39.23) -- (95.91,77.84) -- (96.13,39.01) -- cycle ;
\draw  [draw opacity=0][fill={rgb, 255:red, 80; green, 227; blue, 194 }  ,fill opacity=0.6 ] (117.47,17.33) -- (153.15,16.17) -- (133.71,39.18) -- (134,39.23) -- (124.48,39.18) -- (95.75,39.05) -- cycle ;
\draw  [draw opacity=0][fill={rgb, 255:red, 245; green, 166; blue, 35 }  ,fill opacity=0.6 ] (58.11,18.42) -- (58.11,39.5) -- (58.38,77.86) -- (20.53,115.48) -- (20.73,18.42) -- cycle ;
\draw  [draw opacity=0][fill={rgb, 255:red, 248; green, 231; blue, 28 }  ,fill opacity=0.8 ] (96.68,77.15) -- (58.63,115.74) -- (58.38,77.41) -- cycle ;
\draw  [draw opacity=0][fill={rgb, 255:red, 155; green, 155; blue, 155 }  ,fill opacity=0.6 ] (219.32,18.4) -- (258.05,17.95) -- (257.85,39.03) -- (219.58,77.84) -- (219.32,39.48) -- cycle ;
\draw  [draw opacity=0][fill={rgb, 255:red, 248; green, 231; blue, 28 }  ,fill opacity=0.8 ] (219.26,77.98) -- (257.4,39.48) -- (257.56,77.82) -- cycle ;
\draw  (167.3,115.46) -- (311.65,115.46)(181.73,16.15) -- (181.73,126.49) (304.65,110.46) -- (311.65,115.46) -- (304.65,120.46) (176.73,23.15) -- (181.73,16.15) -- (186.73,23.15)  ;
\draw  [draw opacity=0][fill={rgb, 255:red, 242; green, 145; blue, 197 }  ,fill opacity=0.6 ] (278.67,17.76) -- (257.85,39.03) -- (257.64,17.97) -- cycle ;
\draw  [draw opacity=0][fill={rgb, 255:red, 80; green, 227; blue, 194 }  ,fill opacity=0.6 ] (295.2,39.21) -- (257.11,77.82) -- (257.33,38.99) -- cycle ;
\draw  [draw opacity=0][fill={rgb, 255:red, 80; green, 227; blue, 194 }  ,fill opacity=0.6 ] (278.67,17.31) -- (314.35,16.15) -- (294.91,39.16) -- (295.2,39.21) -- (256.95,39.03) -- cycle ;
\draw  [draw opacity=0][fill={rgb, 255:red, 248; green, 231; blue, 28 }  ,fill opacity=0.8 ] (257.88,77.13) -- (219.83,115.72) -- (219.58,77.39) -- cycle ;
\draw  [draw opacity=0][fill={rgb, 255:red, 155; green, 155; blue, 155 }  ,fill opacity=0.6 ] (539.34,19.4) -- (578.07,18.95) -- (577.87,40.03) -- (539.6,78.84) -- (539.34,40.48) -- cycle ;
\draw  [draw opacity=0][fill={rgb, 255:red, 248; green, 231; blue, 28 }  ,fill opacity=0.8 ] (539.28,78.98) -- (577.42,40.48) -- (577.58,78.82) -- cycle ;
\draw  (487.32,116.46) -- (631.67,116.46)(501.76,17.15) -- (501.76,127.49) (624.67,111.46) -- (631.67,116.46) -- (624.67,121.46) (496.76,24.15) -- (501.76,17.15) -- (506.76,24.15)  ;
\draw  [draw opacity=0][fill={rgb, 255:red, 184; green, 233; blue, 134 }  ,fill opacity=0.6 ] (539.87,116.23) -- (577.45,78.45) -- (577.55,116.13) -- cycle ;
\draw  [draw opacity=0][fill={rgb, 255:red, 248; green, 231; blue, 28 }  ,fill opacity=0.8 ] (577.9,78.14) -- (539.86,116.73) -- (539.6,78.39) -- cycle ;
\draw  [draw opacity=0][fill={rgb, 255:red, 248; green, 231; blue, 28 }  ,fill opacity=0.8 ] (377.76,78.23) -- (415.9,39.73) -- (416.06,78.07) -- cycle ;
\draw  (326.25,115.71) -- (470.6,115.71)(340.69,16.4) -- (340.69,126.74) (463.6,110.71) -- (470.6,115.71) -- (463.6,120.71) (335.69,23.4) -- (340.69,16.4) -- (345.69,23.4)  ;
\draw  [draw opacity=0][fill={rgb, 255:red, 80; green, 227; blue, 194 }  ,fill opacity=0.6 ] (454.15,39.46) -- (416.06,78.07) -- (416.29,39.24) -- cycle ;
\draw  [draw opacity=0][fill={rgb, 255:red, 80; green, 227; blue, 194 }  ,fill opacity=0.6 ] (437.63,17.56) -- (473.3,16.4) -- (453.86,39.42) -- (454.15,39.46) -- (415.9,39.28) -- cycle ;
\draw  [draw opacity=0][fill={rgb, 255:red, 248; green, 231; blue, 28 }  ,fill opacity=0.8 ] (416.83,77.39) -- (378.79,115.98) -- (378.53,77.64) -- cycle ;
\draw [color={rgb, 255:red, 0; green, 0; blue, 0 }  ,draw opacity=0.14 ] [dash pattern={on 4.5pt off 4.5pt}]  (17.04,39.18) -- (134,39.18) ;
\draw [color={rgb, 255:red, 0; green, 0; blue, 0 }  ,draw opacity=0.14 ] [dash pattern={on 4.5pt off 4.5pt}]  (17.04,77.37) -- (134,77.37) ;
\draw [color={rgb, 255:red, 0; green, 0; blue, 0 }  ,draw opacity=0.14 ] [dash pattern={on 4.5pt off 4.5pt}]  (178.04,38.43) -- (295,38.43) ;
\draw [color={rgb, 255:red, 0; green, 0; blue, 0 }  ,draw opacity=0.14 ] [dash pattern={on 4.5pt off 4.5pt}]  (178.04,76.63) -- (295,76.63) ;
\draw [color={rgb, 255:red, 0; green, 0; blue, 0 }  ,draw opacity=0.14 ] [dash pattern={on 4.5pt off 4.5pt}]  (336.78,40.04) -- (453.75,40.04) ;
\draw [color={rgb, 255:red, 0; green, 0; blue, 0 }  ,draw opacity=0.14 ] [dash pattern={on 4.5pt off 4.5pt}]  (336.78,78.23) -- (453.75,78.23) ;
\draw [color={rgb, 255:red, 0; green, 0; blue, 0 }  ,draw opacity=0.14 ] [dash pattern={on 4.5pt off 4.5pt}]  (498.53,40.79) -- (615.5,40.79) ;
\draw [color={rgb, 255:red, 0; green, 0; blue, 0 }  ,draw opacity=0.14 ] [dash pattern={on 4.5pt off 4.5pt}]  (498.53,78.98) -- (615.5,78.98) ;
\draw [color={rgb, 255:red, 0; green, 0; blue, 0 }  ,draw opacity=0.14 ] [dash pattern={on 4.5pt off 4.5pt}]  (257.41,11.33) -- (257.41,119.16) ;
\draw [color={rgb, 255:red, 0; green, 0; blue, 0 }  ,draw opacity=0.14 ] [dash pattern={on 4.5pt off 4.5pt}]  (416.16,14.2) -- (416.16,122.04) ;

\draw (58.34,121.11) node   [align=left] {1};
\draw (96.62,120.65) node   [align=left] {2};
\draw (12.11,77.87) node   [align=left] {1};
\draw (12.41,41.24) node   [align=left] {2};
\draw (35.07,55.1) node [anchor=north west][inner sep=0.75pt]    {$1$};
\draw (76.66,97.29) node [anchor=north west][inner sep=0.75pt]    {$1$};
\draw (70.36,71.17) node [anchor=north west][inner sep=0.75pt]    {$2$};
\draw (69.45,28.83) node [anchor=north west][inner sep=0.75pt]    {$2$};
\draw (95.58,17.57) node [anchor=north west][inner sep=0.75pt]    {$1$};
\draw (109.54,28.83) node [anchor=north west][inner sep=0.75pt]    {$1$};
\draw (0.71,-1.57) node [anchor=north west][inner sep=0.75pt]    {$L_{e_{1}}{}$};
\draw (151.07,112.95) node [anchor=north west][inner sep=0.75pt]    {$r$};
\draw (219.54,121.09) node   [align=left] {1};
\draw (257.82,120.63) node   [align=left] {2};
\draw (173.32,77.1) node   [align=left] {1};
\draw (173.61,40.47) node   [align=left] {2};
\draw (231.56,71.15) node [anchor=north west][inner sep=0.75pt]    {$1$};
\draw (230.66,28.81) node [anchor=north west][inner sep=0.75pt]    {$1$};
\draw (256.78,17.55) node [anchor=north west][inner sep=0.75pt]    {$1$};
\draw (270.74,28.81) node [anchor=north west][inner sep=0.75pt]    {$1$};
\draw (161.91,-1.59) node [anchor=north west][inner sep=0.75pt]    {$L_{e_{1}}{}$};
\draw (312.28,112.95) node [anchor=north west][inner sep=0.75pt]    {$r$};
\draw (539.57,122.09) node   [align=left] {1};
\draw (577.85,121.64) node   [align=left] {2};
\draw (492.59,79.6) node   [align=left] {1};
\draw (492.89,40.72) node   [align=left] {2};
\draw (557.89,98.27) node [anchor=north west][inner sep=0.75pt]    {$1$};
\draw (554.58,67.66) node [anchor=north west][inner sep=0.75pt]    {$1$};
\draw (550.68,29.81) node [anchor=north west][inner sep=0.75pt]    {$1$};
\draw (481.93,-1.27) node [anchor=north west][inner sep=0.75pt]    {$L_{e_{1}}{}$};
\draw (632.3,112.95) node [anchor=north west][inner sep=0.75pt]    {$r$};
\draw (378.5,121.34) node   [align=left] {1};
\draw (416.78,120.89) node   [align=left] {2};
\draw (331.52,78.1) node   [align=left] {1};
\draw (331.82,40.72) node   [align=left] {2};
\draw (390.51,71.4) node [anchor=north west][inner sep=0.75pt]    {$1$};
\draw (429.69,29.07) node [anchor=north west][inner sep=0.75pt]    {$1$};
\draw (320.86,-1.34) node [anchor=north west][inner sep=0.75pt]    {$L_{e_{1}}{}$};
\draw (471.23,112.95) node [anchor=north west][inner sep=0.75pt]    {$r$};
\draw (336.07,128.89) node [anchor=north west][inner sep=0.75pt]   [align=left] {Multiplicity=$2k-4$};
\draw (495.97,128.89) node [anchor=north west][inner sep=0.75pt]   [align=left] {Multiplicity=$k^{2} -3k+1$};
\draw (19.71,128.89) node [anchor=north west][inner sep=0.75pt]   [align=left] {Multiplicity=$1$};
\draw (177.71,128.89) node [anchor=north west][inner sep=0.75pt]   [align=left] {Multiplicity=$1$};

\end{tikzpicture}
}
\caption{Indecomposable direct summands of $PH_0((\mathsf{Star_k})^2_{-,-};\mathbb{F})$.  The number on the colored block indicates the dimension of the vector space in the representation.}
\label{H0summand-1}
\end{figure}

Based on our calculation of the rank of $H_i((\hat{\mathsf{Star}}_{k})^2_r)$, where $k\geq 3$ and  $i=0,1$, we calculated the rank of $H_i(\hat{\mathcal{H}}^2_{r,L_{e_1}})$ for $i=0,1$, where $\hat{\mathcal{H}}$ is the metric graph $H$ such that every edge but $e_1$ has length $1$. 

\begin{theorem}\label{Thm6}
\begin{equation}
\rank(H_0(\hat{\mathcal{H}}^2_{r,L_{e_1}}))=\begin{cases} 1, & \mbox{if } r\leq 1 \\ 
6, & \mbox{if } 1<r\leq 2 \mbox{ and } r\leq L_{e_1}+1 \\
12, & \mbox{if } 1<r\leq 2 \mbox{ and } L_{e_1}+1<r \\
2, & \mbox{if } 2<r \mbox{ and } r\leq L_{e_1}+1 \\
8, & \mbox{if } 2<r \mbox{ and } L_{e_1}+1< r\leq L_{e_1}+2 \\
0, & \mbox{else}
\end{cases}
\end{equation}
and
\begin{equation}
\rank(H_1(\hat{\mathcal{H}}^2_{r,L_{e_1}}))=\begin{cases} 3, & \mbox{if } r\leq 1 \mbox{ and } r\leq L_{e_1}\\ 
5, & \mbox{if } r\leq 1 \mbox{ and } L_{e_1}<r\\ 
2, & \mbox{if } 1<r \mbox{ and } L_{e_1}<r\leq L_{e_1}+1 \\
0, & \mbox{else}
\end{cases}
\end{equation}
\end{theorem}

The rank of $H_0(\hat{\mathcal{H}}^2_{r,L_{e_1}})$ for all $r>0$ and $L_{e_1}>0$ is shown in Figure \ref{fig:h0H}. The rank of $H_1(\hat{\mathcal{H}}^2_{r,L_{e_1}})$ for all $r>0$ and $L_{e_1}>0$ is shown in Figure \ref{fig:h1H}.\\

\begin{minipage}{0.5\textwidth}
\centering
\resizebox{6.5cm}{!}{

\tikzset{every picture/.style={line width=0.75pt}} 

\begin{tikzpicture}[x=0.75pt,y=0.75pt,yscale=-1,xscale=1]

\draw  [draw opacity=0][fill={rgb, 255:red, 74; green, 144; blue, 226 }  ,fill opacity=0.6 ] (40.85,313.04) -- (141.57,217.47) -- (141.82,312.77) -- cycle ;
\draw  [draw opacity=0][fill={rgb, 255:red, 184; green, 233; blue, 134 }  ,fill opacity=0.9 ] (142.21,312.78) -- (237.28,217.76) -- (237.53,312.53) -- cycle ;
\draw  [draw opacity=0][fill={rgb, 255:red, 248; green, 231; blue, 28 }  ,fill opacity=0.6 ] (141.5,35) -- (236.5,35.5) -- (238.79,118.89) -- (141.57,217.47) -- (140.9,120.03) -- cycle ;
\draw  [draw opacity=0][fill={rgb, 255:red, 248; green, 231; blue, 28 }  ,fill opacity=0.6 ] (140.76,217.83) -- (237.65,120.03) -- (238.05,217.42) -- cycle ;
\draw  (5.44,312.73) -- (350.5,312.73)(39.95,12) -- (39.95,346.15) (343.5,307.73) -- (350.5,312.73) -- (343.5,317.73) (34.95,19) -- (39.95,12) -- (44.95,19)  ;
\draw  [draw opacity=0][fill={rgb, 255:red, 208; green, 2; blue, 27 }  ,fill opacity=0.5 ] (333.68,120.33) -- (236.91,218.42) -- (237.48,119.78) -- cycle ;
\draw  [draw opacity=0][fill={rgb, 255:red, 208; green, 2; blue, 27 }  ,fill opacity=0.5 ] (344.5,35.5) -- (344.5,109) -- (333.68,120.33) -- (237.65,120.03) -- (236.5,35.5) -- cycle ;
\draw  [draw opacity=0][fill={rgb, 255:red, 74; green, 144; blue, 226 }  ,fill opacity=0.6 ] (141.5,35) -- (140.9,120.03) -- (141.57,217.47) -- (40.85,313.04) -- (40.5,35) -- cycle ;
\draw  [draw opacity=0][fill={rgb, 255:red, 248; green, 231; blue, 28 }  ,fill opacity=0.6 ] (238.86,216.83) -- (142.2,312.37) -- (141.57,217.47) -- cycle ;
\draw  [draw opacity=0][fill={rgb, 255:red, 155; green, 155; blue, 155 }  ,fill opacity=0.6 ] (344.5,109) -- (344.5,204) -- (237.9,311.87) -- (237.28,217.76) -- cycle ;

\draw (141.47,327.34) node   [align=left] {1};
\draw (238.73,326.19) node   [align=left] {2};
\draw (35.35,217.5) node   [align=left] {1};
\draw (34.21,118.53) node   [align=left] {2};
\draw (13.07,17.96) node [anchor=north west][inner sep=0.75pt]    {$L_{e_{1}}{}$};
\draw (355.76,316.52) node [anchor=north west][inner sep=0.75pt]    {$r$};
\draw (89.35,172.5) node   [align=left] {1};
\draw (187.35,152.5) node   [align=left] {6};
\draw (204.35,286.5) node   [align=left] {12};
\draw (288.35,89.5) node   [align=left] {2};
\draw (288.35,210.5) node   [align=left] {8};

\end{tikzpicture}
}
\captionof{figure}{The rank of $H_0(\hat{\mathcal{H}}^2_{r,L_{e_1}})$}
\label{fig:h0H}  
\end{minipage}%
\begin{minipage}{0.5\textwidth}
\centering
\resizebox{6.5cm}{!}{

\tikzset{every picture/.style={line width=0.75pt}} 

\begin{tikzpicture}[x=0.75pt,y=0.75pt,yscale=-1,xscale=1]

\draw  [draw opacity=0][fill={rgb, 255:red, 139; green, 87; blue, 42 }  ,fill opacity=0.6 ] (40.85,313.05) -- (141.96,217.47) -- (142.21,312.78) -- cycle ;
\draw  [draw opacity=0] (142.21,312.78) -- (237.28,217.76) -- (237.53,312.53) -- cycle ;
\draw  (6.34,313.04) -- (351.4,313.04)(40.85,12.31) -- (40.85,346.45) (344.4,308.04) -- (351.4,313.04) -- (344.4,318.04) (35.85,19.31) -- (40.85,12.31) -- (45.85,19.31)  ;
\draw  [draw opacity=0][fill={rgb, 255:red, 80; green, 227; blue, 194 }  ,fill opacity=0.6 ] (141.5,35) -- (140.9,120.03) -- (141.57,217.47) -- (40.85,313.04) -- (40.5,35) -- cycle ;
\draw  [draw opacity=0][fill={rgb, 255:red, 208; green, 2; blue, 27 }  ,fill opacity=0.6 ] (331.51,26.59) -- (333.7,120.75) -- (142.21,312.78) -- (140.93,217.65) -- cycle ;

\draw (141.47,327.34) node   [align=left] {1};
\draw (238.73,326.19) node   [align=left] {2};
\draw (35.35,217.5) node   [align=left] {1};
\draw (34.21,118.53) node   [align=left] {2};
\draw (13.07,17.96) node [anchor=north west][inner sep=0.75pt]    {$L_{e_{1}}{}$};
\draw (355.76,316.52) node [anchor=north west][inner sep=0.75pt]    {$r$};
\draw (90.21,166.53) node   [align=left] {3};
\draw (111.21,286.53) node   [align=left] {5};
\draw (217.21,186.53) node   [align=left] {2};

\end{tikzpicture}
}
\captionof{figure}{The rank of $H_1(\hat{\mathcal{H}}^2_{r,L_{e_1}})$}
\label{fig:h1H}  
\end{minipage}

Next, we construct a persistence module $PH_i(\hat{\mathcal{H}}^2_{-,-})$ for $i=0,1$. We prove the following theorems regarding the interval decomposability of $PH_0(\hat{\mathcal{H}}^2_{-,-};\mathbb{F})$ and $PH_1(\hat{\mathcal{H}}^2_{-,-};\mathbb{F})$.

\begin{theorem}\label{Thm7}
$PH_1(\hat{\mathcal{H}}^2_{-,-})$ is an interval module, hence it is indecomposable.
\end{theorem}

\begin{theorem}\label{Thm8}
$PH_0(\hat{\mathcal{H}}^2_{-,-})$ is interval decomposable, and the indecomposable direct summands of $PH_0(\hat{\mathcal{H}}^2_{-,-};\mathbb{F})$ are given in Figure \ref{H0summand-H}.
\end{theorem}

\begin{figure}[htbp!]
\centering
\resizebox{18cm}{!}{

\tikzset{every picture/.style={line width=0.75pt}} 

\begin{tikzpicture}[x=0.75pt,y=0.75pt,yscale=-1,xscale=1]

\draw  [draw opacity=0][fill={rgb, 255:red, 74; green, 144; blue, 226 }  ,fill opacity=0.6 ] (24.66,179.15) -- (62.13,143.6) -- (62.22,179.05) -- cycle ;
\draw  [draw opacity=0][fill={rgb, 255:red, 184; green, 233; blue, 134 }  ,fill opacity=0.9 ] (62.22,178.7) -- (97.58,143.36) -- (97.67,178.61) -- cycle ;
\draw  [draw opacity=0][fill={rgb, 255:red, 248; green, 231; blue, 28 }  ,fill opacity=0.6 ] (62.22,16.61) -- (97.67,15.87) -- (98.14,106.58) -- (61.98,143.25) -- (61.73,107.01) -- cycle ;
\draw  [draw opacity=0][fill={rgb, 255:red, 248; green, 231; blue, 28 }  ,fill opacity=0.6 ] (61.68,143.38) -- (97.72,107.01) -- (97.87,143.23) -- cycle ;
\draw  (9.49,178.79) -- (159.78,178.79)(24.52,14.4) -- (24.52,197.06) (152.78,173.79) -- (159.78,178.79) -- (152.78,183.79) (19.52,21.4) -- (24.52,14.4) -- (29.52,21.4)  ;
\draw  [draw opacity=0][fill={rgb, 255:red, 208; green, 2; blue, 27 }  ,fill opacity=0.5 ] (142.82,15.52) -- (97.65,106.91) -- (96.73,15.99) -- cycle ;
\draw  [draw opacity=0][fill={rgb, 255:red, 208; green, 2; blue, 27 }  ,fill opacity=0.5 ] (133.44,107.12) -- (97.44,143.6) -- (97.65,106.91) -- cycle ;
\draw  [draw opacity=0][fill={rgb, 255:red, 208; green, 2; blue, 27 }  ,fill opacity=0.5 ] (152.34,52.67) -- (151.59,88.77) -- (133.16,107.08) -- (133.44,107.12) -- (126.49,106.99) -- (97.29,106.95) -- cycle ;
\draw  [draw opacity=0][fill={rgb, 255:red, 74; green, 144; blue, 226 }  ,fill opacity=0.6 ] (62.22,16.61) -- (61.73,107.01) -- (61.98,143.25) -- (24.52,178.79) -- (24.76,16.98) -- cycle ;
\draw  [draw opacity=0][fill={rgb, 255:red, 248; green, 231; blue, 28 }  ,fill opacity=0.6 ] (98.17,143.01) -- (62.21,178.55) -- (61.98,143.25) -- cycle ;
\draw [line width=1.5]    (24.52,178.79) -- (152.34,52.67) ;
\draw [line width=1.5]    (62.22,179.05) -- (151.59,88.77) ;
\draw [line width=1.5]    (62.22,179.05) -- (62.22,16.61) ;
\draw [line width=1.5]    (97.67,178.61) -- (97.67,15.87) ;
\draw [line width=1.5]    (62.22,179.05) -- (143.04,16.24) ;
\draw [line width=1.5]    (24.52,178.79) -- (107.04,15.99) ;
\draw [line width=1.5]    (97.67,178.61) -- (151.59,124.85) ;
\draw  [draw opacity=0][fill={rgb, 255:red, 208; green, 2; blue, 27 }  ,fill opacity=0.5 ] (152.34,52.67) -- (151.59,15.87) -- (143.17,15.8) -- (98.28,106.33) -- cycle ;
\draw  [draw opacity=0][fill={rgb, 255:red, 155; green, 155; blue, 155 }  ,fill opacity=0.6 ] (151.73,88.17) -- (151.73,124.25) -- (97.81,178.01) -- (97.58,143.01) -- cycle ;
\draw  [draw opacity=0][fill={rgb, 255:red, 184; green, 233; blue, 134 }  ,fill opacity=0.9 ] (228.89,178.84) -- (264.25,143.5) -- (264.34,178.75) -- cycle ;
\draw  [draw opacity=0][fill={rgb, 255:red, 248; green, 231; blue, 28 }  ,fill opacity=0.6 ] (228.89,16.75) -- (264.34,16.01) -- (264.81,106.72) -- (228.65,143.39) -- (228.4,107.15) -- cycle ;
\draw  [draw opacity=0][fill={rgb, 255:red, 248; green, 231; blue, 28 }  ,fill opacity=0.6 ] (228.35,143.52) -- (264.38,107.15) -- (264.54,143.37) -- cycle ;
\draw  (176.16,178.93) -- (326.45,178.93)(191.19,14.54) -- (191.19,197.2) (319.45,173.93) -- (326.45,178.93) -- (319.45,183.93) (186.19,21.54) -- (191.19,14.54) -- (196.19,21.54)  ;
\draw  [draw opacity=0][fill={rgb, 255:red, 208; green, 2; blue, 27 }  ,fill opacity=0.5 ] (309.49,15.66) -- (264.32,107.05) -- (263.4,16.13) -- cycle ;
\draw  [draw opacity=0][fill={rgb, 255:red, 208; green, 2; blue, 27 }  ,fill opacity=0.5 ] (300.1,107.26) -- (264.11,143.74) -- (264.32,107.05) -- cycle ;
\draw  [draw opacity=0][fill={rgb, 255:red, 208; green, 2; blue, 27 }  ,fill opacity=0.5 ] (319.01,52.81) -- (318.26,88.91) -- (299.83,107.22) -- (300.1,107.26) -- (293.16,107.13) -- (263.96,107.09) -- cycle ;
\draw  [draw opacity=0][fill={rgb, 255:red, 248; green, 231; blue, 28 }  ,fill opacity=0.6 ] (264.84,143.15) -- (228.88,178.69) -- (228.65,143.39) -- cycle ;
\draw [line width=1.5]    (228.35,143.52) -- (319.01,52.81) ;
\draw [line width=1.5]    (228.89,179.19) -- (318.26,88.91) ;
\draw [line width=1.5]    (228.89,179.19) -- (228.89,16.75) ;
\draw [line width=1.5]    (264.34,178.75) -- (264.34,16.01) ;
\draw [line width=1.5]    (228.89,179.19) -- (309.71,16.38) ;
\draw [line width=1.5]    (228.4,107.15) -- (273.71,16.13) ;
\draw [line width=1.5]    (264.34,178.75) -- (318.26,124.99) ;
\draw  [draw opacity=0][fill={rgb, 255:red, 208; green, 2; blue, 27 }  ,fill opacity=0.5 ] (319.01,52.81) -- (318.26,16.01) -- (309.84,15.94) -- (264.95,106.47) -- cycle ;
\draw  [draw opacity=0][fill={rgb, 255:red, 155; green, 155; blue, 155 }  ,fill opacity=0.6 ] (318.4,88.32) -- (318.4,124.39) -- (264.48,178.15) -- (264.25,143.15) -- cycle ;
\draw  [draw opacity=0][fill={rgb, 255:red, 184; green, 233; blue, 134 }  ,fill opacity=0.9 ] (393.01,178.53) -- (428.37,143.18) -- (428.47,178.43) -- cycle ;
\draw  [draw opacity=0][fill={rgb, 255:red, 248; green, 231; blue, 28 }  ,fill opacity=0.6 ] (393.01,16.44) -- (428.47,15.69) -- (428.94,106.41) -- (392.77,143.08) -- (392.53,106.83) -- cycle ;
\draw  [draw opacity=0][fill={rgb, 255:red, 248; green, 231; blue, 28 }  ,fill opacity=0.6 ] (392.47,143.21) -- (428.51,106.83) -- (428.66,143.06) -- cycle ;
\draw  (340.28,178.62) -- (490.57,178.62)(355.31,14.23) -- (355.31,196.89) (483.57,173.62) -- (490.57,178.62) -- (483.57,183.62) (350.31,21.23) -- (355.31,14.23) -- (360.31,21.23)  ;
\draw  [draw opacity=0][fill={rgb, 255:red, 248; green, 231; blue, 28 }  ,fill opacity=0.6 ] (428.96,142.84) -- (393.01,178.37) -- (392.77,143.08) -- cycle ;
\draw [line width=1.5]    (392.77,143.08) -- (428.94,106.41) ;
\draw [line width=1.5]    (393.01,178.87) -- (428.96,142.84) ;
\draw [line width=1.5]    (393.01,178.87) -- (393.01,16.44) ;
\draw [line width=1.5]    (428.47,178.43) -- (428.47,15.69) ;
\draw [line width=1.5]    (393.01,178.87) -- (428.51,106.83) ;
\draw [line width=1.5]    (392.53,106.83) -- (429.03,33.35) ;
\draw  [draw opacity=0][fill={rgb, 255:red, 184; green, 233; blue, 134 }  ,fill opacity=0.9 ] (559.89,178.11) -- (595.25,142.76) -- (595.34,178.02) -- cycle ;
\draw  (507.16,178.2) -- (657.45,178.2)(522.19,13.81) -- (522.19,196.47) (650.45,173.2) -- (657.45,178.2) -- (650.45,183.2) (517.19,20.81) -- (522.19,13.81) -- (527.19,20.81)  ;
\draw [line width=1.5]    (559.89,178.46) -- (649.26,88.18) ;
\draw [line width=1.5]    (595.34,178.02) -- (595.34,142.42) ;
\draw [line width=1.5]    (595.34,178.02) -- (649.26,124.26) ;
\draw  [draw opacity=0][fill={rgb, 255:red, 155; green, 155; blue, 155 }  ,fill opacity=0.6 ] (649.4,87.58) -- (649.4,123.66) -- (595.48,177.42) -- (595.25,142.42) -- cycle ;

\draw (363.17,201.74) node [anchor=north west][inner sep=0.75pt]   [align=left] {Multiplicity=4};
\draw (539.89,201.74) node [anchor=north west][inner sep=0.75pt]   [align=left] {Multiplicity=6};
\draw (33.84,202.47) node [anchor=north west][inner sep=0.75pt]   [align=left] {Multiplicity=1};
\draw (200.71,202.47) node [anchor=north west][inner sep=0.75pt]   [align=left] {Multiplicity=1};
\draw (3.56,12.33) node [anchor=north west][inner sep=0.75pt]    {$L_{e_{1}}{}$};
\draw (163.29,174.46) node [anchor=north west][inner sep=0.75pt]    {$r$};
\draw (37.54,74.47) node [anchor=north west][inner sep=0.75pt]    {$1$};
\draw (47.59,130.27) node [anchor=north west][inner sep=0.75pt]    {$1$};
\draw (42.38,161.14) node [anchor=north west][inner sep=0.75pt]    {$1$};
\draw (65.07,137.71) node [anchor=north west][inner sep=0.75pt]    {$1$};
\draw (83.66,133.24) node [anchor=north west][inner sep=0.75pt]    {$1$};
\draw (80.87,162.36) node [anchor=north west][inner sep=0.75pt]    {$1$};
\draw (118.26,127.66) node [anchor=north west][inner sep=0.75pt]    {$1$};
\draw (109.7,30.96) node [anchor=north west][inner sep=0.75pt]    {$1$};
\draw (131.65,44.72) node [anchor=north west][inner sep=0.75pt]    {$1$};
\draw (117.51,91.58) node [anchor=north west][inner sep=0.75pt]    {$1$};
\draw (71.75,33.59) node [anchor=north west][inner sep=0.75pt]    {$1$};
\draw (74.26,91.16) node [anchor=north west][inner sep=0.75pt]    {$1$};
\draw (97.18,13.43) node [anchor=north west][inner sep=0.75pt]  [font=\tiny]  {$1$};
\draw (170.22,12.47) node [anchor=north west][inner sep=0.75pt]    {$L_{e_{1}}{}$};
\draw (329.96,174.6) node [anchor=north west][inner sep=0.75pt]    {$r$};
\draw (231.74,137.85) node [anchor=north west][inner sep=0.75pt]    {$1$};
\draw (250.33,133.38) node [anchor=north west][inner sep=0.75pt]    {$1$};
\draw (247.53,162.5) node [anchor=north west][inner sep=0.75pt]    {$1$};
\draw (284.92,127.8) node [anchor=north west][inner sep=0.75pt]    {$1$};
\draw (276.37,31.1) node [anchor=north west][inner sep=0.75pt]    {$1$};
\draw (298.31,44.86) node [anchor=north west][inner sep=0.75pt]    {$1$};
\draw (284.18,91.72) node [anchor=north west][inner sep=0.75pt]    {$1$};
\draw (238.42,33.73) node [anchor=north west][inner sep=0.75pt]    {$1$};
\draw (240.92,91.3) node [anchor=north west][inner sep=0.75pt]    {$1$};
\draw (263.85,13.57) node [anchor=north west][inner sep=0.75pt]  [font=\tiny]  {$1$};
\draw (334.35,12.15) node [anchor=north west][inner sep=0.75pt]    {$L_{e_{1}}{}$};
\draw (494.09,174.28) node [anchor=north west][inner sep=0.75pt]    {$r$};
\draw (395.86,137.53) node [anchor=north west][inner sep=0.75pt]    {$1$};
\draw (414.46,133.07) node [anchor=north west][inner sep=0.75pt]    {$1$};
\draw (411.66,162.18) node [anchor=north west][inner sep=0.75pt]    {$1$};
\draw (402.55,33.42) node [anchor=north west][inner sep=0.75pt]    {$1$};
\draw (405.05,90.99) node [anchor=north west][inner sep=0.75pt]    {$1$};
\draw (501.22,11.74) node [anchor=north west][inner sep=0.75pt]    {$L_{e_{1}}{}$};
\draw (660.96,173.87) node [anchor=north west][inner sep=0.75pt]    {$r$};
\draw (578.53,161.76) node [anchor=north west][inner sep=0.75pt]    {$1$};
\draw (615.92,127.07) node [anchor=north west][inner sep=0.75pt]    {$1$};

\end{tikzpicture}
}
\caption{Indecomposable direct summands of $PH_0(\hat{\mathcal{H}}^2_{-,-};\mathbb{F})$. The number on the colored block indicates the dimension of the vector space in the representation.}
\label{H0summand-H}
\end{figure}

Based on our calculation of the rank of $H_i(\hat{\mathcal{H}}^2_{r,L_{e_1}})$, where $i=0,1$, we use induction to calculate the rank of $H_i((\hat{\mathcal{H}}_{m,n})^2_{r,L_{e_1}})$ for $i=0,1$ and $m,n\geq 3$, where $\hat{\mathcal{H}}_{m,n}$ is the metric graph $\mathcal{H}_{m,n}$ such that every edge but $e_1$ has length $1$. 

\begin{theorem}\label{Thm9}
\begin{equation}
\rank(H_0((\hat{\mathcal{H}}_{m,n})^2_{-,-}))=\begin{cases} 1, & \mbox{if } r\leq 1 \\ 
m(m-3)+n(n-3)+6, & \mbox{if } 1<r\leq 2 \mbox{ and } r\leq L_{e_1}+1 \\
(m+n)(m+n-5)+6, & \mbox{if } 1<r\leq 2 \mbox{ and } L_{e_1}+1<r \\
2, & \mbox{if } 2<r \mbox{ and } r\leq L_{e_1}+1 \\
2(m-1)(n-1), & \mbox{if } 2<r \mbox{ and } L_{e_1}+1< r\leq L_{e_1}+2 \\
0, & \mbox{else}
\end{cases}
\end{equation}
and
\begin{equation}
\rank(H_1((\hat{\mathcal{H}}_{m,n})^2_{-,-}))=\begin{cases} m(m-3)+n(n-3)+3, & \mbox{if } r\leq 1 \mbox{ and } r\leq L_{e_1}\\ 
(m+n)(m+n-7)+11, & \mbox{if } r\leq 1 \mbox{ and } L_{e_1}<r\\ 
2(m-2)(n-2), & \mbox{if } 1<r \mbox{ and } L_{e_1}<r\leq L_{e_1}+1 \\
0, & \mbox{else}
\end{cases}
\end{equation}
\end{theorem}

The rank of $H_0((\hat{\mathcal{H}}_{m,n})^2_{r,L_{e_1}})$ for all $r>0$ and $L_{e_1}>0$ is shown in Figure \ref{fig:h0Hmn}. When $k\geq 4$, the rank of $H_1((\hat{\mathcal{H}}_{m,n})^2_{r,L_{e_1}})$ for all $r>0$ and $L_{e_1}>0$ is shown in Figure \ref{fig:h1Hmn}.\\

\begin{minipage}{0.5\textwidth}
\centering
\resizebox{6.5cm}{!}{

\tikzset{every picture/.style={line width=0.75pt}} 

\begin{tikzpicture}[x=0.75pt,y=0.75pt,yscale=-1,xscale=1]

\draw  [draw opacity=0][fill={rgb, 255:red, 74; green, 144; blue, 226 }  ,fill opacity=0.6 ] (40.85,313.04) -- (141.57,217.47) -- (141.82,312.77) -- cycle ;
\draw  [draw opacity=0][fill={rgb, 255:red, 184; green, 233; blue, 134 }  ,fill opacity=0.9 ] (142.21,312.78) -- (237.28,217.76) -- (237.53,312.53) -- cycle ;
\draw  [draw opacity=0][fill={rgb, 255:red, 248; green, 231; blue, 28 }  ,fill opacity=0.6 ] (141.5,35) -- (236.5,35.5) -- (238.79,118.89) -- (141.57,217.47) -- (140.9,120.03) -- cycle ;
\draw  [draw opacity=0][fill={rgb, 255:red, 248; green, 231; blue, 28 }  ,fill opacity=0.6 ] (140.76,217.83) -- (237.65,120.03) -- (238.05,217.42) -- cycle ;
\draw  (5.44,312.73) -- (350.5,312.73)(39.95,12) -- (39.95,346.15) (343.5,307.73) -- (350.5,312.73) -- (343.5,317.73) (34.95,19) -- (39.95,12) -- (44.95,19)  ;
\draw  [draw opacity=0][fill={rgb, 255:red, 208; green, 2; blue, 27 }  ,fill opacity=0.5 ] (333.68,120.33) -- (236.91,218.42) -- (237.48,119.78) -- cycle ;
\draw  [draw opacity=0][fill={rgb, 255:red, 208; green, 2; blue, 27 }  ,fill opacity=0.5 ] (344.5,35.5) -- (344.5,109) -- (333.68,120.33) -- (237.65,120.03) -- (236.5,35.5) -- cycle ;
\draw  [draw opacity=0][fill={rgb, 255:red, 74; green, 144; blue, 226 }  ,fill opacity=0.6 ] (141.5,35) -- (140.9,120.03) -- (141.57,217.47) -- (40.85,313.04) -- (40.5,35) -- cycle ;
\draw  [draw opacity=0][fill={rgb, 255:red, 248; green, 231; blue, 28 }  ,fill opacity=0.6 ] (238.05,217.18) -- (141.39,312.73) -- (140.76,217.83) -- cycle ;
\draw  [draw opacity=0][fill={rgb, 255:red, 155; green, 155; blue, 155 }  ,fill opacity=0.6 ] (344.5,109) -- (344.5,204) -- (237.9,311.87) -- (237.28,217.76) -- cycle ;

\draw (141.47,327.34) node   [align=left] {1};
\draw (238.73,326.19) node   [align=left] {2};
\draw (35.35,217.5) node   [align=left] {1};
\draw (34.21,118.53) node   [align=left] {2};
\draw (13.07,17.96) node [anchor=north west][inner sep=0.75pt]    {$L_{e_{1}}{}$};
\draw (355.76,316.52) node [anchor=north west][inner sep=0.75pt]    {$r$};
\draw (89.35,172.5) node   [align=left] {1};
\draw (288.35,89.5) node   [align=left] {2};
\draw (138.4,145.81) node [anchor=north west][inner sep=0.75pt]  [font=\scriptsize,rotate=-359.64]  {$ \begin{array}{l}
m(m-3) +n(n-3)\\
+6
\end{array}$};
\draw (252.69,209.4) node [anchor=north west][inner sep=0.75pt]  [font=\scriptsize]  {$2(m-1)(n-1)$};
\draw (155.73,287.4) node [anchor=north west][inner sep=0.75pt]  [font=\scriptsize]  {$ \begin{array}{l}
( m+n)(m+n-5)\\
+6
\end{array}$};

\end{tikzpicture}
}
\captionof{figure}{The rank of $H_0((\hat{\mathcal{H}}_{m,n})^2_{r,L_{e_1}})$}
\label{fig:h0Hmn}  
\end{minipage}%
\begin{minipage}{0.5\textwidth}
\centering
\resizebox{6.5cm}{!}{

\tikzset{every picture/.style={line width=0.75pt}} 

\begin{tikzpicture}[x=0.75pt,y=0.75pt,yscale=-1,xscale=1]

\draw  [draw opacity=0][fill={rgb, 255:red, 139; green, 87; blue, 42 }  ,fill opacity=0.6 ] (40.85,313.05) -- (141.96,217.47) -- (142.21,312.78) -- cycle ;
\draw  [draw opacity=0] (142.21,312.78) -- (237.28,217.76) -- (237.53,312.53) -- cycle ;
\draw  (6.34,313.04) -- (351.4,313.04)(40.85,12.31) -- (40.85,346.45) (344.4,308.04) -- (351.4,313.04) -- (344.4,318.04) (35.85,19.31) -- (40.85,12.31) -- (45.85,19.31)  ;
\draw  [draw opacity=0][fill={rgb, 255:red, 80; green, 227; blue, 194 }  ,fill opacity=0.6 ] (141.5,35) -- (140.9,120.03) -- (141.57,217.47) -- (40.85,313.04) -- (40.5,35) -- cycle ;
\draw  [draw opacity=0][fill={rgb, 255:red, 208; green, 2; blue, 27 }  ,fill opacity=0.6 ] (331.51,26.59) -- (333.7,120.75) -- (142.21,312.78) -- (140.93,217.65) -- cycle ;

\draw (141.47,327.34) node   [align=left] {1};
\draw (238.73,326.19) node   [align=left] {2};
\draw (35.35,217.5) node   [align=left] {1};
\draw (34.21,118.53) node   [align=left] {2};
\draw (13.07,17.96) node [anchor=north west][inner sep=0.75pt]    {$L_{e_{1}}{}$};
\draw (355.76,316.52) node [anchor=north west][inner sep=0.75pt]    {$r$};
\draw (42.73,150.4) node [anchor=north west][inner sep=0.75pt]  [font=\scriptsize]  {$ \begin{array}{l}
m( m-3) +n( n-3)\\
+3
\end{array}$};
\draw (64.73,285.4) node [anchor=north west][inner sep=0.75pt]  [font=\scriptsize,color={rgb, 255:red, 255; green, 255; blue, 255 }  ,opacity=1 ]  {$ \begin{array}{l}
( m+n)( m+n-7)\\
+10
\end{array}$};
\draw (171.5,185.9) node [anchor=north west][inner sep=0.75pt]  [font=\scriptsize,color={rgb, 255:red, 0; green, 0; blue, 0 }  ,opacity=1 ]  {$2mn-4m-4n+7$};

\end{tikzpicture}
}
\captionof{figure}{The rank of $H_1((\hat{\mathcal{H}}_{m,n})^2_{r,L_{e_1}})$}
\label{fig:h1Hmn}  
\end{minipage}

Next, we construct a persistence module $PH_i((\hat{\mathcal{H}}_{m,n})^2_{-,-})$ for $i=0,1$. We prove the following theorems regarding the interval decomposability of $PH_0((\hat{\mathcal{H}}_{m,n})^2_{-,-};\mathbb{F})$ and $PH_1((\hat{\mathcal{H}}_{m,n})^2_{-,-};\mathbb{F})$.

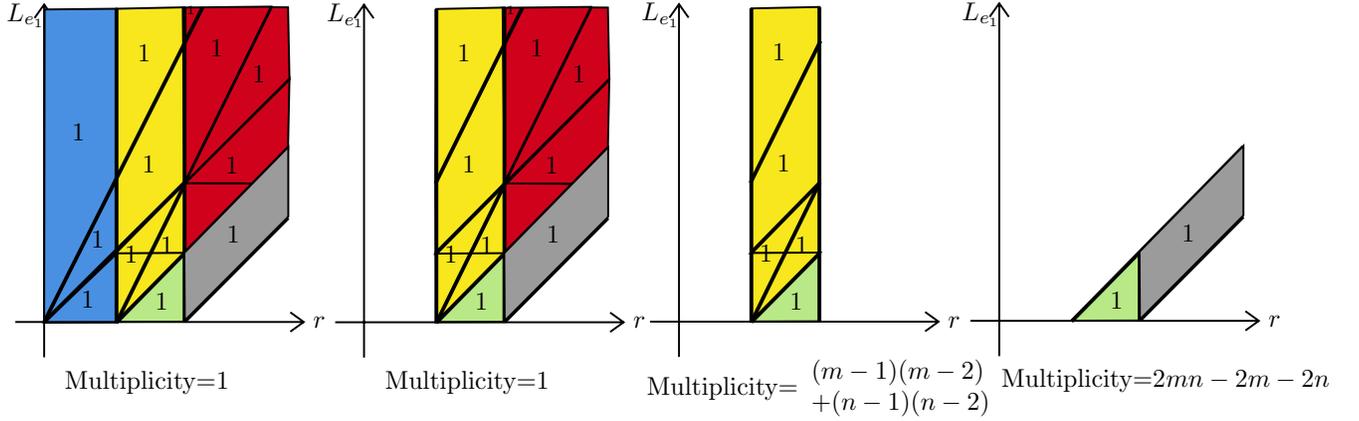
\begin{figure}[htbp!]
\centering
\resizebox{18cm}{!}{

\tikzset{every picture/.style={line width=0.75pt}} 

\begin{tikzpicture}[x=0.75pt,y=0.75pt,yscale=-1,xscale=1]

\draw  [draw opacity=0][fill={rgb, 255:red, 74; green, 144; blue, 226 }  ,fill opacity=0.6 ] (24.66,179.15) -- (62.13,143.6) -- (62.22,179.05) -- cycle ;
\draw  [draw opacity=0][fill={rgb, 255:red, 184; green, 233; blue, 134 }  ,fill opacity=0.9 ] (62.22,178.7) -- (97.58,143.36) -- (97.67,178.61) -- cycle ;
\draw  [draw opacity=0][fill={rgb, 255:red, 248; green, 231; blue, 28 }  ,fill opacity=0.6 ] (62.22,16.61) -- (97.67,15.87) -- (98.14,106.58) -- (61.98,143.25) -- (61.73,107.01) -- cycle ;
\draw  [draw opacity=0][fill={rgb, 255:red, 248; green, 231; blue, 28 }  ,fill opacity=0.6 ] (61.68,143.38) -- (97.72,107.01) -- (97.87,143.23) -- cycle ;
\draw  (9.49,178.79) -- (159.78,178.79)(24.52,14.4) -- (24.52,197.06) (152.78,173.79) -- (159.78,178.79) -- (152.78,183.79) (19.52,21.4) -- (24.52,14.4) -- (29.52,21.4)  ;
\draw  [draw opacity=0][fill={rgb, 255:red, 208; green, 2; blue, 27 }  ,fill opacity=0.5 ] (142.82,15.52) -- (97.65,106.91) -- (96.73,15.99) -- cycle ;
\draw  [draw opacity=0][fill={rgb, 255:red, 208; green, 2; blue, 27 }  ,fill opacity=0.5 ] (133.44,107.12) -- (97.44,143.6) -- (97.65,106.91) -- cycle ;
\draw  [draw opacity=0][fill={rgb, 255:red, 208; green, 2; blue, 27 }  ,fill opacity=0.5 ] (152.34,52.67) -- (151.59,88.77) -- (133.16,107.08) -- (133.44,107.12) -- (126.49,106.99) -- (97.29,106.95) -- cycle ;
\draw  [draw opacity=0][fill={rgb, 255:red, 74; green, 144; blue, 226 }  ,fill opacity=0.6 ] (62.22,16.61) -- (61.73,107.01) -- (61.98,143.25) -- (24.52,178.79) -- (24.76,16.98) -- cycle ;
\draw  [draw opacity=0][fill={rgb, 255:red, 248; green, 231; blue, 28 }  ,fill opacity=0.6 ] (98.17,143.01) -- (62.21,178.55) -- (61.98,143.25) -- cycle ;
\draw [line width=1.5]    (24.52,178.79) -- (152.34,52.67) ;
\draw [line width=1.5]    (62.22,179.05) -- (151.59,88.77) ;
\draw [line width=1.5]    (62.22,179.05) -- (62.22,16.61) ;
\draw [line width=1.5]    (97.67,178.61) -- (97.67,15.87) ;
\draw [line width=1.5]    (62.22,179.05) -- (143.04,16.24) ;
\draw [line width=1.5]    (24.52,178.79) -- (107.04,15.99) ;
\draw [line width=1.5]    (97.67,178.61) -- (151.59,124.85) ;
\draw  [draw opacity=0][fill={rgb, 255:red, 208; green, 2; blue, 27 }  ,fill opacity=0.5 ] (152.34,52.67) -- (151.59,15.87) -- (143.17,15.8) -- (98.28,106.33) -- cycle ;
\draw  [draw opacity=0][fill={rgb, 255:red, 155; green, 155; blue, 155 }  ,fill opacity=0.6 ] (151.73,88.17) -- (151.73,124.25) -- (97.81,178.01) -- (97.58,143.01) -- cycle ;
\draw  [draw opacity=0][fill={rgb, 255:red, 184; green, 233; blue, 134 }  ,fill opacity=0.9 ] (228.89,178.84) -- (264.25,143.5) -- (264.34,178.75) -- cycle ;
\draw  [draw opacity=0][fill={rgb, 255:red, 248; green, 231; blue, 28 }  ,fill opacity=0.6 ] (228.89,16.75) -- (264.34,16.01) -- (264.81,106.72) -- (228.65,143.39) -- (228.4,107.15) -- cycle ;
\draw  [draw opacity=0][fill={rgb, 255:red, 248; green, 231; blue, 28 }  ,fill opacity=0.6 ] (228.35,143.52) -- (264.38,107.15) -- (264.54,143.37) -- cycle ;
\draw  (176.16,178.93) -- (326.45,178.93)(191.19,14.54) -- (191.19,197.2) (319.45,173.93) -- (326.45,178.93) -- (319.45,183.93) (186.19,21.54) -- (191.19,14.54) -- (196.19,21.54)  ;
\draw  [draw opacity=0][fill={rgb, 255:red, 208; green, 2; blue, 27 }  ,fill opacity=0.5 ] (309.49,15.66) -- (264.32,107.05) -- (263.4,16.13) -- cycle ;
\draw  [draw opacity=0][fill={rgb, 255:red, 208; green, 2; blue, 27 }  ,fill opacity=0.5 ] (300.1,107.26) -- (264.11,143.74) -- (264.32,107.05) -- cycle ;
\draw  [draw opacity=0][fill={rgb, 255:red, 208; green, 2; blue, 27 }  ,fill opacity=0.5 ] (319.01,52.81) -- (318.26,88.91) -- (299.83,107.22) -- (300.1,107.26) -- (293.16,107.13) -- (263.96,107.09) -- cycle ;
\draw  [draw opacity=0][fill={rgb, 255:red, 248; green, 231; blue, 28 }  ,fill opacity=0.6 ] (264.84,143.15) -- (228.88,178.69) -- (228.65,143.39) -- cycle ;
\draw [line width=1.5]    (228.35,143.52) -- (319.01,52.81) ;
\draw [line width=1.5]    (228.89,179.19) -- (318.26,88.91) ;
\draw [line width=1.5]    (228.89,179.19) -- (228.89,16.75) ;
\draw [line width=1.5]    (264.34,178.75) -- (264.34,16.01) ;
\draw [line width=1.5]    (228.89,179.19) -- (309.71,16.38) ;
\draw [line width=1.5]    (228.4,107.15) -- (273.71,16.13) ;
\draw [line width=1.5]    (264.34,178.75) -- (318.26,124.99) ;
\draw  [draw opacity=0][fill={rgb, 255:red, 208; green, 2; blue, 27 }  ,fill opacity=0.5 ] (319.01,52.81) -- (318.26,16.01) -- (309.84,15.94) -- (264.95,106.47) -- cycle ;
\draw  [draw opacity=0][fill={rgb, 255:red, 155; green, 155; blue, 155 }  ,fill opacity=0.6 ] (318.4,88.32) -- (318.4,124.39) -- (264.48,178.15) -- (264.25,143.15) -- cycle ;
\draw  [draw opacity=0][fill={rgb, 255:red, 184; green, 233; blue, 134 }  ,fill opacity=0.9 ] (393.01,178.53) -- (428.37,143.18) -- (428.47,178.43) -- cycle ;
\draw  [draw opacity=0][fill={rgb, 255:red, 248; green, 231; blue, 28 }  ,fill opacity=0.6 ] (393.01,16.44) -- (428.47,15.69) -- (428.94,106.41) -- (392.77,143.08) -- (392.53,106.83) -- cycle ;
\draw  [draw opacity=0][fill={rgb, 255:red, 248; green, 231; blue, 28 }  ,fill opacity=0.6 ] (392.47,143.21) -- (428.51,106.83) -- (428.66,143.06) -- cycle ;
\draw  (340.28,178.62) -- (490.57,178.62)(355.31,14.23) -- (355.31,196.89) (483.57,173.62) -- (490.57,178.62) -- (483.57,183.62) (350.31,21.23) -- (355.31,14.23) -- (360.31,21.23)  ;
\draw  [draw opacity=0][fill={rgb, 255:red, 248; green, 231; blue, 28 }  ,fill opacity=0.6 ] (428.96,142.84) -- (393.01,178.37) -- (392.77,143.08) -- cycle ;
\draw [line width=1.5]    (392.77,143.08) -- (428.94,106.41) ;
\draw [line width=1.5]    (393.01,178.87) -- (428.96,142.84) ;
\draw [line width=1.5]    (393.01,178.87) -- (393.01,16.44) ;
\draw [line width=1.5]    (428.47,178.43) -- (428.47,15.69) ;
\draw [line width=1.5]    (393.01,178.87) -- (428.51,106.83) ;
\draw [line width=1.5]    (392.53,106.83) -- (429.03,33.35) ;
\draw  [draw opacity=0][fill={rgb, 255:red, 184; green, 233; blue, 134 }  ,fill opacity=0.9 ] (559.89,178.11) -- (595.25,142.76) -- (595.34,178.02) -- cycle ;
\draw  (507.16,178.2) -- (657.45,178.2)(522.19,13.81) -- (522.19,196.47) (650.45,173.2) -- (657.45,178.2) -- (650.45,183.2) (517.19,20.81) -- (522.19,13.81) -- (527.19,20.81)  ;
\draw [line width=1.5]    (559.89,178.46) -- (649.26,88.18) ;
\draw [line width=1.5]    (595.34,178.02) -- (595.34,142.42) ;
\draw [line width=1.5]    (595.34,178.02) -- (649.26,124.26) ;
\draw  [draw opacity=0][fill={rgb, 255:red, 155; green, 155; blue, 155 }  ,fill opacity=0.6 ] (649.4,87.58) -- (649.4,123.66) -- (595.48,177.42) -- (595.25,142.42) -- cycle ;

\draw (337.17,195.74) node [anchor=north west][inner sep=0.75pt]   [align=left] {Multiplicity=$\begin{array}{{>{\displaystyle}l}}
( m-1)( m-2)\\
+( n-1)( n-2)
\end{array}$};
\draw (521.89,201.74) node [anchor=north west][inner sep=0.75pt]   [align=left] {Multiplicity=$2mn-2m-2n$};
\draw (33.84,202.47) node [anchor=north west][inner sep=0.75pt]   [align=left] {Multiplicity=$1$};
\draw (200.71,202.47) node [anchor=north west][inner sep=0.75pt]   [align=left] {Multiplicity=$1$};
\draw (3.56,12.33) node [anchor=north west][inner sep=0.75pt]    {$L_{e_{1}}{}$};
\draw (163.29,174.46) node [anchor=north west][inner sep=0.75pt]    {$r$};
\draw (37.54,74.47) node [anchor=north west][inner sep=0.75pt]    {$1$};
\draw (47.59,130.27) node [anchor=north west][inner sep=0.75pt]    {$1$};
\draw (42.38,161.14) node [anchor=north west][inner sep=0.75pt]    {$1$};
\draw (65.07,137.71) node [anchor=north west][inner sep=0.75pt]    {$1$};
\draw (83.66,133.24) node [anchor=north west][inner sep=0.75pt]    {$1$};
\draw (80.87,162.36) node [anchor=north west][inner sep=0.75pt]    {$1$};
\draw (118.26,127.66) node [anchor=north west][inner sep=0.75pt]    {$1$};
\draw (109.7,30.96) node [anchor=north west][inner sep=0.75pt]    {$1$};
\draw (131.65,44.72) node [anchor=north west][inner sep=0.75pt]    {$1$};
\draw (117.51,91.58) node [anchor=north west][inner sep=0.75pt]    {$1$};
\draw (71.75,33.59) node [anchor=north west][inner sep=0.75pt]    {$1$};
\draw (74.26,91.16) node [anchor=north west][inner sep=0.75pt]    {$1$};
\draw (97.18,13.43) node [anchor=north west][inner sep=0.75pt]  [font=\tiny]  {$1$};
\draw (170.22,12.47) node [anchor=north west][inner sep=0.75pt]    {$L_{e_{1}}{}$};
\draw (329.96,174.6) node [anchor=north west][inner sep=0.75pt]    {$r$};
\draw (231.74,137.85) node [anchor=north west][inner sep=0.75pt]    {$1$};
\draw (250.33,133.38) node [anchor=north west][inner sep=0.75pt]    {$1$};
\draw (247.53,162.5) node [anchor=north west][inner sep=0.75pt]    {$1$};
\draw (284.92,127.8) node [anchor=north west][inner sep=0.75pt]    {$1$};
\draw (276.37,31.1) node [anchor=north west][inner sep=0.75pt]    {$1$};
\draw (298.31,44.86) node [anchor=north west][inner sep=0.75pt]    {$1$};
\draw (284.18,91.72) node [anchor=north west][inner sep=0.75pt]    {$1$};
\draw (238.42,33.73) node [anchor=north west][inner sep=0.75pt]    {$1$};
\draw (240.92,91.3) node [anchor=north west][inner sep=0.75pt]    {$1$};
\draw (263.85,13.57) node [anchor=north west][inner sep=0.75pt]  [font=\tiny]  {$1$};
\draw (334.35,12.15) node [anchor=north west][inner sep=0.75pt]    {$L_{e_{1}}{}$};
\draw (494.09,174.28) node [anchor=north west][inner sep=0.75pt]    {$r$};
\draw (395.86,137.53) node [anchor=north west][inner sep=0.75pt]    {$1$};
\draw (414.46,133.07) node [anchor=north west][inner sep=0.75pt]    {$1$};
\draw (411.66,162.18) node [anchor=north west][inner sep=0.75pt]    {$1$};
\draw (402.55,33.42) node [anchor=north west][inner sep=0.75pt]    {$1$};
\draw (405.05,90.99) node [anchor=north west][inner sep=0.75pt]    {$1$};
\draw (501.22,11.74) node [anchor=north west][inner sep=0.75pt]    {$L_{e_{1}}{}$};
\draw (660.96,173.87) node [anchor=north west][inner sep=0.75pt]    {$r$};
\draw (578.53,161.76) node [anchor=north west][inner sep=0.75pt]    {$1$};
\draw (615.92,127.07) node [anchor=north west][inner sep=0.75pt]    {$1$};

\end{tikzpicture}
}
\caption{Indecomposable direct summands of $PH_0((\hat{\mathcal{H}}_{m,n})^2_{-,-};\mathbb{F})$. The number on the colored block indicates the dimension of the vector space in the representation.}
\label{H0summand-Hmn}
\end{figure}

\begin{theorem}\label{Thm10}
$PH_i((\hat{\mathcal{H}}_{m,n})^2_{-,-})$ is interval decomposable, where $i=0,1$. The indecomposable direct summands of $PH_0((\hat{\mathcal{H}}_{m,n})^2_{-,-};\mathbb{F})$ are given in Figure \ref{H0summand-Hmn}.
\end{theorem}

\paragraph{Structure of the paper}
Section \ref{prelim} contains the preliminaries and notations of categories used in this paper. Section \ref{conf:star3} provides the calculation of the second configuration spaces of star graphs with restraint parameter $r\in \mathbb{R}_{>0}$ and edge length vector $\mathbf{L}$. Section \ref{chap4:y2} proves $PH_0(Y^2_{-,-};\mathbb{F})$ is not interval decomposable and describes all indecomposable direct summand of $PH_0(Y^2_{-,-};\mathbb{F})$. In addition, section \ref{chap4:y2} proves $PH_1(Y^2_{-,-};\mathbb{F})$ is interval decomposable. Section \ref{decompo:star3} proves $PH_0((\mathsf{Star_k})^2_{-,-};\mathbb{F})$ is not interval decomposable and describes all indecomposable direct summand of $PH_0((\mathsf{Star_k})^2_{-,-};\mathbb{F})$. Section \ref{conf:H} provides the calculation of the second configuration spaces of H graphs with restraint parameter $r\in \mathbb{R}_{>0}$ and edge length vector $\mathbf{L}$. Section \ref{decompo:H33} proves $PH_i(\hat{\mathcal{H}}^2_{-,-};\mathbb{F})$ is interval decomposable for $i=0,1$. The indecomposable direct summands of $PH_0(\hat{\mathcal{H}}^2_{-,-};\mathbb{F})$ are also provided in Section \ref{decompo:H33}. Section \ref{genconf:H} provides the calculation of the second configuration spaces of the generalized H graphs with restraint parameter $r\in \mathbb{R}_{>0}$ and edge length vector $\mathbf{L}$. Section \ref{decompo:Hmn} proves $PH_i((\hat{\mathcal{H}}_{m,n})_{-,-}^2;\mathbb{F})$ is interval decomposable for $i=0,1$. The indecomposable direct summands of $PH_0((\hat{\mathcal{H}}_{m,n})_{-,-}^2;\mathbb{F})$ are also provided in Section \ref{decompo:Hmn}.

\paragraph{Acknowledgements}
I would like to thank my advisor Dr.~Murad~\"Ozayd\i n, for his guidance and helpful discussions during this project.

%% file: chapters/section1.2.tex
\subsection{Parametric Polytopes}\label{ParaPoly}

A parametric polytope $c_{\mathbf{b}}$ in $\mathbb{R}^n$ is determined by $Ax\leq \mathbf{b}$ where $A$ is an $m\times n$ matrix and $\mathbf{b}$ is a $m\times 1$ column vector representing the constraints. Each feasible solution of $Ax\leq \mathbf{b}$ corresponds to a point of $c_{\mathbf{b}}$. Note that $c_{\mathbf{b}}$ is a CW-complex, and the cell structure is also encoded in $Ax\leq \mathbf{b}$. Let $\boldsymbol{\beta}=\{\lambda_1,\dots,\lambda_n\}$. Then $A_{\boldsymbol{\beta}}$ is the submatrix of $A$ with rows $\lambda\in\boldsymbol{\beta}$ and $\mathbf{b}_{\boldsymbol{\beta}}=(b_{\lambda})_{\lambda\in\boldsymbol{\beta}}$. A potential 0-dimensional cell of $c_{\mathbf{b}}$ is given by $A_{\boldsymbol{\beta}}^{-1}\mathbf{b}_{\boldsymbol{\beta}}$ for some $\boldsymbol{\beta}$ (where $A_{\boldsymbol{\beta}}^{-1}$ exists) and such potential $0$-dimensional cell is indeed an $0$-dimensional cell of $c_{\mathbf{b}}$ if $A_{\boldsymbol{\beta}}^{-1}\mathbf{b}_{\boldsymbol{\beta}}$ is a solution of $Ax\leq \mathbf{b}$. We label each 0-cell by $$S_{\mathbf{b}}:=\{\beta\mid A_{\beta}\mbox{ is a solution of } Ax\leq \mathbf{b}\}$$

We assign a label $v_{\mathbf{b}}$ to each vertex of the parametric polytope where $v_{\mathbf{b}}=\{\lambda\mid A_{\lambda}A^{-1}\mathbf{b}= b_\lambda\}$. Let $T_b$ denote the type of $c_b$: $T_{\mathbf{b}}=\{v_{\mathbf{b}}(\beta)\mid \beta\in S_{\mathbf{b}}\}$. Dover and \"Ozayd\i n \cite{dover2013homeomorphism} showed the following lemmas:
\begin{lemma}
The family of polytopes $\{c_{\mathbf{b}}\}_{\mathbf{b}\in\mathbb{R^n}}$ has only finitely many combinatorial types.
\end{lemma}

\begin{lemma}
Let $\mathbf{b}'\in\mathbb{R}^n$. If $T_{\mathbf{b}'}\neq\emptyset$, then the set of all $\mathbf{b}\in\mathbb{R}^n$ with $T_{\mathbf{b}}=T_{\mathbf{b}'}$ is convex.
\end{lemma}

Let $X$ be a finite metric graph. The cellular structure of $X^n_{\mathbf{r}}$ can be characterized by parametric polytopes \cite{dover2013homeomorphism}. Each maximal cell of $X^n_{\mathbf{r}}$ consists of the solutions of inequality system:\\

For $1\leq i\leq n$:
\begin{equation}\label{pp1}
0\leq  x_i\leq  L_{e_i}
\end{equation}

For $1\leq i<j\leq n$:
\begin{equation}\label{pp2}
\begin{aligned}
x_i+x_j+\delta(a_i,a_j)&\geq r_{ij};\\
L_{e_i}-x_i+x_j+\delta(b_i,a_j)&\geq r_{ij};\\
x_i+L_{e_j}-x_j+\delta(a_i,b_j)&\geq r_{ij};\\
L_{e_i}-x_i+L_{e_j}-x_j+\delta({b_i,b_j})&\geq r_{ij}
\end{aligned}
\end{equation}

If there exists $1\leq i,j\leq n$ such that $e_i=e_j$, with further assumption that $x_i\leq x_j$, 
\begin{equation}\label{pp3}
\begin{aligned}
x_i-x_j&\leq 0;\\
x_i-x_j&\leq -r_{ij};\\
L_{e_i}-x_i+x_j+\delta(a_i,b_i)&\geq r_{ij}
\end{aligned}
\end{equation}

Let $X:=X_{\mathbf{L}}$ be a metric graph with edge length vector $\mathbf{L}$. Define an equivalence relation on the parameter space as follows: $\mathbf{r}\sim\mathbf{s}$ if and only if $X^n_{\mathbf{r}}$ and $X^n_{\mathbf{s}}$ have the same combinatorial type. It is known that $X^n_{\mathbf{r}}$ is homeomorphic to $X^n_{\mathbf{s}}$ if $\mathbf{r}\sim\mathbf{s}$.\cite{dover2013homeomorphism} $\mathbf{r}$ is called \textbf{critical} if the combinatorial type (face poset) of the cell of $(X_{\mathbf{L}})^n_{\mathbf{r}}$ changes at $\mathbf{r}$. When $X$ is a metric tree, $\mathbf{r}$ is critical only when one of the constraints in (\ref{pp1}) and one of the constraints in (\ref{pp2}) or (\ref{pp3}) are binding. By considering the critical $\mathbf{r}$ among all possible $\mathbf{L}$, we obtain critical hyperplanes in the parameter space. Those hyperplanes subdivide the parameter space into a disjoint union of subsets which we call ``chambers". Note that $X^n_{\mathbf{r}}$ and $X^n_{\mathbf{s}}$ have the same combinatorial type if $\mathbf{r}$ and $\mathbf{r}$ are contained in the same chamber.

\begin{theorem}[Dover-\"Ozayd\i n, 2013\cite{dover2013homeomorphism}]
For a metric graph $X$ with $E$ edges, the number of homotopy types for the family of spaces $\{X^{2}_{r}\}_{r\geq 0}$ is bounded above by $\frac{9}{2}E^{2}-\frac{5}{2}E+1$, and the number of isotopy types is bounded above by $9E^{2}-5E$.
\end{theorem}

\begin{theorem}[Dover-\"Ozayd\i n, 2013\cite{dover2013homeomorphism}]
Let $n$ be a fixed positive integer and $\Omega$ be a $d$-dimensional affine subspace of $\mathbb{R}^{{n}\choose {2}}$. The number of isotopy types of $\{X^n_{\mathbf{r}}\}_{\mathbf{r}\in\Omega}$ is bounded above by a polynomial of degree $nd$ in the number of edges of $X$.
\end{theorem}

When the norm of parameter $\mathbf{r}$ is sufficiently small, $X^n_{\mathbf{r}}$ is homotopy equivalent to the $n$-th configuration space $X^{\underline{n}}$:

\begin{proposition}
Let $X$ be a connected regular metric graph, and let $L_e$ denote the edge length of $e\in E$ where $E$ is the set of edges of $X$. Define $a:=\min\limits_{e\in E} L_{e}$. If $\lvert\lvert \mathbf{r}\rvert\rvert_{\infty}<\dfrac{a}{n}$, then the inclusion $X^{n}_{\mathbf{r}}\hookrightarrow X^{\underline{n}}$ is a homotopy equivalence.
\end{proposition}
\begin{proof}
If $\mathbf{r}\leq \mathbf{s}$ and $\lvert\lvert \mathbf{r}\rvert\rvert_{\infty}\leq \lvert\lvert \mathbf{s}\rvert\rvert_{\infty}<\dfrac{a}{n}$, then $\mathbf{r}$ and $\mathbf{s}$ belong to the same chamber in the parameter space. Note that each chamber is convex, the maps $X^{n}_{\mathbf{r}}\rightarrow X^{n}_{\mathbf{s}}\hookrightarrow X^{n}_{\mathbf{r}}$ and $X^{n}_{\mathbf{s}}\hookrightarrow X^{n}_{\mathbf{r}}\rightarrow X^{n}_{\mathbf{s}}$ are isotopic to the identity maps. Thus $ \pi_{k}(X^{n}_{\mathbf{s}})\rightarrow \pi_{k}(X^{n}_{\mathbf{r}})$ is an isomorphism for all $k\in\mathbb{N}$. Note that $\underset{}{\varinjlim}X^{n}_{\mathbf{r}}\cong X^{\underline{n}}$, we update the colimit cone with the leg maps $i_{\mathbf{r}}: X^n_{\mathbf{r}}\hookrightarrow X^{\underline{n}}$, for all $\mathbf{r}$. Note that, by Theorem \ref{directlim_exact}, $\underset{}{\varinjlim}$ is exact, therefore, 
$$\underset{}{\varinjlim}\pi_{k}(X^{n}_{\mathbf{r}})\cong\pi_{k}( \underset{}{\varinjlim}X^{n}_{\mathbf{r}})=\pi_{k}(X^{\underline{n}})$$
By Lemma \ref{lemma1-directed}, there exists $\mathbf{r}$ such that the induced map $\pi_k{i_{\mathbf{r}}}:\pi_{k}(X^{n}_{\mathbf{r}})\rightarrow \pi_{k}(X^{\underline{n}})$ is an isomorphism.
By Whitehead's theorem, we conclude that the inclusion $X^{n}_{\mathbf{r}}\hookrightarrow X^{\underline{n}}$ is a homotopy equivalence.
\end{proof}

%% file: chapters/section1.3.tex
\subsection{Persistence Modules}\label{Pers:chap3}
Let $\mathbb{F}$ be a field and $(P,\leq)$ be a poset. Let $\cat{C}$ be a category. $M:(P,\leq)\rightarrow \cat{C}$ is called a \textbf{($P$-indexed) persistence object}\index{persistence object}. $M$ is called a \textbf{($P$-indexed) filtration}\index{(filtration} when $\cat{C}=\cat{Top}$ and $M(i\leq j)$ is an inclusion for all $i\leq j\in P$. $M$ is called a \textbf{($P$-indexed) persistence module}\index{persistence module} when $\cat{C}=\cat{Vect}_{\mathbb{F}}$. In other words, the ($P$-indexed) persistence module $M$ consists of a family of $\mathbb{F}$-vector spaces $\{Mt\mid t\in P\}$ and a doubly-indexed family of linear maps $\{\rho_{st}: Ms\rightarrow Mt\mid s\leq t\}$ where $\rho_{tu}\rho_{st}=\rho_{su}$ for any $s\leq t\leq u$ in $P$ and $\rho_{ss}=\mathrm{id}_{Ms}$ for all $s\in P$. Alternatively, $M$ can be viewed as a poset representation over the poset $(P,\leq)$. We use $\cat{Vect}_{\mathbb{F}}^{(P,\leq)}$ to denote the category of persistence modules over the poset $(P,\leq)$. Note that $\cat{Vect}_{\mathbb{F}}^{(P,\leq)}$ has a full subcategory (denoted by $\cat{vect}_{\mathbb{F}}^{(P,\leq)}$) consists of point-wise finite dimensional persistence modules\footnote{A persistence module $M$ is called point-wise finite dimensional if $M(i)$ is a finite dimensional vector space for all $i\in (P,\leq)$.}. When $(P,\leq)=(\mathbb{R}^n,\leq)$ with the product order, where $(a_1,\dots, a_n)\leq (b_1,\dots, b_n)\in(\mathbb{R}^n,\leq)$ if and only if $a_i\leq b_i$ for all $i=1,\dots,n$, the objects of $\cat{Vect}_{\mathbb{F}}^{(\mathbb{R}^n,\leq)}$ is called \textbf{$n$-parameter persistence modules}.\\

We define the \textbf{$i$-th persistent homology functor} $PH_i(-):\cat{Top}^{(P,\leq)}\rightarrow \cat{Ab}^{(P,\leq)}$ as follows: 
\begin{itemize}
\item At the object level, $PH_i(-)$ sends each $P$-indexed filtration $M$ to a $P$-indexed persistence abelian group $H_i(M)$, where $H_i(M): (P,\leq)\rightarrow \cat{Vect}_{\mathbb{F}}$ sends each $p\in P$ to $H_i(Mp)$ and each arrow $p\leq q$ to the group homomorphism $H_i(Mp)\rightarrow H_i(Mq)$ induced by the inclusion map $Mp\hookrightarrow Mq$ in $\cat{Top}$;
\item An the morphism level, $PH_i(-)$ sends each natural transformation $\alpha: M\Rightarrow N$ to $H_i(\alpha): H_i(M)\Rightarrow H_i(N)$. The naturality of $H_i(\alpha)$ is clear because $H_i(-)$ is a functor hence it preserves commutative diagrams. 
\end{itemize}

In other words, for each $M\in\ob \cat{Top}^{(P,\leq)}$, $PH_i(M)$ can be obtained by post-composing $H_i(-)$ with $M$. The functoriality of $PH_i(-)$ is clear.\\



For example, $PH_i(Y^2_{-,-};\mathbb{F})$ is the $i$-th persistence module corresponding to the bifiltration of the second configuration spaces of the metric Y graph with restraint parameters $r$ and $L_{e_1}$. $PH_i(Y^2_{-,-};\mathbb{F})$ can be viewed as a composition of functors $Y^2_{-,-}:(\mathbb{R}_{>0},\leq)^{\op}\times (\mathbb{R}_{>0},\leq)\rightarrow \cat{Top}$ and $H_i(-;\mathbb{F})$. More generally, when $X$ is a finite metric graph, $PH_i(X^2_{-,-};\mathbb{F})$ is the $i$-th persistence module corresponding to the bifiltration of the second configuration spaces of $X$ with restraint parameters $r$ and $L_{e_1}$. $PH_i(X^2_{-,-};\mathbb{F})$ can be viewed as a composition of functors $X^2_{-,-}: (\mathbb{R}_{>0},\leq)^{\op}\times (\mathbb{R}_{>0},\leq)\rightarrow \cat{Top}$ and $H_i(-;\mathbb{F})$.\\

Let $(P,\leq)$ be a poset and $\cat{I}$ be a subset of $P$. $\cat{I}$ is \textbf{convex} if for any $a\leq b\leq c$ in $P$, $a,c\in\cat{I}$ implies $b\in\cat{I}$. $\cat{I}$ is connected if it is connected as a category. $\cat{I}$ is an \textbf{interval} if it is connected and convex.

\begin{definition}\index{interval module}
Let $(P,\leq)$ be a poset and $\cat{I}\subseteq (P,\leq)$ is an interval. Define the \textbf{interval module} $\mathbb{F}\cat{I}$ as 
$$\mathbb{F}\cat{I}t = \begin{cases} \mathbb{F}, & \mbox{if } t\in \cat{I};\\ 0, & \mbox{if } t\notin \cat{I}\end{cases} \mbox{\quad and \quad} \hom(\mathbb{F}\cat{I}t, \mathbb{F}\cat{I}s) = \begin{cases} \{id_{\mathbb{F}}\}, & \mbox{if } t\leq s\in \cat{I};\\ 0, & \mbox{else} \end{cases} $$
When $(P,\leq)=(\mathbb{R}^n,\leq)$ and $n\geq 2$, we call $\mathbb{F}\cat{I}$ the \textbf{polytope module}\index{polytope module} over $\cat{I}$. 
\end{definition}

\begin{definition}
Let $M\in \cat{Vect}_{\mathbb{F}}^{(P,\leq)}$. $M$ is \textbf{decomposable}\index{decomposable} if there exists non-trivial subrepresentation $N$ and $N'$ such that $Mt\cong Nt\oplus N't$ for all $t\in P$. We say $M$ is \textbf{indecomposable}\index{indecomposable} if it is not decomposable.
\end{definition}

\begin{lemma}\label{thinpoly}
Interval modules are thin and indecomposable.
\end{lemma}
\begin{proof}
Let $M$ be an interval module. The support of $M$ is an interval of $(P,\leq)$. Consider the endomorphism ring of $M$. We want to show $\End{M}\cong\mathbb{F}$. Let $f\in\End(M)$. Note that for every $p\in (P,\leq)$, $f_p: \mathbb{F}\rightarrow \mathbb{F}$ is a linear transformation, hence $f_p(x)=c_fx$ for some $c_f\in\mathbb{F}$. Let $p'$ be a point in the support of $M$. Note that there exists a zigzag path from $p$ to $p'$ because the support of $M$ is an interval. Since $f$ is a morphism between two representations and $M$ is a polytope module, we are forced to have $f_q(x)=c_fx$. Define $\Phi:\End(M)\rightarrow\mathbb{F}$ by $\Phi(f)=c_f$ for each $f\in \End(M)$. It is clear that $\Phi$ is bijective, so we only need to show $\Phi$ is a ring homomorphism. Note that $\Phi(f+g)=c_f+c_g=\Phi(f)+\Phi(g)$ and $\Phi(f\circ g)=c_f\cdot g=\Phi(f)\Phi(g)$, hence $\Phi$ is a ring homomorphism. 
\end{proof}

%% file: chapters/section1.6.tex
\subsection{Mayer-Vietoris Spectral Sequences}
In this section, we give an introduction to the Mayer-Vietoris spectral sequences. The set up for the general spectral sequences can be found in any standard textbooks on homological algebra, for example, \cite{weibel1995introduction} and \cite{rotman2009introduction}. We are going to use $\cat{A}$ to denote an abelian category.\\ 

A \textbf{(homology) spectral sequence}\index{spectral sequence} in $\cat{A}$ is a collection of objects $\{E^r_{pq}\}$ of $\cat{A}$, for all $p,q\in\mathbb{Z}$ and $r\geq a$ (for a fixed $a$), together with morphisms $d^r_{pq}: E^r_{pq}\rightarrow E^r_{p-r,q+r-1}$ such that $d^r_{pq}\circ d^r_{p-r,q+r-1}=0$, satisfying 
$$E^{r+1}_{pq}\cong\frac{\ker d^r_{pq}}{\Ima d^r_{p+r,q-r+1}}$$
For each $r$, the collection of $E^r_{pq}$ is called the \textbf{$r$-th page} of the spectral sequence. $\{E^r_{pq}\}$ is said to be \textbf{bounded} if for each $n:=p+q$ (called the total degree of $E^r_{pq}$), there exists only finitely many non-zero term $E^a_{pq}$.\\ 

When $E^r_{pq}=0$ unless $p\geq 0$ and $q\geq 0$, $\{E^r_{pq}\}$ is called a \textbf{first quadrant spectral sequence}. When $r$ is sufficiently large, $d^r_{pq}: E^r_{pq}\rightarrow E^r_{p-r,q+r-1}$ is the zero map since $E^r_{pq}$=0 or $E^r_{p-r,q+r-1}=0$. Hence $E^r_{pq}=E^{r+1}_{pq}$ for large $r$. Note that every first quadrant spectral sequence is bounded. 

\begin{definition} 
Let $\{E^r_{pq}\}$ be a bounded spectral sequence. Then $\{E^r_{pq}\}$ \textbf{converges} to $H$ (denoted by $E^r_{pq}\Rightarrow H$) if there exists a collection of objects $H_n\in\ob\cat{A}$ and a filtration 
$$0=F_{s}H_n \subseteq \cdots \subseteq F_{p-1}H_n \subseteq F_{p}H_n\subseteq\cdots \subseteq F_{t}H_n=H_n$$
for each $n$ and $$E^{\infty}_{pq}\cong \frac{F_{p}H_n}{F_{p-1}H_n}=\frac{F_{p}H_{p+q}}{F_{p-1}H_{p+q}}$$
for each $p,q\in\mathbb{Z}$ where $p+q=n$.
\end{definition}
We will heavily use the next proposition in this paper.

\begin{proposition}
Let $\{E^r_{pq}\}$ be a bounded spectral sequence where $E^2_{pq}=0$ unless $p=0$ or $p=1$. If $E^r_{pq}\Rightarrow H$, then there exists a short exact sequence 
$$0\rightarrow E^2_{0n}\rightarrow H_n\rightarrow E^2_{1,n-1}\rightarrow 0$$
for each $n$.
\end{proposition}
\begin{proof}
Note that $E^r_{pq}\Rightarrow H$, there exists a filtration of $H$ such that for each $n$,
$$0=F_{s}H_n \subseteq \cdots \subseteq F_{p-1}H_n \subseteq F_{p}H_n\subseteq\cdots \subseteq F_{t}H_n=H_n$$

such that $E^\infty_{pq}\cong \dfrac{F_{p}H_{p+q}}{F_{p-1}H_{p+q}}$. For each $p,q$, $d^2_{pq}: E^2_{pq}\rightarrow E^r_{p-2,q+1}$ is a zero morphism because $E^2_{pq}=0$ unless $p=0$ or $p=1$. Therefore, the $E^2$ page has no no-trivial arrows. Hence $E^2=E^\infty$, and $$E^2_{pq}\cong \frac{F_{p}H_{p+q}}{F_{p-1}H_{p+q}}$$
Note that $n=p+q$, hence 
\begin{itemize}
\item when $p\neq 0$ and $p\neq 1$, $E^2_{pq}=0=\dfrac{F_{p}H_{n}}{F_{p-1}H_{n}}$. Hence $F_{p}H_{n}=F_{p-1}H_{n}$ for all $p\neq 0$ and $p\neq 1$, and the filtration becomes to
$$0=F_sH_n=\cdots=F_{-2}H_n=F_{-1}H_n\subseteq F_{0}H_n\subseteq F_{1}H_n = F_{2}H_n=\cdots=F_{t}H_n=H_n$$
\item when $p=0$, $E^2_{pq}=E^2_{0q}\cong\dfrac{F_{0}H_{n}}{F_{-1}H_{n}}$. Since $F_{-1}H_{n}=0$, we obtain $E^2_{0q}\cong F_{0}H_{0+q}$;
\item when $p=1$, $E^2_{pq}=E^2_{1q}\cong\dfrac{F_{1}H_{n}}{F_{1-1}H_{n}}=\dfrac{H_{n}}{F_{0}H_{n}}$.
\end{itemize}
Note that there is a short exact sequence 
$$0\rightarrow F_0H_n\hookrightarrow H_n\rightarrow H_n/F_0H_n\rightarrow 0$$
Therefore, 
$$0\rightarrow E^2_{0,n}\rightarrow H_n\rightarrow E^2_{1,n-1}\rightarrow 0$$
\end{proof}

\begin{definition}
A \textbf{double complex}\index{double complex} is an ordered triple $(M,d',d'')$, where $M$ is a bigraded module and $d',d''$ are differential maps of bidegree $(-1,0)$ and $(0,-1)$ such that $d'\circ d''+d''\circ d'=0$. 
\end{definition}
In other words, for all $p,q\in\mathbb{Z}$, $$d'_{pq}:M_{p,q}\rightarrow M_{p-1, q}$$ $$d''_{pq}:M_{p,q}\rightarrow M_{p,q-1}$$
\begin{definition}
A \textbf{total complex} of a double complex $(M,d',d'')$ is a complex where its $n$-th term is $\bigoplus\limits_{p+q=n}M_{p,q}$. We use $\mathsf{Tot}(M)_n$ to denote the term $\bigoplus\limits_{p+q=n}M_{p,q}$.
\end{definition}

\begin{lemma}
$\{\mathsf{Tot}(M)_n\}$ is a chain complex, where $d=\sum\limits_{p+q=n} d'_{pq}+d''_{pq}$. 
\end{lemma}

Let $X$ be a CW-complex. Consider a cover $U:=\{X_i\}_{i\in I}$ of $X$ where $I$ is a totally ordered set and $X_i$ is a subcomplex of $X$ for all $i\in I$. Note that for $i_1<i_2<\cdots<i_n$, there is an inclusion map $$C(X_{i_1}\oplus \cdots \oplus X_{i_n})\rightarrow C(X_{i_1}\oplus\cdots \oplus \widehat{X_{i_k}}\oplus\cdots\oplus X_{i_n})$$ where $\widehat{X_{i_k}}$ means we remove the $k$-th direct summand from the given expression. Therefore, it induces a unique map 
$$\bigoplus\limits_{i_1<\cdots<i_n}C(X_{i_1}\oplus \cdots \oplus X_{i_n})\rightarrow \bigoplus\limits_{j_1<\cdots<j_{n-1}}C(X_{j_1}\oplus\cdots\oplus X_{j_{n-1}})$$

Therefore, we obtain a sequence 
\begin{equation}\label{exactseq1}
\begin{aligned}
\cdots\rightarrow \bigoplus\limits_{i_1<\cdots<i_n}C(X_{i_1}\oplus \cdots \oplus X_{i_n})\rightarrow &\bigoplus\limits_{j_1<\cdots<j_{n-1}}C(X_{j_1}\oplus\cdots\oplus X_{j_{n-1}})\rightarrow\cdots\\
&\rightarrow \bigoplus\limits_{i}C(X_{i})\rightarrow C(X)\rightarrow 0
\end{aligned}
\end{equation}
\begin{lemma}
The sequence given in (\ref{exactseq1}) is exact, and $E_{pq}\Rightarrow H_{p+q}(X)$ where $$E_{pq}^2=H_p(E_{pq}^1)=H_p(\bigoplus\limits_{\sigma\in \mathsf{Nerve}(U)} H_q(X_\sigma))$$
\end{lemma}
\begin{proof}
See \cite{brown2012cohomology} pp 166-167.
\end{proof}

%% file: chapters/simplification.tex
\section{A Reduction of $PH_{i}(X^2_{-,-};\mathbb{F})$}
Let $X$ be a finite metric graph. In this section, we construct a persistence module $N$ for each persistence module such that $N$ is naturally isometric to $PH_i(X^2_{-,-};\mathbb{F})$ and the restriction of $N$ to any chamber in the parameter space of $X^2_{r, L_{e_1}}$ is a constant functor. 

\begin{proposition}\label{prop-sim-1}
Let $(J,\leq)$ be a connected poset and $M:(J,\leq)\rightarrow \cat{vect}_\mathbb{F}$ be a persistence module where $M(i\leq j): Mi\rightarrow Mj$ is an isomorphism for all $i\leq j\in P$. Then there exists a persistence module $N:(J,\leq)\rightarrow \cat{vect}_\mathbb{F}$ such that $M\cong N$ and $Ni\rightarrow Nj$ is the identity map for all $i\leq j\in P$.
\end{proposition}
\begin{proof}
Let $\underset{J}{\varinjlim}Mi$ and $\sigma_i:Mi\rightarrow \underset{J}{\varinjlim}Mi$ be the leg map for all $i\in J$. Since for all $y\in \underset{J}{\varinjlim}Mi$, there exists $i\in J$ such that $\sigma_i(x)=y$. Since $M(i\leq j)$ is an isomorphism and $(J,\leq)$ is connected, $\sigma_i$ is surjective. On the other hand, assume there exists $x\in Mi$ such that $\sigma_i(x)=0$. Then there exists $j\geq i\in J$ such that $M(i\leq j)(x)=0$. Because $M(i\leq j)$ is an isomorphism, $x=0$. Therefore, $\sigma_i$ is injective, and $\underset{J}{\varinjlim}Mi\cong Mi$. Since $Mi\cong Mj$ for all $i,j\in J$, $\underset{J}{\varinjlim}Mi\cong Mj$ for all $j\in J$. Note that $\underset{J}{\varinjlim}Mi$ can be viewed as a constant functor $\underset{J}{\varinjlim}Mi:(J,\leq)\rightarrow \cat{vect}_\mathbb{F}$ where $(\underset{J}{\varinjlim}Mi)j=\underset{J}{\varinjlim}Mi$ for all $j\in J$ and $(\underset{J}{\varinjlim}Mi)(j\leq k): (\underset{J}{\varinjlim}Mi)j\rightarrow (\underset{J}{\varinjlim}Mi)k$ is the identity morphism. It is clear that $\underset{J}{\varinjlim}Mi\cong M$. Define $N=\underset{J}{\varinjlim}Mi$. 
\end{proof}

\begin{proposition}\label{prop-sim-2}
Let $(J,\leq)$ be a connected poset and $M:(J,\leq)\rightarrow \cat{vect}_\mathbb{F}$ be a persistence module. Let $(J_1,\leq)$ and $(J_2,\leq)$ be disjoint intervals that contained in $(J,\leq)$ such that $M(i\leq j): Mi\rightarrow Mj$ is an isomorphism for all $i\leq j\in J_k$, $k=1,2$. If there exists $i\in J_1$ and $j\in J_2$ such that $i\leq j$, then there exists a well-defined morphism $\phi_{J_1J_2}:\underset{J_1}{\varinjlim}Mi\rightarrow \underset{J_2}{\varinjlim}Mj$ such that the following diagram commutes, where $\sigma_i$ is the leg map of the colimit cone $\restr{M}{J_1}\Rightarrow \underset{J_1}{\varinjlim}Mi$ and $\tau_j$ is the leg map of the colimit cone $\restr{M}{J_2}\Rightarrow \underset{J_2}{\varinjlim}Mj$
$$
\begin{tikzcd}
  M_i \arrow[r, "M(i\leq j)"] \arrow[d, "\sigma_i"']
    & M_j \arrow[d, "\tau_j"] \\
  \underset{J_1}{\varinjlim}Mi \arrow[r, dashed, "\phi_{J_1J_2}"' ]
&\underset{J_2}{\varinjlim}Mj \end{tikzcd}
$$
\end{proposition}

\begin{proof}
By Proposition \ref{prop-sim-1}, both $\sigma_i$ and $\tau_j$ are isomorphisms. Define $\phi_{J_1J_2}=\tau_j\circ M(i\leq j)\circ \sigma_i^{-1}$. Now we are going to show $\phi_{J_1J_2}$ does not depend on the choice of $\sigma_i$ and $\tau_j$.\\ 

Assume there exists $i'\leq i\in J_1$. Note that 
\begin{equation}
\begin{aligned}
  \tau_j\circ M(i'\leq j)\circ \sigma_{i'}^{-1}&=\tau_j\circ M(i\leq j)\circ M(i'\leq i)\circ \sigma_{i'}^{-1}\\
  &=\tau_j\circ M(i\leq j)\circ \sigma_i^{-1}
\end{aligned}
\end{equation}

On the other hand, assume there exists $i'\in J_1$ such that $i\leq i'\leq j$. Note that 
\begin{equation}
\begin{aligned}
  \tau_j\circ M(i'\leq j)\circ \sigma_{i'}^{-1}&=\tau_j\circ M(i\leq j)\circ M(i\leq i')^{-1}\circ \sigma_{i'}^{-1}\\
  &=\tau_j\circ M(i\leq j)\circ \sigma_i^{-1}
\end{aligned}
\end{equation}

Hence $\phi_{J_1J_2}$ does not depend on the choice of $\sigma_i$ for all $i\in J_1$.

Assume there exists $j\leq j'\in J_2$. Note that 
\begin{equation}
\begin{aligned}
  \tau_{j'}\circ M(i\leq j')\circ \sigma_i^{-1}&=\tau_{j'}\circ M(j\leq j')\circ M(i\leq j)\circ \sigma_i^{-1}\\
  &=\tau_j\circ M(i\leq j)\circ \sigma_i^{-1}
\end{aligned}
\end{equation}

On the other hand, assume there exists $j'\in J_2$ such that $i\leq j'\leq j$. Note that 
\begin{equation}
\begin{aligned}
  \tau_{j'}\circ M(i\leq j')\circ \sigma_i^{-1}&=\tau_{j}\circ M(j'\leq j)\circ M(i\leq j')\circ \sigma_i^{-1}\\
  &=\tau_j\circ M(i\leq j)\circ \sigma_i^{-1}
\end{aligned}
\end{equation}

Hence $\phi_{J_1J_2}$ does not depend on the choice of $\tau_j$ for all $j\in J_2$.
\end{proof}

Let $(X,\delta)$ be a finite metric graph where every edge of $X$ but $e_1$ has length $1$. Since the number of hyperplanes in the parameter space of $X^2_{r,L_{e_1}}$ is finite, the number of chambers in the hyperplane arrangement of $X^2_{r,L_{e_1}}$ is also finite. For each chamber $J$ of the hyperplane arrangement, one can calculate $\underset{J}{\varinjlim}H_i(X^2_{r,L_{e_1}})$. Proposition \ref{prop-sim-1} impies that $\restr{PH_i(X^2_{-,-};\mathbb{F})}{J}\cong \underset{J}{\varinjlim}H_i(X^2_{r,L_{e_1}})$. Using Proposition \ref{prop-sim-2}, we construct a persistence module $N:(\mathbb{R},\leq)^{\op}\times (\mathbb{R},\leq)\rightarrow \cat{vect}_{\mathbb{F}}$ as follows:
\begin{itemize}
    \item For every chamber $J$, $Nj=\underset{J}{\varinjlim}H_i(X^2_{r,L_{e_1}})$ for all $j\in J$;
    \item For all $j_1\leq j_2$, let $J_1$ be the chamber contains $j_1$ and $J_2$ be the chamber contains $j_2$. Then $N(j_1\leq j_2): Nj_1\rightarrow Nj_2$ is
    $$N(j_1\leq j_2)=\begin{cases} id_{\underset{J_1}{\varinjlim}H_i(X^2_{r,L_{e_1}})}, & \mbox{if } J_1=J_2 \\\phi_{J_1J_2} , & \mbox{else } \end{cases}$$
\end{itemize}

Apply proposition \ref{prop-sim-2}, we immediately obtain:

\begin{theorem}\label{thm-sim-1}
Let $(X,\delta)$ be a finite metric graph, and assume every edge of $X$ but $e_1$ has length $1$. Then $$PH_i(X^2_{-,-};\mathbb{F})\cong N$$
\end{theorem}

Note that the functor category $\cat{vect}_{\mathbb{F}}^{(\mathbb{R}_{>0},\leq)^{\op}\times(\mathbb{R}_{>0},\leq)}$ is a Krull-Schmidt category\cite{botnan2020decomposition}. Therefore, to find the in\-decomposable direct summands of $PH_i(X^2_{-,-};\mathbb{F})$, it suffices to find the indecomposable direct summands of $N$. 


%% file: chapters/chapter4-star-nofac.tex
\section{Configuration Spaces of the Star Graphs}\label{conf:star3}
In Section \ref{conf:Y}, we calculated the second configuration space of the Y graph with restraint parameters $r$ and $L_{e_1}$. In Section \ref{star_Hi'}, we calculated the second configuration space of the star graph with restraint parameters $r$ and $L_{e_1}$.

\subsection{Y graph}\label{conf:Y}
Let $Y_{L_{e_1}}$ be a metric graph with the shape of Y where two edges of $Y_{L_{e_1}}$ (labeled by $e_2$ and $e_3$) have length $1$ (i.e., $L_{e_2}=L_{e_3}=1$) and the length of the other edge (labeled by $e_1$) is $L_{e_1}$. Let $r$ denote the minimum distance required between the two robots and let $Y^2_{r,L_{e_1}}$ denote the second configuration space of $Y_{L_{e_1}}$ with restraint parameter $r$. Define $x$ to be the distance between the red robot and the central vertex of Y and $y$ to be the distance between the green robot and the central vertex of Y. The systems of inequalities (\ref{9}) and (\ref{9'}) characterize the cells of $Y^2_{r,L_{e_1}}$.\\

\noindent\begin{minipage}{0.3\linewidth}
\centering
\resizebox{3cm}{!}{
\tikzset{every picture/.style={line width=0.75pt}} 

\begin{tikzpicture}[x=0.75pt,y=0.75pt,yscale=-1,xscale=1]

\draw [line width=2.25]    (147.42,58.58) -- (102.5,103.5) ;
\draw  [fill={rgb, 255:red, 184; green, 233; blue, 150 }  ,fill opacity=1 ] (127,72.75) .. controls (127,69.3) and (129.8,66.5) .. (133.25,66.5) .. controls (136.7,66.5) and (139.5,69.3) .. (139.5,72.75) .. controls (139.5,76.2) and (136.7,79) .. (133.25,79) .. controls (129.8,79) and (127,76.2) .. (127,72.75) -- cycle ;
\draw [line width=2.25]    (67,68.5) -- (102.5,103.5) ;
\draw [line width=2.25]    (102.5,103.5) -- (102.5,161.58) ;
\draw  [fill={rgb, 255:red, 255; green, 16; blue, 0 }  ,fill opacity=1 ] (79,87.25) .. controls (79,83.8) and (81.8,81) .. (85.25,81) .. controls (88.7,81) and (91.5,83.8) .. (91.5,87.25) .. controls (91.5,90.7) and (88.7,93.5) .. (85.25,93.5) .. controls (81.8,93.5) and (79,90.7) .. (79,87.25) -- cycle ;
\draw [line width=2.25]  [dash pattern={on 2.53pt off 3.02pt}]  (48,49.5) -- (67,68.5) ;

\draw (45,68) node [anchor=north west][inner sep=0.75pt]   [align=left] {$e_1$};
\draw (109,136) node [anchor=north west][inner sep=0.75pt]   [align=left] {$e_2$};
\draw (152,53) node [anchor=north west][inner sep=0.75pt]   [align=left] {$e_3$};
\draw (72,93) node [anchor=north west][inner sep=0.75pt]  [color={rgb, 255:red, 255; green, 16; blue, 0}  ,opacity=1 ] [align=left] {$x$};
\draw (138.25,75.75) node [anchor=north west][inner sep=0.75pt]  [color={rgb, 255:red, 126; green, 201; blue, 100 }  ,opacity=1 ] [align=left] {$y$};
\end{tikzpicture}
}
\captionof{figure}{}
\label{Y2}
\end{minipage}%
\noindent\begin{minipage}{0.3\linewidth}
\begin{equation}\label{9}
\begin{aligned}
0&\leq x\leq L_{e_i}\\
0&\leq y\leq L_{e_j}\\
x&+y\geq r
\end{aligned}
\end{equation}
\end{minipage}%
\begin{minipage}{0.3\linewidth}
\begin{equation}\label{9'}
\begin{aligned}
0&\leq x\leq L_{e_i}\\
0&\leq y\leq L_{e_i}\\
\lvert x&-y\rvert\geq r
\end{aligned}
\end{equation}
\end{minipage}

The critical hyperplanes in the parameter space of $Y^2_{r,L_{e_1}}$ are $r=1$, $r=2$, $r=L_{e_1}$, and $r=L_{e_1}+1$. Moreover, the hyperplane arrangement in the parameter space of $Y^2_{r,L_{e_1}}$ is shown in figure \ref{para0}. Note that $Y^2_{r_1,l_1}$ and $Y^2_{r_2,l_2}$ are homotopy equivalent (or even stronger, isotopy equivalent) when the points $(r_1,l_1)$ and $(r_2,l_2)$ are contained in the same chamber.\\

\noindent\begin{minipage}[c]{0.5\textwidth}
\centering
\resizebox{7cm}{!}{

\tikzset{every picture/.style={line width=0.75pt}} 

\begin{tikzpicture}[x=0.75pt,y=0.75pt,yscale=-1,xscale=1]

\draw  [draw opacity=0][fill={rgb, 255:red, 155; green, 155; blue, 155 }  ,fill opacity=0.4 ] (138.9,52.47) -- (237.3,51.32) -- (236.79,104.89) -- (139.57,203.47) -- (138.9,106.03) -- cycle ;
\draw  [draw opacity=0][fill={rgb, 255:red, 155; green, 155; blue, 155 }  ,fill opacity=0.4 ] (138.76,203.83) -- (235.65,106.03) -- (236.05,203.42) -- cycle ;
\draw  (6.75,299.04) -- (373.46,299.04)(43.42,46.75) -- (43.42,327.07) (366.46,294.04) -- (373.46,299.04) -- (366.46,304.04) (38.42,53.75) -- (43.42,46.75) -- (48.42,53.75)  ;
\draw  [draw opacity=0][fill={rgb, 255:red, 155; green, 155; blue, 155 }  ,fill opacity=0.4 ] (43.42,299.04) -- (139.57,203.47) -- (139.82,298.78) -- cycle ;
\draw  [draw opacity=0][fill={rgb, 255:red, 155; green, 155; blue, 155 }  ,fill opacity=0.4 ] (139.82,298.78) -- (235.28,202.82) -- (235.53,298.53) -- cycle ;
\draw  [draw opacity=0][fill={rgb, 255:red, 155; green, 155; blue, 155 }  ,fill opacity=0.4 ] (289.7,50.84) -- (236.79,104.89) -- (236.25,51.39) -- cycle ;
\draw  [draw opacity=0][fill={rgb, 255:red, 155; green, 155; blue, 155 }  ,fill opacity=0.4 ] (331.68,106.33) -- (234.91,204.42) -- (235.48,105.78) -- cycle ;
\draw  [draw opacity=0][fill={rgb, 255:red, 155; green, 155; blue, 155 }  ,fill opacity=0.4 ] (289.7,50.7) -- (386.49,51.53) -- (330.94,106.22) -- (331.68,106.33) -- (313,106) -- (234.5,105.89) -- cycle ;
\draw  [draw opacity=0][fill={rgb, 255:red, 155; green, 155; blue, 155 }  ,fill opacity=0.4 ] (138.9,52.47) -- (138.9,106.03) -- (139.57,203.47) -- (43.42,299.04) -- (43.93,52.47) -- cycle ;
\draw  [draw opacity=0][fill={rgb, 255:red, 155; green, 155; blue, 155 }  ,fill opacity=0.4 ] (236.88,203.18) -- (140.21,298.72) -- (139.58,203.82) -- cycle ;
\draw [line width=1.5]    (43.42,299.04) -- (289.7,50.84) ;
\draw [line width=1.5]    (140.21,299.72) -- (386.49,51.53) ;
\draw [line width=1.5]    (140.21,299.72) -- (138.9,52.47) ;
\draw [line width=1.5]    (235.53,298.53) -- (234.22,51.28) ;

\draw (139.47,313.34) node   [align=left] {1};
\draw (236.73,312.19) node   [align=left] {2};
\draw (35.35,203.5) node   [align=left] {1};
\draw (34.21,108.53) node   [align=left] {2};
\draw (17.07,17.96) node [anchor=north west][inner sep=0.75pt]    {$L_{e_{1}}$};
\draw (376.76,302.52) node [anchor=north west][inner sep=0.75pt]    {$r$};

\end{tikzpicture}
}

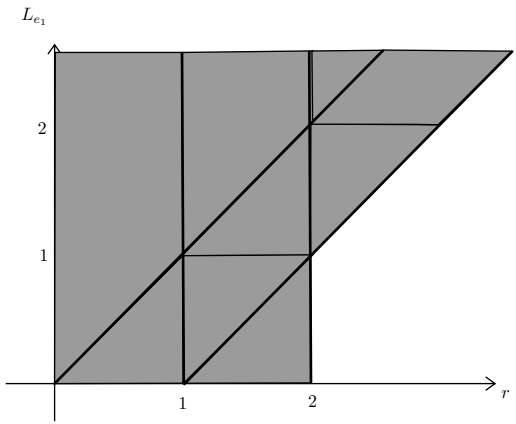
\captionof{figure}{The hyperplane arrangement associated to $Y^2_{r,L_{e_1}}$}
\label{para0}
\end{minipage}
\noindent\begin{minipage}[c]{0.5\textwidth}
\centering
\resizebox{7cm}{!}{

\tikzset{every picture/.style={line width=0.75pt}} 

\begin{tikzpicture}[x=0.75pt,y=0.75pt,yscale=-1,xscale=1]

\draw  [draw opacity=0][fill={rgb, 255:red, 155; green, 155; blue, 155 }  ,fill opacity=1 ] (138.9,52.47) -- (237.3,51.32) -- (236.79,104.89) -- (139.57,203.47) -- (138.9,106.03) -- cycle ;
\draw  [draw opacity=0][fill={rgb, 255:red, 248; green, 231; blue, 28 }  ,fill opacity=1 ] (138.76,203.83) -- (235.65,106.03) -- (236.05,203.42) -- cycle ;
\draw  (6.75,299.04) -- (373.46,299.04)(43.42,46.75) -- (43.42,327.07) (366.46,294.04) -- (373.46,299.04) -- (366.46,304.04) (38.42,53.75) -- (43.42,46.75) -- (48.42,53.75)  ;
\draw  [draw opacity=0][fill={rgb, 255:red, 245; green, 166; blue, 35 }  ,fill opacity=1 ] (43.42,299.04) -- (139.57,203.47) -- (139.82,298.78) -- cycle ;
\draw  [draw opacity=0][fill={rgb, 255:red, 245; green, 166; blue, 35 }  ,fill opacity=1 ] (140.81,298.68) -- (236.27,202.72) -- (236.52,298.43) -- cycle ;
\draw  [draw opacity=0][fill={rgb, 255:red, 245; green, 166; blue, 35 }  ,fill opacity=1 ] (289.7,50.84) -- (236.79,104.89) -- (236.25,51.39) -- cycle ;
\draw  [draw opacity=0][fill={rgb, 255:red, 155; green, 155; blue, 155 }  ,fill opacity=1 ] (331.68,105.33) -- (234.91,203.42) -- (235.48,104.78) -- cycle ;
\draw  [draw opacity=0][fill={rgb, 255:red, 155; green, 155; blue, 155 }  ,fill opacity=1 ] (289.7,50.7) -- (380.32,47.75) -- (330.94,106.22) -- (331.68,106.33) -- (307.51,106.22) -- (234.5,105.89) -- cycle ;
\draw  [draw opacity=0][fill={rgb, 255:red, 80; green, 227; blue, 194 }  ,fill opacity=1 ] (138.9,52.47) -- (138.9,106.03) -- (139.57,203.47) -- (43.42,299.04) -- (43.93,52.47) -- cycle ;
\draw  [draw opacity=0][fill={rgb, 255:red, 248; green, 231; blue, 28 }  ,fill opacity=1 ] (236.86,201.68) -- (140.21,299.72) -- (139.57,202.33) -- cycle ;
\draw  [fill={rgb, 255:red, 248; green, 231; blue, 28 }  ,fill opacity=1 ] (305.53,202.83) .. controls (305.53,199.51) and (308.21,196.83) .. (311.53,196.83) .. controls (314.84,196.83) and (317.53,199.51) .. (317.53,202.83) .. controls (317.53,206.14) and (314.84,208.83) .. (311.53,208.83) .. controls (308.21,208.83) and (305.53,206.14) .. (305.53,202.83) -- cycle ;
\draw  [fill={rgb, 255:red, 245; green, 166; blue, 35 }  ,fill opacity=1 ] (306.5,243.83) .. controls (306.5,240.51) and (309.19,237.83) .. (312.5,237.83) .. controls (315.81,237.83) and (318.5,240.51) .. (318.5,243.83) .. controls (318.5,247.14) and (315.81,249.83) .. (312.5,249.83) .. controls (309.19,249.83) and (306.5,247.14) .. (306.5,243.83) -- cycle ;
\draw  [fill={rgb, 255:red, 80; green, 227; blue, 194 }  ,fill opacity=1 ] (306.38,223.83) .. controls (306.38,220.51) and (309.06,217.83) .. (312.38,217.83) .. controls (315.69,217.83) and (318.38,220.51) .. (318.38,223.83) .. controls (318.38,227.14) and (315.69,229.83) .. (312.38,229.83) .. controls (309.06,229.83) and (306.38,227.14) .. (306.38,223.83) -- cycle ;
\draw  [fill={rgb, 255:red, 155; green, 155; blue, 155 }  ,fill opacity=1 ] (307.38,265.83) .. controls (307.38,262.51) and (310.06,259.83) .. (313.38,259.83) .. controls (316.69,259.83) and (319.38,262.51) .. (319.38,265.83) .. controls (319.38,269.14) and (316.69,271.83) .. (313.38,271.83) .. controls (310.06,271.83) and (307.38,269.14) .. (307.38,265.83) -- cycle ;

\draw (139.47,313.34) node   [align=left] {1};
\draw (236.73,312.19) node   [align=left] {2};
\draw (35.35,203.5) node   [align=left] {1};
\draw (34.21,108.53) node   [align=left] {2};
\draw (17.07,17.96) node [anchor=north west][inner sep=0.75pt]    {$L_{e_{1}}$};
\draw (376.76,302.52) node [anchor=north west][inner sep=0.75pt]    {$r$};
\draw (351,204.27) node   [align=left] {6 points};
\draw (352,244.27) node   [align=left] {2 points};
\draw (349,224.27) node   [align=left] {1 circle};
\draw (353,266.27) node   [align=left] {4 points};

\end{tikzpicture}
}

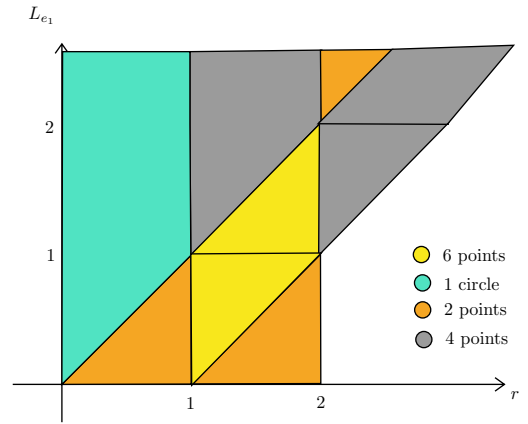
\captionof{figure}{Homotopy type of $Y^2_{r,L_{e_1}}$}
\label{chap1:para}
\end{minipage}

\subsection{$\mathsf{Star}_{k}$, $k\geq 4$}\label{star_Hi'}
Let $\hat{\mathsf{Star}}_{k}$ be the metric graph where every edge but $e_1$ has length $1$. We begin this section with an observation.
\begin{proposition}
Let $k\geq 4$ and $k\in\mathbb{Z}$. Then $(\hat{\mathsf{Star}}_{k})^2_{r,{L_{e_1}}}$ and $Y^2_{r,L_{e_1}}:=(\hat{\mathsf{Star}}_{3})^2_{r,L_{e_1}}$ have the same critical hyperplanes. 
\end{proposition}
\begin{proof}
Let $L_{e_i}$ denote the edge length of $e_i$. Hence $$L_{e_i} = \begin{cases} L_{e_1}, & \mbox{if } i=1 \\ 1, & \mbox{if } i\in\{2,\dots, n\} \end{cases}
$$

Define $x$ to be the distance between the red robot and the central vertex of the star graph and $y$ to be the distance between the green robot and the central vertex of the star graph. The following inequalities (\ref{star1})-(\ref{star1''}) gives the parametric polytope of $(\hat{\mathsf{Star}}_{k})^2_{r,L_{e_1}}$:\\

For $i,j\in\{1,\dots, k\}$,
\begin{equation}\label{star1}
\begin{aligned}
0&\leq x\leq L_{e_i}\\
0&\leq y\leq L_{e_j}\\
\end{aligned}
\end{equation}

For $i,j\in\{1,\dots, k\}$ and $i\neq j$,
\begin{equation}
\begin{aligned}
x&+y\geq r
\end{aligned}
\end{equation}



For $i=j \in\{1,\dots, k\}$,
\begin{equation}\label{star1''}
\begin{aligned}
\lvert x&-y\rvert\geq r
\end{aligned}
\end{equation}

Hence the critical hyperplanes are $r=1$, $r=2$, $r=L_{e_1}$, and $r=L_{e_1}+1$. Therefore, $(\hat{\mathsf{Star}}_{k})^2_{r,{L_{e_1}}}$ and $Y^2_{r,L_{e_1}}$ have the same critical hyperplanes. 
\end{proof}

Let $\mathsf{Star}_{k}$ denote the metric graph where every edge has length $1$.
\begin{proposition}
Let $r$ be a positive number and $k\geq 3$. Then
$$H_0((\mathsf{Star}_k)_r^2) \cong \begin{cases} \mathbb{Z}, & \mbox{if } 0<r\leq 1 \\ \mathbb{Z}^{k^2-k}, & \mbox{if } 1<r\leq 2\\ 0, & \mbox{if } 2<r\end{cases} \mbox{\quad and\quad} H_1((\mathsf{Star}_k)_r^2) \cong \begin{cases} \mathbb{Z}^{k(k-3)+1}, & \mbox{if } 0<r\leq 1 \\ 0, & \mbox{if } 1<r\end{cases} $$
\end{proposition}
\begin{proof}
We run induction on $k$. 
\begin{description}
\item[When $k=3$ ] Consider the following cover of $(\mathrm{Star}_3)_r^2$:
\begin{equation}
\begin{aligned}
U_{11}&=(e_1)^2_r\\
U_{12}&=e_1\times \mathsf{Star}_{2}-\{(x,y)\in e_1\times \mathsf{Star}_{2}\mid \delta(x,y)<r\}\\
U_{21}&=\mathsf{Star}_{2}\times e_1-\{(x,y)\in \mathsf{Star}_{2}\times e_1\mid \delta(x,y)<r\}\\
U_{22}&=(\mathsf{Star}_{2})_r^2
\end{aligned}
\end{equation}
Note that $(e_1)^2_r\simeq  \begin{cases} \{\ast_1,\ast_2\}, & \mbox{if } 0<r\leq 1 \\ \emptyset, & \mbox{if } 1<r\end{cases}$ where $\{\ast_1,\ast_2\}$ is a subspace of $(e_1)^2_r$ consists of two points. On the other hand, $$e_1\times \mathsf{Star}_{2}-\{(x,y)\in e_1\times \mathsf{Star}_{2}\mid \delta(x,y)<r\}\simeq  \begin{cases} \{\ast\} , & \mbox{if } 0<r\leq 1 \\ \{\ast_1,\ast_2\}, & \mbox{if } 1<r\leq 2\\ \emptyset, & \mbox{if } 2<r\end{cases}$$
Similarly,  $$\mathsf{Star}_{2}\times e_1-\{(x,y)\in \mathsf{Star}_{2}\times e_1 \mid \delta(x,y)<r\}\simeq  \begin{cases} \{\ast\} , & \mbox{if } 0<r\leq 1 \\ \{\ast_1,\ast_2\}, & \mbox{if } 1<r\leq 2\\ \emptyset, & \mbox{if } 2<r\end{cases}$$
In addition, note that $$H_0((\mathsf{Star}_2)^2_r)\cong \begin{cases} \mathbb{Z}^2, & \mbox{if } 0<r\leq 2 \\ 0, & \mbox{if } 2<r\end{cases}$$

Now let's consider the intersections of $U_{ij}$. Note that 
 \begin{equation}
\begin{aligned}
U_{11}\cap U_{12}&\simeq U_{11}\cap U_{21}\simeq \begin{cases} \{\ast\} , & \mbox{if } 0<r\leq 1 \\ \emptyset, & \mbox{if } 1<r\end{cases}\\
U_{22}\cap U_{12}&\simeq U_{22}\cap U_{21}\simeq \begin{cases} \{\ast_1,\ast_{2}\} , & \mbox{if } 0<r\leq 1 \\ \emptyset, & \mbox{if } 1<r\end{cases}\\
U_{11}\cap U_{22}&=U_{12}\cap U_{21}=\emptyset
\end{aligned}
\end{equation}
Note that there are natural inclusions 
\begin{equation}\label{star3:eq}
\begin{aligned}
U_{11}\cap U_{12}&\hookrightarrow U_{11}\hookleftarrow U_{11}\cap U_{21}\\
U_{11}\cap U_{12}&\hookrightarrow U_{12}\hookleftarrow U_{22}\cap U_{12}\\
U_{11}\cap U_{21}&\hookrightarrow U_{21}\hookleftarrow U_{22}\cap U_{21}\\
U_{22}\cap U_{12}&\hookrightarrow U_{22}\hookleftarrow U_{22}\cap U_{21}
\end{aligned}
\end{equation}

 \paragraph{When $0<r\leq 1$:}
The $E^1$ page of the Mayer-Vietoris spectral sequence is 
\begin{figure}[htbp!]
\centering
\begin{tikzpicture}
  \matrix (m) [matrix of math nodes,
    nodes in empty cells,nodes={minimum width=5ex,
     minimum height=2.5ex,outer sep=-5pt},
    column sep=1ex,row sep=1ex]{
          q     &      &     &  &  \\
          2     & \vdots & 0 & 0 &\\   
          1     &  H_1((\mathsf{Star}_{2})^2_r) &  0  & 0 \\
          0     &  \mathbb{Z}^4\oplus H_0((\mathsf{Star}_{2})^2_r)  & \mathbb{Z}^{6} & 0\\
    \quad\strut &   0  &  1  &  2 & p\strut \\};
    \draw[-stealth] (m-4-3.west) -- (m-4-2.east) node [above,pos=0.5] {$d^1$};
\draw[thick] (m-1-1.east) -- (m-5-1.east) ;
\draw[thick] (m-5-1.north) -- (m-5-5.north) ;
\end{tikzpicture}
\end{figure}

Let $a$ be the generator of $H_0(U_{11}\cap U_{12})$ and $\hat{a}$ be the generator of $H_0(U_{11}\cap U_{21})$. Let $b_1,b_{2}$ be the generators of $H_0(U_{22}\cap U_{12})$ and $\hat{b}_1,\hat{b}_{2}$ be the generators of $H_0(U_{22}\cap U_{21})$. Note that $U_{11}$ has two path components, while $U_{11}\cap U_{12}$ and $U_{11}\cap U_{21}$ lie in different path components. Choose a generator (denoted by $c_1$) from the path component of $U_{11}$ containing $U_{11}\cap U_{12}$, and choose a generator (denoted by $c_2$) from the path component of $U_{11}$ containing $U_{11}\cap U_{21}$. Let $e$ be the generator of $H_0(U_{12})$ and $\hat{e}$ be the generator of $H_0(U_{21})$ such that $$d^1(a)=c_1+e \mbox{\quad and\quad} d^1(\hat{a})=c_2+\hat{e}$$
On the other hand, let $f_1, f_2$ be the generators of $H_0(U_{22})\cong\mathbb{Z}^2$, where $f_1$ and $f_2$ lie in different path components of $U_{22}$, then $$d^1(b_1)=e+f_1, \quad d^1(b_2)=e+f_2,\quad d^1(\hat{b}_1)=\hat{e}+f_1, \quad d^1(\hat{b}_2)=\hat{e}+f_2$$
Under the given generators, $d^1$ can be represented by the following matrix
$$\left[\begin{array}{cccccc}1 & 0 & 0 & 0 & 0 & 0 \\0 & 1 & 0 & 0 & 0 & 0 \\1 & 0 & 1 & 1 & 0 & 0 \\0 & 1 & 0 & 0 & 1 & 1 \\0 & 0 & 1 & 0 & 1 & 0 \\0 & 0 & 0 & 1 & 0 & 1\end{array}\right]$$
The rank of this matrix is $5$, and the diagonal elements of the Smith Normal form of the above matrix are all equal to $1$. Therefore, $\Ima d^1 \cong\mathbb{Z}^5$ and $\ker d^1\cong \mathbb{Z}^{1}$. Hence we obtain the $E^2$ page of the Mayer-Vietoris spectral sequence: 
\begin{figure}[htbp!]
\centering
\begin{tikzpicture}
  \matrix (m) [matrix of math nodes,
    nodes in empty cells,nodes={minimum width=5ex,
     minimum height=2.5ex,outer sep=-5pt},
    column sep=1ex,row sep=1ex]{
          q     &      &     &  &  \\
          2     & \vdots & 0 & 0 &\\   
          1     &  H_1((\mathsf{Star}_{2})^2_r) &  0  & 0 \\
          0     &  \mathbb{Z}  & \mathbb{Z} & 0\\
    \quad\strut &   0  &  1  &  2 & p\strut \\};
\draw[thick] (m-1-1.east) -- (m-5-1.east) ;
\draw[thick] (m-5-1.north) -- (m-5-5.north) ;
\end{tikzpicture}
\end{figure}

Since there is no non-trivial arrow on the $E^2$ page, $E^2=E^{\infty}$. Therefore, $H_0((\mathsf{Star}_{3})^2_r)\cong\mathbb{Z}$ and $H_1((\mathsf{Star}_{3})^2_r)\cong H_1((\mathsf{Star}_{2})^2_r)\oplus\mathbb{Z}\cong\mathbb{Z}$.

\paragraph{When $1<r\leq 2$:}
The $E^1$ page of the Mayer-Vietoris spectral sequence is 
\begin{figure}[htbp!]
\centering
\begin{tikzpicture}
  \matrix (m) [matrix of math nodes,
    nodes in empty cells,nodes={minimum width=5ex,
     minimum height=2.5ex,outer sep=-5pt},
    column sep=1ex,row sep=1ex]{
          q     &      &     &  &  \\
          2     & \vdots & 0 & 0 &\\   
          1     &  H_1((\mathsf{Star}_{2})^2_r) &  0  & 0 \\
          0     &  \mathbb{Z}^6 & 0 & 0\\
    \quad\strut &   0  &  1  &  2 & p\strut \\};
\draw[thick] (m-1-1.east) -- (m-5-1.east) ;
\draw[thick] (m-5-1.north) -- (m-5-5.north) ;
\end{tikzpicture}
\end{figure}

Since there is no non-trivial arrow on the $E^1$ page, $E^1=E^{\infty}$. Therefore, $H_0((\mathsf{Star}_{3})^2_r)\cong\mathbb{Z}^6$ and $H_1((\mathsf{Star}_{3})^2_r)=0$.

In conclusion, 
\begin{equation}
 H_0((\mathsf{Star}_3)^2_r)\cong\begin{cases} \mathbb{Z} , & \mbox{if } 0<r\leq 1 \\
 \mathbb{Z}^{6}, & \mbox{if } 1<r\leq 2\\ 
0, & \mbox{if } 2<r\end{cases}
\end{equation}
and 
\begin{equation}
 H_1((\mathsf{Star}_3)^2_r)\cong\begin{cases} \mathbb{Z} , & \mbox{if } 0<r\leq 1 \\
0, & \mbox{if } 1<r\end{cases}
\end{equation}

\item[When $k\geq 4$ ]
Consider the following cover of $(\mathrm{Star}_k)_r^2$:
\begin{equation}
\begin{aligned}
U_{11}&=(e_1)^2_r\\
U_{12}&=e_1\times \mathsf{Star}_{k-1}-\{(x,y)\in e_1\times \mathsf{Star}_{k-1}\mid \delta(x,y)<r\}\\
U_{21}&=\mathsf{Star}_{k-1}\times e_1-\{(x,y)\in \mathsf{Star}_{k-1}\times e_1\mid \delta(x,y)<r\}\\
U_{22}&=(\mathsf{Star}_{k-1})_r^2
\end{aligned}
\end{equation}
Note that $(e_1)^2_r\simeq  \begin{cases} \{\ast_1,\ast_2\}, & \mbox{if } 0<r\leq 1 \\ \emptyset, & \mbox{if } 1<r\end{cases}$ where $\{\ast_1,\ast_2\}$ is a subspace of $(e_1)^2_r$ consists of two points. On the other hand, $$e_1\times \mathsf{Star}_{k-1}-\{(x,y)\in e_1\times \mathsf{Star}_{k-1}\mid \delta(x,y)<r\}\simeq  \begin{cases} \{\ast\} , & \mbox{if } 0<r\leq 1 \\ \{\ast_1,\ast_2,\dots, \ast_{k-1}\}, & \mbox{if } 1<r\leq 2\\ \emptyset, & \mbox{if } 2<r\end{cases}$$
 where $\{\ast_1,\ast_2, \dots, \ast_{k-1}\}$ is a subspace of $(e_1)^2_r$ consists of $(k-1)$ points.\\

Now let's consider the intersections of $U_{ij}$. Note that 
 \begin{equation}
\begin{aligned}
U_{11}\cap U_{12}&\simeq U_{11}\cap U_{21}\simeq  \begin{cases} \{\ast\} , & \mbox{if } 0<r\leq 1 \\ \emptyset, & \mbox{if } 1<r\end{cases}\\
U_{22}\cap U_{12}&\simeq U_{22}\cap U_{21}\simeq  \begin{cases} \{\ast_1,\dots, \ast_{k-1}\} , & \mbox{if } 0<r\leq 1 \\ \emptyset, & \mbox{if } 1<r\end{cases}\\
U_{11}\cap U_{22}&=U_{12}\cap U_{21}=\emptyset
\end{aligned}
\end{equation}
Note that there are natural inclusions 
\begin{equation}\label{star:eq}
\begin{aligned}
U_{11}\cap U_{12}&\hookrightarrow U_{11}\hookleftarrow U_{11}\cap U_{21}\\
U_{11}\cap U_{12}&\hookrightarrow U_{12}\hookleftarrow U_{22}\cap U_{12}\\
U_{11}\cap U_{21}&\hookrightarrow U_{21}\hookleftarrow U_{22}\cap U_{21}\\
U_{22}\cap U_{12}&\hookrightarrow U_{22}\hookleftarrow U_{22}\cap U_{21}
\end{aligned}
\end{equation}
 \paragraph{When $0<r\leq 1$:}
The $E^1$ page of the Mayer-Vietoris spectral sequence is 
\begin{figure}[htbp!]
\centering
\begin{tikzpicture}
  \matrix (m) [matrix of math nodes,
    nodes in empty cells,nodes={minimum width=5ex,
     minimum height=2.5ex,outer sep=-5pt},
    column sep=1ex,row sep=1ex]{
          q     &      &     &  &  \\
          2     & \vdots & 0 & 0 &\\   
          1     &  H_1((\mathsf{Star}_{k-1})^2_r) &  0  & 0 \\
          0     &  \mathbb{Z}^4\oplus H_0((\mathsf{Star}_{k-1})^2_r)  & \mathbb{Z}^{2k} & 0\\
    \quad\strut &   0  &  1  &  2 & p\strut \\};
    \draw[-stealth] (m-4-3.west) -- (m-4-2.east) node [above,pos=0.5] {$d^1$};
\draw[thick] (m-1-1.east) -- (m-5-1.east) ;
\draw[thick] (m-5-1.north) -- (m-5-5.north) ;
\end{tikzpicture}
\end{figure}

Let $a$ be the generator of $H_0(U_{11}\cap U_{12})$ and $\hat{a}$ be the generator of $H_0(U_{11}\cap U_{21})$. Let $b_1,\dots, b_{k-1}$ be the generators of $H_0(U_{22}\cap U_{12})$ and $\hat{b}_1,\dots, \hat{b}_{k-1}$ be the generators of $H_0(U_{22}\cap U_{21})$. Note that $U_{11}$ has two path components, while $U_{11}\cap U_{12}$ and $U_{11}\cap U_{21}$ lie in different path components. Choose a generator (denoted by $c_1$) from the path component of $U_{11}$ containing $U_{11}\cap U_{12}$, and choose a generator (denoted by $c_2$) from the path component of $U_{11}$ containing $U_{11}\cap U_{21}$. Let $e$ be the generator of $H_0(U_{12})$ and $\hat{e}$ be the generator of $H_0(U_{21})$ such that $$d^1(a)=c_1+e \mbox{\quad and \quad} d^1(\hat{a})=c_2+\hat{e}$$
Let $f$ be the generator of $H_0(U_{22})\cong\mathbb{Z}$, then for all $i=1,\dots, k-1$, $$d^1(b_i)=e+f\mbox{\quad and \quad}d^1(\hat{b}_i)=\hat{e}+f$$
Under the given generators, $d^1$ can be represented by the following matrix
$$\left[\begin{array}{cccccc}1 & 0 & 0 & 0 & \cdots & 0 \\0 & 1 & 0 & 0 & \cdots & 0 \\1 & 0 & 1 & 0 & \cdots & 0 \\0 & 1 & 0 & 1 & \cdots & 1 \\0 & 0 & 1 & 1 & \cdots & 1\end{array}\right]$$
where the third and the fourth columns are repeated for $(k-1)$ times. The rank of this matrix is $4$, and the diagonal elements of the Smith Normal form of the above matrix are all equal to $1$. Therefore, $\Ima d^1 \cong\mathbb{Z}^4$ and $\ker d^1\cong \mathbb{Z}^{2k-4}$. Hence we obtain the $E^2$ page of the Mayer-Vietoris spectral sequence: 
\begin{figure}[htbp!]
\centering
\begin{tikzpicture}
  \matrix (m) [matrix of math nodes,
    nodes in empty cells,nodes={minimum width=5ex,
     minimum height=2.5ex,outer sep=-5pt},
    column sep=1ex,row sep=1ex]{
          q     &      &     &  &  \\
          2     & \vdots & 0 & 0 &\\   
          1     &  H_1((\mathsf{Star}_{k-1})^2_r) &  0  & 0 \\
          0     &  \mathbb{Z}  & \mathbb{Z}^{2k-4} & 0\\
    \quad\strut &   0  &  1  &  2 & p\strut \\};
\draw[thick] (m-1-1.east) -- (m-5-1.east) ;
\draw[thick] (m-5-1.north) -- (m-5-5.north) ;
\end{tikzpicture}
\end{figure}

Since there is no non-trivial arrow on the $E^2$ page, $E^2=E^{\infty}$. Therefore, we obtain $H_0((\mathsf{Star}_{k})^2_r)\cong\mathbb{Z}$ and the following sequence is a short exact sequence:
\begin{equation}\label{ses:1}
0\rightarrow H_1((\mathsf{Star}_{k-1})^2_r)\rightarrow H_1((\mathsf{Star}_{k})^2_r)\rightarrow \mathbb{Z}^{2k-4}\rightarrow 0
\end{equation}
Note that $\mathbb{Z}^{2k-4}$ is free, the short exact sequence (\ref{ses:1}) splits, hence 
\begin{equation}
\begin{aligned}
H_1((\mathsf{Star}_{k})^2_r)&\cong H_1((\mathsf{Star}_{k-1})^2_r)\oplus \mathbb{Z}^{2k-4}\cong \mathbb{Z}^{2k-4} \oplus \mathbb{Z}^{2(k-1)-4}\oplus \cdots \oplus \mathbb{Z}^{2(4)-4}\oplus \mathbb{Z}\cong \mathbb{Z}^{k(k-3)+1}
\end{aligned}
\end{equation}
 \paragraph{When $1<r\leq 2$:}
The $E^1$ page of the Mayer-Vietoris spectral sequence is 
\begin{figure}[htbp!]
\centering
\begin{tikzpicture}
  \matrix (m) [matrix of math nodes,
    nodes in empty cells,nodes={minimum width=5ex,
     minimum height=2.5ex,outer sep=-5pt},
    column sep=1ex,row sep=1ex]{
          q     &      &     &  &  \\
          2     & \vdots & 0 & 0 &\\   
          1     &  H_1((\mathsf{Star}_{k-1})^2_r) &  0  & 0 \\
          0     &  \mathbb{Z}^{2k-2}\oplus H_0((\mathsf{Star}_{k-1})^2_r)  & 0 & 0\\
    \quad\strut &   0  &  1  &  2 & p\strut \\};
    \draw[-stealth] (m-4-3.west) -- (m-4-2.east) node [above,pos=0.5] {$d^1$};
\draw[thick] (m-1-1.east) -- (m-5-1.east) ;
\draw[thick] (m-5-1.north) -- (m-5-5.north) ;
\end{tikzpicture}
\end{figure}

Since there is no non-trivial arrow on the $E^1$ page, $E^1=E^{\infty}$. Therefore, 
\begin{equation}
\begin{aligned}
H_0((\mathsf{Star}_{k})^2_r)&\cong H_0((\mathsf{Star}_{k-1})^2_r)\oplus\mathbb{Z}^{2k-2}\cong \mathbb{Z}^{2k-2}\oplus \mathbb{Z}^{2(k-1)-2}\oplus \cdots \mathbb{Z}^{2(4)-2}\oplus \mathbb{Z}^6\cong \mathbb{Z}^{k^2-k}  
\end{aligned}
\end{equation}
and
$$H_1((\mathsf{Star}_{k})^2_r)\cong H_1((\mathsf{Star}_{k-1})^2_r)$$
\end{description}
In conclusion, 
\begin{equation}
 H_0((\mathsf{Star}_k)^2_r)\cong\begin{cases} \mathbb{Z} , & \mbox{if } 0<r\leq 1 \\
 \mathbb{Z}^{k^2-k}, & \mbox{if } 1<r\leq 2\\ 
0, & \mbox{if } 2<r\end{cases}
\end{equation}
and 
\begin{equation}
 H_1((\mathsf{Star}_k)^2_r)\cong\begin{cases} \mathbb{Z}^{k(k-3)+1} , & \mbox{if } 0<r\leq 1 \\
0, & \mbox{if } 1<r\end{cases}
\end{equation}
\end{proof}

Now we are ready to compute $H_i((\hat{\mathsf{Star}}_k)_r^2)$. Consider the following cover of $(\hat{\mathsf{Star}}_k)_r^2$:
\begin{equation}\label{cover:star}
\begin{aligned}
U_{11}&=(e_1)^2_r\\
U_{12}&=e_1\times \mathsf{Star}_{k-1}-\{(x,y)\in e_1\times \mathsf{Star}_{k-1}\mid \delta(x,y)<r\}\\
U_{21}&=\mathsf{Star}_{k-1}\times e_1-\{(x,y)\in \mathsf{Star}_{k-1}\times e_1\mid \delta(x,y)<r\}\\
U_{22}&=(\mathsf{Star}_{k-1})_r^2
\end{aligned}
\end{equation}

Note that there are natural inclusions 
\begin{equation}\label{star_k:eq}
\begin{aligned}
U_{11}\cap U_{12}&\hookrightarrow U_{11}\hookleftarrow U_{11}\cap U_{21}\\
U_{11}\cap U_{12}&\hookrightarrow U_{12}\hookleftarrow U_{22}\cap U_{12}\\
U_{11}\cap U_{21}&\hookrightarrow U_{21}\hookleftarrow U_{22}\cap U_{21}\\
U_{22}\cap U_{12}&\hookrightarrow U_{22}\hookleftarrow U_{22}\cap U_{21}
\end{aligned}
\end{equation}

\begin{proposition}
Let $L_{e_1}$ be a positive number. If $r\leq L_{e_1}$ and $r\leq 1$ . Then
$$H_0((\hat{\mathsf{Star}}_k)_r^2) \cong \mathbb{Z}  \mbox{\quad and\quad} H_1((\hat{\mathsf{Star}}_k)_r^2)\cong \mathbb{Z}^{k(k-3)+1}$$
\end{proposition}
\begin{proof}
Consider the cover of $(\hat{\mathsf{Star}}_k)_r^2$ given in Equation (\ref{cover:star}). Note that 
 \begin{equation}
\begin{aligned}
U_{11}&=(e_1)^2_r\simeq  \{\ast_1,\ast_2\}\\
U_{12}&=e_1\times \mathsf{Star}_{k-1}-\{(x,y)\in e_1\times \mathsf{Star}_{k-1}\mid \delta(x,y)<r\}\simeq  \{\ast\}\\
U_{21}&=\mathsf{Star}_{k-1}\times e_1-\{(x,y)\in \mathsf{Star}_{k-1}\times e_1\mid \delta(x,y)<r\}\simeq  \{\ast\}
\end{aligned}
\end{equation}

 Now let's consider the intersections of $U_{ij}$. Note that 
 \begin{equation}
\begin{aligned}
U_{11}\cap U_{12}&\simeq\{\ast\}\simeq U_{11}\cap U_{21}\\
U_{22}\cap U_{12}&\simeq  \{\ast_1,\dots, \ast_{k-1}\}\simeq U_{22}\cap U_{21}\\
U_{11}\cap U_{22}&=\emptyset=U_{12}\cap U_{21}
\end{aligned}
\end{equation}

The $E^1$ page of the Mayer-Vietoris spectral sequence is\\

\begin{figure}[htbp!]
\centering
\begin{tikzpicture}
  \matrix (m) [matrix of math nodes,
    nodes in empty cells,nodes={minimum width=5ex,
     minimum height=2.5ex,outer sep=-5pt},
    column sep=1ex,row sep=1ex]{
          q     &      &     &  &  \\
          2     & \vdots & 0 & 0 &\\   
          1     &  H_1((\mathsf{Star}_{k-1})^2_r) &  0  & 0 \\
          0     &  \mathbb{Z}^4\oplus H_0((\mathsf{Star}_{k-1})^2_r)  & \mathbb{Z}^{2k} & 0\\
    \quad\strut &   0  &  1  &  2 & p\strut \\};
    \draw[-stealth] (m-4-3.west) -- (m-4-2.east) node [above,pos=0.5] {$d^1$};
\draw[thick] (m-1-1.east) -- (m-5-1.east) ;
\draw[thick] (m-5-1.north) -- (m-5-5.north) ;
\end{tikzpicture}
\end{figure}

Recall that $d^1$ is induced by inclusions (\ref{star_k:eq}). Let $a$ be the generator of $H_0(U_{11}\cap U_{12})$ and $\hat{a}$ be the generator of $H_0(U_{11}\cap U_{21})$. Let $b_1,\dots, b_{k-1}$ be the generators of $H_0(U_{22}\cap U_{12})$ and $\hat{b}_1,\dots, \hat{b}_{k-1}$ be the generators of $H_0(U_{22}\cap U_{21})$. Note that $U_{11}$ has two path components, while $U_{11}\cap U_{12}$ and $U_{11}\cap U_{21}$ lie in different path components. Choose a generator (denoted by $c_1$) from the path component of $U_{11}$ containing $U_{11}\cap U_{12}$, and choose a generator (denoted by $c_2$) from the path component of $U_{11}$ containing $U_{11}\cap U_{21}$. Let $e$ be the generator of $H_0(U_{12})$ and $\hat{e}$ be the generator of $H_0(U_{21})$ such that $$d^1(a)=c_1+e\mbox{\quad and \quad}d^1(\hat{a})=c_2+\hat{e}$$

Let $f$ be the generator of $H_0(U_{22})\cong\mathbb{Z}$, then for all $i=1,\dots, k-1$, $$d^1(b_i)=e+f\mbox{\quad and \quad}d^1(\hat{b}_i)=\hat{e}+f$$

Under the given generators, $d^1$ can be represented by the following matrix
$$\left[\begin{array}{cccccc}1 & 0 & 0 & 0 & \cdots & 0 \\0 & 1 & 0 & 0 & \cdots & 0 \\1 & 0 & 1 & 0 & \cdots & 0 \\0 & 1 & 0 & 1 & \cdots & 1 \\0 & 0 & 1 & 1 & \cdots & 1\end{array}\right]$$
where the third and the fourth columns are repeated for $(k-1)$ times. The rank of this matrix is $4$, and the diagonal elements of the Smith Normal form of the above matrix are all equal to $1$. Therefore, $\Ima d^1 \cong\mathbb{Z}^4$ and $\ker d^1\cong \mathbb{Z}^{2k-4}$. Hence we obtain the $E^2$ page of the Mayer-Vietoris spectral sequence: 
\begin{figure}[htbp!]
\centering
\begin{tikzpicture}
  \matrix (m) [matrix of math nodes,
    nodes in empty cells,nodes={minimum width=5ex,
     minimum height=2.5ex,outer sep=-5pt},
    column sep=1ex,row sep=1ex]{
          q     &      &     &  &  \\
          2     & \vdots & 0 & 0 &\\   
          1     &  H_1((\mathsf{Star}_{k-1})^2_r) &  0  & 0 \\
          0     &  \mathbb{Z}  & \mathbb{Z}^{2k-4} & 0\\
    \quad\strut &   0  &  1  &  2 & p\strut \\};
\draw[thick] (m-1-1.east) -- (m-5-1.east) ;
\draw[thick] (m-5-1.north) -- (m-5-5.north) ;
\end{tikzpicture}
\end{figure}

Since there is no non-trivial arrow on the $E^2$ page, $E^2=E^{\infty}$. Therefore, $$H_1((\hat{\mathsf{Star}}_{k})^2_r)\cong H_1((\mathsf{Star}_{k-1})^2_r)\oplus\mathbb{Z}^{2k-4}\cong \mathbb{Z}^{k(k-3)+1}$$ 
$$H_0((\hat{\mathsf{Star}}_{k})^2_r)\cong \mathbb{Z}$$
\end{proof}

\begin{proposition}
Let $L_{e_1}$ be a positive number. If $r>L_{e_1}$ and $r\leq 1$. Then
$$H_0((\hat{\mathsf{Star}}_k)_r^2) \cong \mathbb{Z}  \mbox{\quad and\quad} H_1((\hat{\mathsf{Star}}_k)_r^2)\cong \mathbb{Z}^{(k-1)(k-4)+1}$$
\end{proposition}
\begin{proof}
Consider the cover of $(\hat{\mathsf{Star}}_k)_r^2$ given in Equation (\ref{cover:star}). Note that 
 \begin{equation}
\begin{aligned}
U_{11}&=(e_1)^2_r=\emptyset\\
U_{12}&=e_1\times \mathsf{Star}_{k-1}-\{(x,y)\in e_1\times \mathsf{Star}_{k-1}\mid \delta(x,y)<r\}\simeq  \{\ast_1,\dots,\ast_{k-1}\}\\
U_{21}&=\mathsf{Star}_{k-1}\times e_1-\{(x,y)\in \mathsf{Star}_{k-1}\times e_1\mid \delta(x,y)<r\}\simeq  \{\ast_1, \dots,\ast_{k-1}\}
\end{aligned}
\end{equation}

 Now let's consider the intersections of $U_{ij}$. Note that 
 \begin{equation}
\begin{aligned}
U_{11}\cap U_{12}&=\emptyset=U_{11}\cap U_{21}\\
U_{22}\cap U_{12}&\simeq  \{\ast_1,\dots, \ast_{k-1}\}\simeq U_{22}\cap U_{21}\\
U_{11}\cap U_{22}&=\emptyset=U_{12}\cap U_{21}
\end{aligned}
\end{equation}

The $E^1$ page of the Mayer-Vietoris spectral sequence is 
\begin{figure}[htbp!]
\centering
\begin{tikzpicture}
  \matrix (m) [matrix of math nodes,
    nodes in empty cells,nodes={minimum width=5ex,
     minimum height=2.5ex,outer sep=-5pt},
    column sep=1.5ex,row sep=1ex]{
          q     &      &     &  &  \\
          2     & \vdots & 0 & 0 &\\   
          1     &  H_1((\mathsf{Star}_{k-1})^2_r) &  0  & 0 \\
          0     &  \mathbb{Z}^{2k-2}\oplus H_0((\mathsf{Star}_{k-1})^2_r)  & \mathbb{Z}^{2k-2} & 0\\
    \quad\strut &   0  &  1  &  2 & p\strut \\};
    \draw[-stealth] (m-4-3.west) -- (m-4-2.east) node [above,pos=0.5] {$d^1$};
\draw[thick] (m-1-1.east) -- (m-5-1.east) ;
\draw[thick] (m-5-1.north) -- (m-5-5.north) ;
\end{tikzpicture}
\end{figure}

Recall that $d^1$ is induced by inclusions (\ref{star_k:eq}). Let $b_1,\dots, b_{k-1}$ be the generators of $H_0(U_{22}\cap U_{12})$ and $\hat{b}_1,\dots, \hat{b}_{k-1}$ be the generators of $H_0(U_{22}\cap U_{21})$. Let $e, \dots, e_{k-1}$ be the generators of $H_0(U_{12})$ and $\hat{e}_1,\dots, \hat{e}_{k-1}$ be the generators of $H_0(U_{21})$. Let $f$ be the generator of $H_0(U_{22})\cong\mathbb{Z}$, then for all $i=1,\dots, k-1$, $$d^1(b_i)=e_i+f\mbox{\quad and \quad}d^1(\hat{b}_i)=\hat{e}_i+f$$

Under the given generators, $d^1$ can be represented by the following matrix
$$\left[\begin{array}{ccccc}1 & 0 & \cdots & 0 & 0 \\0 & 1 & \cdots & 0 & 0  \\0 & 0 & \ddots & 0 & 0  \\0 & 0 & \cdots & 1 & 0  \\0 & 0 & \cdots & 0 & 1  \\1 & 1 & \cdots & 1 & 1 \end{array}\right]$$
The rank of this matrix is $2k-2$, and the diagonal elements of the Smith Normal form of the above matrix are all equal to $1$. Therefore, $\Ima d^1 \cong\mathbb{Z}^{2k-2}$ and $\ker d^1=0$. Hence we obtain the $E^2$ page of the Mayer-Vietoris spectral sequence: \\

\begin{figure}[htbp!]
\centering
\begin{tikzpicture}
  \matrix (m) [matrix of math nodes,
    nodes in empty cells,nodes={minimum width=5ex,
     minimum height=2.5ex,outer sep=-5pt},
    column sep=1ex,row sep=1ex]{
          q     &      &     &  &  \\
          2     & \vdots & 0 & 0 &\\   
          1     &  H_1((\mathsf{Star}_{k-1})^2_r) &  0  & 0 \\
          0     &  \mathbb{Z}  &0 & 0\\
    \quad\strut &   0  &  1  &  2 & p\strut \\};
\draw[thick] (m-1-1.east) -- (m-5-1.east) ;
\draw[thick] (m-5-1.north) -- (m-5-5.north) ;
\end{tikzpicture}
\end{figure}

Since there is no non-trivial arrow on the $E^2$ page, $E^2=E^{\infty}$. Therefore, $H_0((\hat{\mathsf{Star}}_{k})^2_r)\cong \mathbb{Z}$ and
$H_1((\hat{\mathsf{Star}}_{k})^2_r)\cong H_1((\mathsf{Star}_{k-1})^2_r)\cong \mathbb{Z}^{(k-1)(k-4)+1}$.
\end{proof}

\begin{proposition}
Let $L_{e_1}$ be a positive number. If $r\leq L_{e_1}$ and $1<r\leq 2$ . Then
$$H_0((\hat{\mathsf{Star}}_k)_r^2) \cong \mathbb{Z}^{k^2-3k+4}  \mbox{\quad and\quad} H_1((\hat{\mathsf{Star}}_k)_r^2)=0$$
\end{proposition}
\begin{proof}
Consider the cover of $(\hat{\mathsf{Star}}_k)_r^2$ given in Equation (\ref{cover:star}). Note that 
 \begin{equation}
\begin{aligned}
U_{11}&=(e_1)^2_r=\{\ast_1,\ast_2\}\\
U_{12}&=e_1\times \mathsf{Star}_{k-1}-\{(x,y)\in e_1\times \mathsf{Star}_{k-1}\mid \delta(x,y)<r\}\simeq  \{\ast\}\\
U_{21}&=\mathsf{Star}_{k-1}\times e_1-\{(x,y)\in \mathsf{Star}_{k-1}\times e_1\mid \delta(x,y)<r\}\simeq  \{\ast\}
\end{aligned}
\end{equation}

 Now let's consider the intersections of $U_{ij}$. Note that 
 \begin{equation}
\begin{aligned}
U_{11}\cap U_{12}&\simeq \{\ast\}\simeq U_{11}\cap U_{21}\\
U_{22}\cap U_{12}&=\emptyset= U_{22}\cap U_{21}\\
U_{11}\cap U_{22}&=\emptyset=U_{12}\cap U_{21}
\end{aligned}
\end{equation}

The $E^1$ page of the Mayer-Vietoris spectral sequence is 
\begin{figure}[htbp!]
\centering
\begin{tikzpicture}
  \matrix (m) [matrix of math nodes,
    nodes in empty cells,nodes={minimum width=5ex,
     minimum height=2.5ex,outer sep=-5pt},
    column sep=1ex,row sep=1ex]{
          q     &      &     &  &  \\
          2     & \vdots & 0 & 0 &\\   
          1     &  H_1((\mathsf{Star}_{k-1})^2_r) &  0  & 0 \\
          0     &  \mathbb{Z}^{4}\oplus H_0((\mathsf{Star}_{k-1})^2_r)  & \mathbb{Z}^{2} & 0\\
    \quad\strut &   0  &  1  &  2 & p\strut \\};
    \draw[-stealth] (m-4-3.west) -- (m-4-2.east) node [above,pos=0.5] {$d^1$};
\draw[thick] (m-1-1.east) -- (m-5-1.east) ;
\draw[thick] (m-5-1.north) -- (m-5-5.north) ;
\end{tikzpicture}
\end{figure}

Recall that $d^1$ is induced by inclusions (\ref{star_k:eq}). Let $a$ be the generator of $H_0(U_{11}\cap U_{12})$ and $\hat{a}$ be the generator of $H_0(U_{11}\cap U_{21})$. Let $e$ be the generator of $H_0(U_{12})$ and $\hat{e}$ be the generator of $H_0(U_{21})$. On the other hand, let $f_1, \dots, f_{k^2-k}$ be the generators of $H_0(U_{22})\cong\mathbb{Z}^{k^2-k}$, then for all $i=1,\dots, k-1$, $$d^1(a)=c_1+e\mbox{\quad and \quad}d^1(\hat{b}_i)=c_2+\hat{e}$$

Under the given generators, $d^1$ can be represented by the following matrix
$$\left[\begin{array}{cc}1 & 0 \\0 & 1 \\1 & 0 \\0 & 1 \\0 & 0 \\\vdots & \vdots \\0 & 0\end{array}\right]$$
The rank of this matrix is $2$, and the diagonal elements of the Smith Normal form of the above matrix are all equal to $1$. Therefore, $\Ima d^1 \cong\mathbb{Z}^2$ and $\ker d^1=0$. Hence we obtain the $E^2$ page of the Mayer-Vietoris spectral sequence: \\

\begin{figure}[htbp!]
\centering
\begin{tikzpicture}
  \matrix (m) [matrix of math nodes,
    nodes in empty cells,nodes={minimum width=5ex,
     minimum height=2.5ex,outer sep=-5pt},
    column sep=1ex,row sep=1ex]{
          q     &      &     &  &  \\
          2     & \vdots & 0 & 0 &\\   
          1     &  H_1((\mathsf{Star}_{k-1})^2_r) &  0  & 0 \\
          0     &  \mathbb{Z}^{k^2-3k+4}  &0 & 0\\
    \quad\strut &   0  &  1  &  2 & p\strut \\};
\draw[thick] (m-1-1.east) -- (m-5-1.east) ;
\draw[thick] (m-5-1.north) -- (m-5-5.north) ;
\end{tikzpicture}
\end{figure}

Since there is no non-trivial arrow on the $E^2$ page, $E^2=E^{\infty}$. Therefore, $H_1((\hat{\mathsf{Star}}_{k})^2_r)\cong H_1((\mathsf{Star}_{k-1})^2_r)=0$ and $H_0((\hat{\mathsf{Star}}_{k})^2_r)\cong \mathbb{Z}^{k^2-3k+4}$.
\end{proof}

\begin{proposition}
Let $L_{e_1}$ be a positive number. If $L_{e_1}<r\leq L_{e_1}+1$ and $1<r\leq 2$ . Then
$$H_0((\hat{\mathsf{Star}}_k)_r^2) \cong \mathbb{Z}^{k^2-k}  \mbox{\quad and\quad} H_1((\hat{\mathsf{Star}}_k)_r^2)=0$$
\end{proposition}
\begin{proof}
Consider the cover of $(\hat{\mathsf{Star}}_k)_r^2$ given in Equation (\ref{cover:star}). Note that 
 \begin{equation}
\begin{aligned}
U_{11}&=(e_1)^2_r=\emptyset\\
U_{12}&=e_1\times \mathsf{Star}_{k-1}-\{(x,y)\in e_1\times \mathsf{Star}_{k-1}\mid \delta(x,y)<r\}\simeq  \{\ast_1, \dots, \ast_{k-1}\}\\
U_{21}&=\mathsf{Star}_{k-1}\times e_1-\{(x,y)\in \mathsf{Star}_{k-1}\times e_1\mid \delta(x,y)<r\}\simeq  \{\ast_1, \dots, \ast_{k-1}\}
\end{aligned}
\end{equation}

In addition, Note that $U_{ij}\cap U_{i'j'}=\emptyset$ when $i\neq i'$ or $j\neq j'$. Therefore, the $E^1$ page of the Mayer-Vietoris spectral sequence is 
\begin{figure}[htbp!]
\centering
\begin{tikzpicture}
  \matrix (m) [matrix of math nodes,
    nodes in empty cells,nodes={minimum width=5ex,
     minimum height=2.5ex,outer sep=-5pt},
    column sep=1ex,row sep=1ex]{
          q     &      &     &  &  \\
          2     & \vdots & 0 & 0 &\\   
          1     &  H_1((\mathsf{Star}_{k-1})^2_r) &  0  & 0 \\
          0     &  \mathbb{Z}^{2k-2}\oplus H_0((\mathsf{Star}_{k-1})^2_r)  & 0 & 0\\
    \quad\strut &   0  &  1  &  2 & p\strut \\};
\draw[thick] (m-1-1.east) -- (m-5-1.east) ;
\draw[thick] (m-5-1.north) -- (m-5-5.north) ;
\end{tikzpicture}
\end{figure}

Since there is no non-trivial arrow on the $E^1$ page, $E^1=E^{\infty}$. Therefore, $H_1((\hat{\mathsf{Star}}_{k})^2_r)\cong H_1((\mathsf{Star}_{k-1})^2_r)=0$ and $H_0((\hat{\mathsf{Star}}_{k})^2_r)\cong \mathbb{Z}^{k^2-k}$. 
\end{proof}

\begin{proposition}
Let $L_{e_1}$ be a positive number. If $r>L_{e_1}+1$ and $1<r\leq 2$ . Then
$$H_0((\hat{\mathsf{Star}}_k)_r^2) \cong \mathbb{Z}^{k^2-3k+2}  \mbox{\quad and\quad} H_1((\hat{\mathsf{Star}}_k)_r^2)=0$$
\end{proposition}
\begin{proof}
Consider the cover of $(\hat{\mathsf{Star}}_k)_r^2$ given in Equation (\ref{cover:star}). Note that 
 \begin{equation}
\begin{aligned}
U_{11}&=(e_1)^2_r=\emptyset\\
U_{12}&=e_1\times \mathsf{Star}_{k-1}-\{(x,y)\in e_1\times \mathsf{Star}_{k-1}\mid \delta(x,y)<r\}=\emptyset\\
U_{21}&=\mathsf{Star}_{k-1}\times e_1-\{(x,y)\in \mathsf{Star}_{k-1}\times e_1\mid \delta(x,y)<r\}=\emptyset
\end{aligned}
\end{equation}

In addition, $U_{ij}\cap U_{i'j'}=\emptyset$ when $i\neq i'$ or $j\neq j'$. Therefore, the $E^1$ page of the Mayer-Vietoris spectral sequence is \\
\begin{figure}[htbp!]
\centering
\begin{tikzpicture}
  \matrix (m) [matrix of math nodes,
    nodes in empty cells,nodes={minimum width=5ex,
     minimum height=2.5ex,outer sep=-5pt},
    column sep=1ex,row sep=1ex]{
          q     &      &     &  &  \\
          2     & \vdots & 0 & 0 &\\   
          1     &  H_1((\mathsf{Star}_{k-1})^2_r) &  0  & 0 \\
          0     &  H_0((\mathsf{Star}_{k-1})^2_r)  & 0 & 0\\
    \quad\strut &   0  &  1  &  2 & p\strut \\};
\draw[thick] (m-1-1.east) -- (m-5-1.east) ;
\draw[thick] (m-5-1.north) -- (m-5-5.north) ;
\end{tikzpicture}
\end{figure}

Since there is no non-trivial arrow on the $E^1$ page, $E^1=E^{\infty}$. Therefore, $H_1((\hat{\mathsf{Star}}_{k})^2_r)\cong H_1((\mathsf{Star}_{k-1})^2_r)=0$ and $H_0((\hat{\mathsf{Star}}_{k})^2_r)\cong \mathbb{Z}^{k^2-3k+2}$.
\end{proof}

\begin{proposition}
Let $L_{e_1}$ be a positive number. If $r\leq L_{e_1}$ and $r>2$ . Then
$$H_0((\hat{\mathsf{Star}}_k)_r^2) \cong \mathbb{Z}^{2}  \mbox{\quad and\quad} H_1((\hat{\mathsf{Star}}_k)_r^2)=0$$
\end{proposition}
\begin{proof}
Consider the cover of $(\hat{\mathsf{Star}}_k)_r^2$ given in Equation (\ref{cover:star}). Note that
 \begin{equation}
\begin{aligned}
U_{11}&=(e_1)^2_r\simeq\{\ast_1,\ast_2\}\\
U_{12}&=e_1\times \mathsf{Star}_{k-1}-\{(x,y)\in e_1\times \mathsf{Star}_{k-1}\mid \delta(x,y)<r\}\simeq\{\ast\}\\
U_{21}&=\mathsf{Star}_{k-1}\times e_1-\{(x,y)\in \mathsf{Star}_{k-1}\times e_1\mid \delta(x,y)<r\}\simeq \{\ast\}
\end{aligned}
\end{equation}

Now let's consider the intersections of $U_{ij}$. Note that 
 \begin{equation}
\begin{aligned}
U_{11}\cap U_{12}&\simeq\{\ast\} \simeq U_{11}\cap U_{21}\\
U_{22}\cap U_{12}&=\emptyset= U_{22}\cap U_{21}\\
U_{11}\cap U_{22}&=\emptyset=U_{12}\cap U_{21}
\end{aligned}
\end{equation}

The $E^1$ page of the Mayer-Vietoris spectral sequence is 
\begin{figure}[htbp!]
\centering
\begin{tikzpicture}
  \matrix (m) [matrix of math nodes,
    nodes in empty cells,nodes={minimum width=5ex,
     minimum height=2.5ex,outer sep=-5pt},
    column sep=1ex,row sep=1ex]{
          q     &      &     &  &  \\
          2     & \vdots & 0 & 0 &\\   
          1     &  H_1((\mathsf{Star}_{k-1})^2_r) &  0  & 0 \\
          0     &  \mathbb{Z}^4\oplus H_0((\mathsf{Star}_{k-1})^2_r)  & \mathbb{Z}^2 & 0\\
    \quad\strut &   0  &  1  &  2 & p\strut \\};
    \draw[-stealth] (m-4-3.west) -- (m-4-2.east) node [above,pos=0.5] {$d^1$};
\draw[thick] (m-1-1.east) -- (m-5-1.east) ;
\draw[thick] (m-5-1.north) -- (m-5-5.north) ;
\end{tikzpicture}
\end{figure}

Recall that $d^1$ is induced by inclusions (\ref{star_k:eq}). Let $a$ be the generator of $H_0(U_{11}\cap U_{12})$ and $\hat{a}$ be the generator of $H_0(U_{11}\cap U_{21})$. Let $c_1, c_2$ be the generators of $H_0(U_{11})$, where $c_1$ and $c_2$ lie in different path components of $U_{11}$. Let $e$ be the generator of $H_0(U_{12})$ and $\hat{e}$ be the generator of $H_0(U_{21})$. Note that $$d^1(a)=c_1+e \mbox{\quad and\quad} d^1(\hat{a})=c_2+\hat{e}$$

Under the given generators, $d^1$ can be represented by the following matrix
$$\left[\begin{array}{cc}1 & 0 \\0 & 1 \\1 & 0 \\0 & 1\end{array}\right]$$
The rank of this matrix is $2$, and the diagonal elements of the Smith Normal form of the above matrix are all equal to $1$. Therefore, $\Ima d^1 \cong\mathbb{Z}^2$ and $\ker d^1=0$. Hence we obtain the $E^2$ page of the Mayer-Vietoris spectral sequence: \\

\begin{figure}[htbp!]
\centering
\begin{tikzpicture}
  \matrix (m) [matrix of math nodes,
    nodes in empty cells,nodes={minimum width=5ex,
     minimum height=2.5ex,outer sep=-5pt},
    column sep=1ex,row sep=1ex]{
          q     &      &     &  &  \\
          2     & \vdots & 0 & 0 &\\   
          1     &  H_1((\mathsf{Star}_{k-1})^2_r) &  0  & 0 \\
          0     &  \mathbb{Z}^{2}  &0 & 0\\
    \quad\strut &   0  &  1  &  2 & p\strut \\};
\draw[thick] (m-1-1.east) -- (m-5-1.east) ;
\draw[thick] (m-5-1.north) -- (m-5-5.north) ;
\end{tikzpicture}
\end{figure}

Since there is no non-trivial arrow on the $E^2$ page, $E^2=E^{\infty}$. Therefore, $H_1((\hat{\mathsf{Star}}_{k})^2_r)\cong H_1((\mathsf{Star}_{k-1})^2_r)=0$ and $H_0((\hat{\mathsf{Star}}_{k})^2_r)\cong \mathbb{Z}^{2}$. 
\end{proof}

\begin{proposition}
Let $L_{e_1}$ be a positive number. If $L_{e_1}<r\leq L_{e_1}+1$ and $r>2$ . Then
$$H_0((\hat{\mathsf{Star}}_k)_r^2) \cong \mathbb{Z}^{2k-2}  \mbox{\quad and\quad} H_1((\hat{\mathsf{Star}}_k)_r^2)=0$$
\end{proposition}
\begin{proof}
Consider the cover of $(\hat{\mathsf{Star}}_k)_r^2$ given in Equation (\ref{cover:star}). Note that $U_{11}=(e_1)^2_r=\emptyset$. Moreover, $$U_{12}=e_1\times \mathsf{Star}_{k-1}-\{(x,y)\in e_1\times \mathsf{Star}_{k-1}\mid \delta(x,y)<r\}\simeq\{\ast_1, \dots, \ast_{k-1}\}$$ and $$U_{21}=\mathsf{Star}_{k-1}\times e_1-\{(x,y)\in \mathsf{Star}_{k-1}\times e_1\mid \delta(x,y)<r\}\simeq \{\ast_1, \dots, \ast_{k-1}\}$$ 

Note that $U_{ij}\cap U_{i'j'}=\emptyset$ when $i\neq i'$ or $j\neq j'$. Therefore, the $E^1$ page of the Mayer-Vietoris spectral sequence is 
\begin{figure}[htbp!]
\centering
\begin{tikzpicture}
  \matrix (m) [matrix of math nodes,
    nodes in empty cells,nodes={minimum width=5ex,
     minimum height=2.5ex,outer sep=-5pt},
    column sep=1ex,row sep=1ex]{
          q     &      &     &  &  \\
          2     & \vdots & 0 & 0 &\\   
          1     &  H_1((\mathsf{Star}_{k-1})^2_r) &  0  & 0 \\
          0     &  \mathbb{Z}^{2k-2}\oplus H_0((\mathsf{Star}_{k-1})^2_r)  &0 & 0\\
    \quad\strut &   0  &  1  &  2 & p\strut \\};
\draw[thick] (m-1-1.east) -- (m-5-1.east) ;
\draw[thick] (m-5-1.north) -- (m-5-5.north) ;
\end{tikzpicture}
\end{figure}

Since there is no non-trivial arrow on the $E^1$ page, $E^1=E^{\infty}$. Therefore, $H_1((\hat{\mathsf{Star}}_{k})^2_r)\cong H_1((\mathsf{Star}_{k-1})^2_r)=0$ and $H_0((\hat{\mathsf{Star}}_{k})^2_r)\cong \mathbb{Z}^{2k-2}$. 
\end{proof}

In summary, we proved Theorem \ref{Thm1}.

\section{Decomposition of $PH_i{((\mathsf{Star}_k)^2_{-,-};\mathbb{F})}$}
In section \ref{chap4:y2}, we decompose the persistence module $PH_0(Y^2_{-,-};\mathbb{F})$ into a direct sum of indecomposables. In Section \ref{decompo:star3}, we decompose the persistence module $PH_0((\mathsf{Star_k})^2_{r,L_{e_1}};\mathbb{F})$ into a direct sum of indecomposables. Moreover, we prove that $PH_1((\mathsf{Star_k})^2_{r,L_{e_1}};\mathbb{F})$ is interval decomposable.

\subsection{Decomposition of $PH_i(Y^2_{-,-};\mathbb{F})$}\label{chap4:y2}
Note that $PH_1(Y^2_{-,-};\mathbb{F})$ is an interval module because the support of $PH_1(Y^2_{-,-};\mathbb{F})$ is an interval, we immediately obtain Theorem \ref{Thm2'}.\\

Now we are going to decompose $PH_0(Y^2_{-,-};\mathbb{F})$. We associate $PH_0(Y^2_{-,-};\mathbb{F})$ with a persistence module $N$, where $N$ is constructed using Proposition \ref{prop-sim-2}. We use the following diagram to represent $N$, where each vector space in the diagram represents a constant functor indexed by a chamber of the parameter space of $PH_0(Y^2_{-,-};\mathbb{F})$.\\

\begin{minipage}{0.5\textwidth}
\begin{equation*}
N=\begin{tikzcd}
  \mathbb{F} & \mathbb{F}^{4}\arrow{l}[swap]{} &\mathbb{F}^2\arrow{l}[swap]{}\\
  \mathbb{F}^2 \arrow{u}{}		  &\mathbb{F}^{6} \arrow{u}{}\arrow{l}[swap]{}                       & \mathbb{F}^{4}\arrow{u}[swap]{}\arrow{l}[swap]{}\\ 
  							     		  &\mathbb{F}^{2} \arrow{u}{} 			         & \\ 
\end{tikzcd}
\end{equation*}
\end{minipage}%
\begin{minipage}{0.45\textwidth}
\centering
\resizebox{5cm}{!}{

\tikzset{every picture/.style={line width=0.75pt}} 

\begin{tikzpicture}[x=0.75pt,y=0.75pt,yscale=-1,xscale=1]

\draw  [draw opacity=0][fill={rgb, 255:red, 0; green, 0; blue, 0 }  ,fill opacity=0.32 ] (105.58,83.34) -- (152.5,130) -- (146.5,147) -- (83.3,147) -- (19.64,83.34) -- cycle ;
\draw  [draw opacity=0][fill={rgb, 255:red, 248; green, 231; blue, 28 }  ,fill opacity=0.63 ] (146.5,146) -- (146.5,146) -- (83.3,146) -- (83.5,106) -- (146.5,146) -- cycle ;
\draw  [draw opacity=0][fill={rgb, 255:red, 0; green, 0; blue, 0 }  ,fill opacity=0.32 ] (218.57,80.16) -- (304.51,80.16) -- (253.7,146) -- (190.5,146) -- (182.5,129) -- cycle ;
\draw  [draw opacity=0][fill={rgb, 255:red, 0; green, 0; blue, 0 }  ,fill opacity=0.32 ] (154.91,16.5) -- (218.57,80.16) -- (182.5,129) -- (150.9,129) -- (105.58,83.34) -- cycle ;
\draw  [draw opacity=0][fill={rgb, 255:red, 0; green, 0; blue, 0 }  ,fill opacity=0.32 ] (253.7,146) -- (318.83,210.66) -- (232.89,210.66) -- (181.5,162) -- (190.5,146) -- cycle ;
\draw  [draw opacity=0][fill={rgb, 255:red, 0; green, 0; blue, 0 }  ,fill opacity=0.32 ] (232.89,210.66) -- (183.56,277.5) -- (119.9,213.84) -- (155.97,165) -- (181.5,162) -- cycle ;
\draw  [draw opacity=0][fill={rgb, 255:red, 0; green, 0; blue, 0 }  ,fill opacity=0.32 ] (146.5,147) -- (156.5,163) -- (119.9,213.84) -- (25.27,213.84) -- (84.77,148) -- cycle ;
\draw  [draw opacity=0][fill={rgb, 255:red, 248; green, 231; blue, 28 }  ,fill opacity=0.63 ] (180.9,128) -- (180.9,128) -- (218.5,45) -- (218.57,80.16) -- (182.5,129) -- cycle ;
\draw [color={rgb, 255:red, 0; green, 0; blue, 0 }  ,draw opacity=1 ][line width=1.5]    (105.58,83.34) -- (152.5,130) ;
\draw [color={rgb, 255:red, 0; green, 0; blue, 0 }  ,draw opacity=1 ][line width=1.5]    (218.57,80.16) -- (182.5,129) ;
\draw [color={rgb, 255:red, 208; green, 2; blue, 27 }  ,draw opacity=1 ][line width=1.5]    (83.3,147) -- (146.5,147) ;
\draw [line width=1.5]    (19.64,83.34) -- (83.3,147) ;
\draw [color={rgb, 255:red, 208; green, 2; blue, 27 }  ,draw opacity=1 ][line width=2.25]    (19.64,83.34) -- (105.58,83.34) ;
\draw [line width=1.5]    (83.3,147) -- (25.27,213.84) ;
\draw [color={rgb, 255:red, 208; green, 2; blue, 27 }  ,draw opacity=1 ][line width=1.5]    (25.27,213.84) -- (119.9,213.84) ;
\draw [line width=1.5]    (304.51,80.16) -- (255.17,147) ;
\draw [color={rgb, 255:red, 208; green, 2; blue, 27 }  ,draw opacity=1 ][line width=2.25]    (218.57,80.16) -- (304.51,80.16) ;
\draw [color={rgb, 255:red, 208; green, 2; blue, 27 }  ,draw opacity=1 ][line width=1.5]    (232.89,210.66) -- (318.83,210.66) ;
\draw [line width=1.5]    (255.17,147) -- (318.83,210.66) ;
\draw [color={rgb, 255:red, 245; green, 166; blue, 35 }  ,draw opacity=1 ][line width=1.5]    (83.5,107) -- (83.3,147) ;
\draw [line width=1.5]    (119.9,213.84) -- (183.56,277.5) ;
\draw [line width=1.5]    (232.89,210.66) -- (183.56,277.5) ;
\draw [line width=1.5]    (154.91,16.5) -- (105.58,83.34) ;
\draw [line width=1.5]    (154.91,16.5) -- (218.57,80.16) ;
\draw [color={rgb, 255:red, 245; green, 166; blue, 35 }  ,draw opacity=1 ][line width=1.5]    (104.5,44) -- (105.58,83.34) ;
\draw [color={rgb, 255:red, 245; green, 166; blue, 35 }  ,draw opacity=1 ][line width=1.5]    (118.5,170) -- (119.9,213.84) ;
\draw [color={rgb, 255:red, 245; green, 166; blue, 35 }  ,draw opacity=1 ][line width=1.5]    (218.5,45) -- (218.57,80.16) ;
\draw [color={rgb, 255:red, 245; green, 166; blue, 35 }  ,draw opacity=1 ][line width=1.5]    (254.5,107) -- (255.17,147) ;
\draw [color={rgb, 255:red, 245; green, 166; blue, 35 }  ,draw opacity=1 ][line width=1.5]    (233.5,169) -- (232.89,210.66) ;
\draw [color={rgb, 255:red, 80; green, 227; blue, 194 }  ,draw opacity=1 ][line width=2.25]    (181.5,162) -- (190.5,146) ;
\draw [color={rgb, 255:red, 80; green, 227; blue, 194 }  ,draw opacity=1 ][line width=2.25]    (182.5,129) -- (190.5,146) ;
\draw [color={rgb, 255:red, 80; green, 227; blue, 194 }  ,draw opacity=1 ][line width=2.25]    (182.5,129) -- (152.5,130) ;
\draw [color={rgb, 255:red, 80; green, 227; blue, 194 }  ,draw opacity=1 ][line width=2.25]    (152.5,130) -- (146.5,147) ;
\draw [color={rgb, 255:red, 80; green, 227; blue, 194 }  ,draw opacity=1 ][line width=2.25]    (156.5,163) -- (146.5,147) ;
\draw [color={rgb, 255:red, 80; green, 227; blue, 194 }  ,draw opacity=1 ][line width=2.25]    (156.5,163) -- (181.5,162) ;
\draw [color={rgb, 255:red, 25; green, 2; blue, 208 }  ,draw opacity=1 ][line width=2.25]    (152.5,130) -- (146.5,147) ;
\draw [color={rgb, 255:red, 25; green, 2; blue, 208 }  ,draw opacity=1 ][line width=2.25]    (146.5,147) -- (156.5,163) ;
\draw [color={rgb, 255:red, 25; green, 2; blue, 208 }  ,draw opacity=1 ][line width=2.25]    (156.5,163) -- (181.5,162) ;
\draw [color={rgb, 255:red, 25; green, 2; blue, 208 }  ,draw opacity=1 ][line width=2.25]    (181.5,162) -- (190.5,146) ;
\draw [color={rgb, 255:red, 25; green, 2; blue, 208 }  ,draw opacity=1 ][line width=2.25]    (182.5,129) -- (190.5,146) ;
\draw [color={rgb, 255:red, 25; green, 2; blue, 208 }  ,draw opacity=1 ][line width=2.25]    (152.5,130) -- (182.5,129) ;
\draw [color={rgb, 255:red, 245; green, 166; blue, 35 }  ,draw opacity=1 ][line width=1.5]    (104.5,44) -- (152.5,130) ;
\draw [color={rgb, 255:red, 245; green, 166; blue, 35 }  ,draw opacity=1 ][line width=1.5]    (218.5,45) -- (182.5,129) ;
\draw [color={rgb, 255:red, 245; green, 166; blue, 35 }  ,draw opacity=1 ][line width=1.5]    (254.5,107) -- (190.5,146) ;
\draw [color={rgb, 255:red, 245; green, 166; blue, 35 }  ,draw opacity=1 ][line width=1.5]    (233.5,169) -- (181.5,162) ;
\draw [color={rgb, 255:red, 245; green, 166; blue, 35 }  ,draw opacity=1 ][line width=1.5]    (118.5,170) -- (156.5,163) ;
\draw [color={rgb, 255:red, 245; green, 166; blue, 35 }  ,draw opacity=1 ][line width=1.5]    (83.5,107) -- (146.5,147) ;
\draw [color={rgb, 255:red, 0; green, 0; blue, 0 }  ,draw opacity=1 ][line width=1.5]    (181.5,162) -- (232.89,210.66) ;
\draw [color={rgb, 255:red, 208; green, 2; blue, 27 }  ,draw opacity=1 ][line width=1.5]    (190.5,146) -- (201.5,146) -- (253.7,146) ;
\draw [color={rgb, 255:red, 0; green, 0; blue, 0 }  ,draw opacity=1 ][line width=1.5]    (155.97,165) -- (119.9,213.84) ;
\draw  [draw opacity=0][fill={rgb, 255:red, 248; green, 231; blue, 28 }  ,fill opacity=0.63 ] (150.9,129) -- (150.9,129) -- (105.58,83.34) -- (104.5,44) -- (152.5,130) -- cycle ;
\draw  [draw opacity=0][fill={rgb, 255:red, 248; green, 231; blue, 28 }  ,fill opacity=0.63 ] (190.5,145) -- (190.5,145) -- (253.7,145) -- (254.5,106) -- (190.5,145) -- cycle ;
\draw  [draw opacity=0][fill={rgb, 255:red, 248; green, 231; blue, 28 }  ,fill opacity=0.63 ] (156.5,163) -- (157.03,161) -- (119.9,213.84) -- (118.5,170) -- (156.5,163) -- cycle ;
\draw  [draw opacity=0][fill={rgb, 255:red, 248; green, 231; blue, 28 }  ,fill opacity=0.63 ] (181.5,162) -- (181.5,162) -- (232.89,210.66) -- (233.5,169) -- (181.5,162) -- cycle ;
\draw [color={rgb, 255:red, 0; green, 0; blue, 0 }  ,draw opacity=1 ]   (25.27,213.84) ;
\draw [shift={(25.27,213.84)}, rotate = 0] [color={rgb, 255:red, 0; green, 0; blue, 0 }  ,draw opacity=1 ][fill={rgb, 255:red, 0; green, 0; blue, 0 }  ,fill opacity=1 ][line width=0.75]      (0, 0) circle [x radius= 6.7, y radius= 6.7]   ;
\draw [color={rgb, 255:red, 0; green, 0; blue, 0 }  ,draw opacity=1 ]   (19.64,83.34) ;
\draw [shift={(19.64,83.34)}, rotate = 0] [color={rgb, 255:red, 0; green, 0; blue, 0 }  ,draw opacity=1 ][fill={rgb, 255:red, 0; green, 0; blue, 0 }  ,fill opacity=1 ][line width=0.75]      (0, 0) circle [x radius= 6.7, y radius= 6.7]   ;
\draw [color={rgb, 255:red, 0; green, 0; blue, 0 }  ,draw opacity=1 ]   (154.91,16.5) ;
\draw [shift={(154.91,16.5)}, rotate = 0] [color={rgb, 255:red, 0; green, 0; blue, 0 }  ,draw opacity=1 ][fill={rgb, 255:red, 0; green, 0; blue, 0 }  ,fill opacity=1 ][line width=0.75]      (0, 0) circle [x radius= 6.7, y radius= 6.7]   ;
\draw [color={rgb, 255:red, 0; green, 0; blue, 0 }  ,draw opacity=1 ]   (318.83,210.66) ;
\draw [shift={(318.83,210.66)}, rotate = 0] [color={rgb, 255:red, 0; green, 0; blue, 0 }  ,draw opacity=1 ][fill={rgb, 255:red, 0; green, 0; blue, 0 }  ,fill opacity=1 ][line width=0.75]      (0, 0) circle [x radius= 6.7, y radius= 6.7]   ;
\draw [color={rgb, 255:red, 0; green, 0; blue, 0 }  ,draw opacity=1 ]   (304.51,80.16) ;
\draw [shift={(304.51,80.16)}, rotate = 0] [color={rgb, 255:red, 0; green, 0; blue, 0 }  ,draw opacity=1 ][fill={rgb, 255:red, 0; green, 0; blue, 0 }  ,fill opacity=1 ][line width=0.75]      (0, 0) circle [x radius= 6.7, y radius= 6.7]   ;
\draw [color={rgb, 255:red, 0; green, 0; blue, 0 }  ,draw opacity=1 ]   (183.56,277.5) ;
\draw [shift={(183.56,277.5)}, rotate = 0] [color={rgb, 255:red, 0; green, 0; blue, 0 }  ,draw opacity=1 ][fill={rgb, 255:red, 0; green, 0; blue, 0 }  ,fill opacity=1 ][line width=0.75]      (0, 0) circle [x radius= 6.7, y radius= 6.7]   ;

\draw (13.64,220.24) node [anchor=north west][inner sep=0.75pt] [font=\large]   {$c_{1}$};
\draw (13.64,90.74) node [anchor=north west][inner sep=0.75pt] [font=\large]   {$c_{2}$};
\draw (167,2.4) node [anchor=north west][inner sep=0.75pt] [font=\large]   {$c_{3}$};
\draw (306.51,83.56) node [anchor=north west][inner sep=0.75pt] [font=\large]   {$c_{4}$};
\draw (320.83,220.24) node [anchor=north west][inner sep=0.75pt] [font=\large]   {$c_{5}$};
\draw (194,272.4) node [anchor=north west][inner sep=0.75pt] [font=\large]   {$c_{6}$};

\end{tikzpicture}
}

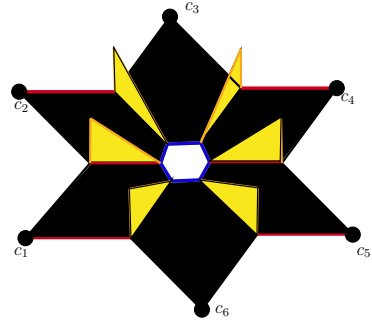
\captionof{figure}{The representatives of the basis of $PH_0(Y^2_{r,L};\mathbb{F})$. The edge $e_1$ is highlighted in red.}
\label{fig:star-h0geo-select}  
\end{minipage}\\

We choose a representative for each basis element of every vector space in the representation $N$ according to the geometry of $Y^2_{r,L_{e_1}}$. Figure \ref{fig:star-h0geo-select} provides the representatives we choose, and (\ref{Y-h0gen}) gives the basis we choose for each vector space of the representation $N$. Under the given bases for each vector space, the behavior of the morphisms can be described as follows: $\alpha$ maps every basis element of $<c_1, c_3, c_4, c_6>$  to $c_1$, $\gamma$ maps $c_1, c_5, c_6$ to $c_1$ and maps $c_2, c_3, c_4$ to $c_2$, $\lambda$ maps $c_3, c_6$ to $c_1$, $\beta$ maps $c_1, c_2$ to $c_1$, maps $c_4, c_5$ to $c_4$, maps $c_3$ to $c_3$, and maps $c_6$ to $c_6$, and $\epsilon$ maps $c_1, c_2$ to $c_1$ and maps $c_4, c_5$ to $c_4$. Moreover, each unlabeled map is an inclusion map. By an abuse of notation, we use $PH_0(Y^2_{-,-};\mathbb{F})$ to denote the $P$-indexed persistence module given by (\ref{Y-h0gen}).\\

\begin{equation}\label{Y-h0gen}
N=\begin{tikzcd}
  <c_1> & <c_1, c_3, c_4, c_6> \arrow{l}[swap]{\alpha} &<c_1, c_4>\arrow[hookrightarrow]{l}[swap]{}\\
  <c_3,c_6> \arrow{u}{\lambda}		  &<c_1, c_2, c_3,c_4,c_5,c_6> \arrow{u}{\beta}\arrow{l}[swap]{\gamma}                       & <c_1, c_2, c_4, c_5>\arrow{u}[swap]{\epsilon}\arrow[hookrightarrow]{l}[swap]{}\\ 
  							     		  &<c_3, c_6> \arrow[hookrightarrow]{u}{} 			         & \\ 
\end{tikzcd}
\end{equation}

\begin{proposition}\label{Y-h0gen-N}
$PH_0(Y^2_{-,-};\mathbb{F})\cong M_1\oplus M_2 \oplus E_1\oplus E_2\oplus F$ where\\

\begin{minipage}{0.3\textwidth}
   $$M_1=\begin{tikzcd}
  \mathbb{F} & \mathbb{F}^{2}\arrow{l}[swap]{\spmat{1 \\1}} &\mathbb{F}\arrow{l}[swap]{\spmat{1 \\0}}\\
  \mathbb{F} \arrow{u}{}		  &\mathbb{F}^{2} \arrow{u}{\spmat{1 & 0\\0 & 1}}\arrow{l}{\spmat{1 \\1}}                       & \mathbb{F}\arrow{u}[swap]{}\arrow{l}[swap]{\spmat{1 \\0}}\\ 
  							     		  &\mathbb{F} \arrow{u}{\spmat{0 \\1}} 			         & \\ 
    \end{tikzcd}
    $$ 
\end{minipage}%
\begin{minipage}{0.3\textwidth}
  $$M_2=\begin{tikzcd}
  0& \mathbb{F}\arrow{l}[swap]{} &\mathbb{F}\arrow{l}[swap]{}\\
 0 \arrow{u}{}		  &\mathbb{F} \arrow{u}{}\arrow{l}[swap]{}                       & \mathbb{F}\arrow{u}[swap]{}\arrow{l}[swap]{}\\ 
  							     		  &0 \arrow{u}{} 			         & \\ 
    \end{tikzcd}
    $$  
\end{minipage}%
\begin{minipage}{0.3\textwidth}
  $$E_1=\begin{tikzcd}
  0& 0\arrow{l}[swap]{} &0\arrow{l}[swap]{}\\
 0 \arrow{u}{}		  &\mathbb{F} \arrow{u}{}\arrow{l}[swap]{}                       & \mathbb{F}\arrow{u}[swap]{}\arrow{l}[swap]{}\\ 
  							     		  &0 \arrow{u}{} 			         & \\ 
    \end{tikzcd}$$  
\end{minipage}

\begin{minipage}{0.3\textwidth}
 $$E_2= \begin{tikzcd}
  0& 0\arrow{l}[swap]{} &0\arrow{l}[swap]{}\\
0 \arrow{u}{}		  &\mathbb{F} \arrow{u}{}\arrow{l}[swap]{}                       & \mathbb{F}\arrow{u}[swap]{}\arrow{l}[swap]{}\\ 
  							     		  &0 \arrow{u}{} 			         & \\ 
    \end{tikzcd}$$   
\end{minipage}%
\begin{minipage}{0.3\textwidth}
 $$F=\begin{tikzcd}
  0& \mathbb{F}\arrow{l}[swap]{} &0\arrow{l}[swap]{}\\
 \mathbb{F} \arrow{u}{}		  &\mathbb{F} \arrow{u}{}\arrow{l}[swap]{}                       & 0\arrow{u}[swap]{}\arrow{l}[swap]{}\\ 
  							     		  &\mathbb{F} \arrow{u}{} 			         & \\ 
\end{tikzcd}$$   
\end{minipage}

\end{proposition}

\begin{proof}
We choose a basis for each vector space $N$: 
\begin{equation}\label{Y-h0gen'}
\begin{tikzcd}
  <c_1> & <c_1, c_3-c_6, c_4-c_1, c_6> \arrow{l}[swap]{\alpha} &<c_1, c_4-c_1>\arrow[hookrightarrow]{l}[swap]{}\\
  <c_3-c_6,c_6> \arrow{u}{\lambda}		  &\left\langle\parbox{4cm}{$c_1, c_2-c_1+c_5-c_4,\\c_3-c_6, c_4-c_1, c_5-c_4, c_6$}\right\rangle  \arrow{u}{\beta}\arrow{l}[swap]{\gamma}                       &\left\langle\parbox{3cm}{$c_1, c_2-c_1+c_5-c_4,\\ c_4-c_1, c_5-c_4$}\right\rangle\arrow{u}[swap]{\epsilon}\arrow[hookrightarrow]{l}[swap]{}\\ 
  							     		  &<c_3-c_6, c_6> \arrow[hookrightarrow]{u}{} 			         & \\ 
\end{tikzcd}
\end{equation}
Define
\begin{equation}\label{Y-h0gen1}
M_1=\begin{tikzcd}
  <c_1> & <c_1, c_6> \arrow{l}[swap]{\spmat{1 \\1}} &<c_1>\arrow[hookrightarrow]{l}[swap]{}\\
  <c_6> \arrow{u}{\id}		  &<c_1,c_6> \arrow{u}{\spmat{1 & 0\\0 & 1}}\arrow{l}{\spmat{1 \\1}}                       & <c_1>\arrow{u}[swap]{\id}\arrow[hookrightarrow]{l}[swap]{}\\ 
  							     		  &<c_6> \arrow[hookrightarrow]{u}{} 			         & \\ 
\end{tikzcd}
\cong \begin{tikzcd}
  \mathbb{F} & \mathbb{F}^{2}\arrow{l}[swap]{\spmat{1 \\1}} &\mathbb{F}\arrow{l}[swap]{}\\
  \mathbb{F} \arrow{u}{}		  &\mathbb{F}^{2} \arrow{u}{\spmat{1 & 0\\0 & 1}}\arrow{l}{\spmat{1 \\1}}                       & \mathbb{F}\arrow{u}[swap]{}\arrow{l}[swap]{}\\ 
  							     		  &\mathbb{F} \arrow{u}{} 			         & \\ 
\end{tikzcd}
\end{equation}

\begin{equation}\label{Y-h0gen2}
M_2=\begin{tikzcd}
  0 & <c_4-c_1> \arrow{l}[swap]{} &<c_4-c_1>\arrow[hookrightarrow]{l}[swap]{}\\
  0 \arrow{u}{}		  &<c_4-c_1> \arrow{u}{\id}\arrow{l}[swap]{}                       & <c_4-c_1>\arrow{u}[swap]{\id}\arrow[hookrightarrow]{l}[swap]{}\\ 
  							     		  &0 \arrow[hookrightarrow]{u}{} 			         & \\ 
\end{tikzcd}\cong \begin{tikzcd}
  0& \mathbb{F}\arrow{l}[swap]{} &\mathbb{F}\arrow{l}[swap]{}\\
 0 \arrow{u}{}		  &\mathbb{F} \arrow{u}{}\arrow{l}[swap]{}                       & \mathbb{F}\arrow{u}[swap]{}\arrow{l}[swap]{}\\ 
  							     		  &0 \arrow{u}{} 			         & \\ 
\end{tikzcd}
\end{equation}

\begin{equation}\label{Y-h0gen3}
E_1=\begin{tikzcd}
  0 & 0 \arrow{l}[swap]{} &0\arrow[hookrightarrow]{l}[swap]{}\\
  0 \arrow{u}{}		  &<c_2-c_1+c_5-c_4> \arrow{u}{}\arrow{l}[swap]{}                       & <c_2-c_1+c_5-c_4>\arrow{u}[swap]{}\arrow[hookrightarrow]{l}[swap]{}\\ 
  							     		  &0 \arrow[hookrightarrow]{u}{} 			         & \\ 
\end{tikzcd}\cong \begin{tikzcd}
  0& 0\arrow{l}[swap]{} &0\arrow{l}[swap]{}\\
 0 \arrow{u}{}		  &\mathbb{F} \arrow{u}{}\arrow{l}[swap]{}                       & \mathbb{F}\arrow{u}[swap]{}\arrow{l}[swap]{}\\ 
  							     		  &0 \arrow{u}{} 			         & \\ 
\end{tikzcd}
\end{equation}

\begin{equation}\label{Y-h0gen4}
E_2=\begin{tikzcd}
  0 & 0 \arrow{l}[swap]{} &0\arrow[hookrightarrow]{l}[swap]{}\\
0 \arrow{u}{}		  &<c_5-c_4> \arrow{u}{}\arrow{l}[swap]{}                       & <c_5-c_4>\arrow{u}[swap]{}\arrow[hookrightarrow]{l}[swap]{}\\ 
  							     		  &0 \arrow[hookrightarrow]{u}{} 			         & \\ 
\end{tikzcd}\cong \begin{tikzcd}
  0& 0\arrow{l}[swap]{} &0\arrow{l}[swap]{}\\
0 \arrow{u}{}		  &\mathbb{F} \arrow{u}{}\arrow{l}[swap]{}                       & \mathbb{F}\arrow{u}[swap]{}\arrow{l}[swap]{}\\ 
  							     		  &0 \arrow{u}{} 			         & \\ 
\end{tikzcd}
\end{equation}

\begin{equation}\label{Y-h0gen5}
F=\begin{tikzcd}
  0 & <c_6-c_3> \arrow{l}[swap]{} &0\arrow[hookrightarrow]{l}[swap]{}\\
  <c_6-c_3> \arrow{u}{}		  &<c_6-c_3> \arrow{u}{\id}\arrow{l}[swap]{\id}                       &0\arrow{u}[swap]{}\arrow[hookrightarrow]{l}[swap]{}\\ 
  							     		  &<c_6-c_3> \arrow[hookrightarrow]{u}{\id} 			         & \\ 
\end{tikzcd}\cong \begin{tikzcd}
  0& \mathbb{F}\arrow{l}[swap]{} &0\arrow{l}[swap]{}\\
 \mathbb{F} \arrow{u}{}		  &\mathbb{F} \arrow{u}{}\arrow{l}[swap]{}                       & 0\arrow{u}[swap]{}\arrow{l}[swap]{}\\ 
  							     		  &\mathbb{F} \arrow{u}{} 			         & \\ 
\end{tikzcd}
\end{equation}

Note that $M_2$, $E_1$, $E_2$, and $F$ are interval modules, hence by Lemma \ref{thinpoly}, they are indecomposable. Now we are going to show that $M_1$ is indecomposable. By the Fitting's lemma, it suffices to show that $\End(M_1)$ contains no idempotent except $0$ and $\id$. Since any linear transformation $\mathbb{F}\rightarrow \mathbb{F}$ is a scalar multiplication, we use $x:\mathbb{F}\rightarrow \mathbb{F}$ to denote the linear transformation that sends $f$ to $xf$ for all $f\in\mathbb{F}$. Let $\phi\in \End(M_1)$. Note that $\phi$ consists of the following data: $x,y,z: \mathbb{F}\rightarrow \mathbb{F}$ and $A:\mathbb{F}^2\rightarrow \mathbb{F}^2$ ($A= \left[\begin{array}{cc}a & b \\c & d\end{array}\right]$ under the given basis) such that 
\begin{equation}\label{Yx}
x\circ \left[\begin{array}{c}1 \\0\end{array}\right]=A\circ\left[\begin{array}{c}1 \\0\end{array}\right]
\end{equation}
\begin{equation}\label{Yy}
y\circ [1\quad 1]=A\circ\left[\begin{array}{c}1 \\0\end{array}\right]
\end{equation}
\begin{equation}\label{Yz}
z\circ \left[\begin{array}{c}0 \\1\end{array}\right]=A\circ\left[\begin{array}{c}0 \\1\end{array}\right]
\end{equation}

Equation (\ref{Yx}) implies \begin{equation}\label{Yx'}
\left[\begin{array}{c}x \\0\end{array}\right]=\left[\begin{array}{c}a \\c\end{array}\right]
\end{equation}

Hence $a=x$ and $c=0$. Update $A=\left[\begin{array}{cc}x & b \\0 & d\end{array}\right]$.\\

Equation (\ref{Yy}) implies \begin{equation}\label{Yy'}
[x\quad b+d]=[y\quad y]
\end{equation}

Hence $x=y$ and $b+d=y$. Update $A=\left[\begin{array}{cc}x & b \\0 & x-b\end{array}\right]$.\\

Equation (\ref{Yz}) implies \begin{equation}\label{Yz'}
\left[\begin{array}{c}b \\x-b\end{array}\right]=\left[\begin{array}{c}0\\z\end{array}\right]
\end{equation}

Hence $b=0$ and $x=z$. Update $A=\left[\begin{array}{cc}x & 0 \\0 & x\end{array}\right]$. Therefore, $\End(M_1)\cong\mathbb{F}$. Since $\mathbb{F}$ is a field, it contains no idempotent except $0$ and $\id$. Thus $M_1$ is indecomposable.

\end{proof}

\subsection{Decomposition of $PH_i((\mathsf{Star_k})^2_{-,-};\mathbb{F})$, when $k\geq 4$}\label{decompo:star3}

Let $k\geq 4$. In Section \ref{star_Hi'}, we computed $H_i((\hat{\mathsf{Star_k}})^2_{r})$ for all $r, L_{e_1}\in\mathbb{R}_{>0}$ and $i=0,1$. Recall that $\hat{\mathsf{Star_k}}$ is a metric star graph where the length of all but one edge (denoted by $e_1$) is $1$, and the length of edge $e_1$ is $L_{e_1}$. In this section, we define $$(\mathsf{Star}_{k})^2_{r,L_{e_1}}:=(\hat{\mathsf{Star}}_{k})^2_r$$
to emphasize that $(\hat{\mathsf{Star_k}})^2_{r}$ is determined by two parameters $r$ and $L_{e_1}$.\\

To decompose $PH_0((\mathsf{Star_k})^2_{r,L_{e_1}};\mathbb{F})$, we associate $PH_0((\mathsf{Star_k})^2_{r,L_{e_1}};\mathbb{F})$ with a persistence module $N$, where $N$ is constructed using Proposition \ref{prop-sim-2}. We use the following diagram to represent $N$, where each vector space in the diagram represents a constant functor indexed by a chamber of the parameter space of $PH_0((\mathsf{Star_k})^2_{r,L_{e_1}};\mathbb{F})$. For every vector space in the representation $N$, we inductively choose a representative for each basis element according to the geometry of  $(\mathsf{Star_k})^2_{r,L_{e_1}}$: when $k=4$, $e_i$ is the arrangement where the red robot is at the leaf of edge $e_i$ while the green robot is at the leaf of edge $e_1$ for all $i=2,3,4$. Let $\hat{e_i}$ to be the arrangement where the red robot is at the leaf of edge $e_1$ while the green robot is at the leaf of edge $e_i$ for all $i=2,3,4$. Define $f_1, \dots, f_6$ to be the representatives $c_1, \dots, c_6$ in Figure \ref{fig:star-h0geo-select}. Now we assume we have a collection of representatives $\{f_1,\dots, f_{k^2-3k+2}\}$ for the basis of $PH_0((\mathsf{Star_{k-1}})^2_{r,L_{e_1}};\mathbb{F})$. Let $e_i$ be the arrangement where the red robot is at the leaf of edge $e_i$ while the green robot is at the leaf of edge $e_1$ for all $i=2, \dots, k$, and $\hat{e_i}$ is the dual of the above arrangement. The basis we choose for every vector space of the representation $N$ is given in (\ref{h0gen}). \\

\begin{equation}\label{h0gen}
N=\begin{aligned}
&\begin{tikzcd}
  \mathbb{F} & \mathbb{F}^{k^2-3k+4} \arrow{l}[swap]{} &\mathbb{F}^2\arrow{l}[swap]{}\\
  \mathbb{F} \arrow{u}{}		  &\mathbb{F}^{k^2-k} \arrow{u}{}\arrow{l}[swap]{}                       & \mathbb{F}^{2k-2}\arrow{u}[swap]{}\arrow{l}[swap]{}\\ 
  							     		  &\mathbb{F}^{k^2-3k+2} \arrow{u}{} 			         & \\ 
\end{tikzcd}=\begin{tikzcd}
  <f> & <e_1, \hat{e}_1, f_1, \dots, f_{k^2-3k+2}> \arrow{l}[swap]{\alpha} &<e_1, \hat{e}_1>\arrow[hookrightarrow]{l}[swap]{}\\
  <f> \arrow{u}{}		  &\left\langle\parbox{3.5cm}{$e_1, \dots, e_{k-1}, \hat{e}_1,\dots, \hat{e}_{k-1}, \\
  f_1, \dots, f_{k^2-3k+2}$}\right\rangle \arrow{u}{\beta}\arrow{l}[swap]{\gamma}                       & \left\langle\parbox{1.75cm}{$e_1, \dots, e_{k-1},\\ \hat{e}_1,\dots, \hat{e}_{k-1}$}\right\rangle\arrow{u}[swap]{\epsilon}\arrow[hookrightarrow]{l}[swap]{}\\ 
  							     		  &<f_1, \dots, f_{k^2-3k+2}> \arrow[hookrightarrow]{u}{} 			         & \\ 
\end{tikzcd}
\end{aligned}
\end{equation}

Under the given bases for each vector space of $N$, the morphisms in $N$ can be described as follows: $\alpha$ and $\gamma$ maps every basis element to $f$, $\beta$ maps $e_i$ to $e_1$, maps $\hat{e}_i$ to $\hat{e}_1$ for all $i=1,\dots, k-1$ and maps $f_j$ to $f_j$ itself for all $j=1,\dots, k^2-3k+2$, and $\epsilon$ maps $e_i$ to $e_1$ and maps $\hat{e}_i$ to $\hat{e}_1$ for all $i=1,\dots, k-1$. All other maps are inclusions.\\

Similarly, we associate $PH_1((\mathsf{Star_k})^2_{-,-};\mathbb{F})$ with a persistence module $N'$ which can be constructed using Proposition \ref{prop-sim-2}.

\begin{equation}\label{h1gen}
N'=\begin{tikzcd}
   \mathbb{F}^{k(k-3)+1} & 0 \arrow{l}[swap]{} &0\arrow{l}[swap]{}\\
  \mathbb{F}^{k^2-5k+5} \arrow{u}{}		  &0 \arrow{u}{}\arrow{l}[swap]{}                       & 0\arrow{u}[swap]{}\arrow{l}[swap]{}\\ 
  							     		  &0 \arrow{u}{} 			         & \\ 
\end{tikzcd}
\end{equation}

Now we prove Theorem \ref{Thm2}.
\begin{proof}

Let $N'_1=\begin{tikzcd}
   \mathbb{F} & 0 \arrow{l}[swap]{} &0\arrow{l}[swap]{}\\
 0 \arrow{u}{}		  &0 \arrow{u}{}\arrow{l}[swap]{}                       & 0\arrow{u}[swap]{}\arrow{l}[swap]{}\\ 
  							     		  &0 \arrow{u}{} 			         & \\ 
\end{tikzcd}$, $N'_2=\begin{tikzcd}
   0 & 0 \arrow{l}[swap]{} &0\arrow{l}[swap]{}\\
 \mathbb{F} \arrow{u}{}		  &0 \arrow{u}{}\arrow{l}[swap]{}                       & 0\arrow{u}[swap]{}\arrow{l}[swap]{}\\ 
  							     		  &0 \arrow{u}{} 			         & \\ 
\end{tikzcd}$ and $N'_3=\begin{tikzcd}
   \mathbb{F} & 0 \arrow{l}[swap]{} &0\arrow{l}[swap]{}\\
 \mathbb{F} \arrow{u}{\id}		  &0 \arrow{u}{}\arrow{l}[swap]{}                       & 0\arrow{u}[swap]{}\arrow{l}[swap]{}\\ 
  							     		  &0 \arrow{u}{} 			         & \\ 
\end{tikzcd}$. Note that $N'_1$ and $N'_2$ are indecomposable because they are simple representations. On the other hand, note that $N'_3$ is an interval module, hence it is indecomposable.\\

Assume there exists $m,n>0$ such that $M:=\begin{tikzcd}
   \mathbb{F}^{m} & 0 \arrow{l}[swap]{} &0\arrow{l}[swap]{}\\
  \mathbb{F}^{n} \arrow{u}{\alpha}		  &0 \arrow{u}{}\arrow{l}[swap]{}                       & 0\arrow{u}[swap]{}\arrow{l}[swap]{}\\ 
  							     		  &0 \arrow{u}{} 			         & \\ 
\end{tikzcd}$
is indecomposable. Note that the commutativity of the diagram does not provide any non-trivial constraint to the endomorphism ring $\End(M)$; therefore, $\End(M)\cong \End(\mathbb{F}^n\xrightarrow{\alpha} \mathbb{F}^m)$, and $M$ is indecomposable implies $\mathbb{F}^n\xrightarrow{\alpha} \mathbb{F}^m$ indecomposable. Note that $\mathbb{F}^n\xrightarrow{\alpha} \mathbb{F}^m$ is an indecomposable $A_2$-quiver and $m,n>0$, we must have $m=n=1$. Therefore, every indecomposable subrepresentation of $PH_1((\mathsf{Star_k})^2_{-,-};\mathbb{F})$ is isomorphic to $N'_1$, $N'_2$, or $N'_3$. 
\end{proof}

Theorem \ref{Thm3} is a consequence of this proposition. 
\begin{proposition}
$PH_0((\hat{\mathsf{Star_k}})^2_{-,-};\mathbb{F})\cong M_1\oplus M_2 \oplus (\bigoplus\limits_{i=2}^{k-1} E_i) \oplus (\bigoplus\limits_{i=2}^{k-1} \hat{E}_i)\oplus (\bigoplus\limits_{j=2}^{k^2-3k+2} F_j)$ where

\begin{minipage}{0.3\textwidth}
  $$M_1=\begin{tikzcd}
  \mathbb{F} & \mathbb{F}^2 \arrow{l}[swap]{\left[1\quad 1\right]}& \mathbb{F} \arrow[l, swap, "\spmat{1 \\0}"] \\
  \mathbb{F} \arrow{u}{\id}		  &\mathbb{F}^2 \arrow{u}{\id}\arrow{l}[swap]{\left[1\quad 1\right]}                       & \mathbb{F}\arrow{u}[swap]{\id}\arrow[l,"\spmat{1 \\0}"]\\ 
  							     		  &\mathbb{F}\arrow[u,"\spmat{0 \\1}"]			         & \\ 
\end{tikzcd}
$$
\end{minipage}%
\begin{minipage}{0.3\textwidth}
  $$M_2=\begin{tikzcd}
  0 & \mathbb{F} \arrow{l}{}& \mathbb{F} \arrow{l}[swap]{\id} \\
 0	 \arrow{u}	  &\mathbb{F} \arrow{u}{\id}\arrow{l}[swap]{}                       & \mathbb{F}\arrow{u}[swap]{\id}\arrow{l}[swap]{\id}\\ 
  							     		  &0 \arrow{u}			         & \\ 
\end{tikzcd}
$$
\end{minipage}%
\begin{minipage}{0.3\textwidth}
  $$E_i=\begin{tikzcd}
  0 & 0 \arrow{l}{}& 0 \arrow{l}\\
 0	 \arrow{u}	  &\mathbb{F} \arrow{u}\arrow{l}[swap]{}                       & \mathbb{F}\arrow{u}\arrow{l}[swap]{\id}\\ 
  							     		  &0 \arrow{u}			         & \\ 
\end{tikzcd}
$$
\end{minipage}

\begin{minipage}{0.3\textwidth}
  $$\hat{E}_i=\begin{tikzcd}
  0 & 0 \arrow{l}{}& 0 \arrow{l}\\
 0	 \arrow{u}	  &\mathbb{F} \arrow{u}\arrow{l}[swap]{}                       & \mathbb{F}\arrow{u}\arrow{l}[swap]{\id}\\ 
  							     		  &0 \arrow{u}			         & \\ 
\end{tikzcd}
$$
\end{minipage}%
\begin{minipage}{0.3\textwidth}
$$F_j=\begin{tikzcd}
  0 & \mathbb{F} \arrow{l}{}& 0\arrow{l}[swap]{} \\
 0	 \arrow{u}	  &\mathbb{F} \arrow{u}{\id}\arrow{l}[swap]{}                       & 0\arrow{u}[swap]{}\arrow{l}[swap]{}\\ 
  							     		  & \mathbb{F}\arrow{u}{\id}			         & \\ 
\end{tikzcd}
$$
\end{minipage}

for all $i=2, \dots, k-1$ and $j=2, \dots, k^2-3k+2$.
\end{proposition}
\begin{proof}
The basis we choose for each vector space in the representation $N$ provided in the diagram below:\\ 

$$
N=\begin{tikzcd}
  <f> &\left\langle\parbox{6cm}{\centering $ e_1, \hat{e}_1-e_1, f_1, f_2-f_1, \dots, f_{k^2-3k+2}-f_1$}\right\rangle \arrow{l}[swap]{\alpha} &<e_1, \hat{e}_1-e_1>\arrow[hookrightarrow]{l}[swap]{}\\
  <f> \arrow{u}{}		  &\left\langle\parbox{5cm}{\centering $e_1, e_2-e_1,\dots,e_{k-1}-e_1, \hat{e}_1-e_1, \newline\hat{e}_2-\hat{e}_1, \dots, \hat{e}_{k-1}-\hat{e}_1, f_1, f_2-f_1, \newline f_3-f_1, \dots, f_{k^2-3k+2}-f_1 $}\right\rangle\arrow{u}{\beta}\arrow{l}[swap]{\gamma}                       & \left\langle\parbox{4.5cm}{\centering $e_1, e_2-e_1,\dots,e_{k-1}-e_1,  \newline\hat{e}_1-e_1,\hat{e}_2-\hat{e}_1, \dots, \hat{e}_{k-1}-\hat{e}_1$}\right\rangle\arrow{u}[swap]{\epsilon}\arrow[hookrightarrow]{l}[swap]{}\\ 
  							     		  &<f_1, f_2-f_1,\dots, f_{k^2-3k+2}-f_1> \arrow[hookrightarrow]{u}{} 			         & \\ 
\end{tikzcd}
$$
Let 
\begin{equation}\label{m1}
\begin{aligned}
M_1:&=\begin{tikzcd}
  <f_1> & <e_1, f_1> \arrow{l}[swap]{\alpha} &<e_1>\arrow[hookrightarrow]{l}[swap]{}\\
  <f_1> \arrow{u}{}		  &<e_1, f_1>\arrow{u}{\beta}\arrow{l}[swap]{\gamma}                       & <e_1>\arrow{u}[swap]{\epsilon}\arrow[hookrightarrow]{l}[swap]{}\\ 
  							     		  &<f_1> \arrow[hookrightarrow]{u}{} 			         & \\ 
\end{tikzcd}\cong \begin{tikzcd}
  \mathbb{F} & \mathbb{F}^2 \arrow{l}[swap]{\left[1\quad 1\right]}& \mathbb{F} \arrow[l, swap, "\spmat{1 \\0}"] \\
  \mathbb{F} \arrow{u}{\id}		  &\mathbb{F}^2 \arrow{u}{\id}\arrow{l}[swap]{\left[1\quad 1\right]}                       & \mathbb{F}\arrow{u}[swap]{\id}\arrow[l,"\spmat{1 \\0}"]\\ 
  							     		  &\mathbb{F}\arrow[u,"\spmat{0 \\1}"]			         & \\ 
\end{tikzcd}
\end{aligned}
\end{equation}

\begin{equation}\label{m2}
\begin{aligned}
M_2:&=\begin{tikzcd}
  0 & <\hat{e}_1-e_1> \arrow{l}[swap]{\alpha} &<\hat{e}_1-e_1>\arrow[hookrightarrow]{l}[swap]{}\\
  0 \arrow{u}{}		  &<\hat{e}_1-e_1>\arrow{u}{\beta}\arrow{l}[swap]{\gamma}                       & <\hat{e}_1-e_1>\arrow{u}[swap]{\epsilon}\arrow[hookrightarrow]{l}[swap]{}\\ 
  							     		  &0 \arrow[hookrightarrow]{u}{} 			         & \\ 
\end{tikzcd}\cong \begin{tikzcd}
  0 & \mathbb{F} \arrow{l}{}& \mathbb{F} \arrow{l}[swap]{\id} \\
 0	 \arrow{u}	  &\mathbb{F} \arrow{u}{\id}\arrow{l}[swap]{}                       & \mathbb{F}\arrow{u}[swap]{\id}\arrow{l}[swap]{\id}\\ 
  							     		  &0 \arrow{u}			         & \\ 
\end{tikzcd}
\end{aligned}
\end{equation}
For $i=2,\dots, k-1$
\begin{equation}\label{Ei}
\begin{aligned}
E_i:&=\begin{tikzcd}
  0 & 0\arrow{l}[swap]{\alpha} &0\arrow[hookrightarrow]{l}[swap]{}\\
  0 \arrow{u}{}		  &<e_i-e_1>\arrow{u}{\beta}\arrow{l}[swap]{\gamma}                       & <e_i-e_1>\arrow{u}[swap]{\epsilon}\arrow[hookrightarrow]{l}[swap]{}\\ 
  							     		  &0 \arrow[hookrightarrow]{u}{} 			         & \\ 
\end{tikzcd}\cong \begin{tikzcd}
  0 & 0 \arrow{l}{}& 0 \arrow{l}[swap]{} \\
 0	 \arrow{u}	  &\mathbb{F} \arrow{u}{}\arrow{l}[swap]{}                       & \mathbb{F}\arrow{u}[swap]{}\arrow{l}[swap]{\id}\\ 
  							     		  &0 \arrow{u}			         & \\ 
\end{tikzcd}
\end{aligned}
\end{equation}

\begin{equation}\label{Ei-hat}
\begin{aligned}
\hat{E}_i:&=\begin{tikzcd}
  0 & 0\arrow{l}[swap]{\alpha} &0\arrow[hookrightarrow]{l}[swap]{}\\
  0 \arrow{u}{}		  &<\hat{e}_i-\hat{e}_1>\arrow{u}{\beta}\arrow{l}[swap]{\gamma}                       & <\hat{e}_i-\hat{e}_1>\arrow{u}[swap]{\epsilon}\arrow[hookrightarrow]{l}[swap]{}\\ 
  							     		  &0 \arrow[hookrightarrow]{u}{} 			         & \\ 
\end{tikzcd}\cong \begin{tikzcd}
  0 & 0 \arrow{l}{}& 0 \arrow{l}[swap]{} \\
 0	 \arrow{u}	  &\mathbb{F} \arrow{u}{}\arrow{l}[swap]{}                       & \mathbb{F}\arrow{u}[swap]{}\arrow{l}[swap]{\id}\\ 
  							     		  &0 \arrow{u}			         & \\ 
\end{tikzcd}
\end{aligned}
\end{equation}

For $j=2,\dots, k^2-3k+2$
\begin{equation}\label{Ei}
\begin{aligned}
F_j:&=\begin{tikzcd}
  0 & <f_j-f_1>\arrow{l}[swap]{\alpha} &0\arrow[hookrightarrow]{l}[swap]{}\\
  0 \arrow{u}{}		  &<f_j-f_1>\arrow{u}{\beta}\arrow{l}[swap]{\gamma}                       & 0\arrow{u}[swap]{\epsilon}\arrow[hookrightarrow]{l}[swap]{}\\ 
  							     		  &<f_j-f_1> \arrow[hookrightarrow]{u}{} 			         & \\ 
\end{tikzcd}\cong \begin{tikzcd}
  0 & \mathbb{F} \arrow{l}{}& 0 \arrow{l}[swap]{} \\
 0	 \arrow{u}	  &\mathbb{F} \arrow{u}{\id}\arrow{l}[swap]{}                       & 0\arrow{u}[swap]{}\arrow{l}[swap]{}\\ 
  							     		  &\mathbb{F} \arrow{u}{\id}			         & \\ 
\end{tikzcd}
\end{aligned}
\end{equation}

Note that $M_2$, $E_i$, $\hat{E}_i$, and $F_j$ (where $i=2,\dots, k-1$ and $j=2,\dots, k^2-3k+2$) are interval modules, hence by Lemma \ref{thinpoly}, they are indecomposable. The indecomposability of $M_1$ can be proved using a similar argument in Proposition \ref{Y-h0gen-N}.
\end{proof}

%% file: chapters/chapter4-H.tex
\section{Configuration Spaces of the H Graph}\label{conf:H}\index{H graph}
In this section, we calculate the second configuration spaces of the H graph with restraint parameters $r$ and $L_{e_1}$. To assign each edge of the H Graph an orientation that matches the orientation of the $Y$ graph, we subdivide the edge $e_1$ of the H graph into two edges $e_6$ and $e_7$, as shown in Figure \ref{fig:Hgraph'}. We denote the geometric realization of the resulting graph by $\mathcal{H}$ and use $\hat{\mathcal{H}}$ to denote the metric graph $\mathcal{H}$ where the length of the bridge is $L_{e_1}$ and the length of any other edge is $1$. \\

\begin{figure}[htbp!]
     \centering
     \begin{subfigure}[b]{0.45\textwidth}
         \centering
         \resizebox{4cm}{!}{

\tikzset{every picture/.style={line width=0.75pt}} 

\begin{tikzpicture}[x=0.75pt,y=0.75pt,yscale=-1,xscale=1]

\draw    (40,32.7) -- (83.99,76.69) ;
\draw [shift={(83.99,76.69)}, rotate = 45] [color={rgb, 255:red, 0; green, 0; blue, 0 }  ][fill={rgb, 255:red, 0; green, 0; blue, 0 }  ][line width=0.75]      (0, 0) circle [x radius= 3.35, y radius= 3.35]   ;
\draw [shift={(40,32.7)}, rotate = 45] [color={rgb, 255:red, 0; green, 0; blue, 0 }  ][fill={rgb, 255:red, 0; green, 0; blue, 0 }  ][line width=0.75]      (0, 0) circle [x radius= 3.35, y radius= 3.35]   ;
\draw    (83.99,76.69) -- (41.63,119.05) ;
\draw [shift={(41.63,119.05)}, rotate = 135] [color={rgb, 255:red, 0; green, 0; blue, 0 }  ][fill={rgb, 255:red, 0; green, 0; blue, 0 }  ][line width=0.75]      (0, 0) circle [x radius= 3.35, y radius= 3.35]   ;
\draw [shift={(83.99,76.69)}, rotate = 135] [color={rgb, 255:red, 0; green, 0; blue, 0 }  ][fill={rgb, 255:red, 0; green, 0; blue, 0 }  ][line width=0.75]      (0, 0) circle [x radius= 3.35, y radius= 3.35]   ;
\draw    (145.9,76.69) -- (189.4,33.19) ;
\draw [shift={(189.4,33.19)}, rotate = 315] [color={rgb, 255:red, 0; green, 0; blue, 0 }  ][fill={rgb, 255:red, 0; green, 0; blue, 0 }  ][line width=0.75]      (0, 0) circle [x radius= 3.35, y radius= 3.35]   ;
\draw [shift={(145.9,76.69)}, rotate = 315] [color={rgb, 255:red, 0; green, 0; blue, 0 }  ][fill={rgb, 255:red, 0; green, 0; blue, 0 }  ][line width=0.75]      (0, 0) circle [x radius= 3.35, y radius= 3.35]   ;
\draw    (145.9,76.69) -- (189.89,120.68) ;
\draw [shift={(189.89,120.68)}, rotate = 45] [color={rgb, 255:red, 0; green, 0; blue, 0 }  ][fill={rgb, 255:red, 0; green, 0; blue, 0 }  ][line width=0.75]      (0, 0) circle [x radius= 3.35, y radius= 3.35]   ;
\draw [shift={(145.9,76.69)}, rotate = 45] [color={rgb, 255:red, 0; green, 0; blue, 0 }  ][fill={rgb, 255:red, 0; green, 0; blue, 0 }  ][line width=0.75]      (0, 0) circle [x radius= 3.35, y radius= 3.35]   ;
\draw    (83.99,76.69) -- (145.9,76.69) ;
\draw [shift={(145.9,76.69)}, rotate = 0] [color={rgb, 255:red, 0; green, 0; blue, 0 }  ][fill={rgb, 255:red, 0; green, 0; blue, 0 }  ][line width=0.75]      (0, 0) circle [x radius= 3.35, y radius= 3.35]   ;
\draw [shift={(83.99,76.69)}, rotate = 0] [color={rgb, 255:red, 0; green, 0; blue, 0 }  ][fill={rgb, 255:red, 0; green, 0; blue, 0 }  ][line width=0.75]      (0, 0) circle [x radius= 3.35, y radius= 3.35]   ;

\draw (106,81.4) node [anchor=north west][inner sep=0.75pt]    {$e_{1}$};
\draw (58,31.4) node [anchor=north west][inner sep=0.75pt]    {$e_{2}$};
\draw (155,31.4) node [anchor=north west][inner sep=0.75pt]    {$e_{4}$};
\draw (58,101) node [anchor=north west][inner sep=0.75pt]    {$e_{3}$};
\draw (156,101) node [anchor=north west][inner sep=0.75pt]    {$e_{5}$};

\end{tikzpicture}
}
\caption{}
\label{fig:Hgraph}
     \end{subfigure}
     \hfill
     \begin{subfigure}[b]{0.45\textwidth}
         \centering
        \resizebox{3.5cm}{!}{

\tikzset{every picture/.style={line width=0.75pt}} 

\begin{tikzpicture}[x=0.75pt,y=0.75pt,yscale=-1,xscale=1]

\draw    (40,32.7) -- (83.99,76.69) ;
\draw [shift={(83.99,76.69)}, rotate = 45] [color={rgb, 255:red, 0; green, 0; blue, 0 }  ][fill={rgb, 255:red, 0; green, 0; blue, 0 }  ][line width=0.75]      (0, 0) circle [x radius= 3.35, y radius= 3.35]   ;
\draw [shift={(57.4,50.1)}, rotate = 45] [fill={rgb, 255:red, 0; green, 0; blue, 0 }  ][line width=0.08]  [draw opacity=0] (8.93,-4.29) -- (0,0) -- (8.93,4.29) -- cycle    ;
\draw [shift={(40,32.7)}, rotate = 45] [color={rgb, 255:red, 0; green, 0; blue, 0 }  ][fill={rgb, 255:red, 0; green, 0; blue, 0 }  ][line width=0.75]      (0, 0) circle [x radius= 3.35, y radius= 3.35]   ;
\draw    (83.99,76.69) -- (41.63,119.05) ;
\draw [shift={(41.63,119.05)}, rotate = 135] [color={rgb, 255:red, 0; green, 0; blue, 0 }  ][fill={rgb, 255:red, 0; green, 0; blue, 0 }  ][line width=0.75]      (0, 0) circle [x radius= 3.35, y radius= 3.35]   ;
\draw [shift={(59.27,101.4)}, rotate = 315] [fill={rgb, 255:red, 0; green, 0; blue, 0 }  ][line width=0.08]  [draw opacity=0] (8.93,-4.29) -- (0,0) -- (8.93,4.29) -- cycle    ;
\draw [shift={(83.99,76.69)}, rotate = 135] [color={rgb, 255:red, 0; green, 0; blue, 0 }  ][fill={rgb, 255:red, 0; green, 0; blue, 0 }  ][line width=0.75]      (0, 0) circle [x radius= 3.35, y radius= 3.35]   ;
\draw    (145.9,76.69) -- (189.4,33.19) ;
\draw [shift={(189.4,33.19)}, rotate = 315] [color={rgb, 255:red, 0; green, 0; blue, 0 }  ][fill={rgb, 255:red, 0; green, 0; blue, 0 }  ][line width=0.75]      (0, 0) circle [x radius= 3.35, y radius= 3.35]   ;
\draw [shift={(171.18,51.4)}, rotate = 135] [fill={rgb, 255:red, 0; green, 0; blue, 0 }  ][line width=0.08]  [draw opacity=0] (8.93,-4.29) -- (0,0) -- (8.93,4.29) -- cycle    ;
\draw [shift={(145.9,76.69)}, rotate = 315] [color={rgb, 255:red, 0; green, 0; blue, 0 }  ][fill={rgb, 255:red, 0; green, 0; blue, 0 }  ][line width=0.75]      (0, 0) circle [x radius= 3.35, y radius= 3.35]   ;
\draw    (145.9,76.69) -- (189.89,120.68) ;
\draw [shift={(189.89,120.68)}, rotate = 45] [color={rgb, 255:red, 0; green, 0; blue, 0 }  ][fill={rgb, 255:red, 0; green, 0; blue, 0 }  ][line width=0.75]      (0, 0) circle [x radius= 3.35, y radius= 3.35]   ;
\draw [shift={(171.43,102.22)}, rotate = 225] [fill={rgb, 255:red, 0; green, 0; blue, 0 }  ][line width=0.08]  [draw opacity=0] (8.93,-4.29) -- (0,0) -- (8.93,4.29) -- cycle    ;
\draw [shift={(145.9,76.69)}, rotate = 45] [color={rgb, 255:red, 0; green, 0; blue, 0 }  ][fill={rgb, 255:red, 0; green, 0; blue, 0 }  ][line width=0.75]      (0, 0) circle [x radius= 3.35, y radius= 3.35]   ;
\draw    (83.99,76.69) -- (115,76.69) ;
\draw [shift={(115,76.69)}, rotate = 0] [color={rgb, 255:red, 0; green, 0; blue, 0 }  ][fill={rgb, 255:red, 0; green, 0; blue, 0 }  ][line width=0.75]      (0, 0) circle [x radius= 3.35, y radius= 3.35]   ;
\draw [shift={(104.49,76.69)}, rotate = 180] [fill={rgb, 255:red, 0; green, 0; blue, 0 }  ][line width=0.08]  [draw opacity=0] (8.93,-4.29) -- (0,0) -- (8.93,4.29) -- cycle    ;
\draw [shift={(83.99,76.69)}, rotate = 0] [color={rgb, 255:red, 0; green, 0; blue, 0 }  ][fill={rgb, 255:red, 0; green, 0; blue, 0 }  ][line width=0.75]      (0, 0) circle [x radius= 3.35, y radius= 3.35]   ;
\draw    (115,76.69) -- (146.01,76.69) ;
\draw [shift={(146.01,76.69)}, rotate = 0] [color={rgb, 255:red, 0; green, 0; blue, 0 }  ][fill={rgb, 255:red, 0; green, 0; blue, 0 }  ][line width=0.75]      (0, 0) circle [x radius= 3.35, y radius= 3.35]   ;
\draw [shift={(124.01,76.69)}, rotate = 0] [fill={rgb, 255:red, 0; green, 0; blue, 0 }  ][line width=0.08]  [draw opacity=0] (8.93,-4.29) -- (0,0) -- (8.93,4.29) -- cycle    ;
\draw [shift={(115,76.69)}, rotate = 0] [color={rgb, 255:red, 0; green, 0; blue, 0 }  ][fill={rgb, 255:red, 0; green, 0; blue, 0 }  ][line width=0.75]      (0, 0) circle [x radius= 3.35, y radius= 3.35]   ;

\draw (92,83) node [anchor=north west][inner sep=0.75pt]    {$e_{6}$};
\draw (58,31.4) node [anchor=north west][inner sep=0.75pt]    {$e_{2}$};
\draw (155,31.4) node [anchor=north west][inner sep=0.75pt]    {$e_{4}$};
\draw (60,103) node [anchor=north west][inner sep=0.75pt]    {$e_{3}$};
\draw (154,103) node [anchor=north west][inner sep=0.75pt]    {$e_{5}$};
\draw (123,83) node [anchor=north west][inner sep=0.75pt]    {$e_{7}$};

\end{tikzpicture}
}
\caption{}
\label{fig:Hgraph'}     
\end{subfigure}
        \caption{}
\end{figure}

Let $L_{e_i}$ denote the length of the $e_i$. Note that $L_{e_6}+L_{e_7}=L_{e_1}$ and $$L_{e_i} = \begin{cases} L_{e_i}, & \mbox{if } i=6,7 \\ 1, & \mbox{if } i\in\{2,3,4,5\} \end{cases}
$$
The parametric polytope of $\mathcal{H}^2_{r,\vec{L}}$ is given by the following inequalities, where $\vec{L}=(L_{e_2},\dots, L_{e_7})$:\\

For $i,j\in\{2,3,4,5,6,7\}$,
\begin{equation}\label{eq:H-01}
\begin{aligned}
0\leq x&\leq L_{e_i}\\
0\leq y&\leq L_{e_j}\\
\end{aligned}
\end{equation}

For $i,j\in\{2,3\}$ or $i,j\in\{4,5\}$, and $i\neq j$,
\begin{equation}
\begin{aligned}
x+y&\geq r
\end{aligned}
\end{equation}


For $i\in\{2,3\}$ and $j=6$ (Respectively, $j\in\{2,3\}$, $i=6$),
\begin{equation}\label{eq:H-02}
\begin{aligned}
x+y&\geq r
\end{aligned}
\end{equation}


For $i\in\{2,3\}$, $j=7$ (Respectively, $i\in\{4,5\}$, $j=6$),
\begin{equation}\label{eq:H-03''}
\begin{aligned}
x+L_{e_1}-y&\geq r
\end{aligned}
\end{equation}

For $j\in\{2,3\}$, $i=7$ (Respectively, $j\in\{4,5\}$, $i=6$),
\begin{equation}\label{eq:H-03‘}
\begin{aligned}
-x+L_{e_1}+y&\geq r
\end{aligned}
\end{equation}

For $i,j\in\{6,7\}$, $i\neq j$,
\begin{equation}\label{eq:H-04}
\begin{aligned}
-x+L_{e_1}-y&\geq r
\end{aligned}
\end{equation}

For $i\in\{2,3\}$, $j\in\{4,5\}$ (Respectively, $i\in\{4,5\}$, $j\in\{2,3\}$),
\begin{equation}\label{eq:H-05}
\begin{aligned}
x+L_{e_1}+y&\geq r
\end{aligned}
\end{equation}


For $i\in\{4,5\}$, $j=7$ (Respectively, $j\in\{4,5\}$, $i=7$),
\begin{equation}\label{eq:H-05''}
\begin{aligned}
x+y&\geq r
\end{aligned}
\end{equation}





For $i,j\in\{2,3,4,5,6,7\}$ and $i=j$,
\begin{equation}\label{eq:H-06'}
\begin{aligned}
\lvert x-y\rvert&\geq r
\end{aligned}
\end{equation}

When $\vec{L}=(1,\dots,1)$ and $L_{e_6}=L_{e_7}=\frac{1}{2}L_{e_1}$, the critical hyperplanes in the parameter space are: $r=1$, $r=2$,$r=\frac{1}{2}L_{e_1}$, $r=1+\frac{1}{2}L_{e_1}$, $r=1+L_{e_1}$, $r=L_{e_1}$, and $r=2+L_{e_1}$.\\

To calculate $H_i((\hat{\mathcal{H}})_{r,L_{e_1}}^2)$, we use the following cover of $(\hat{\mathcal{H}})_{r,L_{e_1}}^2$:
\begin{equation}\label{cover:Hgraph}
\begin{aligned}
U_{11}&=(\hat{\mathsf{Star}}_{3})^2_r\\
U_{12}&=\hat{\mathsf{Star}}_{3}\times\hat{\mathsf{Star}}_{3}-\{(x,y)\in \hat{\mathsf{Star}}_{3}\times \hat{\mathsf{Star}}_{3}\mid \delta(x,y)<r\}\\
U_{21}&=\hat{\mathsf{Star}}_{3}\times\hat{\mathsf{Star}}_{3}-\{(y,x)\in \hat{\mathsf{Star}}_{3}\times \hat{\mathsf{Star}}_{3}\mid \delta(x,y)<r\}\\
U_{22}&=(\hat{\mathsf{Star}}_{3})_r^2
\end{aligned}
\end{equation}
Note that there are natural inclusions:\\
\begin{equation}\label{H:eq0}
\begin{aligned}
U_{11}\cap U_{12}&\hookrightarrow U_{11}\hookleftarrow U_{11}\cap U_{21}\\
U_{11}\cap U_{12}&\hookrightarrow U_{12}\hookleftarrow U_{22}\cap U_{12}\\
U_{11}\cap U_{21}&\hookrightarrow U_{21}\hookleftarrow U_{22}\cap U_{21}\\
U_{22}\cap U_{12}&\hookrightarrow U_{22}\hookleftarrow U_{22}\cap U_{21}
\end{aligned}
\end{equation}

\begin{lemma}\label{ch4:lemmaH}
$\vec{L}=(1,\dots,1)$ and $L_{e_6}=L_{e_7}=\frac{1}{2}L_{e_1}$. Then 
$$H_0(\hat{\mathsf{Star}}_{3}\times\hat{\mathsf{Star}}_{3}-\{(x,y)\in \hat{\mathsf{Star}}_{3}\times \hat{\mathsf{Star}}_{3}\mid \delta(x,y)<r\})\cong  \begin{cases} \mathbb{Z}, & \mbox{if } 0<r\leq L_{e_1}+1 \\ \mathbb{Z}^4, & \mbox{if } L_{e_1}+1 <r\leq L_{e_1}+2 \\
0, & \mbox{else} \end{cases}$$
and 
$$H_1(\hat{\mathsf{Star}}_{3}\times\hat{\mathsf{Star}}_{3}-\{(x,y)\in \hat{\mathsf{Star}}_{3}\times \hat{\mathsf{Star}}_{3}\mid \delta(x,y)<r\})\cong  \begin{cases} \mathbb{Z}, & \mbox{if } L_{e_1}<r\leq L_{e_1}+1 \\0, & \mbox{else} \end{cases}$$
\end{lemma}

\begin{proof}
Consider the following cover of $\hat{\mathsf{Star}}_{3}\times\hat{\mathsf{Star}}_{3}-\{(x,y)\in \hat{\mathsf{Star}}_{3}\times \hat{\mathsf{Star}}_{3}\mid \delta(x,y)<r\}$:
\begin{equation}
\begin{aligned}
V_1&=\hat{\mathsf{Star}}_{3}\times e_7-\{(x,y)\in \hat{\mathsf{Star}}_{3}\times e_7\mid \delta(x,y)<r\}\\
V_2&=\hat{\mathsf{Star}}_{3}\times (e_4\vee e_5)-\{(x,y)\in \hat{\mathsf{Star}}_{3}\times (e_4\vee e_5)\mid \delta(x,y)<r\}
\end{aligned}
\end{equation}

Note that
\begin{equation}
\begin{aligned}
V_1&\simeq \begin{cases} \{\ast\}, & \mbox{if }0<r\leq L_{e_1} \\ \{\ast_1,\ast_2\}, & \mbox{if} L_{e_1}<r\leq L_{e_1}+1\\ \emptyset & \mbox{else}\end{cases}\\
V_2&\simeq \begin{cases} \{\ast\}, & \mbox{if }0<r\leq L_{e_1} \\ S^1, & \mbox{if} L_{e_1}<r\leq L_{e_1}+1\\ \{\ast_1,\ast_2, \ast_3,\ast_4\} , & \mbox{if} L_{e_1}+1<r\leq L_{e_1}+2\\ \emptyset & \mbox{else}\end{cases}\\
\end{aligned}
\end{equation}
and 
\begin{equation}
\begin{aligned}
V_1\cap V_2&\simeq \begin{cases} \{\ast\}, & \mbox{if }0<r\leq L_{e_1} \\ \{\ast_1,\ast_2\}, & \mbox{if} L_{e_1}<r\leq L_{e_1}+1\\ \emptyset & \mbox{else}\end{cases}
\end{aligned}
\end{equation}

Note that $d^1$ is induced from the inclusions $$V_{1}\cap V_{2}\hookrightarrow V_{1}$$ $$V_{1}\cap V_{2}\hookrightarrow V_{2}$$

\paragraph{When $0<r\leq L_{e_1}$:} The $E^1$ page of the Mayer-Vietoris spectral sequence is \\
\begin{figure}[htbp!]
\centering
\begin{tikzpicture}
  \matrix (m) [matrix of math nodes,
    nodes in empty cells,nodes={minimum width=5ex,
     minimum height=2.5ex,outer sep=-5pt},
    column sep=1ex,row sep=1ex]{
          q     &      &     &  &  \\
          2     & \vdots & \vdots  & \vdots  &\\   
          1     &  0 &  0  & 0 & \\
          0     &  \mathbb{Z}^2 & \mathbb{Z} & 0&\\
    \quad\strut &   0  &  1  &  2 & p\strut \\};
    \draw[-stealth] (m-4-3.west) -- (m-4-2.east) node [above,pos=0.5] {$d^1$};
\draw[thick] (m-1-1.east) -- (m-5-1.east) ;
\draw[thick] (m-5-1.north) -- (m-5-5.north) ;
\end{tikzpicture}
\end{figure}

Let $a$ be the generator of $H_0(V_{1}\cap V_{2})$. To understand the behavior of $d^1$, we choose a collection of generators for every abelian group in the $E^1$ page of the spectral sequence. Let $b$ be a generator of $H_0(V_{1})$ and $c$ be a generator of $H_0(V_{2})$ such that $d^1(a)=b+c$. Thus $\dim \Ima d^1=1$ and $\dim \ker d^1=0$, and we obtain the $E^2$ page of the spectral sequence: 
\begin{figure}[htbp!]
\centering
\begin{tikzpicture}
  \matrix (m) [matrix of math nodes,
    nodes in empty cells,nodes={minimum width=5ex,
     minimum height=2.5ex,outer sep=-5pt},
    column sep=1ex,row sep=1ex]{
          q     &      &     &  &  \\
          2     & \vdots & \vdots  & \vdots  &\\   
          1     &  0 &  0  & 0 &\\
          0     &  \mathbb{Z} & 0& 0&\\
    \quad\strut &   0  &  1  &  2 & p\strut \\};
\draw[thick] (m-1-1.east) -- (m-5-1.east) ;
\draw[thick] (m-5-1.north) -- (m-5-5.north) ;
\end{tikzpicture}
\end{figure}

Since there is no non-trivial arrow on the $E^2$ page, $E^2=E^{\infty}$. Therefore, 
$$H_0(\hat{\mathsf{Star}}_{3}\times\hat{\mathsf{Star}}_{3}-\{(x,y)\in \hat{\mathsf{Star}}_{3}\times \hat{\mathsf{Star}}_{3}\mid \delta(x,y)<r\})\cong\mathbb{Z}$$
and
$$H_1(\hat{\mathsf{Star}}_{3}\times\hat{\mathsf{Star}}_{3}-\{(x,y)\in \hat{\mathsf{Star}}_{3}\times \hat{\mathsf{Star}}_{3}\mid \delta(x,y)<r\})=0$$

\paragraph{When $L_{e_1}<r\leq L_{e_1}+1$:} The $E^1$ page of the Mayer-Vietoris spectral sequence is \\
\begin{figure}[htbp!]
\centering
\begin{tikzpicture}
  \matrix (m) [matrix of math nodes,
    nodes in empty cells,nodes={minimum width=5ex,
     minimum height=2.5ex,outer sep=-5pt},
    column sep=1ex,row sep=1ex]{
          q     &      &     &  &  \\
          3    &   \vdots   & \vdots    &\vdots  &  \\
          2     & 0& 0  &0  &\\   
          1     &   \mathbb{Z} &  0  & 0 &\\
          0     &  \mathbb{Z}^3 & \mathbb{Z}^2 & 0&\\
    \quad\strut &   0  &  1  &  2 & p\strut \\};
    \draw[-stealth] (m-5-3.west) -- (m-5-2.east) node [above,pos=0.5] {$d^1$};
\draw[thick] (m-1-1.east) -- (m-6-1.east) ;
\draw[thick] (m-6-1.north) -- (m-6-5.north) ;
\end{tikzpicture}
\end{figure}

To understand the behavior of $d^1$, we choose a collection of generators for every abelian group in the $E^1$ page of the spectral sequence. Let $a_1, a_2$ be the generators of $H_0(V_{1}\cap V_{2})$. Choose generators $b_1$ and $b_2$ for $H_0(V_{1})$ and a generator $c$ for $H_0(V_{2})$ such that $b_1$ lies in the path component that contains $a_1$, and $b_2$ lies in the path component that contains $a_2$. Hence $$d^1(a_1)=b_1+c \mbox{\quad and\quad}d^1(a_2)=b_2+c$$

Under the given generators, $d^1$ can be represented by the following matrix 
$$\left[\begin{array}{cc}1 & 0 \\0 & 1 \\1 & 1\end{array}\right]$$

Run the row reduction, we obtain:
$$\left[\begin{array}{cc}1 & 0 \\0 & 1 \\0 & 0\end{array}\right]$$

The rank of the matrix is $2$, and the diagonal elements of the Smith Normal form of the above matrix are all equal to $1$. Thus $\dim \Ima d^1=2$ and $\dim \ker d^1=1$, and we obtain the $E^2$ page of the spectral sequence: 
\begin{figure}[htbp!]
\centering
\begin{tikzpicture}
  \matrix (m) [matrix of math nodes,
    nodes in empty cells,nodes={minimum width=5ex,
     minimum height=2.5ex,outer sep=-5pt},
    column sep=1ex,row sep=1ex]{
          q     &      &     &  &  \\
          3     & \vdots & \vdots  & \vdots  &\\   
          2     & 0& 0  & 0  &\\  
          1     &  \mathbb{Z} &  0  & 0 & \\
          0     &  \mathbb{Z} & 0& 0&\\
    \quad\strut &   0  &  1  &  2 & p\strut \\};
\draw[thick] (m-1-1.east) -- (m-6-1.east) ;
\draw[thick] (m-6-1.north) -- (m-6-5.north) ;
\end{tikzpicture}
\end{figure}

Since there is no non-trivial arrow on the $E^2$ page, $E^2=E^{\infty}$. Therefore, 
$$H_0(\hat{\mathsf{Star}}_{3}\times\hat{\mathsf{Star}}_{3}-\{(x,y)\in \hat{\mathsf{Star}}_{3}\times \hat{\mathsf{Star}}_{3}\mid \delta(x,y)<r\})\cong\mathbb{Z}$$
and
$$H_1(\hat{\mathsf{Star}}_{3}\times\hat{\mathsf{Star}}_{3}-\{(x,y)\in \hat{\mathsf{Star}}_{3}\times \hat{\mathsf{Star}}_{3}\mid \delta(x,y)<r\})\cong \mathbb{Z}$$

\paragraph{When $L_{e_1}+1<r\leq L_{e_1}+2$:} The $E^1$ page of the Mayer-Vietoris spectral sequence is 
\begin{figure}[htbp!]
\centering
\begin{tikzpicture}
  \matrix (m) [matrix of math nodes,
    nodes in empty cells,nodes={minimum width=5ex,
     minimum height=2.5ex,outer sep=-5pt},
    column sep=1ex,row sep=1ex]{
          q     &      &     &  &  \\
          2     & \vdots & \vdots  & \vdots  &\\   
          1     &  0 &  0  & 0 \\
          0     &  \mathbb{Z}^4 & 0 & 0\\
    \quad\strut &   0  &  1  &  2 & p\strut \\};
\draw[thick] (m-1-1.east) -- (m-5-1.east) ;
\draw[thick] (m-5-1.north) -- (m-5-5.north) ;
\end{tikzpicture}
\end{figure}

Since there is no non-trivial arrow on the $E^1$ page, $E^1=E^{\infty}$. Therefore, 
$$H_0(\hat{\mathsf{Star}}_{3}\times\hat{\mathsf{Star}}_{3}-\{(x,y)\in \hat{\mathsf{Star}}_{3}\times \hat{\mathsf{Star}}_{3}\mid \delta(x,y)<r\})\cong\mathbb{Z}^4$$
and
$$H_1(\hat{\mathsf{Star}}_{3}\times\hat{\mathsf{Star}}_{3}-\{(x,y)\in \hat{\mathsf{Star}}_{3}\times \hat{\mathsf{Star}}_{3}\mid \delta(x,y)<r\})=0$$
\end{proof}

\begin{proposition}
Let $L_{e_1}$ be a positive number. If $r\leq \frac{1}{2}L_{e_1}$ and $r\leq 1$, then
$$H_0((\hat{\mathcal{H}})_{r,L_{e_1}}^2)\cong \mathbb{Z}  \mbox{\quad and\quad} H_1((\hat{\mathcal{H}})_{r,L_{e_1}}^2)\cong \mathbb{Z}^{3}$$
\end{proposition}
\begin{proof}
Consider the cover of $(\hat{\mathcal{H}})_{r,L_{e_1}}^2$ given in Equation (\ref{cover:Hgraph}). The $E^1$ page of the Mayer-Vietoris spectral sequence is 
\begin{figure}[htbp!]
\centering
\begin{tikzpicture}
  \matrix (m) [matrix of math nodes,
    nodes in empty cells,nodes={minimum width=5ex,
     minimum height=2.5ex,outer sep=-5pt},
    column sep=1ex,row sep=1ex]{
          q     &      &     &  &  \\
          3    &   \vdots   & \vdots    &\vdots  &  \\
          2     & 0& 0  &0  &\\   
          1     &   \mathbb{Z}^2 &  0  & 0 &\\
          0     &  \mathbb{Z}^4 & \mathbb{Z}^4 & 0&\\
    \quad\strut &   0  &  1  &  2 & p\strut \\};
    \draw[-stealth] (m-5-3.west) -- (m-5-2.east) node [above,pos=0.5] {$d^1$};
\draw[thick] (m-1-1.east) -- (m-6-1.east) ;
\draw[thick] (m-6-1.north) -- (m-6-5.north) ;
\end{tikzpicture}
\end{figure}

To understand the behavior of $d^1$, we choose a collection of generators for every abelian group in the $E^1$ page of the spectral sequence. Let $a$ be the generator of $H_0(U_{11}\cap U_{12})$ and $\hat{a}$ be the generator of $H_0(U_{11}\cap U_{21})$. Let $b$ be the generator of $H_0(U_{22}\cap U_{12})$ and $\hat{b}$ be the generator of $H_0(U_{22}\cap U_{21})$. Choose a generator (denoted by $c$) for $H_0(U_{11})$,  (denoted by $\hat{c}$) for $H_0(U_{22})$, and let $e$ be the generator of $H_0(U_{12})$ and $\hat{e}$ be the generator of $H_0(U_{21})$ such that
\begin{equation}
\begin{aligned}
d^1(a)&=c+e \mbox{\quad and\quad} d^1(\hat{a})=\hat{c}+\hat{e}\\
d^1(b)&=c+e \mbox{\quad and\quad} d^1(\hat{b})=\hat{c}+\hat{e}
\end{aligned}
\end{equation}

Under the given generators, $d^1$ can be represented by the following matrix 
$$\left[\begin{array}{cccc}1 & 1 & 0 & 0 \\0 & 0 & 1 & 1 \\1 & 0 & 1 & 0 \\0 & 1 & 0 & 1\end{array}\right]$$

Run the row reduction, we obtain:
$$\left[\begin{array}{cccc}1 & 0 & 0 & -1 \\0 & 0 & 1 & 1 \\0 & 0 & 0 & 0 \\0 & 1 & 0 & 1\end{array}\right]$$

The rank of this matrix is $3$, and the diagonal elements of the Smith Normal form of the above matrix are all equal to $1$. Therefore, $\Ima d^1 \cong\mathbb{Z}^3$ and $\ker d^1\cong \mathbb{Z}$. Hence, we obtain the $E^2$ page of the Mayer-Vietoris spectral sequence: 
\begin{figure}[htbp!]
\centering
\begin{tikzpicture}
  \matrix (m) [matrix of math nodes,
    nodes in empty cells,nodes={minimum width=5ex,
     minimum height=2.5ex,outer sep=-5pt},
    column sep=1ex,row sep=1ex]{
          q     &      &     &  &  \\
          3    &   \vdots   & \vdots    &\vdots  &  \\
          2     & 0& 0  &0  &\\   
          1     &   \mathbb{Z}^2 &  0  & 0 &\\
          0     &  \mathbb{Z} & \mathbb{Z} & 0&\\
    \quad\strut &   0  &  1  &  2 & p\strut \\};
   \draw[thick] (m-1-1.east) -- (m-6-1.east) ;
\draw[thick] (m-6-1.north) -- (m-6-5.north) ;
\end{tikzpicture}
\end{figure}

Since there is no non-trivial arrow on the $E^2$ page, $E^2=E^{\infty}$. Therefore, $$H_0((\hat{\mathcal{H}})_{r,L_{e_1}}^2)\cong \mathbb{Z}  \mbox{\quad and\quad} H_1((\hat{\mathcal{H}})_{r,L_{e_1}}^2)\cong \mathbb{Z}^{3}$$
\end{proof}

\begin{proposition}
Let $L_{e_1}$ be a positive number. If $\frac{1}{2}L_{e_1}<r\leq L_{e_1}$ and $r\leq 1$, then
$$H_0((\hat{\mathcal{H}})_{r,L_{e_1}}^2)\cong \mathbb{Z}  \mbox{\quad and\quad} H_1((\hat{\mathcal{H}})_{r,L_{e_1}}^2)\cong \mathbb{Z}^{3}$$
\end{proposition}
\begin{proof}
Consider the cover of $(\hat{\mathcal{H}})_{r,L_{e_1}}^2$ given in Equation (\ref{cover:Hgraph}). The $E^1$ page of the Mayer-Vietoris spectral sequence is\\ 

\begin{figure}[htbp!]
\centering
\begin{tikzpicture}
  \matrix (m) [matrix of math nodes,
    nodes in empty cells,nodes={minimum width=5ex,
     minimum height=2.5ex,outer sep=-5pt},
    column sep=1ex,row sep=1ex]{
          q     &      &     &  &  \\
          3    &   \vdots   & \vdots    &\vdots  &  \\
          2     & 0& 0  &0  &\\   
          1     &   0 &  0  & 0 &\\
          0     &  \mathbb{Z}^6 & \mathbb{Z}^8 & 0&\\
    \quad\strut &   0  &  1  &  2 & p\strut \\};
    \draw[-stealth] (m-5-3.west) -- (m-5-2.east) node [above,pos=0.5] {$d^1$};
\draw[thick] (m-1-1.east) -- (m-6-1.east) ;
\draw[thick] (m-6-1.north) -- (m-6-5.north) ;
\end{tikzpicture}
\end{figure}

To understand the behavior of $d^1$, we choose a collection of generators for every abelian group in the $E^1$ page of the spectral sequence. Let $a_1, a_2$ be the generators of $H_0(U_{11}\cap U_{12})$ such that $a_1$ and $a_2$ lie in different path components of $U_{11}\cap U_{12}$. Similarly, let $\hat{a}_1,\hat{a}_2$ be the generators of $H_0(U_{11}\cap U_{21})$, let $b_1, b_2$ be the generators of $H_0(U_{22}\cap U_{12})$ and $\hat{b}_1, \hat{b}_2$ be the generators of $H_0(U_{22}\cap U_{21})$ such that the generators of each abelian group lie in different path components. Choose generators (denoted by $c_3$ and $c_6$) for $H_0(U_{11})$,  (denoted by $\hat{c}_3$ and $\hat{c}_6$) for $H_0(U_{22})$, and let $e$ be the generator of $H_0(U_{12})$ and $\hat{e}$ be the generator of $H_0(U_{21})$ such that
\begin{equation}
\begin{aligned}
d^1(a_1)&=c_3+e \mbox{\quad and\quad} d^1(a_2)=c_6+e\\ 
d^1(\hat{a}_1)&=c_3+\hat{e} \mbox{\quad and\quad}d^1(\hat{a}_2)=c_6+\hat{e}\\
d^1(b_1)&=\hat{c}_3+e \mbox{\quad and\quad}d^1(b_2)=\hat{c}_6+e\\ 
d^1(\hat{b}_1)&=\hat{c}_3+\hat{e} \mbox{\quad and\quad} d^1(\hat{b}_2)=\hat{c}_6+\hat{e}\\
\end{aligned}
\end{equation}

Under the given generators, $d^1$ can be represented by the following matrix 
$$\left[\begin{array}{cccccccc}1 & 0 & 1 & 0 & 0 & 0 & 0 & 0 \\0 & 1 & 0 & 1 & 0 & 0 & 0 & 0 \\0 & 0 & 0 & 0 & 1 & 0 & 1 & 0 \\0 & 0 & 0 & 0 & 0 & 1 & 0 & 1 \\1 & 1 & 0 & 0 & 1 & 1 & 0 & 0 \\0 & 0 & 1 & 1 & 0 & 0 & 1 & 1\end{array}\right]$$

Run the row reduction, we obtain:\\
$$\left[\begin{array}{cccccccc}1 & 0 & 1 & 0 & 0 & 0 & 0 & 0 \\0 & 1 & 0 & 1 & 0 & 0 & 0 & 0 \\0 & 0 & 1 & 1 & 0 & 0 & 1 & 1 \\0 & 0 & 0 & 0 & 1 & 0 & 1 & 0 \\0 & 0 & 0 & 0 & 0 & 1 & 0 & 1 \\0 & 0 & 0 & 0 & 0 & 0 & 0 & 0\end{array}\right]$$

The rank of this matrix is $5$, and the diagonal elements of the Smith Normal form of the above matrix are all equal to $1$. Therefore, $\Ima d^1 \cong\mathbb{Z}^5$ and $\ker d^1\cong \mathbb{Z}^{3}$. Hence, the $E^2$ page of the Mayer-Vietoris spectral sequence is:
\begin{figure}[htbp!]
\centering
\begin{tikzpicture}
  \matrix (m) [matrix of math nodes,
    nodes in empty cells,nodes={minimum width=5ex,
     minimum height=2.5ex,outer sep=-5pt},
    column sep=1ex,row sep=1ex]{
          q     &      &     &  &  \\
          3    &   \vdots   & \vdots    &\vdots  &  \\
          2     & 0& 0  &0  &\\   
          1     &   0 &  0  & 0 &\\
          0     &  \mathbb{Z} & \mathbb{Z}^3 & 0&\\
    \quad\strut &   0  &  1  &  2 & p\strut \\};
   \draw[thick] (m-1-1.east) -- (m-6-1.east) ;
\draw[thick] (m-6-1.north) -- (m-6-5.north) ;
\end{tikzpicture}
\end{figure}

\newpage
Since there is no non-trivial arrow on the $E^2$ page, $E^2=E^{\infty}$. Therefore, $$H_0((\hat{\mathcal{H}})_{r,L_{e_1}}^2)\cong \mathbb{Z}  \mbox{\quad and\quad} H_1((\hat{\mathcal{H}})_{r,L_{e_1}}^2)\cong \mathbb{Z}^{3}$$
\end{proof}

\begin{proposition}
Let $L_{e_1}$ be a positive number. If $L_{e_1}<r\leq 1$, then
$$H_0((\hat{\mathcal{H}})_{r,L_{e_1}}^2)\cong \mathbb{Z}  \mbox{\quad and\quad} H_1((\hat{\mathcal{H}})_{r,L_{e_1}}^2)\cong \mathbb{Z}^{5}$$
\end{proposition}
\begin{proof}
Consider the cover of $(\hat{\mathcal{H}})_{r,L_{e_1}}^2$ given in Equation (\ref{cover:Hgraph}). The $E^1$ page of the Mayer-Vietoris spectral sequence is \\

\begin{figure}[htbp!]
\centering
\begin{tikzpicture}
  \matrix (m) [matrix of math nodes,
    nodes in empty cells,nodes={minimum width=5ex,
     minimum height=2.5ex,outer sep=-5pt},
    column sep=1ex,row sep=1ex]{
          q     &      &     &  &  \\
          3    &   \vdots   & \vdots    &\vdots  &  \\
          2     & 0& 0  &0  &\\   
          1     &   \mathbb{Z}^2 &  0  & 0 &\\
          0     &  \mathbb{Z}^6 & \mathbb{Z}^8 & 0&\\
    \quad\strut &   0  &  1  &  2 & p\strut \\};
    \draw[-stealth] (m-5-3.west) -- (m-5-2.east) node [above,pos=0.5] {$d^1$};
\draw[thick] (m-1-1.east) -- (m-6-1.east) ;
\draw[thick] (m-6-1.north) -- (m-6-5.north) ;
\end{tikzpicture}
\end{figure}

To understand the behavior of $d^1$, we choose a collection of generators for every abelian group in the $E^1$ page of the spectral sequence. Let $a_1, a_2$ be the generators of $H_0(U_{11}\cap U_{12})$ such that $a_1$ and $a_2$ lie in different path components of $U_{11}\cap U_{12}$. Similarly, let $\hat{a}_1,\hat{a}_2$ be the generators of $H_0(U_{11}\cap U_{21})$ such that $\hat{a}_1$ and $\hat{a}_2$ lie in different path components of $U_{11}\cap U_{21}$. Let $b_1, b_2$ be the generators of $H_0(U_{22}\cap U_{12})$ and $\hat{b}_1, \hat{b}_2$ be the generators of $H_0(U_{22}\cap U_{21})$ such that the generators of each abelian group lie in different path components. Choose generators for $H_0(U_{11})$ (denoted by $c_3$ and $c_6$) and $H_0(U_{22})$ (denoted by $\hat{c}_3$ and $\hat{c}_6$). Let $e$ be the generator of $H_0(U_{12})$ and $\hat{e}$ be the generator of $H_0(U_{21})$ such that
\begin{equation}
\begin{aligned}
d^1(a_1)&=c_3+e \mbox{\quad and\quad} d^1(\hat{a}_1)=c_3+\hat{e}\\
d^1(a_2)&=c_6+e \mbox{\quad and\quad} d^1(\hat{a}_2)=c_6+\hat{e}\\
d^1(b_1)&=\hat{c}_3+e \mbox{\quad and\quad}d^1(\hat{b}_1)=\hat{c}_3+\hat{e}\\
d^1(b_2)&=\hat{c}_6+e \mbox{\quad and\quad}d^1(\hat{b}_2)=\hat{c}_6+\hat{e}\\
\end{aligned}
\end{equation}

Under the given generators, $d^1$ can be represented by the following matrix 
$$\left[\begin{array}{cccccccc}1 & 0 & 1 & 0 & 0 & 0 & 0 & 0 \\0 & 1 & 0 & 1 & 0 & 0 & 0 & 0 \\0 & 0 & 0 & 0 & 1 & 0 & 1 & 0 \\0 & 0 & 0 & 0 & 0 & 1 & 0 & 1 \\1 & 1 & 0 & 0 & 1 & 1 & 0 & 0 \\0 & 0 & 1 & 1 & 0 & 0 & 1 & 1\end{array}\right]$$

Run the row reduction, we obtain:\\
$$\left[\begin{array}{cccccccc}1 & 0 & 1 & 0 & 0 & 0 & 0 & 0 \\0 & 1 & 0 & 1 & 0 & 0 & 0 & 0 \\0 & 0 & 1 & 1 & 0 & 0 & 1 & 1 \\0 & 0 & 0 & 0 & 1 & 0 & 1 & 0 \\0 & 0 & 0 & 0 & 0 & 1 & 0 & 1 \\0 & 0 & 0 & 0 & 0 & 0 & 0 & 0\end{array}\right]$$

The rank of this matrix is $5$, and the diagonal elements of the Smith Normal form of the above matrix are all equal to $1$. Therefore, $\Ima d^1 \cong\mathbb{Z}^5$ and $\ker d^1\cong \mathbb{Z}^{3}$. Hence, the $E^2$ page of the Mayer-Vietoris spectral sequence is
\begin{figure}[htbp!]
\centering
\begin{tikzpicture}
  \matrix (m) [matrix of math nodes,
    nodes in empty cells,nodes={minimum width=5ex,
     minimum height=2.5ex,outer sep=-5pt},
    column sep=1ex,row sep=1ex]{
          q     &      &     &  &  \\
          3    &   \vdots   & \vdots    &\vdots  &  \\
          2     & 0& 0  &0  &\\   
          1     &   \mathbb{Z}^2 &  0  & 0 &\\
          0     &  \mathbb{Z} & \mathbb{Z}^3 & 0&\\
    \quad\strut &   0  &  1  &  2 & p\strut \\};
   \draw[thick] (m-1-1.east) -- (m-6-1.east) ;
\draw[thick] (m-6-1.north) -- (m-6-5.north) ;
\end{tikzpicture}
\end{figure}

Since there is no non-trivial arrow on the $E^2$ page, $E^2=E^{\infty}$. Therefore, $H_0((\hat{\mathcal{H}})_{r,L_{e_1}}^2)\cong \mathbb{Z}$ and $H_1((\hat{\mathcal{H}})_{r,L_{e_1}}^2)\cong \mathbb{Z}^{5}$.
\end{proof}

\begin{proposition}
Let $L_{e_1}$ be a positive number. If $1<r\leq 2$ and $r\leq \frac{1}{2}L_{e_1}$, then
$$H_0((\hat{\mathcal{H}})_{r,L_{e_1}}^2)\cong \mathbb{Z}^6  \mbox{\quad and\quad} H_1((\hat{\mathcal{H}})_{r,L_{e_1}}^2)=0$$
\end{proposition}
\begin{proof}
Consider the cover of $(\hat{\mathcal{H}})_{r,L_{e_1}}^2$ given in Equation (\ref{cover:Hgraph}). The $E^1$ page of the Mayer-Vietoris spectral sequence is \\
\begin{figure}[htbp!]
\centering
\begin{tikzpicture}
  \matrix (m) [matrix of math nodes,
    nodes in empty cells,nodes={minimum width=5ex,
     minimum height=2.5ex,outer sep=-5pt},
    column sep=1ex,row sep=1ex]{
          q     &      &     &  &  \\
          3    &   \vdots   & \vdots    &\vdots  &  \\
          2     & 0& 0  &0  &\\   
          1     &   0 &  0  & 0 &\\
          0     &  \mathbb{Z}^{10} & \mathbb{Z}^4 & 0&\\
    \quad\strut &   0  &  1  &  2 & p\strut \\};
    \draw[-stealth] (m-5-3.west) -- (m-5-2.east) node [above,pos=0.5] {$d^1$};
\draw[thick] (m-1-1.east) -- (m-6-1.east) ;
\draw[thick] (m-6-1.north) -- (m-6-5.north) ;
\end{tikzpicture}
\end{figure}

To understand the behavior of $d^1$, we choose a collection of generators for every abelian group in the $E^1$ page of the spectral sequence. Let $a$ be the generator of $H_0(U_{11}\cap U_{12})$ and $\hat{a}$ be the generator of $H_0(U_{11}\cap U_{21})$. Let $b$ be the generator of $H_0(U_{22}\cap U_{12})$ and $\hat{b}$ be the generator of $H_0(U_{22}\cap U_{21})$. Choose generators (denoted by $c_1$, $c_3$, $c_4$ and $c_6$) for $H_0(U_{11})$ such that no two generators lie in the same path component of $U_{11}$. Similarly, choose generators (denoted by $\hat{c}_1$, $\hat{c}_3$, $\hat{c}_4$ and $\hat{c}_6$) for $H_0(U_{22})$ such that no two generators lie in the same path component of $U_{22}$. Moreover, let $e$ be the generator of $H_0(U_{12})$ and $\hat{e}$ be the generator of $H_0(U_{21})$ such that
\begin{equation}
\begin{aligned}
d^1(a)&=c_1+e \mbox{\quad and\quad} d^1(\hat{a})=c_4+\hat{e}\\
d^1(b)&=\hat{c}_1+e \mbox{\quad and\quad} d^1(\hat{b})=\hat{c}_4+\hat{e}
\end{aligned}
\end{equation}

Under the given generators, $d^1$ can be represented by the following matrix 
$$\left[\begin{array}{cccc}1 & 0 & 0 & 0 \\0 & 1 & 0 & 0 \\0 & 0 & 1 & 0 \\0 & 0 & 0 & 1 \\1 & 0 & 1 & 0 \\0 & 1 & 0 & 1\end{array}\right]$$

Run the row reduction, we obtain:
$$\left[\begin{array}{cccc}1 & 0 & 0 & 0 \\0 & 1 & 0 & 0 \\0 & 0 & 1 & 0 \\0 & 0 & 0 & 1 \\0 & 0 & 0 & 0 \\0 & 0 & 0 & 0\end{array}\right]$$

The rank of this matrix is $4$, and the diagonal elements of the Smith Normal form of the above matrix are all equal to $1$. Therefore, $\Ima d^1 \cong\mathbb{Z}^4$ and $\ker d^1=0$. Hence, we obtain the $E^2$ page of the Mayer-Vietoris spectral sequence: 
\begin{figure}[htbp!]
\centering
\begin{tikzpicture}
  \matrix (m) [matrix of math nodes,
    nodes in empty cells,nodes={minimum width=5ex,
     minimum height=2.5ex,outer sep=-5pt},
    column sep=1ex,row sep=1ex]{
          q     &      &     &  &  \\
          3    &   \vdots   & \vdots    &\vdots  &  \\
          2     & 0& 0  &0  &\\   
          1     &   0 &  0  & 0 &\\
          0     &  \mathbb{Z}^6 & 0 & 0&\\
    \quad\strut &   0  &  1  &  2 & p\strut \\};
   \draw[thick] (m-1-1.east) -- (m-6-1.east) ;
\draw[thick] (m-6-1.north) -- (m-6-5.north) ;
\end{tikzpicture}
\end{figure}

Since there is no non-trivial arrow on the $E^2$ page, $E^2=E^{\infty}$. Therefore, $$H_0((\hat{\mathcal{H}})_{r,L_{e_1}}^2)\cong \mathbb{Z}^6$$ and
$$H_1((\hat{\mathcal{H}})_{r,L_{e_1}}^2)=0$$
\end{proof}

\begin{proposition}\label{H:prop23}
Let $L_{e_1}$ be a positive number. If $1<r\leq 2$ and $\frac{1}{2}L_{e_1}<r\leq L_{e_1}$, then
$$H_0((\hat{\mathcal{H}})_{r,L_{e_1}}^2)\cong \mathbb{Z}^6  \mbox{\quad and\quad} H_1((\hat{\mathcal{H}})_{r,L_{e_1}}^2)=0$$
\end{proposition}
\begin{proof}
Consider the cover of $(\hat{\mathcal{H}})_{r,L_{e_1}}^2$ given in Equation (\ref{cover:Hgraph}). The $E^1$ page of the Mayer-Vietoris spectral sequence is \\
\begin{figure}[htbp!]
\centering
\begin{tikzpicture}
  \matrix (m) [matrix of math nodes,
    nodes in empty cells,nodes={minimum width=5ex,
     minimum height=2.5ex,outer sep=-5pt},
    column sep=1ex,row sep=1ex]{
          q     &      &     &  &  \\
          3    &   \vdots   & \vdots    &\vdots  &  \\
          2     & 0& 0  &0  &\\   
          1     &   0 &  0  & 0 &\\
          0     &  \mathbb{Z}^{14} & \mathbb{Z}^8 & 0&\\
    \quad\strut &   0  &  1  &  2 & p\strut \\};
    \draw[-stealth] (m-5-3.west) -- (m-5-2.east) node [above,pos=0.5] {$d^1$};
\draw[thick] (m-1-1.east) -- (m-6-1.east) ;
\draw[thick] (m-6-1.north) -- (m-6-5.north) ;
\end{tikzpicture}
\end{figure}

To understand the behavior of $d^1$, we choose a collection of generators for every abelian group in the $E^1$ page of the spectral sequence. Let $a_1, a_2$ be the generators of $H_0(U_{11}\cap U_{12})$ and $\hat{a}_1, \hat{a}_2$ be the generators of $H_0(U_{11}\cap U_{21})$. Let $b_1, b_2$ be the generators of $H_0(U_{22}\cap U_{12})$ and $\hat{b}_1, \hat{b}_2$ be the generators of $H_0(U_{22}\cap U_{21})$. Choose generators (denoted by $c_1, \dots, c_6$) for $H_0(U_{11})$ such that no two generators lie in the same path component of $U_{11}$. Similarly, choose generators (denoted by $\hat{c}_1, \dots, \hat{c}_6$) for $H_0(U_{22})$ such that no two generators lie in the same path component of $U_{22}$. Moreover, let $e$ be the generator of $H_0(U_{12})$ and $\hat{e}$ be the generator of $H_0(U_{21})$ such that
\begin{equation}
\begin{aligned}
d^1(a_1)&=c_1+e \mbox{\quad and\quad} d^1(\hat{a}_1)=c_4+\hat{e}\\
d^1(a_2)&=c_2+e \mbox{\quad and\quad} d^1(\hat{a}_2)=c_5+\hat{e}\\
d^1(b_1)&=\hat{c}_1+e \mbox{\quad and\quad} d^1(\hat{b}_1)=\hat{c}_4+\hat{e}\\
d^1(b_2)&=\hat{c}_2+e \mbox{\quad and\quad} d^1(\hat{b}_2)=\hat{c}_5+\hat{e}
\end{aligned}
\end{equation}

Under the given generators, $d^1$ can be represented by the following matrix 
$$\left[\begin{array}{cccccccc}1 & 0 & 0 & 0 & 0 & 0 & 0 & 0 \\0 & 1 & 0 & 0 & 0 & 0 & 0 & 0 \\0 & 0 & 1 & 0 & 0 & 0 & 0 & 0 \\0 & 0 & 0 & 1 & 0 & 0 & 0 & 0 \\0 & 0 & 0 & 0 & 1 & 0 & 0 & 0 \\0 & 0 & 0 & 0 & 0 & 1 & 0 & 0 \\0 & 0 & 0 & 0 & 0 & 0 & 1 & 0 \\0 & 0 & 0 & 0 & 0 & 0 & 0 & 1 \\1 & 1 & 0 & 0 & 1 & 1 & 0 & 0 \\0 & 0 & 1 & 1 & 0 & 0 & 1 & 1\end{array}\right]$$

Run the row reduction, we obtain:
$$\left[\begin{array}{cccccccc}1 & 0 & 0 & 0 & 0 & 0 & 0 & 0 \\0 & 1 & 0 & 0 & 0 & 0 & 0 & 0 \\0 & 0 & 1 & 0 & 0 & 0 & 0 & 0 \\0 & 0 & 0 & 1 & 0 & 0 & 0 & 0 \\0 & 0 & 0 & 0 & 1 & 0 & 0 & 0 \\0 & 0 & 0 & 0 & 0 & 1 & 0 & 0 \\0 & 0 & 0 & 0 & 0 & 0 & 1 & 0 \\0 & 0 & 0 & 0 & 0 & 0 & 0 & 1 \\0 & 0 & 0 & 0 & 0 & 0 & 0 & 0 \\0 & 0 & 0 & 0 & 0 & 0 & 0 & 0\end{array}\right]$$

The rank of this matrix is $8$, and the diagonal elements of the Smith Normal form of the above matrix are all equal to $1$. Therefore, $\Ima d^1 \cong\mathbb{Z}^8$ and $\ker d^1=0$. Hence, we obtain the $E^2$ page of the Mayer-Vietoris spectral sequence: 
\begin{figure}[htbp!]
\centering
\begin{tikzpicture}
  \matrix (m) [matrix of math nodes,
    nodes in empty cells,nodes={minimum width=5ex,
     minimum height=2.5ex,outer sep=-5pt},
    column sep=1ex,row sep=1ex]{
          q     &      &     &  &  \\
          3    &   \vdots   & \vdots    &\vdots  &  \\
          2     & 0& 0  &0  &\\   
          1     &   0 &  0  & 0 &\\
          0     &  \mathbb{Z}^6 & 0 & 0&\\
    \quad\strut &   0  &  1  &  2 & p\strut \\};
   \draw[thick] (m-1-1.east) -- (m-6-1.east) ;
\draw[thick] (m-6-1.north) -- (m-6-5.north) ;
\end{tikzpicture}
\end{figure}

Since there is no non-trivial arrow on the $E^2$ page, $E^2=E^{\infty}$. Therefore, $$H_0((\hat{\mathcal{H}})_{r,L_{e_1}}^2)\cong \mathbb{Z}^6  \mbox{\quad and\quad} H_1((\hat{\mathcal{H}})_{r,L_{e_1}}^2)=0$$
\end{proof}

\begin{proposition}
Let $L_{e_1}$ be a positive number. If $1<r\leq 2$ and $L_{e_1}<r\leq \frac{1}{2}L_{e_1}+1$, then
$$H_0((\hat{\mathcal{H}})_{r,L_{e_1}}^2)\cong \mathbb{Z}^6  \mbox{\quad and\quad} H_1((\hat{\mathcal{H}})_{r,L_{e_1}}^2)\cong \mathbb{Z}^2$$
\end{proposition}
\begin{proof}
Consider the cover of $(\hat{\mathcal{H}})_{r,L_{e_1}}^2$ given in Equation (\ref{cover:Hgraph}). The $E^1$ page of the Mayer-Vietoris spectral sequence is 
\begin{figure}[htbp!]
\centering
\begin{tikzpicture}
  \matrix (m) [matrix of math nodes,
    nodes in empty cells,nodes={minimum width=5ex,
     minimum height=2.5ex,outer sep=-5pt},
    column sep=1ex,row sep=1ex]{
          q     &      &     &  &  \\
          3    &   \vdots   & \vdots    &\vdots  &  \\
          2     & 0& 0  &0  &\\   
          1     &   \mathbb{Z}^2 &  0  & 0 &\\
          0     &  \mathbb{Z}^{14} & \mathbb{Z}^8 & 0&\\
    \quad\strut &   0  &  1  &  2 & p\strut \\};
    \draw[-stealth] (m-5-3.west) -- (m-5-2.east) node [above,pos=0.5] {$d^1$};
\draw[thick] (m-1-1.east) -- (m-6-1.east) ;
\draw[thick] (m-6-1.north) -- (m-6-5.north) ;
\end{tikzpicture}
\end{figure}

To understand the behavior of $d^1$, we choose a collection of generators for every abelian group in the $E^1$ page of the spectral sequence. Let $a_1, a_2$ be the generators of $H_0(U_{11}\cap U_{12})$ and $\hat{a}_1, \hat{a}_2$ be the generators of $H_0(U_{11}\cap U_{21})$. Let $b_1, b_2$ be the generators of $H_0(U_{22}\cap U_{12})$ and $\hat{b}_1, \hat{b}_2$ be the generators of $H_0(U_{22}\cap U_{21})$. Choose generators (denoted by $c_1, \dots, c_6$) for $H_0(U_{11})$ such that no two generators lie in the same path component of $U_{11}$. Similarly, choose generators (denoted by $\hat{c}_1, \dots, \hat{c}_6$) for $H_0(U_{22})$ such that no two generators lie in the same path component of $U_{22}$. Moreover, let $e$ be the generator of $H_0(U_{12})$ and $\hat{e}$ be the generator of $H_0(U_{21})$ such that
\begin{equation}
\begin{aligned}
d^1(a_1)&=c_1+e \mbox{\quad and\quad} d^1(\hat{a}_1)=c_4+\hat{e}\\
d^1(a_2)&=c_2+e \mbox{\quad and\quad} d^1(\hat{a}_2)=c_5+\hat{e}\\
d^1(b_1)&=\hat{c}_1+e \mbox{\quad and\quad} d^1(\hat{b}_1)=\hat{c}_4+\hat{e}\\
d^1(b_2)&=\hat{c}_2+e \mbox{\quad and\quad} d^1(\hat{b}_2)=\hat{c}_5+\hat{e}
\end{aligned}
\end{equation}

Under the given generators, $d^1$ can be represented by the following matrix 
$$\left[\begin{array}{cccccccc}1 & 0 & 0 & 0 & 0 & 0 & 0 & 0 \\0 & 1 & 0 & 0 & 0 & 0 & 0 & 0 \\0 & 0 & 1 & 0 & 0 & 0 & 0 & 0 \\0 & 0 & 0 & 1 & 0 & 0 & 0 & 0 \\0 & 0 & 0 & 0 & 1 & 0 & 0 & 0 \\0 & 0 & 0 & 0 & 0 & 1 & 0 & 0 \\0 & 0 & 0 & 0 & 0 & 0 & 1 & 0 \\0 & 0 & 0 & 0 & 0 & 0 & 0 & 1 \\1 & 1 & 0 & 0 & 1 & 1 & 0 & 0 \\0 & 0 & 1 & 1 & 0 & 0 & 1 & 1\end{array}\right]$$

Run the row reduction, we obtain:
$$\left[\begin{array}{cccccccc}1 & 0 & 0 & 0 & 0 & 0 & 0 & 0 \\0 & 1 & 0 & 0 & 0 & 0 & 0 & 0 \\0 & 0 & 1 & 0 & 0 & 0 & 0 & 0 \\0 & 0 & 0 & 1 & 0 & 0 & 0 & 0 \\0 & 0 & 0 & 0 & 1 & 0 & 0 & 0 \\0 & 0 & 0 & 0 & 0 & 1 & 0 & 0 \\0 & 0 & 0 & 0 & 0 & 0 & 1 & 0 \\0 & 0 & 0 & 0 & 0 & 0 & 0 & 1 \\0 & 0 & 0 & 0 & 0 & 0 & 0 & 0 \\0 & 0 & 0 & 0 & 0 & 0 & 0 & 0\end{array}\right]$$

The rank of this matrix is $8$, and the diagonal elements of the Smith Normal form of the above matrix are all equal to $1$. Therefore, $\Ima d^1 \cong\mathbb{Z}^8$ and $\ker d^1=0$. Hence, we obtain the $E^2$ page of the Mayer-Vietoris spectral sequence: 
\begin{figure}[htbp!]
\centering
\begin{tikzpicture}
  \matrix (m) [matrix of math nodes,
    nodes in empty cells,nodes={minimum width=5ex,
     minimum height=2.5ex,outer sep=-5pt},
    column sep=1ex,row sep=1ex]{
          q     &      &     &  &  \\
          3    &   \vdots   & \vdots    &\vdots  &  \\
          2     & 0& 0  &0  &\\   
          1     &  \mathbb{Z}^2 &  0  & 0 &\\
          0     &  \mathbb{Z}^6 & 0 & 0&\\
    \quad\strut &   0  &  1  &  2 & p\strut \\};
   \draw[thick] (m-1-1.east) -- (m-6-1.east) ;
\draw[thick] (m-6-1.north) -- (m-6-5.north) ;
\end{tikzpicture}
\end{figure}

Since there is no non-trivial arrow on the $E^2$ page, $E^2=E^{\infty}$. Therefore, $$H_0((\hat{\mathcal{H}})_{r,L_{e_1}}^2)\cong \mathbb{Z}^6  \mbox{\quad and\quad} H_1((\hat{\mathcal{H}})_{r,L_{e_1}}^2)\cong \mathbb{Z}^2$$
\end{proof}

\begin{proposition}
Let $L_{e_1}$ be a positive number. If $1<r\leq 2$ and $\frac{1}{2}L_{e_1}+1<r\leq L_{e_1}+1$ , then
$$H_0((\hat{\mathcal{H}})_{r,L_{e_1}}^2)\cong \mathbb{Z}^6  \mbox{\quad and\quad} H_1((\hat{\mathcal{H}})_{r,L_{e_1}}^2)\cong \mathbb{Z}^2$$
\end{proposition}
\begin{proof}
Consider the cover of $(\hat{\mathcal{H}})_{r,L_{e_1}}^2$ given in Equation (\ref{cover:Hgraph}). The $E^1$ page of the Mayer-Vietoris spectral sequence is \\
\newpage
\begin{figure}[htbp!]
\centering
\begin{tikzpicture}
  \matrix (m) [matrix of math nodes,
    nodes in empty cells,nodes={minimum width=5ex,
     minimum height=2.5ex,outer sep=-5pt},
    column sep=1ex,row sep=1ex]{
          q     &      &     &  &  \\
          3    &   \vdots   & \vdots    &\vdots  &  \\
          2     & 0& 0  &0  &\\   
          1     &   \mathbb{Z}^2 &  0  & 0 &\\
          0     &  \mathbb{Z}^{6} & 0 & 0&\\
    \quad\strut &   0  &  1  &  2 & p\strut \\};
\draw[thick] (m-1-1.east) -- (m-6-1.east) ;
\draw[thick] (m-6-1.north) -- (m-6-5.north) ;
\end{tikzpicture}
\end{figure}

Since there is no non-trivial arrow on the $E^1$ page, $E^1=E^{\infty}$. Therefore, $$H_0((\hat{\mathcal{H}})_{r,L_{e_1}}^2)\cong \mathbb{Z}^6  \mbox{\quad and\quad} H_1((\hat{\mathcal{H}})_{r,L_{e_1}}^2)\cong \mathbb{Z}$$
\end{proof}

\begin{proposition}
Let $L_{e_1}$ be a positive number. If $1<r\leq 2$ and $L_{e_1}+1<r$ , then
$$H_0((\hat{\mathcal{H}})_{r,L_{e_1}}^2)\cong \mathbb{Z}^{12}  \mbox{\quad and\quad} H_1((\hat{\mathcal{H}})_{r,L_{e_1}}^2)=0$$
\end{proposition}

\begin{proof}
Consider the cover of $(\hat{\mathcal{H}})_{r,L_{e_1}}^2$ given in Equation (\ref{cover:Hgraph}). The $E^1$ page of the Mayer-Vietoris spectral sequence is\\ 
\begin{figure}[htbp!]
\centering
\begin{tikzpicture}
  \matrix (m) [matrix of math nodes,
    nodes in empty cells,nodes={minimum width=5ex,
     minimum height=2.5ex,outer sep=-5pt},
    column sep=1ex,row sep=1ex]{
          q     &      &     &  &  \\
          3    &   \vdots   & \vdots    &\vdots  &  \\
          2     & 0& 0  &0  &\\   
          1     &   0&  0  & 0 &\\
          0     &  \mathbb{Z}^{12} & 0 & 0&\\
    \quad\strut &   0  &  1  &  2 & p\strut \\};
\draw[thick] (m-1-1.east) -- (m-6-1.east) ;
\draw[thick] (m-6-1.north) -- (m-6-5.north) ;
\end{tikzpicture}
\end{figure}

Since there is no non-trivial arrow on the $E^1$ page, $E^1=E^{\infty}$. Therefore, $$H_0((\hat{\mathcal{H}})_{r,L_{e_1}}^2)\cong \mathbb{Z}^{12}  \mbox{\quad and\quad} H_1((\hat{\mathcal{H}})_{r,L_{e_1}}^2)=0$$
\end{proof}

\begin{proposition}
Let $L_{e_1}$ be a positive number. If $2<r$ and $r\leq \frac{1}{2}L_{e_1}$, then
$$H_0((\hat{\mathcal{H}})_{r,L_{e_1}}^2)\cong \mathbb{Z}^2  \mbox{\quad and\quad} H_1((\hat{\mathcal{H}})_{r,L_{e_1}}^2)=0$$
\end{proposition}
\begin{proof}
Consider the cover of $(\hat{\mathcal{H}})_{r,L_{e_1}}^2$ given in Equation (\ref{cover:Hgraph}). The $E^1$ page of the Mayer-Vietoris spectral sequence is\\

\newpage
\begin{figure}[htbp!]
\centering
\begin{tikzpicture}
  \matrix (m) [matrix of math nodes,
    nodes in empty cells,nodes={minimum width=5ex,
     minimum height=2.5ex,outer sep=-5pt},
    column sep=1ex,row sep=1ex]{
          q     &      &     &  &  \\
          3    &   \vdots   & \vdots    &\vdots  &  \\
          2     & 0& 0  &0  &\\   
          1     &   0 &  0  & 0 &\\
          0     &  \mathbb{Z}^{6} & \mathbb{Z}^{4} & 0&\\
    \quad\strut &   0  &  1  &  2 & p\strut \\};
    \draw[-stealth] (m-5-3.west) -- (m-5-2.east) node [above,pos=0.5] {$d^1$};
\draw[thick] (m-1-1.east) -- (m-6-1.east) ;
\draw[thick] (m-6-1.north) -- (m-6-5.north) ;
\end{tikzpicture}
\end{figure}

To understand the behavior of $d^1$, we choose a collection of generators for every abelian group in the $E^1$ page of the spectral sequence. Let $a$ be the generator of $H_0(U_{11}\cap U_{12})$ and $\hat{a}$ be the generator of $H_0(U_{11}\cap U_{21})$. Let $b$ be the generator of $H_0(U_{22}\cap U_{12})$ and $\hat{b}_1$ be the generator of $H_0(U_{22}\cap U_{21})$. Choose generators (denoted by $c_1, c_4$) for $H_0(U_{11})$ such that no two generators lie in the same path component of $U_{11}$. Similarly, choose generators (denoted by $\hat{c}_1, \hat{c}_4$) for $H_0(U_{22})$ such that no two generators lie in the same path component of $U_{22}$. Moreover, let $e$ be the generator of $H_0(U_{12})$ and $\hat{e}$ be the generator of $H_0(U_{21})$ such that
\begin{equation}
\begin{aligned}
d^1(a)&=c_1+e \mbox{\quad and\quad} d^1(\hat{a})=c_4+\hat{e}\\
d^1(b)&=\hat{c}_1+e \mbox{\quad and\quad} d^1(\hat{b})=\hat{c}_4+\hat{e}
\end{aligned}
\end{equation}

Under the given generators, $d^1$ can be represented by the following matrix 
$$\left[\begin{array}{cccc}1 & 0 & 0 & 0 \\0 & 1 & 0 & 0 \\0 & 0 & 1 & 0 \\0 & 0 & 0 & 1 \\1 & 0 & 1 & 0 \\0 & 1 & 0 & 1\end{array}\right]$$

Run the row reduction, we obtain:
$$\left[\begin{array}{cccc}1 & 0 & 0 & 0 \\0 & 1 & 0 & 0 \\0 & 0 & 1 & 0 \\0 & 0 & 0 & 1 \\0 & 0 & 0 & 0 \\0 & 0 & 0 & 0\end{array}\right]$$

The rank of this matrix is $4$, and the diagonal elements of the Smith Normal form of the above matrix are all equal to $1$. Therefore, $\Ima d^1 \cong\mathbb{Z}^4$ and $\ker d^1=0$. Hence, we obtain the $E^2$ page of the Mayer-Vietoris spectral sequence: 
\begin{figure}[htbp!]
\centering
\begin{tikzpicture}
  \matrix (m) [matrix of math nodes,
    nodes in empty cells,nodes={minimum width=5ex,
     minimum height=2.5ex,outer sep=-5pt},
    column sep=1ex,row sep=1ex]{
          q     &      &     &  &  \\
          3    &   \vdots   & \vdots    &\vdots  &  \\
          2     & 0& 0  &0  &\\   
          1     &  0 &  0  & 0 &\\
          0     &  \mathbb{Z}^2 & 0 & 0&\\
    \quad\strut &   0  &  1  &  2 & p\strut \\};
   \draw[thick] (m-1-1.east) -- (m-6-1.east) ;
\draw[thick] (m-6-1.north) -- (m-6-5.north) ;
\end{tikzpicture}
\end{figure}

Since there is no non-trivial arrow on the $E^2$ page, $E^2=E^{\infty}$. Therefore, $$H_0((\hat{\mathcal{H}})_{r,L_{e_1}}^2)\cong \mathbb{Z}^2  \mbox{\quad and\quad} H_1((\hat{\mathcal{H}})_{r,L_{e_1}}^2)=0$$
\end{proof}

\begin{proposition}
Let $L_{e_1}$ be a positive number. If $2<r$ and $\frac{1}{2}L_{e_1}<r\leq \frac{1}{2}L_{e_1}+1$, then
$$H_0((\hat{\mathcal{H}})_{r,L_{e_1}}^2)\cong \mathbb{Z}^2  \mbox{\quad and\quad} H_1((\hat{\mathcal{H}})_{r,L_{e_1}}^2)=0$$
\end{proposition}
\begin{proof}
Consider the cover of $(\hat{\mathcal{H}})_{r,L_{e_1}}^2$ given in Equation (\ref{cover:Hgraph}). The $E^1$ page of the Mayer-Vietoris spectral sequence is \\

\begin{figure}[htbp!]
\centering
\begin{tikzpicture}
  \matrix (m) [matrix of math nodes,
    nodes in empty cells,nodes={minimum width=5ex,
     minimum height=2.5ex,outer sep=-5pt},
    column sep=1ex,row sep=1ex]{
          q     &      &     &  &  \\
          3    &   \vdots   & \vdots    &\vdots  &  \\
          2     & 0& 0  &0  &\\   
          1     &  0 &  0  & 0 &\\
          0     &  \mathbb{Z}^{10} & \mathbb{Z}^8 & 0&\\
    \quad\strut &   0  &  1  &  2 & p\strut \\};
    \draw[-stealth] (m-5-3.west) -- (m-5-2.east) node [above,pos=0.5] {$d^1$};
\draw[thick] (m-1-1.east) -- (m-6-1.east) ;
\draw[thick] (m-6-1.north) -- (m-6-5.north) ;
\end{tikzpicture}
\end{figure}

To understand the behavior of $d^1$, we choose a collection of generators for every abelian group in the $E^1$ page of the spectral sequence. Let $a_1, a_2$ be the generators of $H_0(U_{11}\cap U_{12})$ and $\hat{a}_1, \hat{a}_2$ be the generators of $H_0(U_{11}\cap U_{21})$. Let $b_1, b_2$ be the generators of $H_0(U_{22}\cap U_{12})$ and $\hat{b}_1, \hat{b}_2$ be the generators of $H_0(U_{22}\cap U_{21})$. Choose generators (denoted by $c_1, c_2, c_4, c_5$) for $H_0(U_{11})$ such that no two generators lie in the same path component of $U_{11}$. Similarly, choose generators (denoted by $\hat{c}_1, \hat{c}_2, \hat{c}_4, \hat{c}_5$) for $H_0(U_{22})$ such that no two generators lie in the same path component of $U_{22}$. Moreover, let $e$ be the generator of $H_0(U_{12})$ and $\hat{e}$ be the generator of $H_0(U_{21})$ such that
\begin{equation}
\begin{aligned}
d^1(a_1)&=c_1+e \mbox{\quad and\quad} d^1(\hat{a}_1)=c_4+\hat{e}\\
d^1(a_2)&=c_2+e \mbox{\quad and\quad} d^1(\hat{a}_2)=c_5+\hat{e}\\
d^1(b_1)&=\hat{c}_1+e \mbox{\quad and\quad} d^1(\hat{b}_1)=\hat{c}_4+\hat{e}\\
d^1(b_2)&=\hat{c}_2+e \mbox{\quad and\quad} d^1(\hat{b}_2)=\hat{c}_5+\hat{e}
\end{aligned}
\end{equation}

Under the given generators, $d^1$ can be represented by the following matrix 
$$\left[\begin{array}{cccccccc}1 & 0 & 0 & 0 & 0 & 0 & 0 & 0 \\0 & 1 & 0 & 0 & 0 & 0 & 0 & 0 \\0 & 0 & 1 & 0 & 0 & 0 & 0 & 0 \\0 & 0 & 0 & 1 & 0 & 0 & 0 & 0 \\0 & 0 & 0 & 0 & 1 & 0 & 0 & 0 \\0 & 0 & 0 & 0 & 0 & 1 & 0 & 0 \\0 & 0 & 0 & 0 & 0 & 0 & 1 & 0 \\0 & 0 & 0 & 0 & 0 & 0 & 0 & 1 \\1 & 1 & 0 & 0 & 1 & 1 & 0 & 0 \\0 & 0 & 1 & 1 & 0 & 0 & 1 & 1\end{array}\right]$$

Run the row reduction, we obtain:
$$\left[\begin{array}{cccccccc}1 & 0 & 0 & 0 & 0 & 0 & 0 & 0 \\0 & 1 & 0 & 0 & 0 & 0 & 0 & 0 \\0 & 0 & 1 & 0 & 0 & 0 & 0 & 0 \\0 & 0 & 0 & 1 & 0 & 0 & 0 & 0 \\0 & 0 & 0 & 0 & 1 & 0 & 0 & 0 \\0 & 0 & 0 & 0 & 0 & 1 & 0 & 0 \\0 & 0 & 0 & 0 & 0 & 0 & 1 & 0 \\0 & 0 & 0 & 0 & 0 & 0 & 0 & 1 \\0 & 0 & 0 & 0 & 0 & 0 & 0 & 0 \\0 & 0 & 0 & 0 & 0 & 0 & 0 & 0\end{array}\right]$$

The rank of this matrix is $8$, and the diagonal elements of the Smith Normal form of the above matrix are all equal to $1$. Therefore, $\Ima d^1 \cong\mathbb{Z}^8$ and $\ker d^1=0$. Hence, we obtain the $E^2$ page of the Mayer-Vietoris spectral sequence: 
\begin{figure}[htbp!]
\centering
\begin{tikzpicture}
  \matrix (m) [matrix of math nodes,
    nodes in empty cells,nodes={minimum width=5ex,
     minimum height=2.5ex,outer sep=-5pt},
    column sep=1ex,row sep=1ex]{
          q     &      &     &  &  \\
          3    &   \vdots   & \vdots    &\vdots  &  \\
          2     & 0& 0  &0  &\\   
          1     & 0 &  0  & 0 &\\
          0     &  \mathbb{Z}^2 & 0 & 0&\\
    \quad\strut &   0  &  1  &  2 & p\strut \\};
   \draw[thick] (m-1-1.east) -- (m-6-1.east) ;
\draw[thick] (m-6-1.north) -- (m-6-5.north) ;
\end{tikzpicture}
\end{figure}

Since there is no non-trivial arrow on the $E^2$ page, $E^2=E^{\infty}$. Therefore, $$H_0((\hat{\mathcal{H}})_{r,L_{e_1}}^2)\cong \mathbb{Z}^2  \mbox{\quad and\quad} H_1((\hat{\mathcal{H}})_{r,L_{e_1}}^2)=0$$
\end{proof}

\begin{proposition}
Let $L_{e_1}$ be a positive number. If $2<r$ and $\frac{1}{2}L_{e_1}+1<r\leq L_{e_1}$ , then
$$H_0((\hat{\mathcal{H}})_{r,L_{e_1}}^2)\cong \mathbb{Z}^{2}  \mbox{\quad and\quad} H_1((\hat{\mathcal{H}})_{r,L_{e_1}}^2)=0$$
\end{proposition}

\begin{proof}
Consider the cover of $(\hat{\mathcal{H}})_{r,L_{e_1}}^2$ given in Equation (\ref{cover:Hgraph}). Note that when $2<r$ and $\frac{1}{2}L_{e_1}+1<r\leq L_{e_1}$, $U_{11}=\emptyset=U_{22}$. The $E^1$ page of the Mayer-Vietoris spectral sequence is 
\begin{figure}[htbp!]
\centering
\begin{tikzpicture}
  \matrix (m) [matrix of math nodes,
    nodes in empty cells,nodes={minimum width=5ex,
     minimum height=2.5ex,outer sep=-5pt},
    column sep=1ex,row sep=1ex]{
          q     &      &     &  &  \\
          3    &   \vdots   & \vdots    &\vdots  &  \\
          2     & 0& 0  &0  &\\   
          1     &   0&  0  & 0 &\\
          0     &  \mathbb{Z}^{2} & 0 & 0&\\
    \quad\strut &   0  &  1  &  2 & p\strut \\};
\draw[thick] (m-1-1.east) -- (m-6-1.east) ;
\draw[thick] (m-6-1.north) -- (m-6-5.north) ;
\end{tikzpicture}
\end{figure}

Since there is no non-trivial arrow on the $E^1$ page, $E^1=E^{\infty}$. Therefore, $$H_0((\hat{\mathcal{H}})_{r,L_{e_1}}^2)\cong \mathbb{Z}^{2}  \mbox{\quad and\quad} H_1((\hat{\mathcal{H}})_{r,L_{e_1}}^2)=0$$
\end{proof}

\begin{proposition}
Let $L_{e_1}$ be a positive number. If $2<r$ and $L_{e_1}<r\leq L_{e_1}+1$ , then
$$H_0((\hat{\mathcal{H}})_{r,L_{e_1}}^2)\cong \mathbb{Z}^{2}  \mbox{\quad and\quad} H_1((\hat{\mathcal{H}})_{r,L_{e_1}}^2)\cong \mathbb{Z}^2$$
\end{proposition}

\begin{proof}
Consider the cover of $(\hat{\mathcal{H}})_{r,L_{e_1}}^2$ given in Equation (\ref{cover:Hgraph}). Note that when $2<r$ and $L_{e_1}<r\leq L_{e_1}+1$, $U_{11}=\emptyset=U_{22}$. The $E^1$ page of the Mayer-Vietoris spectral sequence is 
\begin{figure}[htbp!]
\centering
\begin{tikzpicture}
  \matrix (m) [matrix of math nodes,
    nodes in empty cells,nodes={minimum width=5ex,
     minimum height=2.5ex,outer sep=-5pt},
    column sep=1ex,row sep=1ex]{
          q     &      &     &  &  \\
          3    &   \vdots   & \vdots    &\vdots  &  \\
          2     & 0& 0  &0  &\\   
          1     &   \mathbb{Z}^2&  0  & 0 &\\
          0     &  \mathbb{Z}^{2} & 0 & 0&\\
    \quad\strut &   0  &  1  &  2 & p\strut \\};
\draw[thick] (m-1-1.east) -- (m-6-1.east) ;
\draw[thick] (m-6-1.north) -- (m-6-5.north) ;
\end{tikzpicture}
\end{figure}

Since there is no non-trivial arrow on the $E^1$ page, $E^1=E^{\infty}$. Therefore, $$H_0((\hat{\mathcal{H}})_{r,L_{e_1}}^2)\cong \mathbb{Z}^{2}  \mbox{\quad and\quad} H_1((\hat{\mathcal{H}})_{r,L_{e_1}}^2)\cong \mathbb{Z}^2$$
\end{proof}

\begin{proposition}
Let $L_{e_1}$ be a positive number. If $2<r$ and $L_{e_1}+1<r\leq L_{e_1}+2$ , then
$$H_0((\hat{\mathcal{H}})_{r,L_{e_1}}^2)\cong \mathbb{Z}^{8}  \mbox{\quad and\quad} H_1((\hat{\mathcal{H}})_{r,L_{e_1}}^2)=0$$
\end{proposition}

\begin{proof}
Consider the cover of $(\hat{\mathcal{H}})_{r,L_{e_1}}^2$ given in Equation (\ref{cover:Hgraph}). Note that when $2<r$ and $L_{e_1}+1<r\leq L_{e_1}+2$, $U_{11}=\emptyset=U_{22}$. The $E^1$ page of the Mayer-Vietoris spectral sequence is 
\begin{figure}[htbp!]
\centering
\begin{tikzpicture}
  \matrix (m) [matrix of math nodes,
    nodes in empty cells,nodes={minimum width=5ex,
     minimum height=2.5ex,outer sep=-5pt},
    column sep=1ex,row sep=1ex]{
          q     &      &     &  &  \\
          3    &   \vdots   & \vdots    &\vdots  &  \\
          2     & 0& 0  &0  &\\   
          1     &   0&  0  & 0 &\\
          0     &  \mathbb{Z}^{8} & 0 & 0&\\
    \quad\strut &   0  &  1  &  2 & p\strut \\};
\draw[thick] (m-1-1.east) -- (m-6-1.east) ;
\draw[thick] (m-6-1.north) -- (m-6-5.north) ;
\end{tikzpicture}
\end{figure}

Since there is no non-trivial arrow on the $E^1$ page, $E^1=E^{\infty}$. Therefore, $$H_0((\hat{\mathcal{H}})_{r,L_{e_1}}^2)\cong \mathbb{Z}^{8}  \mbox{\quad and\quad} H_1((\hat{\mathcal{H}})_{r,L_{e_1}}^2)=0$$
\end{proof}

In summary, The rank of $H_0(\hat{\mathcal{H}}^2_{r,L_{e_1}})$ for all $r>0$ and $L_{e_1}>0$ is shown in Figure \ref{fig:h0H}. The rank of $H_1(\hat{\mathcal{H}}^2_{r,L_{e_1}})$ for all $r>0$ and $L_{e_1}>0$ is shown in Figure \ref{fig:h1H}.

\section{Decomposition of $PH_i(\hat{\mathcal{H}}^2_{-,-};\mathbb{F})$}\label{decompo:H33}


In this section, we are going to decompose $PH_0(\hat{\mathcal{H}}^2_{-,-};\mathbb{F})$ and and $PH_1(\hat{\mathcal{H}}^2_{-,-};\mathbb{F})$. We associate $PH_0(\hat{\mathcal{H}}^2_{-,-};\mathbb{F})$ with a persistence module $N$ which can be constructed using Proposition \ref{prop-sim-2}. We use the following diagram to represent $N$, where each vector space in the diagram represents a constant functor indexed by a chamber of the parameter space of $PH_0(\hat{\mathcal{H}}^2_{-,-};\mathbb{F})$. For every vector space in the representation $N$, we choose a representative for each basis element according to the geometry of  $\hat{\mathcal{H}}^2_{r,L_{e_1}}$. The bases we choose for each vector space of the representation $N$ are given in (\ref{h0genH33}). Note that $e_1\simeq e_2\simeq e_3\simeq e_4$ ($\hat{e}_1\simeq \hat{e}_2\simeq \hat{e}_3\simeq \hat{e}_4$, respectively) for certain parameters $(r,L_{e_1})$. We use $e$ ($\hat{e}$, respectively) to represent the homotopy class.

\begin{equation}\label{h0genH33}
N=\begin{aligned}
&\begin{tikzcd}
  \mathbb{F} & \mathbb{F}^{6} \arrow{l}[swap]{} &\mathbb{F}^2\arrow{l}[swap]{}\\
  \mathbb{F} \arrow{u}{}  &\mathbb{F}^{6} \arrow{u}{}\arrow{l}[swap]{}  & \mathbb{F}^{2}\arrow{u}[swap]{}\arrow{l}[swap]{}\\ 
  	\mathbb{F}\arrow{u}{}	 &\mathbb{F}^{6} \arrow{u}{}\arrow{l}[swap]{}	 & \mathbb{F}^{2}\arrow{u}[swap]{}\\ 
	  	  	 &\mathbb{F}^{6} \arrow{u}{} & \mathbb{F}^{2}\arrow{u}[swap]{}\arrow{l}[swap]{}\\ 
			 &\mathbb{F}^{12} \arrow{u}{} & \mathbb{F}^{8}\arrow{u}[swap]{}\arrow{l}[swap]{}\\ 
\end{tikzcd}
=\begin{tikzcd}
<e>                  & <c_3,c_6,\hat{c}_3,\hat{c}_6, e,\hat{e}> \arrow{l}[swap]{\alpha} 				&<e,\hat{e}>\arrow{l}[swap]{\iota}\\
<e> \arrow{u}{}  &<c_3,c_6,\hat{c}_3,\hat{c}_6, e,\hat{e}> \arrow{u}{}\arrow{l}[swap]{\alpha}  		& <e,\hat{e}>\arrow{u}[swap]{}\arrow{l}[swap]{\iota}\\ 
<e>\arrow{u}{}	 &<c_3,c_6,\hat{c}_3,\hat{c}_6, e,\hat{e}> \arrow{u}{}\arrow{l}[swap]{\alpha}	 	& <e,\hat{e}>\arrow{u}[swap]{}\\ 
	  	  	 &<c_3,c_6,\hat{c}_3,\hat{c}_6, e,\hat{e}> \arrow{u}{} 					& <e,\hat{e}>\arrow{u}[swap]{}\arrow{l}[swap]{\iota}\\ 
			 &\left\langle\parbox{3.2cm}{$c_3,c_6,\hat{c}_3,\hat{c}_6, e_1, \dots, e_4,
			 \newline \hat{e}_1, \dots,\hat{e}_4$}\right\rangle \arrow{u}{\beta} & <e_1, \dots, e_4,\hat{e}_1, \dots,\hat{e}_4>\arrow{u}[swap]{\gamma}\arrow{l}[swap]{\varsigma}\\ 
\end{tikzcd}
\end{aligned}
\end{equation}

\begin{figure}[htbp!]
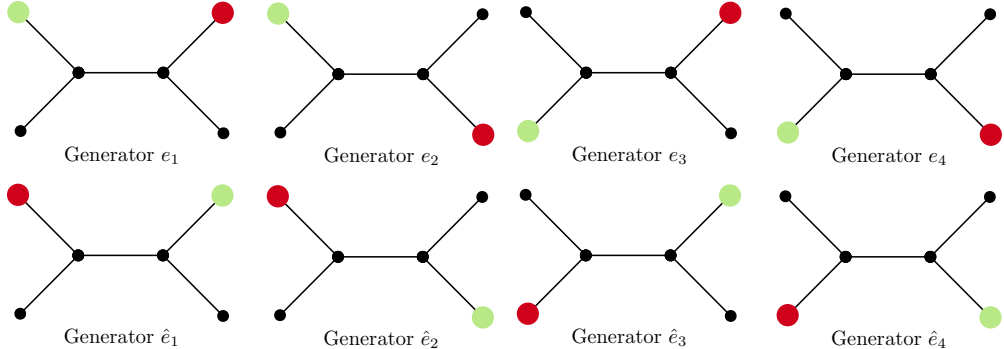

    \centering
    \resizebox{13.5cm}{!}{

\tikzset{every picture/.style={line width=0.75pt}} 


    }
    \caption{Representatives of $e_i$ and $\hat{e}_i$ for $i=1,2,3,4$.}
    \label{fig:H-h0-generator}
\end{figure}

In (\ref{h0genH33}), $c_3$ denote the arrangement where the red robot is at the leaf of edge $e_2$ while the green robot is at the leaf of edge $e_3$, and let $c_6$ denote the arrangement where the red robot is at the leaf of edge $e_4$ while the green robot is at the leaf of edge $e_5$. Respectively, $\hat{c}_3$ represents the arrangement where the red robot is at the leaf of edge $e_3$ while the green robot is at the leaf of edge $e_2$, and $\hat{c}_6$ represents the arrangement where the red robot is at the leaf of edge $e_5$ while the green robot is at the leaf of edge $e_4$. Figure \ref{fig:H-h0-generator} provides the other representatives we choose.

With the given basis for each vector space, the behavior of the morphisms in $N$ can be described as follows:
\begin{itemize}
\item $\alpha$ maps every basis element of $ <c_3,c_6,\hat{c}_3,\hat{c}_6, e,\hat{e}>$ to $e$;
\item $\beta$ maps $e_1, \dots, e_4$ to $e$, maps $\hat{e}_1, \dots,\hat{e}_4$ to $\hat{e}$, maps $c_i$ to $c_i$, and maps $\hat{c}_i$ to $\hat{c}_i$ for $i=3,6$;
\item $\gamma$ maps $e_1, \dots, e_4$ to $e$ and maps $\hat{e}_1, \dots,\hat{e}_4$ to $\hat{e}$;
\item $\iota$ and $\varsigma$ are inclusions;
\item unlabeled vertical maps are the identity maps.
\end{itemize}

Similarly, we associate $PH_1(\hat{\mathcal{H}}^2_{-,-};\mathbb{F})$ with a persistence module $N'$ where $N'$ is constructed using Proposition \ref{prop-sim-2}.

\begin{equation}\label{h1genH33}
\begin{aligned}
N'=\begin{tikzcd}
  \mathbb{F}^3 & 0 \arrow{l}[swap]{} &0\arrow{l}[swap]{}\\
  \mathbb{F}^3 \arrow{u}{}  &0 \arrow{u}{}\arrow{l}[swap]{}  & 0\arrow{u}[swap]{}\arrow{l}[swap]{}\\ 
  	\mathbb{F}^5\arrow{u}{\eta}	 &\mathbb{F}^2 \arrow{u}{}\arrow{l}[swap]{\zeta}	 & 0\arrow{u}[swap]{}\\ 
	  	  	 &\mathbb{F}^2 \arrow{u}{} & \mathbb{F}^2\arrow{u}[swap]{}\arrow{l}[swap]{}\\ 
			 &0 \arrow{u}{} & 0\arrow{u}[swap]{}\arrow{l}[swap]{}\\ 
\end{tikzcd}\\
\end{aligned}
\end{equation}

 Note that in (\ref{h1genH33}), $\zeta$ is the inclusion map and $\eta\circ \zeta=0$. 

\begin{theorem}
$PH_1(\hat{\mathcal{H}}^2_{-,-};\mathbb{F})$ is interval decomposable.
\end{theorem}
\begin{proof}
We denote the restriction of $N'$ to its support by $M$.\\

\begin{minipage}{0.3\textwidth}
\begin{equation*}\label{h1genH33'}
M=\begin{aligned}
\begin{tikzcd}
  \mathbb{F}^3 & 	 &	\\
  \mathbb{F}^3 \arrow{u}{}  &	 & 	\\ 
  	\mathbb{F}^5\arrow{u}{\eta}	 &\mathbb{F}^2 \arrow{l}[swap]{\zeta}	 & 	\\ 
	  	  	 &\mathbb{F}^2 \arrow{u}{} & \mathbb{F}^2\arrow{l}[swap]{}\\ 
\end{tikzcd}
\end{aligned}
\end{equation*}   
\end{minipage}%
\begin{minipage}{0.3\textwidth}
\begin{equation*}\label{h1genH33'1}
\begin{aligned}
M_1:=\begin{tikzcd}
  0 & 	 &	\\
  0\arrow{u}{} &	 & 	\\ 
  	\mathbb{F}^2\arrow{u}{}	 &\mathbb{F}^2 \arrow{l}[swap]{\id}	 & 	\\ 
	  	  	 &\mathbb{F}^2 \arrow{u}{} & \mathbb{F}^2\arrow{l}[swap]{}\\ 
\end{tikzcd}
\end{aligned}
\end{equation*}
\end{minipage}%
\begin{minipage}{0.33\textwidth}
\begin{equation*}\label{h1genH33'2}
\begin{aligned}
M_2:=\begin{tikzcd}
  \mathbb{F}^3 & 	 &	\\
  \mathbb{F}^3 \arrow{u}{}  &	 & 	\\ 
  	\mathbb{F}^3\arrow{u}[swap]{\restr{\eta}{\coker \zeta}}	 &0\arrow{l}[swap]{}	 & 	\\ 
	  	  	 &0\arrow{u}{} & 0\arrow{l}[swap]{}\\ 
\end{tikzcd}
\end{aligned}
\end{equation*}
\end{minipage}

Note that $\eta\circ \zeta=0$, and all unlabeled morphisms are identity maps, see (\ref{h1genH33}). Moreover, $\zeta$ is the inclusion map; hence, $M$ can be decomposed into two subrepresentations: $M_1$ and $M_2$. (See above.) Since $M_1$ is an $A_4$-quiver, and $M_2$ is an $A_3$-quiver, they are interval decomposable. Note that the intervals in the support of $M_1$ and $M_2$ are also intervals in the support of $N'$. Therefore, we can extend each indecomposable representation of $M_i$ (where $i=1,2$) to a subrepresentation of $N'$ by putting $0$ to every vertex of $N'$ which is not in $M_i$ and a trivial morphism between two vertices where at least one vertex is not in the support of $M_i$.

\end{proof}

\begin{theorem}
$PH_0(\hat{\mathcal{H}}^2_{-,-};\mathbb{F})$ is interval decomposable. 
\end{theorem}
\begin{proof}
We choose a basis for each vector space $PH_0(\hat{\mathcal{H}}^2_{r,L_{e_1}};\mathbb{F})$:\\ 
\begin{equation}\label{h0genH33'}
\begin{aligned}
&\begin{tikzcd}
<e>                  & \left\langle\parbox{4cm}{\centering$c_3-e,c_6-e,\hat{c}_3-\hat{e},\hat{c}_6-\hat{e}, e,\hat{e}-e$}\right\rangle \arrow{l}[swap]{\alpha} 				&<e,\hat{e}-e>\arrow{l}[swap]{\iota}\\
<e> \arrow{u}{}  &\left\langle\parbox{4cm}{\centering$c_3-e,c_6-e,\hat{c}_3-\hat{e},\hat{c}_6-\hat{e}, e,\hat{e}-e$}\right\rangle  \arrow{u}{}\arrow{l}[swap]{\alpha}  		& <e,\hat{e}-e>\arrow{u}[swap]{}\arrow{l}[swap]{\iota}\\ 
<e>\arrow{u}{}	 &\left\langle\parbox{4cm}{\centering$c_3-e,c_6-e,\hat{c}_3-\hat{e},\hat{c}_6-\hat{e}, e,\hat{e}-e$}\right\rangle  \arrow{u}{}\arrow{l}[swap]{\alpha}	 	& <e,\hat{e}-e>\arrow{u}[swap]{}\\ 
	  	  	 &\left\langle\parbox{4cm}{\centering$c_3-e,c_6-e,\hat{c}_3-\hat{e},\hat{c}_6-\hat{e}, e,\hat{e}-e$}\right\rangle \arrow{u}{} 					& <e,\hat{e}-e>\arrow{u}[swap]{}\arrow{l}[swap]{\iota}\\ 
			 &\left\langle\parbox{4cm}{\centering$c_3-e_1,c_6-e_1,\hat{c}_3-\hat{e}_1,\hat{c}_6-\hat{e}_1, e_1, e_2-e_1,\dots, e_4-e_1, \hat{e}_1-e_1,\hat{e}_2-\hat{e}_1,\dots,  \hat{e}_4-\hat{e}_1$}\right\rangle \arrow{u}{\beta} & \left\langle\parbox{4cm}{\centering$e_1, e_2-e_1,\dots, e_4-e_1, \hat{e}_1-e_1,\hat{e}_2-\hat{e}_1,\dots,  \hat{e}_4-\hat{e}_1$}\right\rangle \arrow{u}[swap]{\gamma}\arrow{l}[swap]{\varsigma}\\ 
\end{tikzcd}
\end{aligned}
\end{equation}
\begin{equation}\label{n1}
\begin{aligned}
N_1&=\begin{tikzcd}
<e>                  & <e>  \arrow{l}[swap]{\id} 				&<e>  \arrow{l}[swap]{\id}\\
<e> \arrow{u}{}  &<e> \arrow{u}{}\arrow{l}[swap]{\id}  		& <e>  \arrow{u}[swap]{}\arrow{l}[swap]{\id}\\ 
<e>\arrow{u}{}	 &<e>  \arrow{u}{}\arrow{l}[swap]{\id}	 	& <e>  \arrow{u}[swap]{}\\ 
	  	  	 &<e>  \arrow{u}{} 					& <e>  \arrow{u}[swap]{}\arrow{l}[swap]{\id}\\ 
			 &<e_1>  \arrow{u}{\beta} & <e_1>  \arrow{u}[swap]{\gamma}\arrow{l}[swap]{\id}\\ 
\end{tikzcd}\cong\begin{tikzcd}
\mathbb{F}                  & \mathbb{F}  \arrow{l}[swap]{\id} 				&\mathbb{F}  \arrow{l}[swap]{\id}\\
\mathbb{F} \arrow{u}{}  &\mathbb{F} \arrow{u}{}\arrow{l}[swap]{\id}  		& \mathbb{F}  \arrow{u}[swap]{}\arrow{l}[swap]{\id}\\ 
\mathbb{F}\arrow{u}{}	 &\mathbb{F} \arrow{u}{}\arrow{l}[swap]{\id}	 	& \mathbb{F}   \arrow{u}[swap]{\id}\\ 
	  	  	 &\mathbb{F} \arrow{u}{\id} 					& \mathbb{F}   \arrow{u}[swap]{}\arrow{l}[swap]{\id}\\ 
			 &\mathbb{F}\arrow{u}{\beta} 				& \mathbb{F}\arrow{u}[swap]{\gamma}\arrow{l}[swap]{\id}\\ 
\end{tikzcd}
\end{aligned}
\end{equation}

\begin{equation}\label{n2}
\begin{aligned}
N_2&=\begin{tikzcd}
0                  & <\hat{e}-e>  \arrow{l}[swap]{} 				&<\hat{e}-e>   \arrow{l}[swap]{\id}\\
0 \arrow{u}{}  &<\hat{e}-e>  \arrow{u}{}\arrow{l}[swap]{}  		& <\hat{e}-e>   \arrow{u}[swap]{}\arrow{l}[swap]{\id}\\ 
0\arrow{u}{}	 &<\hat{e}-e>   \arrow{u}{}\arrow{l}[swap]{}	 	& <\hat{e}-e>   \arrow{u}[swap]{}\\ 
	  	  	 &<\hat{e}-e>   \arrow{u}{} 					& <\hat{e}-e>   \arrow{u}[swap]{}\arrow{l}[swap]{\id}\\ 
			 &<\hat{e}_1-e_1>  \arrow{u}{\beta} & <\hat{e}_1-e_1>  \arrow{u}[swap]{\gamma}\arrow{l}[swap]{\id}\\ 
\end{tikzcd}\cong\begin{tikzcd}
0                  & \mathbb{F}  \arrow{l}[swap]{} 				&\mathbb{F}  \arrow{l}[swap]{}\\
0 \arrow{u}{}  &\mathbb{F} \arrow{u}{}\arrow{l}[swap]{}  		& \mathbb{F}  \arrow{u}[swap]{}\arrow{l}[swap]{}\\ 
0\arrow{u}{}	 &\mathbb{F} \arrow{u}{}\arrow{l}[swap]{}	 	& \mathbb{F}   \arrow{u}[swap]{}\\ 
	  	  	 &\mathbb{F} \arrow{u}{} 					& \mathbb{F}   \arrow{u}[swap]{}\arrow{l}[swap]{}\\ 
			 &\mathbb{F}\arrow{u}{\beta} 				& \mathbb{F}\arrow{u}[swap]{\gamma}\arrow{l}[swap]{\id}\\ 
\end{tikzcd}
\end{aligned}
\end{equation}

For $i=3,6$
\begin{equation}\label{ne1}
\begin{aligned}
E_i&=\begin{tikzcd}
0                  & <c_i-e>  \arrow{l}[swap]{} 				&0  \arrow{l}[swap]{}\\
0 \arrow{u}{}  &<c_i-e>  \arrow{u}{}\arrow{l}[swap]{}  		& 0  \arrow{u}[swap]{}\arrow{l}[swap]{}\\ 
0\arrow{u}{}	 &<c_i-e>  \arrow{u}{}\arrow{l}[swap]{}	 	& 0   \arrow{u}[swap]{}\\ 
	  	  	 &<c_i-e>   \arrow{u}{} 					& 0   \arrow{u}[swap]{}\arrow{l}[swap]{}\\ 
			 &<c_i-e_1> \arrow{u}{\beta} 				&0  \arrow{u}[swap]{}\arrow{l}[swap]{}\\ 
\end{tikzcd}\cong\begin{tikzcd}
0                  & \mathbb{F}  \arrow{l}[swap]{} 				&0  \arrow{l}[swap]{}\\
0 \arrow{u}{}  &\mathbb{F} \arrow{u}{}\arrow{l}[swap]{}  		& 0  \arrow{u}[swap]{}\arrow{l}[swap]{}\\ 
0\arrow{u}{}	 &\mathbb{F} \arrow{u}{}\arrow{l}[swap]{}	 	& 0   \arrow{u}[swap]{}\\ 
	  	  	 &\mathbb{F} \arrow{u}{} 					& 0   \arrow{u}[swap]{}\arrow{l}[swap]{}\\ 
			 &\mathbb{F}\arrow{u}{\beta} 				& 0\arrow{u}[swap]{}\arrow{l}[swap]{}\\ 
\end{tikzcd}
\end{aligned}
\end{equation}
and 
\begin{equation}\label{ne2}
\begin{aligned}
\hat{E}_i&=\begin{tikzcd}
0                  & <\hat{c}_i-\hat{e}>  \arrow{l}[swap]{} 				&0  \arrow{l}[swap]{}\\
0 \arrow{u}{}  &<\hat{c}_i-\hat{e}>  \arrow{u}{}\arrow{l}[swap]{}  		& 0  \arrow{u}[swap]{}\arrow{l}[swap]{}\\ 
0\arrow{u}{}	 &<\hat{c}_i-\hat{e}> \arrow{u}{}\arrow{l}[swap]{}	 	& 0   \arrow{u}[swap]{}\\ 
	  	  	 &<\hat{c}_i-\hat{e}>  \arrow{u}{} 					& 0   \arrow{u}[swap]{}\arrow{l}[swap]{}\\ 
			 &<\hat{c}_i-\hat{e}_1>\arrow{u}{\beta} 				& 0\arrow{u}[swap]{}\arrow{l}[swap]{}\\ 
\end{tikzcd}\cong\begin{tikzcd}
0                  & \mathbb{F}  \arrow{l}[swap]{} 				&0  \arrow{l}[swap]{}\\
0 \arrow{u}{}  &\mathbb{F} \arrow{u}{}\arrow{l}[swap]{}  		& 0  \arrow{u}[swap]{}\arrow{l}[swap]{}\\ 
0\arrow{u}{}	 &\mathbb{F} \arrow{u}{}\arrow{l}[swap]{}	 	& 0   \arrow{u}[swap]{}\\ 
	  	  	 &\mathbb{F} \arrow{u}{} 					& 0   \arrow{u}[swap]{}\arrow{l}[swap]{}\\ 
			 &\mathbb{F}\arrow{u}{\beta} 				& 0\arrow{u}[swap]{}\arrow{l}[swap]{}\\ 
\end{tikzcd}
\end{aligned}
\end{equation}
For $j=2,3,4$
\begin{equation}\label{nf1}
\begin{aligned}
F_j&=\begin{tikzcd}
0                  & 0  \arrow{l}[swap]{} 				&0  \arrow{l}[swap]{}\\
0 \arrow{u}{}  &0 \arrow{u}{}\arrow{l}[swap]{}  		& 0  \arrow{u}[swap]{}\arrow{l}[swap]{}\\ 
0\arrow{u}{}	 &0 \arrow{u}{}\arrow{l}[swap]{}	 	& 0   \arrow{u}[swap]{}\\ 
	  	  	 &0 \arrow{u}{} 					& 0   \arrow{u}[swap]{}\arrow{l}[swap]{}\\ 
			 &<e_i-e_1>\arrow{u}{} 				& <e_i-e_1>\arrow{u}[swap]{}\arrow{l}[swap]{\id}\\ 
\end{tikzcd}\cong\begin{tikzcd}
0                  & 0  \arrow{l}[swap]{} 				&0  \arrow{l}[swap]{}\\
0 \arrow{u}{}  &0 \arrow{u}{}\arrow{l}[swap]{}  		& 0  \arrow{u}[swap]{}\arrow{l}[swap]{}\\ 
0\arrow{u}{}	 &0 \arrow{u}{}\arrow{l}[swap]{}	 	& 0   \arrow{u}[swap]{}\\ 
	  	  	 &0 \arrow{u}{} 					& 0   \arrow{u}[swap]{}\arrow{l}[swap]{}\\ 
			 &\mathbb{F}\arrow{u}{} 				& \mathbb{F}\arrow{u}[swap]{}\arrow{l}[swap]{\id}\\ 
\end{tikzcd}
\end{aligned}
\end{equation}
and
\begin{equation}\label{nf2}
\begin{aligned}
\hat{F}_j&=\begin{tikzcd}
0                  & 0  \arrow{l}[swap]{} 				&0  \arrow{l}[swap]{}\\
0 \arrow{u}{}  &0 \arrow{u}{}\arrow{l}[swap]{}  		& 0  \arrow{u}[swap]{}\arrow{l}[swap]{}\\ 
0\arrow{u}{}	 &0 \arrow{u}{}\arrow{l}[swap]{}	 	& 0   \arrow{u}[swap]{}\\ 
	  	  	 &0 \arrow{u}{} 					& 0   \arrow{u}[swap]{}\arrow{l}[swap]{}\\ 
			 &<\hat{e}_i-\hat{e}_1>\arrow{u}{} 				& <\hat{e}_i-\hat{e}_1>\arrow{u}[swap]{}\arrow{l}[swap]{\id}\\ 
\end{tikzcd}\cong\begin{tikzcd}
0                  & 0  \arrow{l}[swap]{} 				&0  \arrow{l}[swap]{}\\
0 \arrow{u}{}  &0 \arrow{u}{}\arrow{l}[swap]{}  		& 0  \arrow{u}[swap]{}\arrow{l}[swap]{}\\ 
0\arrow{u}{}	 &0 \arrow{u}{}\arrow{l}[swap]{}	 	& 0   \arrow{u}[swap]{}\\ 
	  	  	 &0 \arrow{u}{} 					& 0   \arrow{u}[swap]{}\arrow{l}[swap]{}\\ 
			 &\mathbb{F}\arrow{u}{} 				& \mathbb{F}\arrow{u}[swap]{}\arrow{l}[swap]{\id}\\ 
\end{tikzcd}
\end{aligned}
\end{equation}

It is clear that $PH_0(\hat{\mathcal{H}}^2_{-,-};\mathbb{F})\cong N_1\oplus N_2 \oplus E_3\oplus E_6 \oplus \hat{E}_3\oplus \hat{E}_6\oplus (\bigoplus\limits_{j=2}^{4} F_j)\oplus (\bigoplus\limits_{j=2}^{4} \hat{F}_j)$. Note that $N_1$, $N_2$, $E_i$, $\hat{E}_i$, $F_j$, and $\hat{F}_j$ (where $i=3,6$ and $j=2,3,4$) are interval modules, hence by Lemma \ref{thinpoly}, they are indecomposable. Therefore, $PH_0(\hat{\mathcal{H}}^2_{-,-};\mathbb{F})$ is interval decomposable.
\end{proof}

\section{Configuration Spaces of the Generalized H Graphs}\label{genconf:H}\index{generalized H graph}
Let $m,n\geq 3$. A generalized H graph, denoted by $\mathcal{H}_{m,n}$, is a tree shown in Figure \ref{fig:genHgraph}. One can obtain $\mathcal{H}_{m,n}$ by concatenating $\mathsf{Star}_m$ and $\mathsf{Star}_n$ at a degree $1$ vertex $x_0$. We want to assign each edge an orientation where the orientation is matched with the orientation that we assigned for the star graph. Hence we subdivide the bridge of $\mathcal{H}_{m,n}$ by introducing a new vertex (and replacing the bridge with two new edges), as shown in Figure \ref{fig:genHgraph'}. We use $f$ to denote the new edge incident to the center of $\mathsf{Star}_m$ and $f'$ to denote the new edge incident to the center of $\mathsf{Star}_n$. By an abuse of notation, we denote the resulting graph by $\mathcal{H}_{m,n}$, and we use $\hat{\mathcal{H}}_{m,n}$ to denote the special $\mathcal{H}$ where the length of the bridge (before the subdivision) is $L_{e_1}$ and the length of any other edge is $1$. In this section, we describe the second configuration spaces of $\hat{\mathcal{H}}_{m,n}$ with restraint parameters $r$ and $L_{e_1}$. 

 \begin{figure}[htbp!]
     \centering
     \begin{subfigure}[b]{0.45\textwidth}
         \centering
         \resizebox{5cm}{!}{

\tikzset{every picture/.style={line width=0.75pt}} 

\begin{tikzpicture}[x=0.75pt,y=0.75pt,yscale=-1,xscale=1]

\draw    (40,32.7) -- (83.99,76.69) ;
\draw [shift={(83.99,76.69)}, rotate = 45] [color={rgb, 255:red, 0; green, 0; blue, 0 }  ][fill={rgb, 255:red, 0; green, 0; blue, 0 }  ][line width=0.75]      (0, 0) circle [x radius= 3.35, y radius= 3.35]   ;
\draw [shift={(40,32.7)}, rotate = 45] [color={rgb, 255:red, 0; green, 0; blue, 0 }  ][fill={rgb, 255:red, 0; green, 0; blue, 0 }  ][line width=0.75]      (0, 0) circle [x radius= 3.35, y radius= 3.35]   ;
\draw    (83.99,76.69) -- (41.63,119.05) ;
\draw [shift={(41.63,119.05)}, rotate = 135] [color={rgb, 255:red, 0; green, 0; blue, 0 }  ][fill={rgb, 255:red, 0; green, 0; blue, 0 }  ][line width=0.75]      (0, 0) circle [x radius= 3.35, y radius= 3.35]   ;
\draw [shift={(83.99,76.69)}, rotate = 135] [color={rgb, 255:red, 0; green, 0; blue, 0 }  ][fill={rgb, 255:red, 0; green, 0; blue, 0 }  ][line width=0.75]      (0, 0) circle [x radius= 3.35, y radius= 3.35]   ;
\draw    (145.9,76.69) -- (189.4,33.19) ;
\draw [shift={(189.4,33.19)}, rotate = 315] [color={rgb, 255:red, 0; green, 0; blue, 0 }  ][fill={rgb, 255:red, 0; green, 0; blue, 0 }  ][line width=0.75]      (0, 0) circle [x radius= 3.35, y radius= 3.35]   ;
\draw [shift={(145.9,76.69)}, rotate = 315] [color={rgb, 255:red, 0; green, 0; blue, 0 }  ][fill={rgb, 255:red, 0; green, 0; blue, 0 }  ][line width=0.75]      (0, 0) circle [x radius= 3.35, y radius= 3.35]   ;
\draw    (145.9,76.69) -- (189.89,120.68) ;
\draw [shift={(189.89,120.68)}, rotate = 45] [color={rgb, 255:red, 0; green, 0; blue, 0 }  ][fill={rgb, 255:red, 0; green, 0; blue, 0 }  ][line width=0.75]      (0, 0) circle [x radius= 3.35, y radius= 3.35]   ;
\draw [shift={(145.9,76.69)}, rotate = 45] [color={rgb, 255:red, 0; green, 0; blue, 0 }  ][fill={rgb, 255:red, 0; green, 0; blue, 0 }  ][line width=0.75]      (0, 0) circle [x radius= 3.35, y radius= 3.35]   ;
\draw    (83.99,76.69) -- (146.01,76.69) ;
\draw [shift={(146.01,76.69)}, rotate = 0] [color={rgb, 255:red, 0; green, 0; blue, 0 }  ][fill={rgb, 255:red, 0; green, 0; blue, 0 }  ][line width=0.75]      (0, 0) circle [x radius= 3.35, y radius= 3.35]   ;
\draw [shift={(83.99,76.69)}, rotate = 0] [color={rgb, 255:red, 0; green, 0; blue, 0 }  ][fill={rgb, 255:red, 0; green, 0; blue, 0 }  ][line width=0.75]      (0, 0) circle [x radius= 3.35, y radius= 3.35]   ;
\draw    (27,57.5) -- (83.99,76.69) ;
\draw [shift={(83.99,76.69)}, rotate = 18.61] [color={rgb, 255:red, 0; green, 0; blue, 0 }  ][fill={rgb, 255:red, 0; green, 0; blue, 0 }  ][line width=0.75]      (0, 0) circle [x radius= 3.35, y radius= 3.35]   ;
\draw [shift={(27,57.5)}, rotate = 18.61] [color={rgb, 255:red, 0; green, 0; blue, 0 }  ][fill={rgb, 255:red, 0; green, 0; blue, 0 }  ][line width=0.75]      (0, 0) circle [x radius= 3.35, y radius= 3.35]   ;
\draw    (146.01,76.69) -- (203,55.5) ;
\draw [shift={(203,55.5)}, rotate = 339.6] [color={rgb, 255:red, 0; green, 0; blue, 0 }  ][fill={rgb, 255:red, 0; green, 0; blue, 0 }  ][line width=0.75]      (0, 0) circle [x radius= 3.35, y radius= 3.35]   ;
\draw [shift={(146.01,76.69)}, rotate = 339.6] [color={rgb, 255:red, 0; green, 0; blue, 0 }  ][fill={rgb, 255:red, 0; green, 0; blue, 0 }  ][line width=0.75]      (0, 0) circle [x radius= 3.35, y radius= 3.35]   ;

\draw (58,31.4) node [anchor=north west][inner sep=0.75pt]    {$e_{2}$};
\draw (155,31.4) node [anchor=north west][inner sep=0.75pt]    {$e_{m+2}$};
\draw (58,102.4) node [anchor=north west][inner sep=0.75pt]    {$e_{m}$};
\draw (153,104.4) node [anchor=north west][inner sep=0.75pt]    {$e_{m+n}$};
\draw (39,64.9) node [anchor=north west][inner sep=0.75pt]    {$e_{3}$};
\draw (37.86,78.63) node [anchor=north west][inner sep=0.75pt]  [rotate=-63.61]  {$\cdots $};
\draw (177,65) node [anchor=north west][inner sep=0.75pt]    {$e_{m+3}$};
\draw (182.61,105.22) node [anchor=north west][inner sep=0.75pt]  [rotate=-284.03]  {$\cdots $};
\draw (107,81.4) node [anchor=north west][inner sep=0.75pt]    {$e_{1}$};

\end{tikzpicture}
}
\caption{}
\label{fig:genHgraph}
     \end{subfigure}
     \hfill
     \begin{subfigure}[b]{0.45\textwidth}
         \centering
        \resizebox{5cm}{!}{

\tikzset{every picture/.style={line width=0.75pt}} 

\begin{tikzpicture}[x=0.75pt,y=0.75pt,yscale=-1,xscale=1]

\draw    (40,32.7) -- (83.99,76.69) ;
\draw [shift={(83.99,76.69)}, rotate = 45] [color={rgb, 255:red, 0; green, 0; blue, 0 }  ][fill={rgb, 255:red, 0; green, 0; blue, 0 }  ][line width=0.75]      (0, 0) circle [x radius= 3.35, y radius= 3.35]   ;
\draw [shift={(57.4,50.1)}, rotate = 45] [fill={rgb, 255:red, 0; green, 0; blue, 0 }  ][line width=0.08]  [draw opacity=0] (8.93,-4.29) -- (0,0) -- (8.93,4.29) -- cycle    ;
\draw [shift={(40,32.7)}, rotate = 45] [color={rgb, 255:red, 0; green, 0; blue, 0 }  ][fill={rgb, 255:red, 0; green, 0; blue, 0 }  ][line width=0.75]      (0, 0) circle [x radius= 3.35, y radius= 3.35]   ;
\draw    (83.99,76.69) -- (41.63,119.05) ;
\draw [shift={(41.63,119.05)}, rotate = 135] [color={rgb, 255:red, 0; green, 0; blue, 0 }  ][fill={rgb, 255:red, 0; green, 0; blue, 0 }  ][line width=0.75]      (0, 0) circle [x radius= 3.35, y radius= 3.35]   ;
\draw [shift={(59.27,101.4)}, rotate = 315] [fill={rgb, 255:red, 0; green, 0; blue, 0 }  ][line width=0.08]  [draw opacity=0] (8.93,-4.29) -- (0,0) -- (8.93,4.29) -- cycle    ;
\draw [shift={(83.99,76.69)}, rotate = 135] [color={rgb, 255:red, 0; green, 0; blue, 0 }  ][fill={rgb, 255:red, 0; green, 0; blue, 0 }  ][line width=0.75]      (0, 0) circle [x radius= 3.35, y radius= 3.35]   ;
\draw    (145.9,76.69) -- (189.4,33.19) ;
\draw [shift={(189.4,33.19)}, rotate = 315] [color={rgb, 255:red, 0; green, 0; blue, 0 }  ][fill={rgb, 255:red, 0; green, 0; blue, 0 }  ][line width=0.75]      (0, 0) circle [x radius= 3.35, y radius= 3.35]   ;
\draw [shift={(171.18,51.4)}, rotate = 135] [fill={rgb, 255:red, 0; green, 0; blue, 0 }  ][line width=0.08]  [draw opacity=0] (8.93,-4.29) -- (0,0) -- (8.93,4.29) -- cycle    ;
\draw [shift={(145.9,76.69)}, rotate = 315] [color={rgb, 255:red, 0; green, 0; blue, 0 }  ][fill={rgb, 255:red, 0; green, 0; blue, 0 }  ][line width=0.75]      (0, 0) circle [x radius= 3.35, y radius= 3.35]   ;
\draw    (145.9,76.69) -- (189.89,120.68) ;
\draw [shift={(189.89,120.68)}, rotate = 45] [color={rgb, 255:red, 0; green, 0; blue, 0 }  ][fill={rgb, 255:red, 0; green, 0; blue, 0 }  ][line width=0.75]      (0, 0) circle [x radius= 3.35, y radius= 3.35]   ;
\draw [shift={(171.43,102.22)}, rotate = 225] [fill={rgb, 255:red, 0; green, 0; blue, 0 }  ][line width=0.08]  [draw opacity=0] (8.93,-4.29) -- (0,0) -- (8.93,4.29) -- cycle    ;
\draw [shift={(145.9,76.69)}, rotate = 45] [color={rgb, 255:red, 0; green, 0; blue, 0 }  ][fill={rgb, 255:red, 0; green, 0; blue, 0 }  ][line width=0.75]      (0, 0) circle [x radius= 3.35, y radius= 3.35]   ;
\draw    (83.99,76.69) -- (115,76.69) ;
\draw [shift={(115,76.69)}, rotate = 0] [color={rgb, 255:red, 0; green, 0; blue, 0 }  ][fill={rgb, 255:red, 0; green, 0; blue, 0 }  ][line width=0.75]      (0, 0) circle [x radius= 3.35, y radius= 3.35]   ;
\draw [shift={(104.49,76.69)}, rotate = 180] [fill={rgb, 255:red, 0; green, 0; blue, 0 }  ][line width=0.08]  [draw opacity=0] (8.93,-4.29) -- (0,0) -- (8.93,4.29) -- cycle    ;
\draw [shift={(83.99,76.69)}, rotate = 0] [color={rgb, 255:red, 0; green, 0; blue, 0 }  ][fill={rgb, 255:red, 0; green, 0; blue, 0 }  ][line width=0.75]      (0, 0) circle [x radius= 3.35, y radius= 3.35]   ;
\draw    (115,76.69) -- (146.01,76.69) ;
\draw [shift={(146.01,76.69)}, rotate = 0] [color={rgb, 255:red, 0; green, 0; blue, 0 }  ][fill={rgb, 255:red, 0; green, 0; blue, 0 }  ][line width=0.75]      (0, 0) circle [x radius= 3.35, y radius= 3.35]   ;
\draw [shift={(124.01,76.69)}, rotate = 0] [fill={rgb, 255:red, 0; green, 0; blue, 0 }  ][line width=0.08]  [draw opacity=0] (8.93,-4.29) -- (0,0) -- (8.93,4.29) -- cycle    ;
\draw [shift={(115,76.69)}, rotate = 0] [color={rgb, 255:red, 0; green, 0; blue, 0 }  ][fill={rgb, 255:red, 0; green, 0; blue, 0 }  ][line width=0.75]      (0, 0) circle [x radius= 3.35, y radius= 3.35]   ;
\draw    (27,57.5) -- (83.99,76.69) ;
\draw [shift={(83.99,76.69)}, rotate = 18.61] [color={rgb, 255:red, 0; green, 0; blue, 0 }  ][fill={rgb, 255:red, 0; green, 0; blue, 0 }  ][line width=0.75]      (0, 0) circle [x radius= 3.35, y radius= 3.35]   ;
\draw [shift={(49.33,65.02)}, rotate = 18.61] [fill={rgb, 255:red, 0; green, 0; blue, 0 }  ][line width=0.08]  [draw opacity=0] (8.93,-4.29) -- (0,0) -- (8.93,4.29) -- cycle    ;
\draw [shift={(27,57.5)}, rotate = 18.61] [color={rgb, 255:red, 0; green, 0; blue, 0 }  ][fill={rgb, 255:red, 0; green, 0; blue, 0 }  ][line width=0.75]      (0, 0) circle [x radius= 3.35, y radius= 3.35]   ;
\draw    (146.01,76.69) -- (203,55.5) ;
\draw [shift={(203,55.5)}, rotate = 339.6] [color={rgb, 255:red, 0; green, 0; blue, 0 }  ][fill={rgb, 255:red, 0; green, 0; blue, 0 }  ][line width=0.75]      (0, 0) circle [x radius= 3.35, y radius= 3.35]   ;
\draw [shift={(179.19,64.35)}, rotate = 159.6] [fill={rgb, 255:red, 0; green, 0; blue, 0 }  ][line width=0.08]  [draw opacity=0] (8.93,-4.29) -- (0,0) -- (8.93,4.29) -- cycle    ;
\draw [shift={(146.01,76.69)}, rotate = 339.6] [color={rgb, 255:red, 0; green, 0; blue, 0 }  ][fill={rgb, 255:red, 0; green, 0; blue, 0 }  ][line width=0.75]      (0, 0) circle [x radius= 3.35, y radius= 3.35]   ;

\draw (92,81.4) node [anchor=north west][inner sep=0.75pt]    {$f$};
\draw (58,31.4) node [anchor=north west][inner sep=0.75pt]    {$e_{2}$};
\draw (155,31.4) node [anchor=north west][inner sep=0.75pt]    {$e_{m+2}$};
\draw (58,102.4) node [anchor=north west][inner sep=0.75pt]    {$e_{m}$};
\draw (153,104.4) node [anchor=north west][inner sep=0.75pt]    {$e_{m+n}$};
\draw (123,82.09) node [anchor=north west][inner sep=0.75pt]    {$f'$};
\draw (39,64.9) node [anchor=north west][inner sep=0.75pt]    {$e_{3}$};
\draw (37.86,78.63) node [anchor=north west][inner sep=0.75pt]  [rotate=-63.61]  {$\cdots $};
\draw (177,65) node [anchor=north west][inner sep=0.75pt]    {$e_{m+3}$};
\draw (182.61,105.22) node [anchor=north west][inner sep=0.75pt]  [rotate=-284.03]  {$\cdots $};

\end{tikzpicture}        
}
\caption{}
\label{fig:genHgraph'}     
\end{subfigure}
        \caption{}
\end{figure}

Let $L_{e_i}$ denote the length of the $e_i$ for $i=2, \dots, m, m+2, \dots m+n$. Let $L_{f}$ denote the length of the $f$ and $L_{f'}$ denote the length of the $f'$. The parametric polytope of $(\mathcal{H}_{m,n})^2_{r,\vec{L}}$ is given by the following inequalities, where $\vec{L}=(L_{e_2},\dots, L_{e_m}, L_{e_{m+2}},\dots, L_{e_{m+n}}, L_{f},L_{f'})$ and $L_{f}+L_{f'}=L_{e_1}$:\\

For $i,j\in\{2,\dots, m\}$, $i\neq j$ (Respectively, $i,j\in\{m+2,\dots, m+n\}$, $i\neq j$):

\begin{equation}\label{eq:gH-01}
\begin{aligned}
0\leq x&\leq L_{e_i}\\
0\leq y&\leq L_{e_j}\\
x+y&\geq r
\end{aligned}
\end{equation}


For $i\in\{2,\dots, m\}$, $j\in\{m+2,\dots, m+n\}$ (Respectively, $i\in\{m+2,\dots, m+n\}$, $j\in\{2,\dots, m\}$):
\begin{equation}\label{eq:gH-05}
\begin{aligned}
0\leq x&\leq L_{e_i}\\
0\leq y&\leq L_{e_j}\\
x+L_{e_1}+y&\geq r
\end{aligned}
\end{equation}


For $i,j\in\{2,\dots, m\}$, $i=j$ (Respectively, $i,j\in\{m+2,\dots, m+n\}$, $i=j$):
\begin{equation}\label{eq:gH-06}
\begin{aligned}
0\leq x&\leq L_{e_i}\\
0\leq y&\leq L_{e_i}\\
\lvert x-y\rvert &\geq r
\end{aligned}
\end{equation}

\begin{minipage}{0.45\textwidth}
In addition,
\begin{equation}\label{eq:gH-04}
\begin{aligned}
0\leq x&\leq L_{f}\\
0\leq y&\leq L_{f'}\\
-x+L_{e_1}-y&\geq r
\end{aligned}
\end{equation}
\end{minipage}%
\begin{minipage}{0.45\textwidth}

\begin{equation}\label{eq:gH-04'}
\begin{aligned}
0\leq x&\leq L_{f'}\\
0\leq y&\leq L_{f}\\
-x+L_{e_1}-y&\geq r
\end{aligned}
\end{equation}
\end{minipage}%

\begin{minipage}{0.45\textwidth}
For $i\in\{2,\dots, m\}$:
\begin{equation}\label{eq:gH-02}
\begin{aligned}
0\leq x&\leq L_{e_i}\\
0\leq y&\leq L_{f}\\
x+y&\geq r
\end{aligned}
\end{equation} 
\end{minipage}%
\begin{minipage}{0.45\textwidth}
For $j\in\{2,\dots, m\}$:
\begin{equation}\label{eq:gH-02'}
\begin{aligned}
0\leq x&\leq L_{f}\\
0\leq y&\leq L_{e_j}\\
x+y&\geq r
\end{aligned}
\end{equation}  
\end{minipage}

\begin{minipage}{0.45\textwidth}
For $i\in\{2,\dots, m\}$:
\begin{equation}\label{eq:gH-03}
\begin{aligned}
0\leq x&\leq L_{e_i}\\
0\leq y&\leq L_{f'}\\
x+L_{e_1}-y&\geq r
\end{aligned}
\end{equation}  
\end{minipage}%
\begin{minipage}{0.45\textwidth}
For $j\in\{2,\dots, m\}$:
\begin{equation}\label{eq:gH-03‘}
\begin{aligned}
0\leq x&\leq L_{f'}\\
0\leq y&\leq L_{e_j}\\
-x+L_{e_1}+y&\geq r
\end{aligned}
\end{equation}    
\end{minipage}

\begin{minipage}{0.45\textwidth}
For $i\in\{m+2,\dots, m+n\}$, $j=f'$:
\begin{equation}\label{eq:gH-05''}
\begin{aligned}
0\leq x&\leq L_{e_i}\\
0\leq y&\leq L_{f'}\\
x+y&\geq r
\end{aligned}
\end{equation}
\end{minipage}%
\begin{minipage}{0.45\textwidth}
For $j\in\{m+2,\dots, m+n\}$, $i=f'$:
\begin{equation}\label{eq:gH-05'''}
\begin{aligned}
0\leq x&\leq L_{f'}\\
0\leq y&\leq L_{e_j}\\
x+y&\geq r
\end{aligned}
\end{equation}     
\end{minipage}

\begin{minipage}{0.45\textwidth}
For $i\in\{m+2,\dots, m+n\}$, $j=f$:
\begin{equation}\label{eq:gH-05''''}
\begin{aligned}
0\leq x&\leq L_{e_i}\\
0\leq y&\leq L_{f}\\
x+L_{e_1}-y&\geq r
\end{aligned}
\end{equation}  
\end{minipage}%
\begin{minipage}{0.45\textwidth}
For $j\in\{m+2,\dots, m+n\}$, $i=f$:
\begin{equation}\label{eq:gH-05'''''}
\begin{aligned}
0\leq x&\leq L_{f}\\
0\leq y&\leq L_{e_j}\\
-x+L_{e_1}+y&\geq r
\end{aligned}
\end{equation}     
\end{minipage}




\begin{minipage}{0.45\textwidth}
For $i=j=f$:
\begin{equation}\label{eq:gH-07}
\begin{aligned}
0\leq x&\leq L_{f}\\
0\leq y&\leq L_{f}\\
\lvert x-y\rvert &\geq r
\end{aligned}
\end{equation}
\end{minipage}%
\begin{minipage}{0.45\textwidth}
For $i=j=f'$:
\begin{equation}\label{eq:gH-08'}
\begin{aligned}
0\leq x&\leq L_{f'}\\
0\leq y&\leq L_{f'}\\
\lvert x-y\rvert &\geq r
\end{aligned}
\end{equation}
\end{minipage}

Compare the inequality systems (\ref{eq:H-01})-(\ref{eq:H-06'}) with inequality systems (\ref{eq:gH-01})-(\ref{eq:gH-08'}) associated to $(\mathcal{H}_{m,n})^2_{r,\vec{L}}$ when $\vec{L}=(1,\dots, 1, 1,\dots, 1, L_{f},L_{f'})$ and $L_{f}=L_{f'}=\frac{1}{2}L_{e_1}$, we conclude that $(\hat{\mathcal{H}}_{m,n})^2_{r,L_{e_1}}$ and $\hat{\mathcal{H}}^2_{r,L_{e_1}}$ have the same type of parametric polytopes. Therefore, $(\hat{\mathcal{H}}_{m,n})^2_{r,L_{e_1}}$ and $\hat{\mathcal{H}}^2_{r,L_{e_1}}$ have the same critical hyperplanes. \\

Consider the following cover of $(\hat{\mathcal{H}}_{m,n})_{r,L_{e_1}}^2$:
\begin{equation}\label{cover:gHgraph}
\begin{aligned}
U_{11}&=(e_m)^2_r\\
U_{12}&=e_m\times \hat{\mathcal{H}}_{m-1,n} -\{(x,y)\in e_m\times \hat{\mathcal{H}}_{m-1,n}\mid \delta(x,y)<r\}\\
U_{21}&=\hat{\mathcal{H}}_{m-1,n}\times e_m -\{(y,x)\in \hat{\mathcal{H}}_{m-1,n}\times e_m\mid \delta(x,y)<r\}\\
U_{22}&=(\hat{\mathcal{H}}_{m-1,n})_r^2
\end{aligned}
\end{equation}

Note that there are natural inclusions:\\

\begin{equation}\label{inc:gH}
\begin{aligned}
U_{11}\cap U_{12}&\hookrightarrow U_{11}\hookleftarrow U_{11}\cap U_{21}\\
U_{11}\cap U_{12}&\hookrightarrow U_{12}\hookleftarrow U_{22}\cap U_{12}\\
U_{11}\cap U_{21}&\hookrightarrow U_{21}\hookleftarrow U_{22}\cap U_{21}\\
U_{22}\cap U_{12}&\hookrightarrow U_{22}\hookleftarrow U_{22}\cap U_{21}
\end{aligned}
\end{equation}

Now we are going to calculate $H_i(\hat{\mathcal{H}}_{m,n})^2_{r,L_{e_1}}$ in $2$ steps. First we calculate $H_i(\hat{\mathcal{H}}_{m,3})^2_{r,L_{e_1}}$ for all $m\geq 3$, then we calculate $H_i(\hat{\mathcal{H}}_{m,n})^2_{r,L_{e_1}}$ for all $n\geq 3$.\\

\paragraph{Step 1: Calculating $H_i(\hat{\mathcal{H}}_{m,3})^2_{r,L_{e_1}}$.}
\begin{lemma}
$\vec{L}=(1,\dots,1,L_{f},L_{f'})\in(\mathbb{R}_{>0})^{m+3}$ where $L_{f}=L_{f'}=\frac{1}{2}L_{e_1}$. Then 
$$H_0(U_{12})\cong  \begin{cases} \mathbb{Z}, & \mbox{if } 0<r\leq 1 \\ 
\mathbb{Z}^{m-1}, & \mbox{if } 1 <r\leq 2  \mbox{\quad and\quad} r\leq 1+L_{e_1}  \\ 
\mathbb{Z}^{m}, & \mbox{if } 1 <r\leq 2  \mbox{\quad and\quad} r >1+L_{e_1}  \\
\mathbb{Z}, & \mbox{if } 2<r \mbox{\quad and\quad} r\leq 1+L_{e_1} \\
\mathbb{Z}^2, & \mbox{if } 2<r \mbox{\quad and\quad} 1+L_{e_1}<r\leq 2+L_{e_1}\\
0, & \mbox{else} \end{cases}$$
and 
$$H_1(U_{12})\cong  0$$
\end{lemma}

\begin{proposition}
Let $L_{e_1}$ be a positive number. If $r\leq L_{e_1}$ and $r\leq 1$, then
$$H_0((\hat{\mathcal{H}}_{m,3})_{r,L_{e_1}}^2)\cong \mathbb{Z}  \mbox{\quad and\quad} H_1((\hat{\mathcal{H}}_{m,3})_{r,L_{e_1}}^2)\cong \mathbb{Z}^{m^2-3m+3}$$
\end{proposition}
\begin{proof}
Consider the cover of $(\hat{\mathcal{H}}_{m,3})_{r,L_{e_1}}^2$ given in Equation (\ref{cover:gHgraph}). The $E^1$ page of the Mayer-Vietoris spectral sequence is 
\begin{figure}[htbp!]
\centering
\begin{tikzpicture}
  \matrix (m) [matrix of math nodes,
    nodes in empty cells,nodes={minimum width=5ex,
     minimum height=2.5ex,outer sep=-5pt},
    column sep=1ex,row sep=1ex]{
          q     &      &     &  &  \\
          3    &   \vdots   & \vdots    &\vdots  &  \\
          2     & 0& 0  &0  &\\   
          1     &   H_1((\hat{\mathcal{H}}_{m-1,3})_{r,L_{e_1}}^2) &  0  & 0 &\\
          0     &  \mathbb{Z}^5 & \mathbb{Z}^{2m} & 0&\\
    \quad\strut &   0  &  1  &  2 & p\strut \\};
    \draw[-stealth] (m-5-3.west) -- (m-5-2.east) node [above,pos=0.5] {$d^1$};
\draw[thick] (m-1-1.east) -- (m-6-1.east) ;
\draw[thick] (m-6-1.north) -- (m-6-5.north) ;
\end{tikzpicture}
\end{figure}

 Let $a$ be the generator of $H_0(U_{11}\cap U_{12})$ and $\hat{a}$ be the generator of $H_0(U_{11}\cap U_{21})$. Let $b_2,\dots, b_{m-1}$ and $b'$ be the generators of $H_0(U_{22}\cap U_{12})$. Let $\hat{b}_2, \dots, \hat{b}_{m-1}$ and $\hat{b}'$ be the generators of $H_0(U_{22}\cap U_{21})$. Choose generators (denoted by $c_1$ and $c_2$) for $H_0(U_{11})$. By induction hypothesis, $H_0(U_{22})\cong\mathbb{Z}$ and we choose a generator (denoted by $f$) for $H_0(U_{22})$. Let $e$ be the generator of $H_0(U_{12})$ and $\hat{e}$ be the generator of $H_0(U_{21})$ such that
 
\begin{equation}
\begin{aligned}
d^1(a)&=c_1+e \mbox{\quad and\quad} d^1(\hat{a})=c_2+\hat{e}\\
d^1(b_i)&=e+f \mbox{\quad and\quad} d^1(\hat{b}_i)=\hat{e}+f, \quad i=2,\dots,m-1\\
d^1(b')&=e+f \mbox{\quad and\quad} d^1(\hat{b}')=\hat{e}+f
\end{aligned}
\end{equation}

Under the given generators, $d^1$ can be represented by the following matrix. 

$$\left[\begin{array}{cccccccc}1 & 0 & 0 & \cdots  & 0 & 0 & \cdots  & 0 \\0 & 1 & 0 & \cdots  & 0 & 0 & \cdots  & 0 \\0 & 0 & 1 & \cdots & 1 & 1 & \cdots & 1 \\1 & 0 & 1 & \cdots & 1 & 0 & \cdots  & 0 \\0 & 1 & 0 & \cdots  & 0 & 1 & \cdots & 1\end{array}\right]$$
The rank of this matrix is $4$, and the diagonal elements of the Smith Normal form of the above matrix are all equal to $1$. Therefore, $\Ima d^1 \cong\mathbb{Z}^4$ and $\ker d^1\cong \mathbb{Z}^{2m-4}$. Hence we obtain the $E^2$ page of the Mayer-Vietoris spectral sequence: \\

\begin{figure}[htbp!]
\centering
\begin{tikzpicture}
  \matrix (m) [matrix of math nodes,
    nodes in empty cells,nodes={minimum width=5ex,
     minimum height=2.5ex,outer sep=-5pt},
    column sep=1ex,row sep=1ex]{
          q     &      &     &  &  \\
          3    &   \vdots   & \vdots    &\vdots  &  \\
          2     & 0& 0  &0  &\\   
          1     &  \mathbb{Z}^{m^2-5m+7}  &  0  & 0 &\\
          0     &  \mathbb{Z} & \mathbb{Z}^{2m-4} & 0&\\
    \quad\strut &   0  &  1  &  2 & p\strut \\};
   \draw[thick] (m-1-1.east) -- (m-6-1.east) ;
\draw[thick] (m-6-1.north) -- (m-6-5.north) ;
\end{tikzpicture}
\end{figure}

Here we applied the induction hypothesis $H_1((\hat{\mathcal{H}}_{m-1,3})_{r,L_{e_1}}^2)\cong \mathbb{Z}^{m^2-5m+7}$. Since there is no non-trivial arrow on the $E^2$ page, $E^2=E^{\infty}$. Therefore, $$H_0((\hat{\mathcal{H}}_{m,3})_{r,L_{e_1}}^2)\cong \mathbb{Z}$$ and  $$H_1((\hat{\mathcal{H}}_{m,3})_{r,L_{e_1}}^2)\cong \mathbb{Z}^{m^2-5m+7}\oplus \mathbb{Z}^{2m-4}\cong\mathbb{Z}^{m^2-3m+3}$$
\end{proof}

\begin{proposition}
Let $L_{e_1}$ be a positive number. If $L_{e_1}<r\leq 1$, then
$$H_0((\hat{\mathcal{H}}_{m,3})_{r,L_{e_1}}^2)\cong \mathbb{Z}  \mbox{\quad and\quad} H_1((\hat{\mathcal{H}}_{m,3})_{r,L_{e_1}}^2)\cong \mathbb{Z}^{m^2-m-1}$$
\end{proposition}
\begin{proof}
Consider the cover of $(\hat{\mathcal{H}}_{m,3})_{r,L_{e_1}}^2$ given in Equation (\ref{cover:gHgraph}). The $E^1$ page of the Mayer-Vietoris spectral sequence is \\

\newpage
\begin{figure}[htbp!]
\centering
\begin{tikzpicture}
  \matrix (m) [matrix of math nodes,
    nodes in empty cells,nodes={minimum width=5ex,
     minimum height=2.5ex,outer sep=-5pt},
    column sep=1ex,row sep=1ex]{
          q     &      &     &  &  \\
          3    &   \vdots   & \vdots    &\vdots  &  \\
          2     & 0& 0  &0  &\\   
          1     &   H_1((\hat{\mathcal{H}}_{m-1,3})_{r,L_{e_1}}^2) &  0  & 0 &\\
          0     &  \mathbb{Z}^5 & \mathbb{Z}^{2m+2} & 0&\\
    \quad\strut &   0  &  1  &  2 & p\strut \\};
    \draw[-stealth] (m-5-3.west) -- (m-5-2.east) node [above,pos=0.5] {$d^1$};
\draw[thick] (m-1-1.east) -- (m-6-1.east) ;
\draw[thick] (m-6-1.north) -- (m-6-5.north) ;
\end{tikzpicture}
\end{figure}

Let $a$ be the generator of $H_0(U_{11}\cap U_{12})$ and $\hat{a}$ be the generator of $H_0(U_{11}\cap U_{21})$. Let $b_2,\dots, b_{m-1}, b'_1$ and $b'_2$ be the generators of $H_0(U_{22}\cap U_{12})$. Let $\hat{b}_2, \dots, \hat{b}_{m-1}, \hat{b}'_1$ and $\hat{b}'_2$ be the generators of $H_0(U_{22}\cap U_{21})$. Choose generators (denoted by $c_1$ and $c_2$) for $H_0(U_{11})$. By induction hypothesis, $H_0(U_{22})\cong\mathbb{Z}$ and we choose a generator (denoted by $f$) for $H_0(U_{22})$. Let $e$ be the generator of $H_0(U_{12})$ and $\hat{e}$ be the generator of $H_0(U_{21})$ such that
\begin{equation}
\begin{aligned}
d^1(a)&=c_1+e \mbox{\quad and\quad} d^1(\hat{a})=c_2+\hat{e}\\
d^1(b_i)&=e+f \mbox{\quad and\quad} d^1(\hat{b}_i)=\hat{e}+f, \quad i=2,\dots,m-1\\
d^1(b'_j)&=e+f \mbox{\quad and\quad} d^1(\hat{b}'_j)=\hat{e}+f, \quad j=1,2
\end{aligned}
\end{equation}

Under the given generators, $d^1$ can be represented by the following matrix 
$$\left[\begin{array}{cccccccc}1 & 0 & 0 & \cdots  & 0 & 0 & \cdots  & 0 \\0 & 1 & 0 & \cdots  & 0 & 0 & \cdots  & 0 \\0 & 0 & 1 & \cdots & 1 & 1 & \cdots & 1 \\1 & 0 & 1 & \cdots & 1 & 0 & \cdots  & 0 \\0 & 1 & 0 & \cdots  & 0 & 1 & \cdots & 1\end{array}\right]$$
The rank of this matrix is $4$, and the diagonal elements of the Smith Normal form of the above matrix are all equal to $1$. Therefore, $\Ima d^1\cong\mathbb{Z}^4$ and $\ker d^1\cong \mathbb{Z}^{2m-2}$. Hence, we obtain the $E^2$ page of the Mayer-Vietoris spectral sequence: 
\begin{figure}[htbp!]
\centering
\begin{tikzpicture}
  \matrix (m) [matrix of math nodes,
    nodes in empty cells,nodes={minimum width=5ex,
     minimum height=2.5ex,outer sep=-5pt},
    column sep=1ex,row sep=1ex]{
          q     &      &     &  &  \\
          3    &   \vdots   & \vdots    &\vdots  &  \\
          2     & 0& 0  &0  &\\   
          1     &  \mathbb{Z}^{(m-1)^2-(m-1)-1}  &  0  & 0 &\\
          0     &  \mathbb{Z} & \mathbb{Z}^{2m-2} & 0&\\
    \quad\strut &   0  &  1  &  2 & p\strut \\};
   \draw[thick] (m-1-1.east) -- (m-6-1.east) ;
\draw[thick] (m-6-1.north) -- (m-6-5.north) ;
\end{tikzpicture}
\end{figure}

Here we applied the induction hypothesis $H_1((\hat{\mathcal{H}}_{m-1,3})_{r,L_{e_1}}^2)\cong \mathbb{Z}^{(m-1)^2-(m-1)-1}$. Since there is no non-trivial arrow on the $E^2$ page, $E^2=E^{\infty}$. Therefore, $$H_0((\hat{\mathcal{H}}_{m,3})_{r,L_{e_1}}^2)\cong \mathbb{Z}$$ and $$H_1((\hat{\mathcal{H}}_{m,3})_{r,L_{e_1}}^2)\cong \mathbb{Z}^{(m-1)^2-(m-1)-2}\oplus \mathbb{Z}^{2m-2}\cong\mathbb{Z}^{m^2-m-1}$$
\end{proof}

\begin{proposition}
Let $L_{e_1}$ be a positive number. If $1<r\leq 2$ and $r\leq L_{e_1}$, then
$$H_0((\hat{\mathcal{H}}_{m,3})_{r,L_{e_1}}^2)\cong \mathbb{Z}^{m^2-3m+6}  \mbox{\quad and\quad} H_1((\hat{\mathcal{H}}_{m,3})_{r,L_{e_1}}^2)=0$$
\end{proposition}

\begin{proof}
Consider the cover of $(\hat{\mathcal{H}}_{m,3})_{r,L_{e_1}}^2$ given in Equation (\ref{cover:gHgraph}). The $E^1$ page of the Mayer-Vietoris spectral sequence is 
\begin{figure}[htbp!]
\centering
\begin{tikzpicture}
  \matrix (m) [matrix of math nodes,
    nodes in empty cells,nodes={minimum width=5ex,
     minimum height=2.5ex,outer sep=-5pt},
    column sep=1ex,row sep=1ex]{
          q     &      &     &  &  \\
          3    &   \vdots   & \vdots    &\vdots  &  \\
          2     & 0& 0  &0  &\\   
          1     &   H_1((\hat{\mathcal{H}}_{m-1,3})_{r,L_{e_1}}^2) &  0  & 0 &\\
          0     &  \mathbb{Z}^{2m-2}\oplus H_0((\hat{\mathcal{H}}_{m-1,3})_{r,L_{e_1}}^2)  & \mathbb{Z}^{2} & 0&\\
    \quad\strut &   0  &  1  &  2 & p\strut \\};
    \draw[-stealth] (m-5-3.west) -- (m-5-2.east) node [above,pos=0.5] {$d^1$};
\draw[thick] (m-1-1.east) -- (m-6-1.east) ;
\draw[thick] (m-6-1.north) -- (m-6-5.north) ;
\end{tikzpicture}
\end{figure}

Since $H_0(U_{11}\cap U_{12})=H_0(U_{11}\cap U_{21})=\emptyset$, we do not  choose generators for their $0$-th homology groups. Let $b$ be the generator of $H_0(U_{22}\cap U_{12})$ and $\hat{b}$ be the generator of $H_0(U_{22}\cap U_{21})$. By induction hypothesis, $H_0(U_{22})\cong\mathbb{Z}^{(m-1)(m-4)+6}$ and we choose a collection of generators (denoted by $f_{ij}$ and $\hat{f}_{ij}$, where $2\leq i<j\leq m-1$ or $m+1\leq i<j\leq m+2$, $f'$, and $\hat{f}'$) for $H_0(U_{22})$. Let $e_2,\dots, e_{m-1}, e'$ be the generators of $H_0(U_{12})$ and $\hat{e}_2,\dots, \hat{e}_{m-1}, \hat{e}'$ be the generators of $H_0(U_{21})$ such that
\begin{equation}
\begin{aligned}
d^1(b)&=e'+f' \mbox{\quad and\quad}
d^1(\hat{b})=\hat{e}'+\hat{f}'\\
\end{aligned}
\end{equation}

Under the given generators, $d^1$ can be represented by the following matrix 
$$\left[\begin{array}{cc}1 & 0 \\0 & 1 \\1 & 0\\0 & 1\end{array}\right]$$
The rank of this matrix is $2$, and the diagonal elements of the Smith Normal form of the above matrix are all equal to $1$. Therefore, $\Ima d^1\cong\mathbb{Z}^{2}$ and $\coker d^1\cong \mathbb{Z}^{m^2-3m+6}$. Hence, we obtain the $E^2$ page of the Mayer-Vietoris spectral sequence: 
\begin{figure}[htbp!]
\centering
\begin{tikzpicture}
  \matrix (m) [matrix of math nodes,
    nodes in empty cells,nodes={minimum width=5ex,
     minimum height=2.5ex,outer sep=-5pt},
    column sep=1ex,row sep=1ex]{
          q     &      &     &  &  \\
          3    &   \vdots   & \vdots    &\vdots  &  \\
          2     & 0& 0  &0  &\\   
          1     &  0 &  0  & 0 &\\
          0     & \mathbb{Z}^{m^2-3m+6}  & 0 & 0&\\
    \quad\strut &   0  &  1  &  2 & p\strut \\};
   \draw[thick] (m-1-1.east) -- (m-6-1.east) ;
\draw[thick] (m-6-1.north) -- (m-6-5.north) ;
\end{tikzpicture}
\end{figure}

Here we applied the induction hypothesis $H_1((\hat{\mathcal{H}}_{m-1,3})_{r,L_{e_1}}^2)=0$. Since there is no non-trivial arrow on the $E^2$ page, $E^2=E^{\infty}$. Therefore, $H_0((\hat{\mathcal{H}}_{m,3})_{r,L_{e_1}}^2)\cong \mathbb{Z}^{m^2-3m+6}$ and $ H_1((\hat{\mathcal{H}}_{m,3})_{r,L_{e_1}}^2)=0$.
\end{proof}

\begin{proposition}
Let $L_{e_1}$ be a positive number. If $1<r\leq 2$ and $L_{e_1}<r\leq L_{e_1}+1$, then
$$H_0((\hat{\mathcal{H}}_{m,3})_{r,L_{e_1}}^2)\cong \mathbb{Z}^{m^2-3m+6}  \mbox{\quad and\quad} H_1((\hat{\mathcal{H}}_{m,3})_{r,L_{e_1}}^2)=\mathbb{Z}^{2m-4}$$
\end{proposition}

\begin{proof}
Consider the cover of $(\hat{\mathcal{H}}_{m,3})_{r,L_{e_1}}^2$ given in Equation (\ref{cover:gHgraph}). The $E^1$ page of the Mayer-Vietoris spectral sequence is 
\begin{figure}[htbp!]
\centering
\begin{tikzpicture}
  \matrix (m) [matrix of math nodes,
    nodes in empty cells,nodes={minimum width=5ex,
     minimum height=2.5ex,outer sep=-5pt},
    column sep=1ex,row sep=1ex]{
          q     &      &     &  &  \\
          3    &   \vdots   & \vdots    &\vdots  &  \\
          2     & 0& 0  &0  &\\   
          1     &   H_1((\hat{\mathcal{H}}_{m-1,3})_{r,L_{e_1}}^2) &  0  & 0 &\\
          0     &  \mathbb{Z}^{2m-2}\oplus H_0((\hat{\mathcal{H}}_{m-1,3})_{r,L_{e_1}}^2)  & \mathbb{Z}^{4} & 0&\\
    \quad\strut &   0  &  1  &  2 & p\strut \\};
    \draw[-stealth] (m-5-3.west) -- (m-5-2.east) node [above,pos=0.5] {$d^1$};
\draw[thick] (m-1-1.east) -- (m-6-1.east) ;
\draw[thick] (m-6-1.north) -- (m-6-5.north) ;
\end{tikzpicture}
\end{figure}

Since $H_0(U_{11}\cap U_{12})=H_0(U_{11}\cap U_{21})=\emptyset$, we do not  choose generators for their $0$-th homology groups. Let $b_1$ and $b_2$ be the generators of $H_0(U_{22}\cap U_{12})$. Let $\hat{b}_1$ and $\hat{b}_2$ be the generators of $H_0(U_{22}\cap U_{21})$. By induction hypothesis, $H_0(U_{22})\cong\mathbb{Z}^{(m-1)(m-4)+6}$ and we choose a collection of generators (denoted by $f_{ij}$ and $\hat{f}_{ij}$, where $2\leq i<j\leq m-1$ or $m+1\leq i<j\leq m+2$, $f'$, and $\hat{f}'$) for $H_0(U_{22})$. Let $e_2,\dots, e_{m-1}, e'$ be the generators of $H_0(U_{12})$ and $\hat{e}_2,\dots, \hat{e}_{m-1}, \hat{e}'$ be the generators of $H_0(U_{21})$ such that for $i=1,2$, 
\begin{equation}
\begin{aligned}
d^1(b_i)&=e'+f' \mbox{\quad and\quad}
d^1(\hat{b}_i)=\hat{e}'+\hat{f}'\\
\end{aligned}
\end{equation}

Under the given generators, $d^1$ can be represented by the following matrix 
$$\left[\begin{array}{cccc}1 & 1 & 0 & 0 \\0 & 0 & 1 & 1 \\1 & 1 & 0 & 0 \\0 & 0& 1 & 1 \end{array}\right]$$
The rank of this matrix is $2$, and the diagonal elements of the Smith Normal form of the above matrix are all equal to $1$. Therefore, $\Ima d^1\cong\mathbb{Z}^{2}$ and $\coker d^1\cong\mathbb{Z}^{2m-2+(m-1)(m-4)+6-2}\cong \mathbb{Z}^{m^2-3m+6}$. Hence, we obtain the $E^2$ page of the Mayer-Vietoris spectral sequence: 
\begin{figure}[htbp!]
\centering
\begin{tikzpicture}
  \matrix (m) [matrix of math nodes,
    nodes in empty cells,nodes={minimum width=5ex,
     minimum height=2.5ex,outer sep=-5pt},
    column sep=1ex,row sep=1ex]{
          q     &      &     &  &  \\
          3    &   \vdots   & \vdots    &\vdots  &  \\
          2     & 0& 0  &0  &\\   
          1     &  \mathbb{Z}^{2(m-1)-4} &  0  & 0 &\\
          0     &  \mathbb{Z}^{m^2-3m+6} & \mathbb{Z}^2 & 0&\\
    \quad\strut &   0  &  1  &  2 & p\strut \\};
   \draw[thick] (m-1-1.east) -- (m-6-1.east) ;
\draw[thick] (m-6-1.north) -- (m-6-5.north) ;
\end{tikzpicture}
\end{figure}

Here we applied the induction hypothesis $H_1((\hat{\mathcal{H}}_{m-1,3})_{r,L_{e_1}}^2)=\mathbb{Z}^{2(m-1)-4}$ and $H_0((\hat{\mathcal{H}}_{m-1,3})_{r,L_{e_1}}^2)=\mathbb{Z}^{(m-1)(m-4)+6}$. Since there is no non-trivial arrow on the $E^2$ page, $E^2=E^{\infty}$. Therefore, $$H_0((\hat{\mathcal{H}}_{m,3})_{r,L_{e_1}}^2)\cong \mathbb{Z}^{m^2-3m+6}$$ and $$H_1((\hat{\mathcal{H}}_{m,3})_{r,L_{e_1}}^2)=\mathbb{Z}^{2m-4}$$
\end{proof}

\begin{proposition}
Let $L_{e_1}$ be a positive number. If $1<r\leq 2$ and $L_{e_1}+1<r\leq L_{e_1}+2$, then
$$H_0((\hat{\mathcal{H}}_{m,3})_{r,L_{e_1}}^2)\cong \mathbb{Z}^{m^2+m}  \mbox{\quad and\quad} H_1((\hat{\mathcal{H}}_{m,3})_{r,L_{e_1}}^2)=0$$
\end{proposition}

\begin{proof}
Consider the cover of $(\hat{\mathcal{H}}_{m,3})_{r,L_{e_1}}^2$ given in Equation (\ref{cover:gHgraph}). The $E^1$ page of the Mayer-Vietoris spectral sequence is 
\begin{figure}[htbp!]
\centering
\begin{tikzpicture}
  \matrix (m) [matrix of math nodes,
    nodes in empty cells,nodes={minimum width=5ex,
     minimum height=2.5ex,outer sep=-5pt},
    column sep=1ex,row sep=1ex]{
          q     &      &     &  &  \\
          3    &   \vdots   & \vdots    &\vdots  &  \\
          2     & 0& 0  &0  &\\   
          1     &   H_1((\hat{\mathcal{H}}_{m-1,3})_{r,L_{e_1}}^2) &  0  & 0 &\\
          0     &  \mathbb{Z}^{2m}\oplus H_0((\hat{\mathcal{H}}_{m-1,3})_{r,L_{e_1}}^2)  & 0 & 0&\\
    \quad\strut &   0  &  1  &  2 & p\strut \\};
\draw[thick] (m-1-1.east) -- (m-6-1.east) ;
\draw[thick] (m-6-1.north) -- (m-6-5.north) ;
\end{tikzpicture}
\end{figure}

By induction hypothesis, $H_0(U_{22})\cong\mathbb{Z}^{(m-1)^2+(m-1)}$. Since there is no non-trivial arrow on the $E^1$ page, $E^1=E^{\infty}$. Therefore, $H_0((\hat{\mathcal{H}}_{m,3})_{r,L_{e_1}}^2)\cong \mathbb{Z}^{m^2+m}$ and $H_1((\hat{\mathcal{H}}_{m,3})_{r,L_{e_1}}^2)=0$.
\end{proof}

\begin{proposition}
Let $L_{e_1}$ be a positive number. If $2<r$ and $r\leq L_{e_1}$, then
$$H_0((\hat{\mathcal{H}}_{m,3})_{r,L_{e_1}}^2)\cong \mathbb{Z}^{2}  \mbox{\quad and\quad} H_1((\hat{\mathcal{H}}_{m,3})_{r,L_{e_1}}^2)=0$$
\end{proposition}

\begin{proof}
Consider the cover of $(\hat{\mathcal{H}}_{m,3})_{r,L_{e_1}}^2$ given in Equation (\ref{cover:gHgraph}). The $E^1$ page of the Mayer-Vietoris spectral sequence is 
\begin{figure}[htbp!]
\centering
\begin{tikzpicture}
  \matrix (m) [matrix of math nodes,
    nodes in empty cells,nodes={minimum width=5ex,
     minimum height=2.5ex,outer sep=-5pt},
    column sep=1ex,row sep=1ex]{
          q     &      &     &  &  \\
          3    &   \vdots   & \vdots    &\vdots  &  \\
          2     & 0& 0  &0  &\\   
          1     &   H_1((\hat{\mathcal{H}}_{m-1,3})_{r,L_{e_1}}^2) &  0  & 0 &\\
          0     &  \mathbb{Z}^{2}\oplus H_0((\hat{\mathcal{H}}_{m-1,3})_{r,L_{e_1}}^2)  & \mathbb{Z}^{2} & 0&\\
    \quad\strut &   0  &  1  &  2 & p\strut \\};
    \draw[-stealth] (m-5-3.west) -- (m-5-2.east) node [above,pos=0.5] {$d^1$};
\draw[thick] (m-1-1.east) -- (m-6-1.east) ;
\draw[thick] (m-6-1.north) -- (m-6-5.north) ;
\end{tikzpicture}
\end{figure}

Since $H_0(U_{11}\cap U_{12})=H_0(U_{11}\cap U_{21})=\emptyset$, we do not  choose generators for their $0$-th homology groups. Let $b$ be the generator of $H_0(U_{22}\cap U_{12})$ and $\hat{b}$ be the generator of $H_0(U_{22}\cap U_{21})$. By induction hypothesis, $H_0(U_{22})\cong\mathbb{Z}^{(m-1)(m-4)+6}$ and we choose a set of generators (denoted by $f'$ and $\hat{f}'$) for $H_0(U_{22})$. Let $e$ be the generator of $H_0(U_{12})$ and $\hat{e}'$ be the generator of $H_0(U_{21})$ such that for $i=1,2$, 
\begin{equation}
\begin{aligned}
d^1(b_i)&=e'+f' \mbox{\quad and\quad}
d^1(\hat{b}_i)=\hat{e}'+\hat{f}'\\
\end{aligned}
\end{equation}

Under the given generators, $d^1$ can be represented by the following matrix 
$$\left[\begin{array}{cc}1 & 0 \\0 & 1  \\1 & 0  \\0 & 1 \end{array}\right]$$
The rank of this matrix is $2$, and the diagonal elements of the Smith Normal form of the above matrix are all equal to $1$. Therefore, $\Ima d^1\cong\mathbb{Z}^{2}$ and $\coker d^1\cong \mathbb{Z}^{2}$. Hence, we obtain the $E^2$ page of the Mayer-Vietoris spectral sequence: 
\begin{figure}[htbp!]
\centering
\begin{tikzpicture}
  \matrix (m) [matrix of math nodes,
    nodes in empty cells,nodes={minimum width=5ex,
     minimum height=2.5ex,outer sep=-5pt},
    column sep=1ex,row sep=1ex]{
          q     &      &     &  &  \\
          3    &   \vdots   & \vdots    &\vdots  &  \\
          2     & 0& 0  &0  &\\   
          1     &  0 &  0  & 0 &\\
          0     &  \mathbb{Z}^{2}  & 0 & 0&\\
    \quad\strut &   0  &  1  &  2 & p\strut \\};
   \draw[thick] (m-1-1.east) -- (m-6-1.east) ;
\draw[thick] (m-6-1.north) -- (m-6-5.north) ;
\end{tikzpicture}
\end{figure}

Here we applied the induction hypothesis $H_1((\hat{\mathcal{H}}_{m-1,3})_{r,L_{e_1}}^2)=0$ and $H_0((\hat{\mathcal{H}}_{m-1,3})_{r,L_{e_1}}^2)=\mathbb{Z}^2$. Since there is no non-trivial arrow on the $E^2$ page, $E^2=E^{\infty}$. Therefore, $H_0((\hat{\mathcal{H}}_{m,3})_{r,L_{e_1}}^2)\cong \mathbb{Z}^{2}$ and $H_1((\hat{\mathcal{H}}_{m,3})_{r,L_{e_1}}^2)=0$.
\end{proof}

\begin{proposition}
Let $L_{e_1}$ be a positive number. If $2<r$ and $L_{e_1}<r\leq L_{e_1}+1$, then
$$H_0((\hat{\mathcal{H}}_{m,3})_{r,L_{e_1}}^2)\cong \mathbb{Z}^{2}  \mbox{\quad and\quad} H_1((\hat{\mathcal{H}}_{m,3})_{r,L_{e_1}}^2)=\mathbb{Z}^{2m-4}$$
\end{proposition}

\begin{proof}
Consider the cover of $(\hat{\mathcal{H}}_{m,3})_{r,L_{e_1}}^2$ given in Equation (\ref{cover:gHgraph}). The $E^1$ page of the Mayer-Vietoris spectral sequence is 
\begin{figure}[htbp!]
\centering
\begin{tikzpicture}
  \matrix (m) [matrix of math nodes,
    nodes in empty cells,nodes={minimum width=5ex,
     minimum height=2.5ex,outer sep=-5pt},
    column sep=1ex,row sep=1ex]{
          q     &      &     &  &  \\
          3    &   \vdots   & \vdots    &\vdots  &  \\
          2     & 0& 0  &0  &\\   
          1     &   H_1((\hat{\mathcal{H}}_{m-1,3})_{r,L_{e_1}}^2) &  0  & 0 &\\
          0     &  \mathbb{Z}^{2}\oplus H_0((\hat{\mathcal{H}}_{m-1,3})_{r,L_{e_1}}^2)  & \mathbb{Z}^{4} & 0&\\
    \quad\strut &   0  &  1  &  2 & p\strut \\};
    \draw[-stealth] (m-5-3.west) -- (m-5-2.east) node [above,pos=0.5] {$d^1$};
\draw[thick] (m-1-1.east) -- (m-6-1.east) ;
\draw[thick] (m-6-1.north) -- (m-6-5.north) ;
\end{tikzpicture}
\end{figure}

Since $H_0(U_{11}\cap U_{12})=H_0(U_{11}\cap U_{21})=\emptyset$, we do not  choose generators for their $0$-th homology groups. Let $b_1$ and $b_2$ be the generators of $H_0(U_{22}\cap U_{12})$. Let $\hat{b}_1$ and $\hat{b}_2$ be the generators of $H_0(U_{22}\cap U_{21})$. By induction hypothesis, $H_0(U_{22})\cong\mathbb{Z}^{2}$ and we choose a set of generators (denoted by $f'$ and $\hat{f}'$) for $H_0(U_{22})$. Let $e$ be the generator of $H_0(U_{12})$ and $\hat{e}'$ be the generator of $H_0(U_{21})$ such that for $i=1,2$, 
\begin{equation}
\begin{aligned}
d^1(b_i)&=e'+f' \mbox{\quad and\quad}
d^1(\hat{b}_i)=\hat{e}'+\hat{f}'\\
\end{aligned}
\end{equation}

Under the given generators, $d^1$ can be represented by the following matrix 
$$\left[\begin{array}{cccc}1 & 1 & 0 & 0 \\0 & 0 & 1 & 1 \\1 & 1 & 0 & 0  \\0 & 0 & 1 & 1 \end{array}\right]$$
The rank of this matrix is $2$, and the diagonal elements of the Smith Normal form of the above matrix are all equal to $1$. Therefore, $\Ima d^1\cong\mathbb{Z}^{2}$ and $\coker d^1\cong \mathbb{Z}^{2}$. Hence, we obtain the $E^2$ page of the Mayer-Vietoris spectral sequence: 
\begin{figure}[htbp!]
\centering
\begin{tikzpicture}
  \matrix (m) [matrix of math nodes,
    nodes in empty cells,nodes={minimum width=5ex,
     minimum height=2.5ex,outer sep=-5pt},
    column sep=1ex,row sep=1ex]{
          q     &      &     &  &  \\
          3    &   \vdots   & \vdots    &\vdots  &  \\
          2     & 0& 0  &0  &\\   
          1     &  \mathbb{Z}^{2(m-1)-4} &  0  & 0 &\\
          0     &  \mathbb{Z}^{2}  & \mathbb{Z}^{2} & 0&\\
    \quad\strut &   0  &  1  &  2 & p\strut \\};
   \draw[thick] (m-1-1.east) -- (m-6-1.east) ;
\draw[thick] (m-6-1.north) -- (m-6-5.north) ;
\end{tikzpicture}
\end{figure}

Here we applied the induction hypothesis $H_1((\hat{\mathcal{H}}_{m-1,3})_{r,L_{e_1}}^2)\cong\mathbb{Z}^{2(m-1)-4}$ and $H_0((\hat{\mathcal{H}}_{m-1,3})_{r,L_{e_1}}^2)=\mathbb{Z}^{2}$. Since there is no non-trivial arrow on the $E^2$ page, $E^2=E^{\infty}$. Therefore, $H_0((\hat{\mathcal{H}}_{m,3})_{r,L_{e_1}}^2)\cong \mathbb{Z}^{2}$ and $H_1((\hat{\mathcal{H}}_{m,3})_{r,L_{e_1}}^2)\cong\mathbb{Z}^{2m-4}$.
\end{proof}

\begin{proposition}
Let $L_{e_1}$ be a positive number. If $2<r$ and $L_{e_1}+1<r\leq L_{e_1}+2$, then
$$H_0((\hat{\mathcal{H}}_{m,3})_{r,L_{e_1}}^2)\cong \mathbb{Z}^{4(m-3)+8}  \mbox{\quad and\quad} H_1((\hat{\mathcal{H}}_{m,3})_{r,L_{e_1}}^2)=0$$
\end{proposition}

\begin{proof}
Consider the cover of $(\hat{\mathcal{H}}_{m,3})_{r,L_{e_1}}^2$ given in Equation (\ref{cover:gHgraph}). The $E^1$ page of the Mayer-Vietoris spectral sequence is 
\begin{figure}[htbp!]
\centering
\begin{tikzpicture}
  \matrix (m) [matrix of math nodes,
    nodes in empty cells,nodes={minimum width=5ex,
     minimum height=2.5ex,outer sep=-5pt},
    column sep=1ex,row sep=1ex]{
          q     &      &     &  &  \\
          3    &   \vdots   & \vdots    &\vdots  &  \\
          2     & 0& 0  &0  &\\   
          1     &   H_1((\hat{\mathcal{H}}_{m-1,3})_{r,L_{e_1}}^2) &  0  & 0 &\\
          0     &  \mathbb{Z}^{4}\oplus H_0((\hat{\mathcal{H}}_{m-1,3})_{r,L_{e_1}}^2)  & 0 & 0&\\
    \quad\strut &   0  &  1  &  2 & p\strut \\};
\draw[thick] (m-1-1.east) -- (m-6-1.east) ;
\draw[thick] (m-6-1.north) -- (m-6-5.north) ;
\end{tikzpicture}
\end{figure}

By induction hypothesis, $H_0(U_{22})\cong\mathbb{Z}^{4(m-4)+8}$ and $H_1(U_{22})=0$. Since there is no non-trivial arrow on the $E^1$ page, $E^1=E^{\infty}$. Therefore, $H_0((\hat{\mathcal{H}}_{m,3})_{r,L_{e_1}}^2)\cong \mathbb{Z}^{4(m-3)+8}$ and $H_1((\hat{\mathcal{H}}_{m,3})_{r,L_{e_1}}^2)=0$.
\end{proof}

\paragraph{Step 2: Calculating $H_i(\hat{\mathcal{H}}_{m,n})^2_{r,L_{e_1}}$ where $n\geq 4$}
\begin{lemma}
$\vec{L}=(1,\dots,1,L_{f},L_{f'})\in(\mathbb{R}_{>0})^{m+3}$ where $L_{f}=L_{f'}=\frac{1}{2}L_{e_1}$. Then 
$$H_0(U_{12})\cong  \begin{cases} \mathbb{Z}, & \mbox{if } 0<r\leq 1 \\ 
\mathbb{Z}^{n-1}, & \mbox{if } 1 <r\leq 2  \mbox{\quad and\quad} r\leq 1+L_{e_1}  \\ 
\mathbb{Z}^{m+n-3}, & \mbox{if } 1 <r\leq 2  \mbox{\quad and\quad} r >1+L_{e_1}  \\
\mathbb{Z}, & \mbox{if } 2<r \mbox{\quad and\quad} r\leq 1+L_{e_1} \\
\mathbb{Z}^{m-1}, & \mbox{if } 2<r \mbox{\quad and\quad} 1+L_{e_1}<r\leq 2+L_{e_1}\\
0, & \mbox{else} \end{cases}$$
and 
$$H_1(U_{12})\cong  0$$
\end{lemma}

\begin{proposition}
Let $L_{e_1}$ be a positive number. If $r\leq L_{e_1}$ and $r\leq 1$, then
$$H_0((\hat{\mathcal{H}}_{m,n})_{r,L_{e_1}}^2)\cong \mathbb{Z}  \mbox{\quad and\quad} H_1((\hat{\mathcal{H}}_{m,3})_{r,L_{e_1}}^2)\cong \mathbb{Z}^{m(m-3)+n(n-3)+3}$$
\end{proposition}
\begin{proof}
Consider the cover of $(\hat{\mathcal{H}}_{m,n})_{r,L_{e_1}}^2$ given in Equation (\ref{cover:gHgraph}). The $E^1$ page of the Mayer-Vietoris spectral sequence is 
\begin{figure}[htbp!]
\centering
\begin{tikzpicture}
  \matrix (m) [matrix of math nodes,
    nodes in empty cells,nodes={minimum width=5ex,
     minimum height=2.5ex,outer sep=-5pt},
    column sep=1ex,row sep=1ex]{
          q     &      &     &  &  \\
          3    &   \vdots   & \vdots    &\vdots  &  \\
          2     & 0& 0  &0  &\\   
          1     &   H_1((\hat{\mathcal{H}}_{m,n-1})_{r,L_{e_1}}^2) &  0  & 0 &\\
          0     &  \mathbb{Z}^4\oplus H_0((\hat{\mathcal{H}}_{m,n-1})_{r,L_{e_1}}^2) & \mathbb{Z}^{2n} & 0&\\
    \quad\strut &   0  &  1  &  2 & p\strut \\};
    \draw[-stealth] (m-5-3.west) -- (m-5-2.east) node [above,pos=0.5] {$d^1$};
\draw[thick] (m-1-1.east) -- (m-6-1.east) ;
\draw[thick] (m-6-1.north) -- (m-6-5.north) ;
\end{tikzpicture}
\end{figure}

Let $a$ be the generator of $H_0(U_{11}\cap U_{12})$ and $\hat{a}$ be the generator of $H_0(U_{11}\cap U_{21})$. Let $b'_2,\dots, b'_{n-1}$ and $b$ be the generators of $H_0(U_{22}\cap U_{12})$. Let $\hat{b}'_2, \dots, \hat{b}'_{n-1}$ and $\hat{b}$ be the generators of $H_0(U_{22}\cap U_{21})$. Choose generators (denoted by $c_1$ and $c_2$) for $H_0(U_{11})$. By induction hypothesis, $H_0(U_{22})\cong\mathbb{Z}$ and we choose a generator (denoted by $f$) for $H_0(U_{22})$. Let $e'$ be the generator of $H_0(U_{12})$ and $\hat{e}'$ be the generator of $H_0(U_{21})$ such that
\begin{equation}
\begin{aligned}
d^1(a)&=c_1+e' \mbox{\quad and\quad} d^1(\hat{a})=c_2+\hat{e}'\\
d^1(b_i)&=e'+f \mbox{\quad and\quad} d^1(\hat{b}_i)=\hat{e}'+f, \quad i=2,\dots,n-1\\
d^1(b')&=e'+f \mbox{\quad and\quad} d^1(\hat{b}')=\hat{e}'+f
\end{aligned}
\end{equation}

Under the given generators, $d^1$ can be represented by the following matrix 
$$\left[\begin{array}{cccccccc}1 & 0 & 0 & \cdots  & 0 & 0 & \cdots  & 0 \\0 & 1 & 0 & \cdots  & 0 & 0 & \cdots  & 0 \\0 & 0 & 1 & \cdots & 1 & 1 & \cdots & 1 \\1 & 0 & 1 & \cdots & 1 & 0 & \cdots  & 0 \\0 & 1 & 0 & \cdots  & 0 & 1 & \cdots & 1\end{array}\right]$$
The rank of this matrix is $4$, and the diagonal elements of the Smith Normal form of the above matrix are all equal to $1$. Therefore, $\Ima d^1 \cong\mathbb{Z}^4$ and $\ker d^1\cong \mathbb{Z}^{2n-4}$. Hence, we obtain the $E^2$ page of the Mayer-Vietoris spectral sequence: 
\begin{figure}[htbp!]
\centering
\begin{tikzpicture}
  \matrix (m) [matrix of math nodes,
    nodes in empty cells,nodes={minimum width=5ex,
     minimum height=2.5ex,outer sep=-5pt},
    column sep=1ex,row sep=1ex]{
          q     &      &     &  &  \\
          3    &   \vdots   & \vdots    &\vdots  &  \\
          2     & 0& 0  &0  &\\   
          1     &  \mathbb{Z}^{(n-1)(n-4)+m(m-3)+3}  &  0  & 0 &\\
          0     &  \mathbb{Z} & \mathbb{Z}^{2n-4} & 0&\\
    \quad\strut &   0  &  1  &  2 & p\strut \\};
   \draw[thick] (m-1-1.east) -- (m-6-1.east) ;
\draw[thick] (m-6-1.north) -- (m-6-5.north) ;
\end{tikzpicture}
\end{figure}

Here we applied the induction hypothesis $H_1((\hat{\mathcal{H}}_{m,n-1})_{r,L_{e_1}}^2)\cong \mathbb{Z}^{(n-1)(n-4)+m(m-3)+3}$. Since there is no non-trivial arrow on the $E^2$ page, $E^2=E^{\infty}$. Therefore, $$H_0((\hat{\mathcal{H}}_{m,n})_{r,L_{e_1}}^2)\cong \mathbb{Z}$$ and $$H_1((\hat{\mathcal{H}}_{m,n})_{r,L_{e_1}}^2)\cong \mathbb{Z}^{m^2-5m+7}\oplus \mathbb{Z}^{2m-4}\cong\mathbb{Z}^{m(m-3)+n(n-3)+3}$$
\end{proof}

\begin{proposition}
Let $L_{e_1}$ be a positive number. If $L_{e_1}<r\leq 1$, then
$$H_0((\hat{\mathcal{H}}_{m,n})_{r,L_{e_1}}^2)\cong \mathbb{Z}  \mbox{\quad and\quad} H_1((\hat{\mathcal{H}}_{m,n})_{r,L_{e_1}}^2)\cong \mathbb{Z}^{(m+n)(m+n-7)+11}$$
\end{proposition}
\begin{proof}
Consider the cover of $(\hat{\mathcal{H}}_{m,n})_{r,L_{e_1}}^2$ given in Equation (\ref{cover:gHgraph}). The $E^1$ page of the Mayer-Vietoris spectral sequence is 
\begin{figure}[htbp!]
\centering
\begin{tikzpicture}
  \matrix (m) [matrix of math nodes,
    nodes in empty cells,nodes={minimum width=5ex,
     minimum height=2.5ex,outer sep=-5pt},
    column sep=1ex,row sep=1ex]{
          q     &      &     &  &  \\
          3    &   \vdots   & \vdots    &\vdots  &  \\
          2     & 0& 0  &0  &\\   
          1     &   H_1((\hat{\mathcal{H}}_{m,n-1})_{r,L_{e_1}}^2) &  0  & 0 &\\
          0     &  \mathbb{Z}^4\oplus H_0((\hat{\mathcal{H}}_{m,n-1})_{r,L_{e_1}}^2) & \mathbb{Z}^{2m+2n-4} & 0&\\
    \quad\strut &   0  &  1  &  2 & p\strut \\};
    \draw[-stealth] (m-5-3.west) -- (m-5-2.east) node [above,pos=0.5] {$d^1$};
\draw[thick] (m-1-1.east) -- (m-6-1.east) ;
\draw[thick] (m-6-1.north) -- (m-6-5.north) ;
\end{tikzpicture}
\end{figure}

Let $a$ be the generator of $H_0(U_{11}\cap U_{12})$ and $\hat{a}$ be the generator of $H_0(U_{11}\cap U_{21})$. Let $b'_2,\dots, b'_{n-1}$, $b_2$,\dots, $b_{m-1}$ and $b_m$ be the generators of $H_0(U_{22}\cap U_{12})$. Let $\hat{b}'_2, \dots, \hat{b}'_{n-1}$,  $\hat{b}_2$,\dots, $\hat{b}_{m-1}$ and $\hat{b}_m$ be the generators of $H_0(U_{22}\cap U_{21})$. Choose generators (denoted by $c_1$ and $c_2$) for $H_0(U_{11})$. By induction hypothesis, $H_0(U_{22})\cong\mathbb{Z}$ and we choose a generator (denoted by $f$) for $H_0(U_{22})$. Let $e'$ be the generator of $H_0(U_{12})$ and $\hat{e}'$ be the generator of $H_0(U_{21})$ such that
\begin{equation}
\begin{aligned}
d^1(a)&=c_1+e' \mbox{\quad and\quad} d^1(\hat{a})=c_2+\hat{e}'\\
d^1(b_i)&=e'+f \mbox{\quad and\quad} d^1(\hat{b}_i)=\hat{e}'+f, \quad i=2,\dots,n-1\\
d^1(b'_j)&=e'+f \mbox{\quad and\quad} d^1(\hat{b}'_j)=\hat{e}'+f, \quad j=2,\dots,m
\end{aligned}
\end{equation}

Under the given generators, $d^1$ can be represented by the following matrix 
$$\left[\begin{array}{cccccccc}1 & 0 & 0 & \cdots  & 0 & 0 & \cdots  & 0 \\0 & 1 & 0 & \cdots  & 0 & 0 & \cdots  & 0 \\0 & 0 & 1 & \cdots & 1 & 1 & \cdots & 1 \\1 & 0 & 1 & \cdots & 1 & 0 & \cdots  & 0 \\0 & 1 & 0 & \cdots  & 0 & 1 & \cdots & 1\end{array}\right]$$
The rank of this matrix is $4$, and the diagonal elements of the Smith Normal form of the above matrix are all equal to $1$. Therefore, $\Ima d^1\cong\mathbb{Z}^4$ and $\ker d^1\cong \mathbb{Z}^{2m+2n-8}$. Hence, we obtain the $E^2$ page of the Mayer-Vietoris spectral sequence: 
\begin{figure}[htbp!]
\centering
\begin{tikzpicture}
  \matrix (m) [matrix of math nodes,
    nodes in empty cells,nodes={minimum width=5ex,
     minimum height=2.5ex,outer sep=-5pt},
    column sep=1ex,row sep=1ex]{
          q     &      &     &  &  \\
          3    &   \vdots   & \vdots    &\vdots  &  \\
          2     & 0& 0  &0  &\\   
          1     &  \mathbb{Z}^{(m+n-1)(m+n-8)+11}  &  0  & 0 &\\
          0     &  \mathbb{Z} & \mathbb{Z}^{2m+2n-8} & 0&\\
    \quad\strut &   0  &  1  &  2 & p\strut \\};
   \draw[thick] (m-1-1.east) -- (m-6-1.east) ;
\draw[thick] (m-6-1.north) -- (m-6-5.north) ;
\end{tikzpicture}
\end{figure}

Here we applied the induction hypothesis $H_1((\hat{\mathcal{H}}_{m,n-1})_{r,L_{e_1}}^2)\cong \mathbb{Z}^{(m+n-1)(m+n-8)+11}$. Since there is no non-trivial arrow on the $E^2$ page, $E^2=E^{\infty}$. Therefore, $$H_0((\hat{\mathcal{H}}_{m,n})_{r,L_{e_1}}^2)\cong \mathbb{Z}$$ and $$ H_1((\hat{\mathcal{H}}_{m,n})_{r,L_{e_1}}^2)\cong \mathbb{Z}^{(m+n-1)(m+n-8)+10}\oplus \mathbb{Z}^{2m+2n-8}\cong\mathbb{Z}^{(m+n)(m+n-7)+11}$$
\end{proof}

\begin{proposition}
Let $L_{e_1}$ be a positive number. If $1<r\leq 2$ and $r\leq L_{e_1}$, then
$$H_0((\hat{\mathcal{H}}_{m,n})_{r,L_{e_1}}^2)\cong \mathbb{Z}^{m(m-3)+n(n-3)+6}  \mbox{\quad and\quad} H_1((\hat{\mathcal{H}}_{m,n})_{r,L_{e_1}}^2)=0$$
\end{proposition}
\begin{proof}
Consider the cover of $(\hat{\mathcal{H}}_{m,n})_{r,L_{e_1}}^2$ given in Equation (\ref{cover:gHgraph}). The $E^1$ page of the Mayer-Vietoris spectral sequence is \\

\begin{figure}[htbp!]
\centering
\begin{tikzpicture}
  \matrix (m) [matrix of math nodes,
    nodes in empty cells,nodes={minimum width=5ex,
     minimum height=2.5ex,outer sep=-5pt},
    column sep=1ex,row sep=1ex]{
          q     &      &     &  &  \\
          3    &   \vdots   & \vdots    &\vdots  &  \\
          2     & 0& 0  &0  &\\   
          1     &   H_1((\hat{\mathcal{H}}_{m,n-1})_{r,L_{e_1}}^2) &  0  & 0 &\\
          0     &  \mathbb{Z}^{2n-2}\oplus H_0((\hat{\mathcal{H}}_{m,n-1})_{r,L_{e_1}}^2) & \mathbb{Z}^{2} & 0&\\
    \quad\strut &   0  &  1  &  2 & p\strut \\};
    \draw[-stealth] (m-5-3.west) -- (m-5-2.east) node [above,pos=0.5] {$d^1$};
\draw[thick] (m-1-1.east) -- (m-6-1.east) ;
\draw[thick] (m-6-1.north) -- (m-6-5.north) ;
\end{tikzpicture}
\end{figure}

Let $b$ be the generator of $H_0(U_{22}\cap U_{12})$ and $\hat{b}$ be the generator of $H_0(U_{22}\cap U_{21})$. By induction hypothesis, $H_0(U_{22})\cong\mathbb{Z}^{m(m-3)+(n-1)(n-4)+6}$ and we choose generators (denoted by $f_{ij}$, $\hat{f}_{ij}$ for $2\leq i<j\leq m$ and $m+2\leq i<j\leq m+n-1$, $f_{2,m+2}$ and $\hat{f}_{2,m+2}$) for $H_0(U_{22})$. Let $e'_2, \dots, e'_{n-1}, e$ be the generators of $H_0(U_{12})$ and $\hat{e}'_1,\dots, \hat{e}'_{n-1}, \hat{e}$ be the generator of $H_0(U_{21})$ such that
\begin{equation}
\begin{aligned}
d^1(b)&=e'+f_{2,m+2} \mbox{\quad and\quad} d^1(\hat{b})=\hat{e}'+\hat{f}_{2,m+2}
\end{aligned}
\end{equation}

Under the given generators, $d^1$ can be represented by the following matrix. 

$$\left[\begin{array}{cc}1 & 0 \\0 & 1 \\1 & 0\\0 & 1\end{array}\right]$$
The rank of this matrix is $2$, and the diagonal elements of the Smith Normal form of the above matrix are all equal to $1$. Therefore, $\Ima d^1\cong\mathbb{Z}^2$ and $\ker d^1=0$. Hence, we obtain the $E^2$ page of the Mayer-Vietoris spectral sequence: 
\begin{figure}[htbp!]
\centering
\begin{tikzpicture}
  \matrix (m) [matrix of math nodes,
    nodes in empty cells,nodes={minimum width=5ex,
     minimum height=2.5ex,outer sep=-5pt},
    column sep=1ex,row sep=1ex]{
          q     &      &     &  &  \\
          3    &   \vdots   & \vdots    &\vdots  &  \\
          2     & 0& 0  &0  &\\   
          1     &  0 &  0  & 0 &\\
          0     &  \mathbb{Z}^{m(m-3)+n(n-3)+6} & 0 & 0&\\
    \quad\strut &   0  &  1  &  2 & p\strut \\};
   \draw[thick] (m-1-1.east) -- (m-6-1.east) ;
\draw[thick] (m-6-1.north) -- (m-6-5.north) ;
\end{tikzpicture}
\end{figure}

Here we applied the induction hypothesis $H_1((\hat{\mathcal{H}}_{m,n-1})_{r,L_{e_1}}^2)=0$. Since there is no non-trivial arrow on the $E^2$ page, $E^2=E^{\infty}$. Therefore, $$H_0((\hat{\mathcal{H}}_{m,n})_{r,L_{e_1}}^2)\cong \mathbb{Z}^{m(m-3)+n(n-3)+6}$$
$$H_1((\hat{\mathcal{H}}_{m,n})_{r,L_{e_1}}^2)=0$$
\end{proof}

\begin{proposition}
Let $L_{e_1}$ be a positive number. If $1<r\leq 2$ and $L_{e_1}<r\leq L_{e_1}+1$, then
$$H_0((\hat{\mathcal{H}}_{m,n})_{r,L_{e_1}}^2)\cong \mathbb{Z}^{m(m-3)+n(n-3)+6}  \mbox{\quad and\quad} H_1((\hat{\mathcal{H}}_{m,n})_{r,L_{e_1}}^2)\cong \mathbb{Z}^{2(m-2)(n-2)} $$
\end{proposition}
\begin{proof}
Consider the cover of $(\hat{\mathcal{H}}_{m,n})_{r,L_{e_1}}^2$ given in Equation (\ref{cover:gHgraph}). The $E^1$ page of the Mayer-Vietoris spectral sequence is\\

\begin{figure}[htbp!]
\centering
\begin{tikzpicture}
  \matrix (m) [matrix of math nodes,
    nodes in empty cells,nodes={minimum width=5ex,
     minimum height=2.5ex,outer sep=-5pt},
    column sep=1ex,row sep=1ex]{
          q     &      &     &  &  \\
          3    &   \vdots   & \vdots    &\vdots  &  \\
          2     & 0& 0  &0  &\\   
          1     &   H_1((\hat{\mathcal{H}}_{m,n-1})_{r,L_{e_1}}^2) &  0  & 0 &\\
          0     &  \mathbb{Z}^{2n-2}\oplus H_0((\hat{\mathcal{H}}_{m,n-1})_{r,L_{e_1}}^2) & \mathbb{Z}^{2m-2} & 0&\\
    \quad\strut &   0  &  1  &  2 & p\strut \\};
    \draw[-stealth] (m-5-3.west) -- (m-5-2.east) node [above,pos=0.5] {$d^1$};
\draw[thick] (m-1-1.east) -- (m-6-1.east) ;
\draw[thick] (m-6-1.north) -- (m-6-5.north) ;
\end{tikzpicture}
\end{figure}

Let $b_2,\dots, b_m$ be the generators of $H_0(U_{22}\cap U_{12})$ and $\hat{b}_2,\dots,\hat{b}_m$ be the generators of $H_0(U_{22}\cap U_{21})$. By induction hypothesis, $H_0(U_{22})\cong\mathbb{Z}^{m(m-3)+(n-1)(n-4)+6}$ and we choose generators (denoted by $f_{ij}$, $\hat{f}_{ij}$ for $2\leq i<j\leq m$ and $m+2\leq i<j\leq m+n-1$, $f_{2,m+2}$ and $\hat{f}_{2,m+2}$) for $H_0(U_{22})$. Let $e'_2, \dots, e'_{n-1}, e$ be the generators of $H_0(U_{12})$ and $\hat{e}'_1,\dots, \hat{e}'_{n-1}, \hat{e}$ be the generator of $H_0(U_{21})$ such that
\begin{equation}
\begin{aligned}
d^1(b)&=e'+f_{2,m+2} \mbox{\quad and\quad} d^1(\hat{b})=\hat{e}'+\hat{f}_{2,m+2}
\end{aligned}
\end{equation}

Under the given generators, $d^1$ can be represented by the following matrix 
$$\left[\begin{array}{cc}1 & 0 \\0 & 1 \\1 & 0\\0 & 1\end{array}\right]$$
The rank of this matrix is $2$, and the diagonal elements of the Smith Normal form of the above matrix are all equal to $1$. Therefore, $\Ima d^1\cong\mathbb{Z}^2$ and $\ker d^1=2m-4$. Hence, we obtain the $E^2$ page of the Mayer-Vietoris spectral sequence: 
\begin{figure}[htbp!]
\centering
\begin{tikzpicture}
  \matrix (m) [matrix of math nodes,
    nodes in empty cells,nodes={minimum width=5ex,
     minimum height=2.5ex,outer sep=-5pt},
    column sep=1ex,row sep=1ex]{
          q     &      &     &  &  \\
          3    &   \vdots   & \vdots    &\vdots  &  \\
          2     & 0& 0  &0  &\\   
          1     &  \mathbb{Z}^{2(m-2)(n-3)} &  0  & 0 &\\
          0     &  \mathbb{Z}^{m(m-3)+n(n-3)+6} & \mathbb{Z}^{2m-4}  & 0&\\
    \quad\strut &   0  &  1  &  2 & p\strut \\};
   \draw[thick] (m-1-1.east) -- (m-6-1.east) ;
\draw[thick] (m-6-1.north) -- (m-6-5.north) ;
\end{tikzpicture}
\end{figure}

Here we applied the induction hypothesis $H_1((\hat{\mathcal{H}}_{m,n-1})_{r,L_{e_1}}^2)\cong\mathbb{Z}^{2(m-2)(n-3)}$. Since there is no non-trivial arrow on the $E^2$ page, $E^2=E^{\infty}$. Therefore, $$H_0((\hat{\mathcal{H}}_{m,n})_{r,L_{e_1}}^2)\cong \mathbb{Z}^{m(m-3)+n(n-3)+6}$$ $$H_1((\hat{\mathcal{H}}_{m,n})_{r,L_{e_1}}^2)\cong\mathbb{Z}^{2(m-2)(n-2)}$$
\end{proof}

\begin{proposition}
Let $L_{e_1}$ be a positive number. If $1<r\leq 2$ and $L_{e_1}+1<r\leq L_{e_1}+2$, then
$$H_0((\hat{\mathcal{H}}_{m,n})_{r,L_{e_1}}^2)\cong \mathbb{Z}^{(m+n)(m+n-5)+6}  \mbox{\quad and\quad} H_1((\hat{\mathcal{H}}_{m,n})_{r,L_{e_1}}^2)=0$$
\end{proposition}
\begin{proof}
Consider the cover of $(\hat{\mathcal{H}}_{m,n})_{r,L_{e_1}}^2$ given in Equation (\ref{cover:gHgraph}). The $E^1$ page of the Mayer-Vietoris spectral sequence is\\

\begin{figure}[htbp!]
\centering
\begin{tikzpicture}
  \matrix (m) [matrix of math nodes,
    nodes in empty cells,nodes={minimum width=5ex,
     minimum height=2.5ex,outer sep=-5pt},
    column sep=1ex,row sep=1ex]{
          q     &      &     &  &  \\
          3    &   \vdots   & \vdots    &\vdots  &  \\
          2     & 0& 0  &0  &\\   
          1     &   H_1((\hat{\mathcal{H}}_{m,n-1})_{r,L_{e_1}}^2) &  0  & 0 &\\
          0     &  \mathbb{Z}^{2m+2n-6}\oplus H_0((\hat{\mathcal{H}}_{m,n-1})_{r,L_{e_1}}^2) & 0 & 0&\\
    \quad\strut &   0  &  1  &  2 & p\strut \\};
\draw[thick] (m-1-1.east) -- (m-6-1.east) ;
\draw[thick] (m-6-1.north) -- (m-6-5.north) ;
\end{tikzpicture}
\end{figure}

Since there is no non-trivial arrow on the $E^1$ page, $E^1=E^{\infty}$. Therefore, $H_0((\hat{\mathcal{H}}_{m,n})_{r,L_{e_1}}^2)\cong \mathbb{Z}^{(m+n)(m+n-5)+6}$ and $H_1((\hat{\mathcal{H}}_{m,n})_{r,L_{e_1}}^2)=0$.
\end{proof}

\begin{proposition}
Let $L_{e_1}$ be a positive number. If $2<r$ and $r\leq L_{e_1}$, then
$$H_0((\hat{\mathcal{H}}_{m,n})_{r,L_{e_1}}^2)\cong \mathbb{Z}^{2}  \mbox{\quad and\quad} H_1((\hat{\mathcal{H}}_{m,n})_{r,L_{e_1}}^2)=0$$
\end{proposition}
\begin{proof}
Consider the cover of $(\hat{\mathcal{H}}_{m,n})_{r,L_{e_1}}^2$ given in Equation (\ref{cover:gHgraph}). The $E^1$ page of the Mayer-Vietoris spectral sequence is 
\begin{figure}[htbp!]
\centering
\begin{tikzpicture}
  \matrix (m) [matrix of math nodes,
    nodes in empty cells,nodes={minimum width=5ex,
     minimum height=2.5ex,outer sep=-5pt},
    column sep=1ex,row sep=1ex]{
          q     &      &     &  &  \\
          3    &   \vdots   & \vdots    &\vdots  &  \\
          2     & 0& 0  &0  &\\   
          1     &   H_1((\hat{\mathcal{H}}_{m,n-1})_{r,L_{e_1}}^2) &  0  & 0 &\\
          0     &  \mathbb{Z}^{2}\oplus H_0((\hat{\mathcal{H}}_{m,n-1})_{r,L_{e_1}}^2) & \mathbb{Z}^{2} & 0&\\
    \quad\strut &   0  &  1  &  2 & p\strut \\};
    \draw[-stealth] (m-5-3.west) -- (m-5-2.east) node [above,pos=0.5] {$d^1$};
\draw[thick] (m-1-1.east) -- (m-6-1.east) ;
\draw[thick] (m-6-1.north) -- (m-6-5.north) ;
\end{tikzpicture}
\end{figure}

Let $b$ be the generator of $H_0(U_{22}\cap U_{12})$ and $\hat{b}$ be the generator of $H_0(U_{22}\cap U_{21})$. By induction hypothesis, $H_0(U_{22})\cong\mathbb{Z}^{2}$ and we choose generators (denoted by $f_{2,m+2}$ and $\hat{f}_{2,m+2}$) for $H_0(U_{22})$. Let $e$ be the generator of $H_0(U_{12})$ and $\hat{e}$ be the generator of $H_0(U_{21})$ such that
\begin{equation}
\begin{aligned}
d^1(b)&=e+f_{2,m+2} \mbox{\quad and\quad} d^1(\hat{b})=\hat{e}+\hat{f}_{2,m+2}
\end{aligned}
\end{equation}

Under the given generators, $d^1$ can be represented by the following matrix 
$$\left[\begin{array}{cc}1 & 0 \\0 & 1 \\1 & 0\\0 & 1\end{array}\right]$$
The rank of this matrix is $2$, and the diagonal elements of the Smith Normal form of the above matrix are all equal to $1$. Therefore, $\Ima d^1\cong\mathbb{Z}^2$ and $\ker d^1=0$. Hence, we obtain the $E^2$ page of the Mayer-Vietoris spectral sequence: 
\begin{figure}[htbp!]
\centering
\begin{tikzpicture}
  \matrix (m) [matrix of math nodes,
    nodes in empty cells,nodes={minimum width=5ex,
     minimum height=2.5ex,outer sep=-5pt},
    column sep=1ex,row sep=1ex]{
          q     &      &     &  &  \\
          3    &   \vdots   & \vdots    &\vdots  &  \\
          2     & 0& 0  &0  &\\   
          1     &  0 &  0  & 0 &\\
          0     &  \mathbb{Z}^{2} & 0 & 0&\\
    \quad\strut &   0  &  1  &  2 & p\strut \\};
   \draw[thick] (m-1-1.east) -- (m-6-1.east) ;
\draw[thick] (m-6-1.north) -- (m-6-5.north) ;
\end{tikzpicture}
\end{figure}

Here we applied the induction hypothesis $H_1((\hat{\mathcal{H}}_{m,n-1})_{r,L_{e_1}}^2)=0$ and $H_0((\hat{\mathcal{H}}_{m,n-1})_{r,L_{e_1}}^2)\cong\mathbb{Z}^2$. Since there is no non-trivial arrow on the $E^2$ page, $E^2=E^{\infty}$. Therefore, $H_0((\hat{\mathcal{H}}_{m,n})_{r,L_{e_1}}^2)\cong \mathbb{Z}^{2}$ and $H_1((\hat{\mathcal{H}}_{m,n})_{r,L_{e_1}}^2)=0$.
\end{proof}

\begin{proposition}
Let $L_{e_1}$ be a positive number. If $2<r$ and $L_{e_1}<r\leq L_{e_1}+1$, then
$$H_0((\hat{\mathcal{H}}_{m,n})_{r,L_{e_1}}^2)\cong \mathbb{Z}^{2}  \mbox{\quad and\quad} H_1((\hat{\mathcal{H}}_{m,n})_{r,L_{e_1}}^2)\cong\mathbb{Z}^{2(m-2)(n-2)}$$
\end{proposition}
\begin{proof}
Consider the cover of $(\hat{\mathcal{H}}_{m,n})_{r,L_{e_1}}^2$ given in Equation (\ref{cover:gHgraph}). The $E^1$ page of the Mayer-Vietoris spectral sequence is 
\begin{figure}[htbp!]
\centering
\begin{tikzpicture}
  \matrix (m) [matrix of math nodes,
    nodes in empty cells,nodes={minimum width=5ex,
     minimum height=2.5ex,outer sep=-5pt},
    column sep=1ex,row sep=1ex]{
          q     &      &     &  &  \\
          3    &   \vdots   & \vdots    &\vdots  &  \\
          2     & 0& 0  &0  &\\   
          1     &   H_1((\hat{\mathcal{H}}_{m,n-1})_{r,L_{e_1}}^2) &  0  & 0 &\\
          0     &  \mathbb{Z}^{2}\oplus H_0((\hat{\mathcal{H}}_{m,n-1})_{r,L_{e_1}}^2) & \mathbb{Z}^{2m-2} & 0&\\
    \quad\strut &   0  &  1  &  2 & p\strut \\};
    \draw[-stealth] (m-5-3.west) -- (m-5-2.east) node [above,pos=0.5] {$d^1$};
\draw[thick] (m-1-1.east) -- (m-6-1.east) ;
\draw[thick] (m-6-1.north) -- (m-6-5.north) ;
\end{tikzpicture}
\end{figure}

Let $b_2, \dots, b_m$ be the generators of $H_0(U_{22}\cap U_{12})$ and $\hat{b}_2,\dots,\hat{b}_m $ be the generators of $H_0(U_{22}\cap U_{21})$. By induction hypothesis, $H_0(U_{22})\cong\mathbb{Z}^{2}$ and we choose generators (denoted by $f_{2,m+2}$ and $\hat{f}_{2,m+2}$) for $H_0(U_{22})$. Let $e$ be the generator of $H_0(U_{12})$ and $\hat{e}$ be the generator of $H_0(U_{21})$ such that
\begin{equation}
\begin{aligned}
d^1(b_i)&=e+f_{2,m+2},  \mbox{\quad and\quad} d^1(\hat{b}_i)=\hat{e}+\hat{f}_{2,m+2}, \quad i=2,\dots, m
\end{aligned}
\end{equation}

Under the given generators, $d^1$ can be represented by the following matrix 
$$\left[\begin{array}{cccccc}1 & \cdots & 1 & 0 & \cdots & 0 \\0 & \cdots & 0 & 1 & \cdots & 1 \\1 & \cdots & 1 & 0 & \cdots & 0 \\0 & \cdots & 0 & 1 & \cdots & 1\end{array}\right]$$
The rank of this matrix is $2$, and the diagonal elements of the Smith Normal form of the above matrix are all equal to $1$. Therefore, $\Ima d^1\cong\mathbb{Z}^2$ and $\ker d^1=2m-4$. Hence, we obtain the $E^2$ page of the Mayer-Vietoris spectral sequence: 
\begin{figure}[htbp!]
\centering
\begin{tikzpicture}
  \matrix (m) [matrix of math nodes,
    nodes in empty cells,nodes={minimum width=5ex,
     minimum height=2.5ex,outer sep=-5pt},
    column sep=1ex,row sep=1ex]{
          q     &      &     &  &  \\
          3    &   \vdots   & \vdots    &\vdots  &  \\
          2     & 0& 0  &0  &\\   
          1     & \mathbb{Z}^{2(m-2)(n-3)}  &  0  & 0 &\\
          0     &  \mathbb{Z}^{2} & \mathbb{Z}^{2m-4} & 0&\\
    \quad\strut &   0  &  1  &  2 & p\strut \\};
   \draw[thick] (m-1-1.east) -- (m-6-1.east) ;
\draw[thick] (m-6-1.north) -- (m-6-5.north) ;
\end{tikzpicture}
\end{figure}

Here we applied the induction hypothesis $H_1((\hat{\mathcal{H}}_{m,n-1})_{r,L_{e_1}}^2)\cong \mathbb{Z}^{2(m-2)(n-3)}$ and $H_0((\hat{\mathcal{H}}_{m,n-1})_{r,L_{e_1}}^2)\cong\mathbb{Z}^2$. Since there is no non-trivial arrow on the $E^2$ page, $E^2=E^{\infty}$. Therefore, $H_0((\hat{\mathcal{H}}_{m,n})_{r,L_{e_1}}^2)\cong \mathbb{Z}^{2}$ and $ H_1((\hat{\mathcal{H}}_{m,n})_{r,L_{e_1}}^2)\cong\mathbb{Z}^{2(m-2)(n-2)}$.
\end{proof}

\begin{proposition}
Let $L_{e_1}$ be a positive number. If $2<r$ and $L_{e_1}+1<r\leq L_{e_1}+2$, then
$$H_0((\hat{\mathcal{H}}_{m,n})_{r,L_{e_1}}^2)\cong \mathbb{Z}^{2(m-1)(n-1)}  \mbox{\quad and\quad} H_1((\hat{\mathcal{H}}_{m,n})_{r,L_{e_1}}^2)=0$$
\end{proposition}
\begin{proof}
Consider the cover of $(\hat{\mathcal{H}}_{m,n})_{r,L_{e_1}}^2$ given in Equation (\ref{cover:gHgraph}). The $E^1$ page of the Mayer-Vietoris spectral sequence is 
\begin{figure}[htbp!]
\centering
\begin{tikzpicture}
  \matrix (m) [matrix of math nodes,
    nodes in empty cells,nodes={minimum width=5ex,
     minimum height=2.5ex,outer sep=-5pt},
    column sep=1ex,row sep=1ex]{
          q     &      &     &  &  \\
          3    &   \vdots   & \vdots    &\vdots  &  \\
          2     & 0& 0  &0  &\\   
          1     &   0 &  0  & 0 &\\
          0     &  \mathbb{Z}^{2m-2}\oplus \mathbb{Z}^{2(m-1)(n-2)} & 0 & 0&\\
    \quad\strut &   0  &  1  &  2 & p\strut \\};
\draw[thick] (m-1-1.east) -- (m-6-1.east) ;
\draw[thick] (m-6-1.north) -- (m-6-5.north) ;
\end{tikzpicture}
\end{figure}

Here we applied the induction hypothesis $H_1((\hat{\mathcal{H}}_{m,n-1})_{r,L_{e_1}}^2)=0$ and $H_0((\hat{\mathcal{H}}_{m,n-1})_{r,L_{e_1}}^2)\cong\mathbb{Z}^{2(m-1)(n-2)}$. Since there is no non-trivial arrow on the $E^1$ page, $E^1=E^{\infty}$. Therefore, $H_0((\hat{\mathcal{H}}_{m,n})_{r,L_{e_1}}^2)\cong \mathbb{Z}^{2(m-1)(n-1)}$ and $ H_1((\hat{\mathcal{H}}_{m,n})_{r,L_{e_1}}^2)=0$.
\end{proof}

In summary, when $m,n\geq 3$, The rank of $H_0((\hat{\mathcal{H}}_{m,n})^2_{r,L_{e_1}})$ for all $r>0$ and $L_{e_1}>0$ is shown in Figure \ref{fig:h0Hmn}. When $k\geq 4$, the rank of $H_1((\hat{\mathcal{H}}_{m,n})^2_{r,L_{e_1}})$ for all $r>0$ and $L_{e_1}>0$ is shown in Figure \ref{fig:h1Hmn}.

\section{Decomposition of $PH_i((\hat{\mathcal{H}}_{m,n})_{-,-}^2;\mathbb{F})$}\label{decompo:Hmn}

In this section, we are going to decompose $PH_i((\hat{\mathcal{H}}_{m,n})_{-,-}^2;\mathbb{F})$ where $i=0,1$. We associate $PH_0((\hat{\mathcal{H}}_{m,n})_{-,-}^2;\mathbb{F})$ with a persistence module $N$ which can be constructed using Proposition \ref{prop-sim-2}. We use the following diagram to represent $N$, where each vector space in the diagram represents a constant functor indexed by a chamber of the parameter space of $PH_0((\hat{\mathcal{H}}_{m,n})_{-,-}^2;\mathbb{F})$.

\begin{equation}\label{h0genHmn1}
N=\begin{tikzcd}
  \mathbb{F} & \mathbb{F}^{m^2+n^2-3m-3n+6} \arrow{l}[swap]{} &\mathbb{F}^2\arrow{l}[swap]{}\\
  \mathbb{F} \arrow{u}{}  &\mathbb{F}^{m^2+n^2-3m-3n+6} \arrow{u}{}\arrow{l}[swap]{}  & \mathbb{F}^{2}\arrow{u}[swap]{}\arrow{l}[swap]{}\\ 
  	\mathbb{F}\arrow{u}{}	 &\mathbb{F}^{m^2+n^2-3m-3n+6} \arrow{u}{}\arrow{l}[swap]{}	 & \mathbb{F}^{2}\arrow{u}[swap]{}\\ 
	  	  	 &\mathbb{F}^{m^2+n^2-3m-3n+6} \arrow{u}{} & \mathbb{F}^{2}\arrow{u}[swap]{}\arrow{l}[swap]{}\\ 
			 &\mathbb{F}^{m^2+n^2+2mn-5m-5n+6} \arrow{u}{} & \mathbb{F}^{2(m-1)(n-1)}\arrow{u}[swap]{}\arrow{l}[swap]{}\\ 
\end{tikzcd}
\end{equation}

Figure \ref{fig:Hmn-h0-generator} provides the configurations of two distinct robots (represented by red and green dots) that correspond to the representatives we choose for the bases of the vector spaces in $N$.\\

\begin{figure}[htbp!]
    \centering
    \resizebox{17.5cm}{!}{

\tikzset{every picture/.style={line width=0.75pt}} 


    }
    \caption{}
    \label{fig:Hmn-h0-generator}
\end{figure}

\begin{equation}\label{h0genHmn2}
N=\begin{tikzcd}
<f>                  & \left\langle\parbox{7cm}{$
                e'_2,\dots, e'_{n-1},\hat{e}'_2,\dots, \hat{e}'_{n-1}, f_{2,m+2},\hat{f}_{2,m+2},\\ f_{ij},\hat{f}_{ij}\mid 2\leq i<j\leq m, m+2\leq i<j\leq m+n-1$}\right\rangle \arrow{l}[swap]{\alpha} 				&<f_{2,m+2},\hat{f}_{2,m+2}>\arrow{l}[swap]{\iota}\\
<f> \arrow{u}{}  &\left\langle\parbox{7cm}{$
                e'_2,\dots, e'_{n-1},\hat{e}'_2,\dots, \hat{e}'_{n-1}, f_{2,m+2},\hat{f}_{2,m+2},\\ f_{ij},\hat{f}_{ij}\mid 2\leq i<j\leq m, m+2\leq i<j\leq m+n-1$}\right\rangle  \arrow{u}{}\arrow{l}[swap]{\alpha}  		& <f_{2,m+2},\hat{f}_{2,m+2}>\arrow{u}[swap]{}\arrow{l}[swap]{\iota}\\ 
<f>\arrow{u}{}	 &\left\langle\parbox{7cm}{$
                e'_2,\dots, e'_{n-1},\hat{e}'_2,\dots, \hat{e}'_{n-1}, f_{2,m+2},\hat{f}_{2,m+2},\\ f_{ij},\hat{f}_{ij}\mid 2\leq i<j\leq m, m+2\leq i<j\leq m+n-1$}\right\rangle  \arrow{u}{}\arrow{l}[swap]{\alpha}	 	& <f_{2,m+2},\hat{f}_{2,m+2}>\arrow{u}[swap]{}\\ 
	  	  	 &\left\langle\parbox{7cm}{$
                e'_2,\dots, e'_{n-1},\hat{e}'_2,\dots, \hat{e}'_{n-1}, f_{2,m+2},\hat{f}_{2,m+2},\\ f_{ij},\hat{f}_{ij}\mid 2\leq i<j\leq m, m+2\leq i<j\leq m+n-1$}\right\rangle  \arrow{u}{} 					& <f_{2,m+2},\hat{f}_{2,m+2}>\arrow{u}[swap]{}\arrow{l}[swap]{\iota}\\ 
			 &\left\langle\parbox{7cm}{$e_2,\dots, e_m,\hat{e}_2,\dots, \hat{e}_m,
                e'_2,\dots, e'_{n-1},\hat{e}'_2,\dots, \hat{e}'_{n-1},\\
                f_{l,m+2},\dots, f_{l,m+n-1},
			 \hat{f}_{l,m+2},\dots, \hat{f}_{l,m+n-1}, \\f_{ij},\hat{f}_{ij}\mid 2\leq i<j\leq m, m+2\leq i<j\leq m+n-1, 2\leq l\leq m$}\right\rangle  \arrow{u}{\beta} & \left\langle\parbox{3cm}{$e_2,\dots, e_m,\hat{e}_2,\dots, \hat{e}_m,\\
                f_{l,m+2},\dots,f_{l,m+n-1},\\
              \hat{f}_{l,m+2},\dots, \hat{f}_{l,m+n-1}\mid 2\leq l\leq m$}\right\rangle\arrow{u}[swap]{\gamma}\arrow{l}[swap]{\varsigma}\\ 
\end{tikzcd}
\end{equation}

With the given basis for each vector space, the behavior of the morphisms in $N$ can be described as follows:\\

\begin{itemize}
\item $\alpha$ maps every basis element of $\mathbb{F}^{m^2+n^2-3m-3n+6}$ to $f$;
\item For all $2\leq l\leq m$, $\beta$ maps $e_2, \dots, e_m, f_{l,m+2}, \dots, f_{l,m+n-1} $ to $f_{2,m+2}$ and maps $\hat{e}_1, \dots \hat{e}_m, \hat{f}_{l,m+2}, \dots, \hat{f}_{l,m+n-1}$ to $\hat{f}_{2,m+2}$. Moreover, $\beta$ maps $e'_k$ to $e'_k$ and maps $\hat{e}'_k$ to $\hat{e}'_k$ for all $2\leq k\leq n-1$. In addition, $\beta$ maps $f_{ij}$ to $f_{ij}$ and $\hat{f}_{ij}$ to $\hat{f}_{ij}$ for all $2\leq i<j\leq m$ and $m+2\leq i<j\leq m+n-1$. 
\item $\gamma$ maps $f_{l,m+k}$ to $f_{2,m+2}$ and maps $\hat{f}_{l,m+k}$ to $\hat{f}_{2,m+2}$ for all $2\leq k\leq n-1$ and $2\leq l\leq m$. Moreover, $\gamma$ maps $e_i$ to $f_{2,m+2}$ and $\hat{e}_i$ to $\hat{f}_{2,m+2}$ for all $i=2,\dots,m$;
\item $\iota$ and $\varsigma$ are inclusions;
\item unlabeled vertical maps are the identity maps.
\end{itemize}

Similarly, we associate $PH_1((\hat{\mathcal{H}}_{-,-})_{r,L_{e_1}}^2;\mathbb{F})$ with a persistence module $N'$, where $N'$ is constructed using Proposition \ref{prop-sim-2}.

\begin{equation}\label{h1genHmn}
\begin{aligned}
N'=\begin{tikzcd}
  \mathbb{F}^{m^2-3m+n^2-3n+3} & 0 \arrow{l}[swap]{} &0\arrow{l}[swap]{}\\
  \mathbb{F}^{m^2-3m+n^2-3n+3}  \arrow{u}{}  &0 \arrow{u}{}\arrow{l}[swap]{}  & 0\arrow{u}[swap]{}\arrow{l}[swap]{}\\ 
  \mathbb{F}^{m^2+n^2+2mn-7m-7n+11} \arrow{u}{\eta}	 &\mathbb{F}^{2(m-2)(n-2)} \arrow{u}{}\arrow{l}[swap]{\zeta}	 & 0\arrow{u}[swap]{}\\ 
	  	  	 &\mathbb{F}^{2(m-2)(n-2)} \arrow{u}{} & \mathbb{F}^{2(m-2)(n-2)}\arrow{u}[swap]{}\arrow{l}[swap]{}\\ 
			 &0 \arrow{u}{} & 0\arrow{u}[swap]{}\arrow{l}[swap]{}\\ 
\end{tikzcd}\\
\end{aligned}
\end{equation}
In (\ref{h1genHmn}), $\zeta$ is the inclusion map, and $\eta\circ \zeta=0$. Moreover, all unlabeled non-trivial vertical maps are identity maps.\\

\begin{theorem}
$PH_1((\hat{\mathcal{H}}_{-,-})_{a,b}^2;\mathbb{F})$ is interval decomposable.
\end{theorem}
\begin{proof}
We denote the restriction of $N'$ to its support by $M$. Note that $\eta\circ \zeta=0$ and the endomorphism of $\mathbb{F}^{2(m-2)(n-2)}$ is the identity map, see (\ref{h1genHmn}). Moreover, $\zeta$ is the inclusion map, hence $M$ can be decomposed into two subrepresentations $M_1$ and $M_2$:\\

\begin{minipage}{0.45\textwidth}
   \begin{equation*}\label{h1genHmn'1}
\begin{aligned}
M_1=\begin{tikzcd}
  0					& 	 					&\\
  0  \arrow{u}{} 		&						&\\ 
 \mathbb{F}^{2(m-2)(n-2)}  \arrow{u}{}	 &\mathbb{F}^{2(m-2)(n-2)} \arrow{l}[swap]{\id}	 &	\\ 
	  	  	 			&\mathbb{F}^{2(m-2)(n-2)} \arrow{u}{} & \mathbb{F}^{2(m-2)(n-2)}\arrow{l}[swap]{}
\end{tikzcd}
\end{aligned}
\end{equation*} 
\end{minipage}%
\begin{minipage}{0.45\textwidth}
  \begin{equation*}\label{h1genHmn'2}
\begin{aligned}
M_2=\begin{tikzcd}
  \mathbb{F}^{m^2-3m+n^2-3n+3} 					& 	 					&\\
  \mathbb{F}^{m^2-3m+n^2-3n+3}  \arrow{u}{} 		&						&\\ 
   \mathbb{F}^{m^2-3m+n^2-3n+3}  \arrow{u}{}	 &0 \arrow{l}[swap]{}	 &	\\ 
	  	  	 			&0 \arrow{u}{} &0\arrow{l}[swap]{}
\end{tikzcd}
\end{aligned}
\end{equation*}  
\end{minipage}
\vspace{1em}

Since $M_1$ is an $A_4$-quiver, and $M_2$ is an $A_3$-quiver, they are interval decomposable. Note that the intervals in the support of $M_1$ and $M_2$ are also intervals in the support of $N'$. Therefore, we can extend each indecomposable representation of $M_i$ (where $i=1,2$) to a subrepresentation of $N'$ by putting $0$ to every vertex of $N'$ which is not in $M_i$ and a trivial morphism between two vertices where at least one vertex is not in the support of $M_i$. Hence, $N'$ is interval decomposable. Note that $PH_1((\hat{\mathcal{H}}_{-,-})_{a,b}^2;\mathbb{F})\cong N'$, we conclude that $PH_1((\hat{\mathcal{H}}_{-,-})_{a,b}^2;\mathbb{F})$ is interval decomposable.
\end{proof}

\begin{theorem}
$PH_0((\hat{\mathcal{H}}_{-,-})_{r,L_{e_1}}^2;\mathbb{F})$ is interval decomposable. 
\end{theorem}
\begin{proof}
Equation (\ref{h0genHmn2}) provides us with the behavior of each arrow with given basis elements. Our goal is to find a new basis of each vector space in $N$ such that each morphism in (\ref{h0genHmn2}) maps every basis to another basis or zero, depending on the geometry of $(\hat{\mathcal{H}}_{m,n})_{r,L_{e_1}}^2$. We choose the basis for each vector space provided in the diagram below:\\ 
\begin{equation}\label{h0genHmn2'}
\begin{tikzcd}
<f>                  & \left\langle\parbox{8cm}{$e'_k-f_{2,m+2}, \hat{e}'_k-\hat{f}_{2,m+2}, f_{2,m+2},\hat{f}_{2,m+2}-f_{2,m+2}, f_{ij}-f_{2,m+2},\hat{f}_{ij}-\hat{f}_{2,m+2}\mid 2\leq k\leq n-1, 2\leq i<j\leq m, m+2\leq i<j\leq m+n-1$}\right\rangle \arrow{l}[swap]{\alpha} 				&<f_{2,m+2},\hat{f}_{2,m+2}-f_{2,m+2}>\arrow{l}[swap]{\iota}\\
<f> \arrow{u}{}  &\left\langle\parbox{8cm}{$e'_k-f_{2,m+2}, \hat{e}'_k-\hat{f}_{2,m+2}, f_{2,m+2},\hat{f}_{2,m+2}-f_{2,m+2}, f_{ij}-f_{2,m+2},\hat{f}_{ij}-\hat{f}_{2,m+2}\mid 2\leq k\leq n-1, 2\leq i<j\leq m, m+2\leq i<j\leq m+n-1$}\right\rangle \arrow{u}{}\arrow{l}[swap]{\alpha}  		& <f_{2,m+2},\hat{f}_{2,m+2}-f_{2,m+2}>\arrow{u}[swap]{}\arrow{l}[swap]{\iota}\\ 
<f>\arrow{u}{}	 &\left\langle\parbox{8cm}{$e'_k-f_{2,m+2}, \hat{e}'_k-\hat{f}_{2,m+2}, f_{2,m+2},\hat{f}_{2,m+2}-f_{2,m+2}, f_{ij}-f_{2,m+2},\hat{f}_{ij}-\hat{f}_{2,m+2}\mid 2\leq k\leq n-1, 2\leq i<j\leq m, m+2\leq i<j\leq m+n-1$}\right\rangle \arrow{u}{}\arrow{l}[swap]{\alpha}	 	& <f_{2,m+2},\hat{f}_{2,m+2}-f_{2,m+2}>\arrow{u}[swap]{}\\ 
	  	  	 &\left\langle\parbox{8cm}{$e'_k-f_{2,m+2}, \hat{e}'_k-\hat{f}_{2,m+2}, f_{2,m+2},\hat{f}_{2,m+2}-f_{2,m+2}, f_{ij}-f_{2,m+2},\hat{f}_{ij}-\hat{f}_{2,m+2}\mid 2\leq k\leq n-1, 2\leq i<j\leq m, m+2\leq i<j\leq m+n-1$}\right\rangle \arrow{u}{} 					& <f_{2,m+2},\hat{f}_{2,m+2}-f_{2,m+2}>\arrow{u}[swap]{}\arrow{l}[swap]{\iota}\\ 
			 &<A>\arrow{u}{\beta} & <B>\arrow{u}[swap]{\gamma}\arrow{l}[swap]{\varsigma}\\ 
\end{tikzcd}
\end{equation}
where $A$ consists of the following elements:\\
\begin{enumerate}
\item $e_r-f_{2,m+2}$ and $\hat{r}_k-\hat{f}_{2,m+2}$ for all $2\leq r\leq m$;
\item $e'_k-f_{2,m+2}$ and $\hat{e}'_k-\hat{f}_{2,m+2}$ for all $2\leq k\leq n-1$;
\item $f_{ij}-f_{2,m+2},\hat{f}_{ij}-\hat{f}_{2,m+2}$ for all $2\leq i<j\leq m$ and $m+2\leq i<j\leq m+n-1$;
\item $f_{2,m+2}, f_{2,m+3}-f_{2,m+2}, \dots, f_{2,m+n-1}-f_{2,m+2}$;
\item $f_{s,m+t}-f_{2,m+t}$ for all $3\leq s\leq m$ and $2\leq t\leq n-1$;
\item $\hat{f}_{l,m+t}-f_{l,m+t}$ for all $2\leq l\leq m$ and $2\leq t\leq n-1$
\end{enumerate}
 and $B$ consists of the following elements:\\
 \begin{enumerate}
 \item $e_r-f_{2,m+2}$ and $\hat{r}_k-\hat{f}_{2,m+2}$ for all $2\leq r\leq m$;
\item $f_{2,m+2}, f_{2,m+3}-f_{2,m+2}, \dots, f_{2,m+n-1}-f_{2,m+2}$;
\item $f_{s,m+t}-f_{2,m+t}$ for all $3\leq s\leq m$ and $2\leq t\leq n-1$;
\item $\hat{f}_{l,m+t}-f_{l,m+t}$ for all $2\leq l\leq m$ and $2\leq t\leq n-1$
\end{enumerate}

Define
\begin{equation}\label{n1'}
\begin{aligned}
N_1&=\begin{tikzcd}
<f>                  & <f_{2,m+2}>  \arrow{l}[swap]{\id} 				&<f_{2,m+2}>  \arrow{l}[swap]{\id}\\
<f> \arrow{u}{}  &<f_{2,m+2}> \arrow{u}{}\arrow{l}[swap]{\id}  		& <f_{2,m+2}>  \arrow{u}[swap]{}\arrow{l}[swap]{\id}\\ 
<f>\arrow{u}{}	 &<f_{2,m+2}>  \arrow{u}{}\arrow{l}[swap]{\id}	 	& <f_{2,m+2}>  \arrow{u}[swap]{}\\ 
	  	  	 &<f_{2,m+2}>  \arrow{u}{} 					& <f_{2,m+2}>  \arrow{u}[swap]{}\arrow{l}[swap]{\id}\\ 
			 &<f_{2,m+2}>  \arrow{u}{\beta} & <f_{2,m+2}>  \arrow{u}[swap]{\gamma}\arrow{l}[swap]{\id}\\ 
\end{tikzcd}\cong\begin{tikzcd}
\mathbb{F}                  & \mathbb{F}  \arrow{l}[swap]{\id} 				&\mathbb{F}  \arrow{l}[swap]{\id}\\
\mathbb{F} \arrow{u}{}  &\mathbb{F} \arrow{u}{}\arrow{l}[swap]{\id}  		& \mathbb{F}  \arrow{u}[swap]{}\arrow{l}[swap]{\id}\\ 
\mathbb{F}\arrow{u}{}	 &\mathbb{F} \arrow{u}{}\arrow{l}[swap]{\id}	 	& \mathbb{F}   \arrow{u}[swap]{\id}\\ 
	  	  	 &\mathbb{F} \arrow{u}{\id} 					& \mathbb{F}   \arrow{u}[swap]{}\arrow{l}[swap]{\id}\\ 
			 &\mathbb{F}\arrow{u}{\beta} 				& \mathbb{F}\arrow{u}[swap]{\gamma}\arrow{l}[swap]{\id}\\ 
\end{tikzcd}
\end{aligned}
\end{equation}

\begin{equation}\label{n2'}
\begin{aligned}
N_2&=\begin{tikzcd}
0                  & <\hat{f}_{2,m+2}-f_{2,m+2}>  \arrow{l}[swap]{} 				&<\hat{f}_{2,m+2}-f_{2,m+2}>   \arrow{l}[swap]{\id}\\
0 \arrow{u}{}  &<\hat{f}_{2,m+2}-f_{2,m+2}>  \arrow{u}{}\arrow{l}[swap]{}  		& <\hat{f}_{2,m+2}-f_{2,m+2}>   \arrow{u}[swap]{}\arrow{l}[swap]{\id}\\ 
0\arrow{u}{}	 &<\hat{f}_{2,m+2}-f_{2,m+2}>   \arrow{u}{}\arrow{l}[swap]{}	 	& <\hat{f}_{2,m+2}-f_{2,m+2}>   \arrow{u}[swap]{}\\ 
	  	  	 &<\hat{f}_{2,m+2}-f_{2,m+2}>   \arrow{u}{} 					& <\hat{f}_{2,m+2}-f_{2,m+2}>   \arrow{u}[swap]{}\arrow{l}[swap]{\id}\\ 
			 &<\hat{f}_{2,m+2}-f_{2,m+2}>  \arrow{u}{\beta} & <\hat{f}_{2,m+2}-f_{2,m+2}>  \arrow{u}[swap]{\gamma}\arrow{l}[swap]{\id}\\ 
\end{tikzcd}\cong\begin{tikzcd}
0                  & \mathbb{F}  \arrow{l}[swap]{} 				&\mathbb{F}  \arrow{l}[swap]{}\\
0 \arrow{u}{}  &\mathbb{F} \arrow{u}{}\arrow{l}[swap]{}  		& \mathbb{F}  \arrow{u}[swap]{}\arrow{l}[swap]{}\\ 
0\arrow{u}{}	 &\mathbb{F} \arrow{u}{}\arrow{l}[swap]{}	 	& \mathbb{F}   \arrow{u}[swap]{}\\ 
	  	  	 &\mathbb{F} \arrow{u}{} 					& \mathbb{F}   \arrow{u}[swap]{}\arrow{l}[swap]{}\\ 
			 &\mathbb{F}\arrow{u}{\beta} 				& \mathbb{F}\arrow{u}[swap]{\gamma}\arrow{l}[swap]{\id}\\ 
\end{tikzcd}
\end{aligned}
\end{equation}

For $2\leq k\leq n-1$
\begin{equation}\label{ne1'}
\begin{aligned}
E_k&=\begin{tikzcd}
0                  & <e'_k-f_{2,m+2}>  \arrow{l}[swap]{} 				&0  \arrow{l}[swap]{}\\
0 \arrow{u}{}  &<e'_k-f_{2,m+2}>  \arrow{u}{}\arrow{l}[swap]{}  		& 0  \arrow{u}[swap]{}\arrow{l}[swap]{}\\ 
0\arrow{u}{}	 &<e'_k-f_{2,m+2}>  \arrow{u}{}\arrow{l}[swap]{}	 	& 0   \arrow{u}[swap]{}\\ 
	  	  	 &<e'_k-f_{2,m+2}>   \arrow{u}{} 					& 0   \arrow{u}[swap]{}\arrow{l}[swap]{}\\ 
			 &<e'_k-f_{2,m+2}> \arrow{u}{\beta} 				&0  \arrow{u}[swap]{}\arrow{l}[swap]{}\\ 
\end{tikzcd}\cong\begin{tikzcd}
0                  & \mathbb{F}  \arrow{l}[swap]{} 				&0  \arrow{l}[swap]{}\\
0 \arrow{u}{}  &\mathbb{F} \arrow{u}{}\arrow{l}[swap]{}  		& 0  \arrow{u}[swap]{}\arrow{l}[swap]{}\\ 
0\arrow{u}{}	 &\mathbb{F} \arrow{u}{}\arrow{l}[swap]{}	 	& 0   \arrow{u}[swap]{}\\ 
	  	  	 &\mathbb{F} \arrow{u}{} 					& 0   \arrow{u}[swap]{}\arrow{l}[swap]{}\\ 
			 &\mathbb{F}\arrow{u}{\beta} 				& 0\arrow{u}[swap]{}\arrow{l}[swap]{}\\ 
\end{tikzcd}
\end{aligned}
\end{equation}
and 
\begin{equation}\label{ne2'}
\begin{aligned}
\hat{E}_k&=\begin{tikzcd}
0                  & <\hat{e}'_k-\hat{f}_{2,m+2}>  \arrow{l}[swap]{} 				&0  \arrow{l}[swap]{}\\
0 \arrow{u}{}  &<\hat{e}'_k-\hat{f}_{2,m+2}>  \arrow{u}{}\arrow{l}[swap]{}  		& 0  \arrow{u}[swap]{}\arrow{l}[swap]{}\\ 
0\arrow{u}{}	 &<\hat{e}'_k-\hat{f}_{2,m+2}> \arrow{u}{}\arrow{l}[swap]{}	 	& 0   \arrow{u}[swap]{}\\ 
	  	  	 &<\hat{e}'_k-\hat{f}_{2,m+2}>  \arrow{u}{} 					& 0   \arrow{u}[swap]{}\arrow{l}[swap]{}\\ 
			 &<\hat{e}'_k-\hat{f}_{2,m+2}>\arrow{u}{\beta} 				& 0\arrow{u}[swap]{}\arrow{l}[swap]{}\\ 
\end{tikzcd}\cong\begin{tikzcd}
0                  & \mathbb{F}  \arrow{l}[swap]{} 				&0  \arrow{l}[swap]{}\\
0 \arrow{u}{}  &\mathbb{F} \arrow{u}{}\arrow{l}[swap]{}  		& 0  \arrow{u}[swap]{}\arrow{l}[swap]{}\\ 
0\arrow{u}{}	 &\mathbb{F} \arrow{u}{}\arrow{l}[swap]{}	 	& 0   \arrow{u}[swap]{}\\ 
	  	  	 &\mathbb{F} \arrow{u}{} 					& 0   \arrow{u}[swap]{}\arrow{l}[swap]{}\\ 
			 &\mathbb{F}\arrow{u}{\beta} 				& 0\arrow{u}[swap]{}\arrow{l}[swap]{}\\ 
\end{tikzcd}
\end{aligned}
\end{equation}

For $ 2\leq i<j\leq m, m+2\leq i<j\leq m+n-1$ 
\begin{equation}\label{ne3}
\begin{aligned}
F_{ij}&=\begin{tikzcd}
0                  & <f_{ij}-f_{2,m+2}>  \arrow{l}[swap]{} 				&0  \arrow{l}[swap]{}\\
0 \arrow{u}{}  &<f_{ij}-f_{2,m+2}>  \arrow{u}{}\arrow{l}[swap]{}  		& 0  \arrow{u}[swap]{}\arrow{l}[swap]{}\\ 
0\arrow{u}{}	 &<f_{ij}-f_{2,m+2}>  \arrow{u}{}\arrow{l}[swap]{}	 	& 0   \arrow{u}[swap]{}\\ 
	  	  	 &<f_{ij}-f_{2,m+2}>   \arrow{u}{} 					& 0   \arrow{u}[swap]{}\arrow{l}[swap]{}\\ 
			 &<f_{ij}-f_{2,m+2}> \arrow{u}{\beta} 				&0  \arrow{u}[swap]{}\arrow{l}[swap]{}\\ 
\end{tikzcd}\cong\begin{tikzcd}
0                  & \mathbb{F}  \arrow{l}[swap]{} 				&0  \arrow{l}[swap]{}\\
0 \arrow{u}{}  &\mathbb{F} \arrow{u}{}\arrow{l}[swap]{}  		& 0  \arrow{u}[swap]{}\arrow{l}[swap]{}\\ 
0\arrow{u}{}	 &\mathbb{F} \arrow{u}{}\arrow{l}[swap]{}	 	& 0   \arrow{u}[swap]{}\\ 
	  	  	 &\mathbb{F} \arrow{u}{} 					& 0   \arrow{u}[swap]{}\arrow{l}[swap]{}\\ 
			 &\mathbb{F}\arrow{u}{\beta} 				& 0\arrow{u}[swap]{}\arrow{l}[swap]{}\\ 
\end{tikzcd}
\end{aligned}
\end{equation}
and 
\begin{equation}\label{ne4}
\begin{aligned}
\hat{F}_{ij}&=\begin{tikzcd}
0                  & <\hat{f}_{ij}-\hat{f}_{2,m+2}>  \arrow{l}[swap]{} 				&0  \arrow{l}[swap]{}\\
0 \arrow{u}{}  &<\hat{f}_{ij}-\hat{f}_{2,m+2}>  \arrow{u}{}\arrow{l}[swap]{}  		& 0  \arrow{u}[swap]{}\arrow{l}[swap]{}\\ 
0\arrow{u}{}	 &<\hat{f}_{ij}-\hat{f}_{2,m+2}> \arrow{u}{}\arrow{l}[swap]{}	 	& 0   \arrow{u}[swap]{}\\ 
	  	  	 &<\hat{f}_{ij}-\hat{f}_{2,m+2}>  \arrow{u}{} 					& 0   \arrow{u}[swap]{}\arrow{l}[swap]{}\\ 
			 &<\hat{f}_{ij}-\hat{f}_{2,m+2}>\arrow{u}{\beta} 				& 0\arrow{u}[swap]{}\arrow{l}[swap]{}\\ 
\end{tikzcd}\cong\begin{tikzcd}
0                  & \mathbb{F}  \arrow{l}[swap]{} 				&0  \arrow{l}[swap]{}\\
0 \arrow{u}{}  &\mathbb{F} \arrow{u}{}\arrow{l}[swap]{}  		& 0  \arrow{u}[swap]{}\arrow{l}[swap]{}\\ 
0\arrow{u}{}	 &\mathbb{F} \arrow{u}{}\arrow{l}[swap]{}	 	& 0   \arrow{u}[swap]{}\\ 
	  	  	 &\mathbb{F} \arrow{u}{} 					& 0   \arrow{u}[swap]{}\arrow{l}[swap]{}\\ 
			 &\mathbb{F}\arrow{u}{\beta} 				& 0\arrow{u}[swap]{}\arrow{l}[swap]{}\\ 
\end{tikzcd}
\end{aligned}
\end{equation}

For all $3\leq t\leq n-1$
\begin{equation}\label{nf1'}
\begin{aligned}
G_{2,t}&=\begin{tikzcd}
0                  & 0  \arrow{l}[swap]{} 				&0  \arrow{l}[swap]{}\\
0 \arrow{u}{}  &0 \arrow{u}{}\arrow{l}[swap]{}  		& 0  \arrow{u}[swap]{}\arrow{l}[swap]{}\\ 
0\arrow{u}{}	 &0 \arrow{u}{}\arrow{l}[swap]{}	 	& 0   \arrow{u}[swap]{}\\ 
	  	  	 &0 \arrow{u}{} 					& 0   \arrow{u}[swap]{}\arrow{l}[swap]{}\\ 
			 &<f_{2,m+t}-f_{2,m+2}>\arrow{u}{} 				& <f_{2,m+t}-f_{2,m+2}>\arrow{u}[swap]{}\arrow{l}[swap]{\id}\\ 
\end{tikzcd}\cong\begin{tikzcd}
0                  & 0  \arrow{l}[swap]{} 				&0  \arrow{l}[swap]{}\\
0 \arrow{u}{}  &0 \arrow{u}{}\arrow{l}[swap]{}  		& 0  \arrow{u}[swap]{}\arrow{l}[swap]{}\\ 
0\arrow{u}{}	 &0 \arrow{u}{}\arrow{l}[swap]{}	 	& 0   \arrow{u}[swap]{}\\ 
	  	  	 &0 \arrow{u}{} 					& 0   \arrow{u}[swap]{}\arrow{l}[swap]{}\\ 
			 &\mathbb{F}\arrow{u}{} 				& \mathbb{F}\arrow{u}[swap]{}\arrow{l}[swap]{\id}\\ 
\end{tikzcd}
\end{aligned}
\end{equation}
and
\begin{equation}\label{nf2'}
\begin{aligned}
\hat{G}_{2,t}&=\begin{tikzcd}
0                  & 0  \arrow{l}[swap]{} 				&0  \arrow{l}[swap]{}\\
0 \arrow{u}{}  &0 \arrow{u}{}\arrow{l}[swap]{}  		& 0  \arrow{u}[swap]{}\arrow{l}[swap]{}\\ 
0\arrow{u}{}	 &0 \arrow{u}{}\arrow{l}[swap]{}	 	& 0   \arrow{u}[swap]{}\\ 
	  	  	 &0 \arrow{u}{} 					& 0   \arrow{u}[swap]{}\arrow{l}[swap]{}\\ 
			 &<\hat{f}_{2,m+t}-f_{2,m+t}>\arrow{u}{} 				& <\hat{f}_{2,m+t}-f_{2,m+t}>\arrow{u}[swap]{}\arrow{l}[swap]{\id}\\ 
\end{tikzcd}\cong
\begin{tikzcd}
0                  & 0  \arrow{l}[swap]{} 				&0  \arrow{l}[swap]{}\\
0 \arrow{u}{}  &0 \arrow{u}{}\arrow{l}[swap]{}  		& 0  \arrow{u}[swap]{}\arrow{l}[swap]{}\\ 
0\arrow{u}{}	 &0 \arrow{u}{}\arrow{l}[swap]{}	 	& 0   \arrow{u}[swap]{}\\ 
	  	  	 &0 \arrow{u}{} 					& 0   \arrow{u}[swap]{}\arrow{l}[swap]{}\\ 
			 &\mathbb{F}\arrow{u}{} 				& \mathbb{F}\arrow{u}[swap]{}\arrow{l}[swap]{\id}\\ 
\end{tikzcd}
\end{aligned}
\end{equation}

For all $2\leq t\leq n-1$ and $3\leq s\leq m$, 
\begin{equation}\label{nf3}
\begin{aligned}
G_{s,t}&=\begin{tikzcd}
0                  & 0  \arrow{l}[swap]{} 				&0  \arrow{l}[swap]{}\\
0 \arrow{u}{}  &0 \arrow{u}{}\arrow{l}[swap]{}  		& 0  \arrow{u}[swap]{}\arrow{l}[swap]{}\\ 
0\arrow{u}{}	 &0 \arrow{u}{}\arrow{l}[swap]{}	 	& 0   \arrow{u}[swap]{}\\ 
	  	  	 &0 \arrow{u}{} 					& 0   \arrow{u}[swap]{}\arrow{l}[swap]{}\\ 
			 &<f_{s,m+t}-f_{2,m+t}>\arrow{u}{} 				& <f_{s,m+t}-f_{2,m+t}>\arrow{u}[swap]{}\arrow{l}[swap]{\id}\\ 
\end{tikzcd}\cong\begin{tikzcd}
0                  & 0  \arrow{l}[swap]{} 				&0  \arrow{l}[swap]{}\\
0 \arrow{u}{}  &0 \arrow{u}{}\arrow{l}[swap]{}  		& 0  \arrow{u}[swap]{}\arrow{l}[swap]{}\\ 
0\arrow{u}{}	 &0 \arrow{u}{}\arrow{l}[swap]{}	 	& 0   \arrow{u}[swap]{}\\ 
	  	  	 &0 \arrow{u}{} 					& 0   \arrow{u}[swap]{}\arrow{l}[swap]{}\\ 
			 &\mathbb{F}\arrow{u}{} 				& \mathbb{F}\arrow{u}[swap]{}\arrow{l}[swap]{\id}\\ 
\end{tikzcd}
\end{aligned}
\end{equation}
and
\begin{equation}\label{nf4}
\begin{aligned}
\hat{G}_{s,t}&=\begin{tikzcd}
0                  & 0  \arrow{l}[swap]{} 				&0  \arrow{l}[swap]{}\\
0 \arrow{u}{}  &0 \arrow{u}{}\arrow{l}[swap]{}  		& 0  \arrow{u}[swap]{}\arrow{l}[swap]{}\\ 
0\arrow{u}{}	 &0 \arrow{u}{}\arrow{l}[swap]{}	 	& 0   \arrow{u}[swap]{}\\ 
	  	  	 &0 \arrow{u}{} 					& 0   \arrow{u}[swap]{}\arrow{l}[swap]{}\\ 
			 &<\hat{f}_{s,m+t}-\hat{f}_{2,m+t}>\arrow{u}{} 				& <\hat{f}_{s,m+t}-\hat{f}_{2,m+t}>\arrow{u}[swap]{}\arrow{l}[swap]{\id}\\ 
\end{tikzcd}\cong\begin{tikzcd}
0                  & 0  \arrow{l}[swap]{} 				&0  \arrow{l}[swap]{}\\
0 \arrow{u}{}  &0 \arrow{u}{}\arrow{l}[swap]{}  		& 0  \arrow{u}[swap]{}\arrow{l}[swap]{}\\ 
0\arrow{u}{}	 &0 \arrow{u}{}\arrow{l}[swap]{}	 	& 0   \arrow{u}[swap]{}\\ 
	  	  	 &0 \arrow{u}{} 					& 0   \arrow{u}[swap]{}\arrow{l}[swap]{}\\ 
			 &\mathbb{F}\arrow{u}{} 				& \mathbb{F}\arrow{u}[swap]{}\arrow{l}[swap]{\id}\\ 
\end{tikzcd}
\end{aligned}
\end{equation}

For all $2\leq r\leq m$, 
\begin{equation}\label{nf5}
\begin{aligned}
K_r&=\begin{tikzcd}
0                  & 0  \arrow{l}[swap]{} 				&0  \arrow{l}[swap]{}\\
0 \arrow{u}{}  &0 \arrow{u}{}\arrow{l}[swap]{}  		& 0  \arrow{u}[swap]{}\arrow{l}[swap]{}\\ 
0\arrow{u}{}	 &0 \arrow{u}{}\arrow{l}[swap]{}	 	& 0   \arrow{u}[swap]{}\\ 
	  	  	 &0 \arrow{u}{} 					& 0   \arrow{u}[swap]{}\arrow{l}[swap]{}\\ 
			 &<e_r-f_{2,m+t}>\arrow{u}{} 				& <e_r-f_{2,m+t}>\arrow{u}[swap]{}\arrow{l}[swap]{\id}\\ 
\end{tikzcd}\cong\begin{tikzcd}
0                  & 0  \arrow{l}[swap]{} 				&0  \arrow{l}[swap]{}\\
0 \arrow{u}{}  &0 \arrow{u}{}\arrow{l}[swap]{}  		& 0  \arrow{u}[swap]{}\arrow{l}[swap]{}\\ 
0\arrow{u}{}	 &0 \arrow{u}{}\arrow{l}[swap]{}	 	& 0   \arrow{u}[swap]{}\\ 
	  	  	 &0 \arrow{u}{} 					& 0   \arrow{u}[swap]{}\arrow{l}[swap]{}\\ 
			 &\mathbb{F}\arrow{u}{} 				& \mathbb{F}\arrow{u}[swap]{}\arrow{l}[swap]{\id}\\ 
\end{tikzcd}
\end{aligned}
\end{equation}
and
\begin{equation}\label{nf6}
\begin{aligned}
\hat{K}_{r}&=\begin{tikzcd}
0                  & 0  \arrow{l}[swap]{} 				&0  \arrow{l}[swap]{}\\
0 \arrow{u}{}  &0 \arrow{u}{}\arrow{l}[swap]{}  		& 0  \arrow{u}[swap]{}\arrow{l}[swap]{}\\ 
0\arrow{u}{}	 &0 \arrow{u}{}\arrow{l}[swap]{}	 	& 0   \arrow{u}[swap]{}\\ 
	  	  	 &0 \arrow{u}{} 					& 0   \arrow{u}[swap]{}\arrow{l}[swap]{}\\ 
			 &<\hat{e}_{r}-\hat{f}_{2,m+t}>\arrow{u}{} 				& <\hat{e}_{r}-\hat{f}_{2,m+t}>\arrow{u}[swap]{}\arrow{l}[swap]{\id}\\ 
\end{tikzcd}\cong\begin{tikzcd}
0                  & 0  \arrow{l}[swap]{} 				&0  \arrow{l}[swap]{}\\
0 \arrow{u}{}  &0 \arrow{u}{}\arrow{l}[swap]{}  		& 0  \arrow{u}[swap]{}\arrow{l}[swap]{}\\ 
0\arrow{u}{}	 &0 \arrow{u}{}\arrow{l}[swap]{}	 	& 0   \arrow{u}[swap]{}\\ 
	  	  	 &0 \arrow{u}{} 					& 0   \arrow{u}[swap]{}\arrow{l}[swap]{}\\ 
			 &\mathbb{F}\arrow{u}{} 				& \mathbb{F}\arrow{u}[swap]{}\arrow{l}[swap]{\id}\\ 
\end{tikzcd}
\end{aligned}
\end{equation}

It is clear that 
\begin{equation}
\begin{aligned}
N \cong N_1&\oplus N_2 \oplus  (\bigoplus\limits_{k=2}^{n-1} E_k) \oplus  (\bigoplus\limits_{k=2}^{n-1} \hat{E}_k)\oplus (\bigoplus\limits_{2\leq i<j\leq m} F_{ij})\oplus(\bigoplus\limits_{2\leq i<j\leq m} \hat{F}_{ij})\\
&\oplus  (\bigoplus\limits_{m+2\leq i<j\leq m+n-1} F_{ij})\oplus (\bigoplus\limits_{m+2\leq i<j\leq m+n-1} \hat{F}_{ij})\\
&\oplus(\bigoplus\limits_{3\leq t\leq n-1} G_{2,t})\oplus (\bigoplus\limits_{3\leq t\leq n-1} \hat{G}_{2,t})\\
&\oplus(\bigoplus\limits_{3\leq s\leq m} \bigoplus\limits_{2\leq t\leq n-1} G_{s,t})\oplus(\bigoplus\limits_{3\leq s\leq m} \bigoplus\limits_{2\leq t\leq n-1} \hat{G}_{s,t})\\
&\oplus(\bigoplus\limits_{2\leq r\leq m} K_{r})\oplus (\bigoplus\limits_{2\leq t\leq m} \hat{K}_{r})
\end{aligned}
\end{equation}

Note that $N_1$, $N_2$, $E_k$, $\hat{E}_k$, $F_{ij}$, $\hat{F}_{ij}$, $G_{2,t}$, $\hat{G}_{2,t}$, $G_{s,t}$, $\hat{G}_{s,t}$, $K_r$, and $\hat{K}_r$ are interval modules, hence by Lemma \ref{thinpoly}, they are indecomposable. Therefore, $N$ is interval decomposable. Since 
$PH_0((\hat{\mathcal{H}}_{m,n})_{-,-}^2;\mathbb{F})\cong N$, we conclude that $PH_0((\hat{\mathcal{H}}_{-,-})_{r,L_{e_1}}^2;\mathbb{F})$ is interval decomposable.
\end{proof}

%% file: chapters/section1.2.5.tex
\section{Limits and Colimits}
In this section, a diagram in $\cat{C}$ (over a category $\cat{J}$, called a shape) is a functor $\mathcal{F}: \cat{J}\rightarrow \cat{C}$, where $\cat{J}$ is the path category of a directed graph. 
\begin{example}[Trivial diagram]
Let $c\in\ob\cat{C}$. A trivial diagram in $\cat{C}$ over a given shape $\cat{J}$ is a functor $c: \cat{J}\rightarrow \cat{C}$ where $c$ sends each vertex of $\cat{J}$ to $c$ and each arrow of $\cat{J}$ to $\id_{c}$.
\end{example}
\begin{definition}\index{cone}
Let $\mathcal{F}: \cat{J}\rightarrow \cat{C}$ be a functor. A \textbf{cone over $\mathcal{F}$} with summit $c$ is a natural transformation $\lambda: c\Rightarrow \mathcal{F}$ and a \textbf{cone under $\mathcal{F}$} with nadir $c$ is a natural transformation $\lambda: \mathcal{F}\Rightarrow c$. 
\end{definition}
\begin{definition}\index{limit}
Let $\mathcal{F}: \cat{J}\rightarrow \cat{C}$ be a functor. The \textbf{limit} of $\mathcal{F}$, denoted by $\lim_{\cat{J}}\mathcal{F}$, is an object of $\cat{C}$ along with a universal cone $\lambda: \lim_{\cat{J}}\mathcal{F}\Rightarrow \mathcal{F}$ satisfying the following condition: for all $\mu: c\Rightarrow \mathcal{F}$, there exists a unique natural transformation $\alpha: c\Rightarrow \lim_{\cat{J}}\mathcal{F}$ such that $\mu_i=\lambda_i\circ\alpha_i$ for all $i\in \ob\cat{J}$.
\end{definition}
Equivalently, we can define $\lim_{\cat{J}}\mathcal{F}$ to be the representative of the functor $\mathsf{Cone}(-, \mathcal{F}):\cat{C}\rightarrow \cat{Set}$, in other words, $$\hom_{\cat{C}}(-,\lim_{\cat{J}}\mathcal{F})\cong \mathsf{Cone}(-, \mathcal{F})$$
Dually, we can define colimits. 
\begin{definition}\index{colimit}
Let $\mathcal{F}: \cat{J}\rightarrow \cat{C}$ be a functor. The \textbf{colimit} of $\mathcal{F}$, denoted by $\colim_{\cat{J}}\mathcal{F}$, is an object of $\cat{C}$ along with a universal cone $\lambda: \mathcal{F}\Rightarrow \lim_{\cat{J}}\mathcal{F}$ satisfying the following condition: for all $\mu: \mathcal{F} \Rightarrow c$, there exists a unique natural transformation $\alpha: \colim_{\cat{J}}\mathcal{F}\Rightarrow c$ such that $\mu_i=\alpha_i \circ\lambda_i$ for all $i\in \ob\cat{J}$.
\end{definition}
\begin{example}[Product and Coproduct]
Let $\cat{J}$ be a discrete category and $\mathcal{F}: \cat{J}\rightarrow \cat{C}$. Assume $\cat{C}$ has limits and colimits over $\cat{J}$. The product\index{product} of $\mathcal{F}_j$, denoted by $\prod\limits_{j\in\cat{J}}\mathcal{F}_j$, is the limit of $\mathcal{F}$. In other words, $\prod\limits_{j\in\cat{J}}\mathcal{F}_j=\lim_{\cat{J}}\mathcal{F}$. Dually, the coproduct\index{coproduct} of $\mathcal{F}_j$, denoted by $\coprod\limits_{j\in\cat{J}}\mathcal{F}_j$, is the colimit of $\mathcal{F}$. In other words, $\coprod\limits_{j\in\cat{J}}\mathcal{F}_j=\colim_{\cat{J}}\mathcal{F}$.
\end{example}

A set $\mathcal{I}$ is called a \textbf{directed set} if it has a preorder with an additional property that every pair of elements of $\mathcal{I}$ has an upper bound. A \textbf{directed system} in $\cat{C}$ over the directed set $\mathcal{I}$ (denoted by $\{X_{i}\}_{i\in\mathcal{I}}$) is a functor $X:\mathcal{I}\rightarrow\cat{C}$. Since $\mathcal{I}$ is a preorder, $X$ has the following properties:
\begin{enumerate}
\item for all $i,j\in\mathcal{I}$, $\lvert \hom_{\cat{C}}(X_i,X_j)\rvert\leq 1$;
\item $X(i\leq i)=\id_{X_i}$ for all $i\in\mathcal{I}$;
\item $X(j\leq k)\circ X(i\leq j)=X(i\leq k)$, for all $i,j,k\in\mathcal{I}$ where $i\leq j\leq k$.
\end{enumerate}
\begin{definition}\index{direct limit}
Let $\{X_{i}\}_{i\in\mathcal{I}}$ be a directed system. The \textbf{direct limit} of $\{X_{i}\}_{i\in\mathcal{I}}$ is the colimit of $X$, in other words, $$\underset{}{\dirlim}X_{i}=\colim_{\mathcal{I}}X$$
\end{definition}
Notation: We denote the (unique) map $X(i\leq j): X_{i}\rightarrow X_{j}$ ($i\leq j$) by $f_{ij}$; we use $p_{i}$ represents the map $X_{i}\rightarrow \underset{}{\dirlim}X_{i}$, which is a leg map of the colimit cone. 

\begin{lemma}\label{lemma1-directed}
Let $\{X_{i}\}_{i\in\mathcal{I}}$ be a directed system where $X_{i}\in\cat{Mod}_{R}$ for all $i\in\mathcal{I}$. Then $L$ is the direct limit of $\{X_{i}\}_{i\in\mathcal{I}}$ iff the following conditions hold:\begin{enumerate}
\item\label{c1} Given $x\in L$, there exists $i\in\mathcal{I}$ and $x_{i}\in X_{i}$ such that $p_{i}(x_{i})=x$. 
\item \label{c2} If there exists $x_{i}\in X_{i}$ for some $i\in\mathcal{I}$ such that $p_{i}(x_{i})=0$, there exists $j\in\mathcal{I}$ with $j\geq i$ such that $f_{ij}(x_{i})=0$.
\end{enumerate}
\end{lemma}
\begin{pflemma}
$(\Rightarrow)$ \ \ We first proof (\ref{c1}). If there is an element of $L$, say $x$, which is not in the image of any $p_{i}$ for $i\in\mathcal{I}$. Let $K:=L-\{x\}$. It is clear that $\{X_{i}\}_{i\in\mathcal{I}}\Rightarrow K$ is a cocone (here we treat $K$ as the constant cocone), hence by the universal property of direct limit, there exists a unique map $L\rightarrow K$. Note that there is an inclusion $K\hookrightarrow L$. Therefore, $K=L$, and it contradicts our assumption that $K:=L-\{x\}$. \\

We now prove (\ref{c2}). Assume (for such $x_{i}$) $f_{ij}(x_{i})\neq 0$ for all $j\geq i$. Let $K=L$ and let $q_{i}$ denote the map $X_{i}\rightarrow K$ such that $q_{i}(x_{i})=a\neq 0$ for some $a\in L-\{0\}$ and $q_{i}(y)=p_{i}(y)$ if $y\neq x_{0}$; define $q_{j}f_{ij}(x_{i})=a$. Then we have a cocone $\{X_{i}\}_{i\in\mathcal{I}}\Rightarrow K$. By the universal property of direct limit, there exists a unique map $\phi: L\rightarrow K$. $\phi$ is an $R$-module homomorphism, hence it has to send $p_{i}(x_{i})=0$ to $0$. Therefore, $$a=q_i(x_i)=\phi\circ p_{i}(x_i)=\phi(0)=0$$ Contradicting the assumption $a\neq 0$.\\

$(\Leftarrow)$ We want to show for any cocone under $\{X_i\}_{i\in I}$ with nadir $M$ with leg maps $(m_i)_{i\in I}$, there exists a unique map $$\phi: L\rightarrow M$$ such that $m_i=\phi\circ p_i$ for all $i\in I$. By the condition (\ref{c1}), we are forced to define $\phi(x)=m_i(x_i)$ where $p_i(x_i)=x$ for some $i\in I$. Suppose there is another $x_j$ in $X_j$ such that $p_j(x_j)=x$, WLOG, we assume $i\leq j$. Notice that $p_j(f_{ij}(x_i))=p_i(x_i)$ hence $p_j(x_j-f_{ij}(x_i))=p_j(x_j)-p_i(x_i)=x-x=0$. Therefore, by condition (\ref{c2}), there exists $k\in I$ with $k\geq j$ such that $f_{jk}(x_j-f_{ij}(x_i))=0$. Hence $m_k f_{jk}(x_j-f_{ij}(x_i))=0$, i.e., $m_j(x_j)-m_i(x_i)=0$, which implies $\phi$ is well-defined.

\end{pflemma}

\begin{theorem}
Let $\{X_{i}\}_{i\in\mathcal{I}}$ be a directed system where $X_{i}$ are open subset of an ambient topological space and $X_{i}\rightarrow X_{j}$ is an inclusion map if $i\leq j$. Then 
$$\underset{}{\dirlim}\pi_{n}(X_{i})\cong\pi_{n}( \underset{}{\dirlim}X_{i})$$
\end{theorem}

\begin{proof}
Given $[\sigma]\in \pi_{n}( \underset{}{\dirlim}X_{i})$, then $\sigma$ is a map from $S^{n}$ to $\underset{}{\dirlim}X_{i}=\bigcup X_{i}$. Since $S^{n}$ is compact, there exists $X_{k}$ such that $\Ima{\sigma}\subseteq X_{k}$. That means $[\sigma]\in \pi_{n}(X_{k})$.

On the other hand, if there exists $[\sigma]\in\pi_{n}(X_{i})$ such that $\sigma\sim 0$ in $\underset{}{\dirlim}X_{i}$. Let $H:S^{n}\times I\rightarrow \underset{}{\dirlim}X_{i}$ be the homotopy. Note that $S^{n}\times I$ is compact, there must exist $X_{k}$ such that $\Ima{H}\subseteq X_{k}$. Hence $\sigma\sim 0$ in $X_{k}$.

By the above lemma, we conclude that $\underset{}{\dirlim}\pi_{n}(X_{i})\cong\pi_{n}( \underset{}{\dirlim}X_{i})$.

\end{proof}

\begin{theorem}\label{directlim_exact}
$\underset{}{\dirlim}-$ is exact.
\end{theorem}
\begin{proof}
Let $\{X_{i}\}_{i\in\mathcal{I}}, \{Y_{i}\}_{i\in\mathcal{I}}, \{Z_{i}\}_{i\in\mathcal{I}}$ be directed systems where for each $i\in\mathcal{I}$, $X_{i}\xrightarrow{f_{i}} Y_{i}\xrightarrow{g_{i}} Z_{i}$ is exact at $Y_{i}$. Denote $\alpha_{ij}:Y_{i}\rightarrow Y_{j}$ and $\beta_{ij}:Z_{i}\rightarrow Z_{j}$. Denote $p_{i}: X_{i}\rightarrow \underset{}{\dirlim} X_{i}$, $q_{i}: Y_{i}\rightarrow \underset{}{\dirlim} Y_{i}$ and  $r_{i}: Z_{i}\rightarrow \underset{}{\dirlim} Z_{i}$ We want to show the following sequence is exact:
$$\underset{}{\dirlim} X_{i}\xrightarrow{F} \underset{}{\dirlim} Y_{i}\xrightarrow{G} \underset{}{\dirlim} Z_{i}$$

First, the existence of $F$ and $G$ is clear by the universal property of direct limits.

Second, notice that $g_{i}\circ f_{i}=0$ for all $i\in\mathcal{I}$ (which means the map $g_{i}\circ f_{i}:X_{i}\rightarrow Z_{i}$ factors through $0$), hence we have $G\circ F=0$.

On the other hand, given $y\in\ker{G}$, we have $G(y)=0$. By the previous lemma, there exists $y_{i}\in Y_{i}$ for some $i\in\mathcal{I}$ such that $q_{i}(y_{i})=y$. Hence we have $(G\circ q_{i})(y_{i})=G(y)=0$. Note that $G\circ q_{i}=r_{i}g_{i}$ we have $r_{i}g_{i}(y_{i})=0$. By the previous lemma, there exists $Z_{j}$ such that $\beta_{ij}(g_{i}(y_{i}))=0$. Note that $\beta_{ij}g_{i}=g_{j}\alpha_{ij}$ hence $g_{j}\alpha_{ij}(y_{i})=0$. Note that $X_{j}\xrightarrow{f_{j}} Y_{j}\xrightarrow{g_{j}} Z_{j}$ is exact, there exists $x_{j}\in X_{j}$ such that $f_{j}(x_{j})=\alpha_{ij}(y_{i})$. Hence we have $$F\circ p_{j}(x_{j})=q_{j}f_{j}(x_{j})=q_{j}\alpha_{ij}(y_{i})=q_{i}(y_{i})=y$$

Hence $\ker{G}\subseteq \Ima{F}$.
\end{proof}